\documentclass[preprint,11pt]{elsarticle}

\usepackage{algpseudocode}
\usepackage{algorithm}
\algdef{SE}{Begin}{End}{\textbf{Begin}}{\textbf{End}}

\usepackage{graphicx}
\usepackage{booktabs}
\usepackage{array}
\usepackage{subcaption}
\usepackage{tikz}
\usetikzlibrary{shapes}
\usetikzlibrary{spy}
\usepackage{amsmath, amsthm, amsfonts, amssymb}
\usepackage[utf8x]{inputenc} 
\usepackage[colorlinks=true,urlcolor=blue,linkcolor=blue]{hyperref}
\usepackage{lmodern} 
\usepackage{color}
\usepackage{url}
\usepackage{multirow}
\usepackage{appendix}
\usepackage{epstopdf} 

\newcommand{\abs}[1]{\left\vert#1\right\vert}

\newcommand{\dn}[1]{\partial_{n}{#1}}
\usepackage{bm}
\newcommand{\vect}[1]{\boldsymbol{\mathbf{#1}}}

\DeclareMathOperator*{\argmin}{arg\,min}
\renewcommand{\div}[1]{\operatorname{div}\left(#1\right)}

\newcommand{\VV}{\vect{\theta}}
\newcommand{\Vn}{\theta_{n}}
\newcommand{\un}{u}
\newcommand{\ud}{w}
\newcommand{\uo}{u_{\omega}} 

\newtheorem{theorem}{Theorem}[section]
\newtheorem{definition}{Definition}   
\newtheorem{remark}{Remark}

\newtheorem{problem}[theorem]{Problem}

\usepackage{parskip}

\journal{a journal for publication}

\begin{document}

\begin{frontmatter}


 
\title{A robust alternating direction method of multipliers numerical scheme in a shape optimization setting for solving geometric inverse problems}

\author[jftradd]{Julius Fergy Tiongson Rabago}
\address[jftradd]{Faculty of Mathematics and Physics, Institute of Science and Engineering\\ Kanazawa University, Kanazawa 920-1192, Japan} 
\ead{jftrabago@gmail.com, rabagojft@se.kanazawa-u.ac.jp}

\author[hadriadd]{Aissam Hadri}
\address[hadriadd]{Universit\'e Ibn Zohr, Facult\'e Polydisciplinaire, Ouarzazate, Morocco} 
\ead{aissamhadri20@gmail.com}
	
\author[afraitesadd]{Lekbir Afraites}
\address[afraitesadd]{EMI, FST, B\'{e}ni-Mellal, Universit\'{e} Sultan Moulay Slimane, Morocco} 
\ead{l.afraites@usms.ma}

\author[hendyadd,hendyadd2]{Ahmed S. Hendy}
\address[hendyadd]{Department of Computational Mathematics and Computer Science, Institute of Natural Sciences and Mathematics, Ural Federal University, 19 Mira St., Ekaterinburg 620002, Russia} 
\address[hendyadd2]{Department of Mathematics, Faculty of Science, Benha University, Benha 13511, Egypt} 
\ead{ahmed.hendy@fsc.bu.edu.eg}

\author[zakyadd,zakyadd2]{ Mahmoud A. Zaky\corref{mycorrespondingauthor}}
\address[zakyadd]{Department of Mathematics and Statistics, College of Science, Imam Mohammad Ibn Saud Islamic University, Riyadh, Saudi Arabia} 
\address[zakyadd2]{Department of Applied Mathematics, National Research Centre, Dokki, Cairo 12622, Egypt} 
\ead{ma.zaky@yahoo.com}

\begin{abstract}
The alternating direction method of multipliers within a shape optimization framework is developed for solving geometric inverse problems, focusing on a cavity identification problem from the perspective of non-destructive testing and evaluation techniques.
The rationale behind this method is to achieve more accurate detection of unknown inclusions with pronounced concavities, emphasizing the aspect of shape optimization.
Several numerical results to illustrate the applicability and efficiency of the method are presented for various shape detection problems. These numerical experiments are conducted in both two- and three-dimensional settings, with a focus on cases involving noise-contaminated data. 
The main finding of the study is that the proposed method significantly outperforms conventional shape optimization methods in reconstructing unknown cavity shapes.
\end{abstract}




\end{frontmatter}


\section{Introduction}\label{sec:introduction}
In this study, we propose a novel application of the alternating direction method of multipliers (ADMM) to shape inverse problems in a shape optimization setting. 
Specifically, we address the problem of identifying a perfectly conducting inclusion inside a larger bounded domain from boundary measurements. 
Let $\Omega$ be a given smooth (at least of Lipschitz class), open, and simply connected bounded set in $\mathbb{R}^{d}$, where $d\in \{2,3\}$.
Inside $\Omega$, we assume the existence of an unknown simply connected subdomain denoted by $\omega$ with a ${\mathcal{C}}^{1,1}$ regular boundary such that ${\partial}{\omega} \cap \partial \Omega = \emptyset$.
We fix $\delta > 0$ and define $\mathcal{O}_{\delta}$ as the set of all open subsets $\omega$ strictly included in $\Omega$, with a ${\mathcal{C}}^{1,1}$ boundary, such that $d(x,\partial \Omega) > \delta$ for all $x\in \omega$, and $\Omega\setminus\overline{\omega}$ is connected.
For some technical purposes, we explicitly assume that $\overline{\omega} \subset {\Omega}_{\delta}$, where ${\Omega}_{\delta}$ is another simply connected subdomain of $\Omega$ with a boundary that is $\mathcal{C}^{\infty}$ and lies within the $\delta/2$-neighborhood of ${\partial}{\Omega}$ (e.g., $\{x\in\Omega \mid d(x,{\partial}{\omega}) > \delta/2\} \subset \Omega_{\delta/2} \subset \{x\in\Omega \mid d(x,{\partial}{\omega}) > \delta/2\} $).
Additionally, we assume that the material composing $\omega$ is perfectly conducting, in contrast to the ``background material" inside $\Omega\setminus \overline{\omega}$, which is assumed to have a constant conductivity $\sigma=1$.
For a given Neumann flux $g$ on the accessible boundary ${\partial}{\Omega}$ and Dirichlet boundary measurement $f$, the inclusion $\omega$ and the electrostatic potential $u$ solve the overdetermined problem:
\begin{equation}
\label{eq:overdetermined_problem}
 		-\Delta u = 0 \ \text{in} \ \Omega\setminus\overline{\omega},\quad
		\dn{u} = g \quad \text{and} \quad
	 	 u = f\ \text{on} \ \partial \Omega,\quad
		u = 0 \ \text{on} \ \partial \omega,
\end{equation}
where $\partial_n:=\partial/\partial n$ stands for the outward normal derivative to $\Omega\setminus\overline{\omega}$.
Here, to be precise, we let $g \in H^{1/2}(\partial \Omega)$ and $f \in H^{3/2}(\partial \Omega)$. 
These aforementioned regularities are assumed for simplicity of discussion, and $f$ belongs to the image of the Neumann-to-Dirichlet map $\Upsilon_{{\partial}{\Omega}}: g \in H^{1/2}({\partial}{\Omega}) \mapsto f = \operatorname{trace}(u) \in H^{3/2}({\partial}{\Omega})$, where $u$ solves the equation \eqref{eq:overdetermined_problem} without the Dirichlet boundary condition.

The inverse geometry problem reads as follows:
\begin{equation}\label{eq:inverse_problem}
	\text{find}~~ \omega\in \mathcal{O}_{\delta}~~ \text{and}~~ u ~~\text{that satisfies the overdetermined system}~~ \eqref{eq:overdetermined_problem}.
\end{equation}
The model problem being examined is a specific case of the more general conductivity reconstruction problem and is severely ill-posed in the sense of Hadamard \cite{EpplerHarbrecht2005}.
The problem appears in many applications such as nondestructive testing of materials and has received extensive attention in the literature, with thorough theoretical and numerical investigations; see, e.g., 
\cite{EpplerHarbrecht2005,Afraites2022,AkdumanKress2002,AlessandriniIsakovPowell1995,AlessandriniDiazValenzuela1996,BourgeoisDarde2010,ChapkoKress2005,HettlichRundell1998,Isakov2006} and references therein. 
For instance, the issue of the existence and uniqueness of the solution to the problem from boundary measurement data has been studied by several authors; see, e.g., \cite{AlessandriniIsakovPowell1995,AlessandriniDiazValenzuela1996,BourgeoisDarde2010,Isakov2006}.
We recall in the following theorem the important identifiability result for this inverse problem which shows that the inclusion $\omega$ (and thus the potential $u$) is unique.
\begin{theorem}[{\cite[Thm. 1.1]{BourgeoisDarde2010}}]\label{thm:identifiability_result}
	The Cauchy pair $(f,g)\neq(0,0)$ uniquely determine $\omega$ and $u$ satisfying \eqref{eq:overdetermined_problem}.
\end{theorem}
The proof of this classical theorem is a slight adaptation of the proof given in \cite[Thm. 5.1, p. 106]{ColtonKress2013} for the same problem with the Helmholtz operator instead of the Laplace operator.
Theorem \ref{thm:identifiability_result} implies that we can reasonably attempt to retrieve the unknown obstacle (cavity or void) $\omega$ and the unknown function $u$ from the Cauchy data $(f,g)$ on ${\partial}{\Omega}$.
In the past decades, a vast literature on numerical approaches has emerged.
Some of these methods are based on the parameterization of obstacles \cite{ColtonKress2013}, while others rely on shape optimization techniques \cite{SokolowskiZolesio1992, DelfourZolesio2011, HenrotPierre2018} or topological gradients, as described, for instance, in \cite{CeaGarreauGuillaumeMasmoudi2000, Bonnet2008}.
In the specific case of the Laplace equation in two spatial dimensions, methods based on conformal mappings have also been explored, as shown in \cite{HaddarKress2005}.
Among the studies employing such methods, the obstacle characterized by homogeneous Dirichlet data is addressed in \cite{HaddarKress2005} and further discussed in \cite{KressRundell2005, Rundell2009}.
Another successful approach involves level set techniques (see, e.g., \cite{BourgeoisDarde2010}), which transform the problem of finding geometry into the problem of locating the zero level-set of a function.
Since their introduction in \cite{OsherSethian1998}, level set techniques have been extensively used in the framework of inverse problems, primarily due to their ability to handle topological changes. 
This is illustrated in the context of the inverse obstacle problem, as exemplified in \cite{Santosa1996,Burger1998,LitmanLesselierSantosa2004,BenAmeurBurgerHackl2004}.
In this study, our primary focus is on identifying cavities using ADMM in a shape optimization setting, with an emphasis on achieving more accurate detection of cavity shapes, especially those with pronounced concavities. 
The proposed method will be implemented numerically using Lagrangian methods.

In the conventional shape optimization approach to solve problem \eqref{eq:inverse_problem}, various formulations can be considered. 
For instance, one can track the Dirichlet data in the $L^{2}$ sense by considering the least-squares functional
\begin{equation}\label{eq:least_squares}
	J(\Omega\setminus\overline{\omega}):=J_{LS}(\Omega\setminus\overline{\omega})=\frac{1}{2}\int_{\partial \Omega} \abs{u - f}^2 {d}{s},
\end{equation}
where $\un:=\un(\Omega\setminus\overline{\omega})$ is subject to the well-posed mixed Dirichlet-Neumann problem
\begin{equation}\label{eq:Neumann_problem}
 		-\Delta \un = 0 \ \text{in} \ \Omega\setminus\overline{\omega},\qquad\quad
	 	\dn{\un} = g \ \text{on} \ \partial \Omega,\qquad\quad
		\un = 0 \ \text{on} \ \partial \omega.	
\end{equation}
One can also track the Neumann data in least-squares sense instead of \eqref{eq:least_squares}.
However, this approach requires more regularity in the states (and adjoint states) and may be impractical in numerical experiments where high regularity of the state variable is not guaranteed.
Another approach is to consider minimizing the energy-gap cost functional
\begin{equation}\label{eq::KV}
	J_{KV}(\Omega\setminus\overline{\omega})=\frac{1}{2}\int_{\Omega\setminus\overline{\omega}}\abs{\nabla (\ud-\un)}^2 \, dx,
\end{equation}
where $\un$ solves \eqref{eq:Neumann_problem} while $\ud:=\ud(\Omega\setminus\overline{\omega})$ is the solution of the Dirichlet problem
\begin{equation}\label{eq:Dirichlet_problem}
 		-\Delta \ud = 0 \ \text{in} \ \Omega\setminus\overline{\omega},\qquad\quad
	 	\ud = f \ \text{on} \ \partial \Omega,\qquad\quad
		\ud = 0 \ \text{on} \ \partial \omega.
\end{equation}
The latter approach -- also known as Kohn-Vogelius method \cite{KohnVogelius1987} -- was first studied in \cite{RocheSokolowski1996}, and then was re-examined in \cite{EpplerHarbrecht2005}.
Here, $J_{KV}$ is positive and vanishes only if $\ud=\un$ (i.e., $\omega$ fits the exact inclusion). 
Meanwhile, another optimization reformulation of \eqref{eq:inverse_problem} via the so-called coupled complex boundary method was recently proposed in \cite{Afraites2022}.
The idea of the method is to transform the overdetermined problem into a complex boundary value problem.
This involves introducing a complex Robin boundary condition that couples the Dirichlet and Neumann boundary conditions on the unknown (free) boundary.
Subsequently, the goal is to optimize the cost function, which is constructed using the imaginary part of the solution across the entire domain; see \cite{Afraites2022} for more details. 
In these investigations, the shape derivatives of the costs were computed and then used in an algorithm to numerically resolve the minimization problems.
The issue of \textit{ill-posedness} of the shape problems was also addressed in the said papers.
A similar study, but for the case of perfectly insulating material, can be found in \cite{AfraitesDambrineEpplerKateb2007}.
In this work, to illustrate the proposed ADMM in the context of a shape inverse problem via shape optimization settings, we will focus on the least-squares approach using \eqref{eq:least_squares}.

The rest of the note is organized as follows. In Section \ref{sec:main_contribution}, we will present the main contribution of the study by discussing how ADMM is applied to solve a new shape optimization formulation of equation \eqref{eq:inverse_problem} with equality or inequality constraints. 
The section begins with the main optimization problem, followed by the formulation of the ADMM algorithm. It then proceeds to discuss the subproblems within the main iteration procedure. 
In Section \ref{sec:numerics}, we provide numerical examples to illustrate the feasibility of the method, highlighting how it outperforms classical shape optimization methods in detecting unknown inclusions with pronounced concavities. 
The numerical experiments are carried out in two-dimensional (2D) and three-dimensional (3D) settings and under noisy data. 
Finally, Section \ref{sec:conclusion} includes a short conclusion and a statement of future work.
\section{Main Contribution}\label{sec:main_contribution}
As alluded to in the Introduction, our main intent is to propose a novel application of ADMM, also known as split Bregman \cite{GoldsteinOsher2009}, to shape identification problems, using \eqref{eq:inverse_problem} as a toy problem. 
ADMM was developed by Glowinski and Marrocco \cite{GlowinskiMarroco1975} and Gabay and Mercier \cite{GabayMercier1976} in the 1970s, with roots dating back to the 1950s. 
It is well-suited for handling convex optimization problems. 
The method takes the form of a decomposition-coordination procedure, leveraging the advantages of algorithms such as dual decomposition, the method of multipliers, Douglas--Rachford splitting, Dykstra's alternating projections, Bregman iterative algorithms for $l_{1}$ problems, proximal methods, and augmented Lagrangian methods \cite{Boydetal2011, GlowinskiOsherYin2016}.
Despite being introduced almost half a century ago, ADMM has gained popularity in recent decades, largely attributed to its efficient applications in various areas of modern technology, including computer vision, image processing, statistical learning, and more. 
For specific applications of ADMM to PDE-constrained optimal control problems, we direct readers, for example, to \cite{NeitzelTroltzsch2012, GlowinskiSongYuan2020,SratiOulmelkAfraitesHadri2023}, and the references therein.
\subsection{Proposed approach and the ADMM algorithm}
In this section, we will demonstrate how ADMM is adapted to our present problem. 
To start, we reformulate our original shape inverse problem \eqref{eq:inverse_problem} into the following shape optimization problem with equality or inequality constraints.
\begin{problem}\label{prob:optimal_shape_problem}
Let $a$ and $b$, $b\geqslant a$, be given fixed constants.
Find the shape $\omega^{\ast}$ in the space of admissible set
\[
	\mathcal{O}_{ad}=\left\{ \omega \mid \text{$\omega \in \mathcal{O}_{\delta}$ and $a\leqslant \uo\leqslant b$ a.e. in $\Omega$ where $\uo$ solves problem \eqref{eq:Neumann_problem}} \right\}
\]
such that
	\begin{equation}\label{eq:control} 
\omega^{\ast} 
	= \argmin_{\omega \in \mathcal{O}_{ad}} J(\Omega\setminus\overline{\omega})
	:= \argmin_{\omega \in \mathcal{O}_{ad}} \left\{ \frac{1}{2} \int_{\partial \Omega} \vert \uo-f\vert^2 \, dx \right\}.
	\end{equation}	
\end{problem}
The proof of the existence of the optimal shape solution to the above shape problem can be addressed rigorously using, for example, the ideas developed in \cite{HKKP2004a, HKKP2003,HaslingerMakinen2003,HenrotPierre2018}, or with the tools furnished in, for instance, \cite{BoulkhemairChakib2007, Boulkhemairetal2008, Boulkhemairetal2013}. 
Note that, since \eqref{eq:Neumann_problem} is uniquely solvable in $H_{{\partial}{\omega},0}^{1}({\Omega}\setminus\overline{\omega}) := \{ \varphi \in H^{1}(\Omega\setminus\overline{\omega}) \mid \varphi = 0 \ \text{on ${\partial}{\omega}$}\}$, one can define the map $\omega \mapsto u_{\omega}$, the graph of which is given by
\[
	\mathcal{F}=\{ (\omega,u) \mid \text{$\omega \in \mathcal{O}_{ad}$, $u_{\omega}$ solves \eqref{eq:Neumann_problem} on $\Omega\setminus\overline{\omega}$} \}.
\]
Problem \ref{prob:optimal_shape_problem} is then equivalent to minimizing $J(\omega, u_{\omega})=J(\Omega\setminus\overline{\omega})$ on $\mathcal{F}$. 
To prove the existence of a solution to this minimization problem, it is necessary to endow the set $\mathcal{F}$ with a topology for which it is compact and then demonstrate that $J$ is lower semi-continuous; for further details, refer to, e.g., \cite[Chap. 4]{HenrotPierre2018}.
While we omit the detailed proof here, a similar approach to the one used in \cite{RabagoAzegami2019b} can be applied, at least for the case of two dimensions. 
For the three-dimensional case, a more general approach using the concept of convergence of sets in the sense of Hausdorff \cite{Holzleitner2001}, combined with the (uniform) cone property (see \cite{Chenais1975}), can be employed. 
Refer also to \cite[Chap. 4]{HenrotPierre2018}.

To solve the above state-constrained shape optimal control problem with an equality or inequality constraint, we will apply ADMM. 
The method allows us to divide the global problem into a series of easily solvable subproblems. 
For the said purpose, we introduce an auxiliary variable $v$ satisfying $v=\uo$ a.e. in $\Omega\setminus\overline{\omega}$.
Then, problem \eqref{eq:control} can be written as follows  
\begin{equation}\label{eq:control_Uad}
	(\omega^{\ast},v^{\ast}) = \argmin_{(\omega,v)\in \mathcal{E}} \left\{J(\Omega\setminus\overline{\omega})+U_{\mathcal{K}}(v)\right\},
\end{equation} 
where the set $\mathcal{K}$ is the closed convex non-empty set of $L^{2}(\Omega\setminus\overline{\omega})$ defined by
\[
	\mathcal{K}=\left\{ v \in L^{2}(\Omega\setminus\overline{\omega}) \mid a\leqslant v\leqslant b \ \ \text{a.e. in $\Omega\setminus\overline{\omega}$} \right\},
\]
while $U_{\mathcal{K}}$ is the indicator functional of the set $\mathcal{K}$; that is, $U_{\mathcal{K}}(v) = 0$ if $v \in K$, otherwise, $U_{\mathcal{K}}(v) = \infty$ if $ v \in L^{2}(\Omega\setminus\overline{\omega})\setminus \mathcal{K}$.
Meanwhile, $\mathcal{E}$ is defined as follows
\[
	\mathcal{E}=\left\{ (\omega,v)\in \mathcal{O}_{ad} \times L^{2}(\Omega\setminus\overline{\omega}) \mid \text{$\uo=v$ a.e. in $\Omega\setminus\overline{\omega}$} \right\}.
\]
To apply ADMM to the control model $\eqref{eq:control_Uad},$ we need to define the augmented Lagrangian functional first.
This is possible since the minimum of problem $\eqref{eq:control_Uad}$ is the saddle point of the following Lagrangian functional 
\begin{equation}\label{eq:augmented_lagrangian}
	\mathcal{L}_{\beta}(\omega,v;\lambda)=J(\Omega\setminus\overline{\omega})+U_{\mathcal{K}}(v)+\frac{\beta}{2}  \int_{\Omega\setminus\overline{\omega}} \vert \uo-v\vert^2 \, dx  +   \int_{\Omega\setminus\overline{\omega}} \lambda  (\uo-v ) \, dx,
 \end{equation}
where $\lambda$ is the Lagrange multiplier and $\beta>0$ is a penalty parameter.

Now, to find a saddle point of the Lagrangian functional $\mathcal{L}$, we will implement an approximation procedure based on ADMM. 
Specifically, starting with initial values $v^{0}, \lambda^{0} \in L^{2}(\Omega\setminus\overline{\omega})$, we will iteratively compute the optimizer of $\mathcal{L}$ for $k = 1, 2, \ldots$ by solving the following sequence of minimization problems:
\begin{align}
	\omega^{k+1}	&= \argmin_{\omega \in \mathcal{O}_{ad}}\mathcal{L}_{\beta}(\omega,v^{k};\lambda^{k}); \label{eq:controle}\tag{SP1} \\
	v^{k+1}		&= \argmin_{v \in L^{2}(\Omega\setminus\overline{\omega})}\mathcal{L}_{\beta}(\omega^{k+1},v;\lambda^{k}); \label{eq:etat}\tag{SP2}\\
	\lambda^{k+1}	&= \lambda^{k}+\beta (u_{\omega^{k+1}}-v^{k+1} ). \label{eq:parametre1}\tag{SP3} 
\end{align}
For the sake of technical simplicity and to streamline certain arguments, we assume that both $\Omega$ and $\omega$ are smooth domains with $\mathcal{C}^{2,1}$ regularity. Furthermore, in this work, we consider a fixed value for the penalty parameter $\beta > 0$ to simplify our discussion.
While it is possible to develop an optimization scheme for $\beta$ within our main algorithm using bilevel optimization \cite{Dempe2020} (see, e.g., \cite{OulmelkAfraitesHadriNachaoui2022}), we have opted to keep $\beta$ fixed. This choice consistently yields good results, as demonstrated in the numerical section of the paper.
Consequently, utilizing the augmented Lagrangian functional given in \eqref{eq:augmented_lagrangian}, we outline the ADMM scheme in Algorithm \ref{algo:ADMM_algorithm}.
\begin{algorithm}
\begin{enumerate}\itemsep0.3em 
	\item \textit{Input} Fix $\beta$, $a$, and $b$, and define the Cauchy pair $(f, g) \in L^{2}({\partial}{\Omega})^{2}$.
	\item \textit{Initialization} Set the initial values $ v^{0},\lambda^{0} \in L^{2}(\Omega\setminus\overline{\omega}) $.
	\item \textit{Iteration} For $k=1,2, \ldots$, compute $\{v^{k},\lambda^{k}\}$ via \eqref{eq:controle}--\eqref{eq:parametre1} by doing the sequence of computations
	\[
		\{v^{k},\lambda^{k}\} \ \stackrel{\eqref{eq:controle}}{\longrightarrow} \ \omega^{k+1} 
		\ \stackrel{\eqref{eq:etat}}{\longrightarrow} \ v^{k+1}
		\ \stackrel{\eqref{eq:parametre1}}{\longrightarrow} \ \lambda^{k+1}.\vspace{-3pt}
	\]
	\item \textit{Stop Test} Repeat \textit{Iteration} until convergence.
\end{enumerate}
\caption{ADMM algorithm for the solution of problem \eqref{eq:control}.}
\label{algo:ADMM_algorithm}  
\end{algorithm}

In the next two subsections, we will decouple the subproblems \eqref{eq:controle} and \eqref{eq:etat} by solving each of these minimization problems separately.
\subsection{Solution of $\omega$-subproblem \eqref{eq:controle}}
We first look for the solution of the first $\omega$-subproblem \eqref{eq:controle} of Algorithm \ref{algo:ADMM_algorithm} where we minimize the augmented Lagrangian functional $\mathcal{L}_{\beta}$ with respect to $\omega$.
The $\omega$-subproblem \eqref{eq:controle} is given as follows
\[
 \omega^{k+1}=\argmin_{\omega \in \mathcal{O}_{ad}} \left\{ J(\Omega\setminus\overline{\omega})+U_{\mathcal{K}}(v^{k})+\frac{\beta}{2} \int_{\Omega\setminus\overline{\omega}} \vert \uo-v^{k}\vert^2 \, dx +   \int_{\Omega\setminus\overline{\omega}} \lambda^{k}  (\uo-v^{k}) \, dx\right\}.
\]
Let us consider the following shape functional
\[
\mathcal{G}^{k}(\Omega\setminus\overline{\omega})=\mathcal{L}_{\beta}(\omega,v^{k};\lambda^{k})=J(\Omega\setminus\overline{\omega})+\frac{\beta}{2} \int_{\Omega\setminus\overline{\omega}} \vert \uo-v^{k}\vert^2 \, dx  +   \int_{\Omega\setminus\overline{\omega}} \lambda^{k}  (\uo-v^{k} ) \, dx,
\]
where $J(\Omega\setminus\overline{\omega})=\displaystyle\frac{1}{2}\int_{\partial \Omega}  \vert \uo-f\vert^2 \, dx$ and $\uo$ solves problem \eqref{eq:Neumann_problem} associated with $\omega$.

The resolution of the $\omega$-subproblem \eqref{eq:controle} requires the shape derivative of $\mathcal{G}^{k}$.
In this regard, we let $\VV$ be a sufficiently smooth vector\footnote{At least $\mathcal{C}^{2,1}$ smooth for our argumentation, but a $\mathcal{C}^{1,1}$ regularity assumption is sufficient.} field on $\mathbb{R}^{d}$ with compact support in ${\Omega}_{\delta}$.
For ease of writing, we write $\Vn = \VV \cdot \vect{n}$ and we denote the set of admissible vector field $\VV$ by $\vect{\Theta}$.
Without further notice, we always assume that $\VV \in \vect{\Theta}$.

To proceed, let us first define the shape derivative of the shape functional $\mathcal{G}$ according to \cite[Sec. 4.3.2, Eq. (3.6), p. 172]{DelfourZolesio2011} in the following definition.
\begin{definition}\label{def:shape_derivative}
Let $\mathcal{O}_{ad}$ denote the set of admissible domains $\Omega\setminus\overline{\omega}$.
The functional $\mathcal{G} : \mathcal{O}_{ad} \to \mathbb{R}$ has a directional \textit{first-order} \textit{Eulerian derivative} at $\Omega\setminus\overline{\omega}$ in the direction of a given deformation field $\VV\in \vect{\Theta}$ if the limit
\begin{equation}
\label{eq:limit_J}
	\lim_{t \searrow0} \frac{{\mathcal{G}(\Omega\setminus\overline{\omega}_{t})-\mathcal{G}(\Omega\setminus\overline{\omega})}}{t} =: {d}\mathcal{G}(\Omega\setminus\overline{\omega})[\VV]
\end{equation}
exists.
The shape functional $\mathcal{G}$ is said to be \textit{shape differentiable} at $\Omega\setminus\overline{\omega}\in\mathcal{O}_{ad}$ in the direction of $\VV\in \vect{\Theta}$ if the map $\VV \mapsto {d}\mathcal{G}(\Omega\setminus\overline{\omega})[\VV]$ is linear and continuous.
In this case, we refer to this map as the \textit{shape gradient} of $\mathcal{G}$.
\end{definition}
\begin{remark}
In the case of the classical shape optimization formulation \eqref{eq:least_squares}, the set of admissible domains $\mathcal{O}_{ad}$ can be defined simply as the set $\mathcal{O}_{\delta}$.
Here, the set of admissible domains $\mathcal{O}_{ad}$ in Problem \ref{prob:optimal_shape_problem} incorporates the additional inequality constraint for the development of the proposed ADMM.
\end{remark}

Let $D_{\omega} \mathcal{G}^{k}(\Omega\setminus\overline{\omega})$ be the first-order shape derivative of $\mathcal{G}^{k}(\cdot)$ at $\omega$ in the direction of the vector field $\VV$.
Then, formally, we have the following computations
\begin{equation}
\begin{aligned}\label{derivative:lagrangian}
	d\mathcal{G}^{k}(\Omega\setminus\overline{\omega})[\VV]
		&= \langle D_{\omega} \mathcal{G}^{k}(\Omega\setminus\overline{\omega}), \VV \rangle_{{\partial}{\omega}}
		= \int_{\partial \omega} D_{\omega} \mathcal{G}^{k}(\Omega\setminus\overline{\omega}) \cdot \VV \, dx \\
		&= \displaystyle \int_{\partial \Omega}\big(\uo-f\big) {{\un}^{\prime}} \, dx
			+ \beta  \int_{\Omega\setminus\overline{\omega}} \big(\uo-v^{k} \big) {{\un}^{\prime}} \, dx 
			+ \frac{\beta}{2} \int_{\Omega\setminus\overline{\omega}} \div{ \big(\uo-v^{k} \big)^{2} \VV} \, dx \\
		&\qquad
			+ \int_{\Omega\setminus\overline{\omega}} \lambda^{k}  {{\un}^{\prime}} \, dx
			+ \int_{\Omega\setminus\overline{\omega}} \div{ \lambda^{k} \big(\uo-v^{k} \big) \VV} \, dx \\
		&=: dJ(\Omega\setminus\overline{\omega})[\VV]
			+ \beta  \int_{\Omega\setminus\overline{\omega}} \big(\uo-v^{k} \big) {{\un}^{\prime}} \, dx 
				+ \frac{\beta}{2} \int_{\partial \omega} \big(v^{k} \big)^2 \Vn  \, ds \\
		&\qquad + \int_{\Omega\setminus\overline{\omega}} \lambda^{k}   {{\un}^{\prime}} \, dx
				- \int_{\partial \omega} \lambda^{k} v^{k} \Vn \, ds.
\end{aligned}
\end{equation}
In above computation, ${\un}^{\prime} := {\un}^{\prime}(\Omega\setminus\overline{\omega})[\VV] = \lim_{t \searrow 0} \frac{1}{t}(u(\Omega\setminus\overline{\omega}_t) - u(\Omega\setminus\overline{\omega}))$ represents the shape derivative of the state variable $u$, which solves the following well-posed PDE system (refer, for example, to \cite{EpplerHarbrecht2005}).
\begin{equation}\label{eq:shape_derivative_of_the_state}
 		-\Delta \un^{\prime} = 0 \ \text{in} \ \Omega\setminus\overline{\omega},\qquad\quad
	 	\dn{\un^{\prime}} = 0 \ \text{on} \ \partial \Omega,\qquad\quad
		\un^{\prime} = -\dn{\un}\Vn \ \text{on} \ \partial \omega.	
\end{equation}
In equation \eqref{derivative:lagrangian}, we have used the fact that $\un = 0$ on $\partial \omega$ and the assumption that $\operatorname{supp}(\VV) \subset {\Omega}_{\delta}$, and so $\VV = \vect{0}$ on $\partial \Omega$.

We emphasize here that the solvability of \eqref{eq:shape_derivative_of_the_state} in $H^{1}(\Omega\setminus\overline{\omega})$ requires additional regularity of the domain and the data. 
Specifically, the existence of a (unique) weak solution $\un^{\prime}\in H^{1}(\Omega\setminus\overline{\omega})$ corresponding to the variational formulation of \eqref{eq:shape_derivative_of_the_state} is guaranteed by the Lax-Milgram lemma, provided that $\omega$ is of class $\mathcal{C}^{2,1}$ and $g \in H^{3/2}({\partial}{\omega})$. 
It is noteworthy that, under these regularity assumptions, the state variable $u$ is not only $H^{1}_{{\partial}{\omega},0}(\Omega\setminus\overline{\omega})$-smooth but is also $H^{3}(\Omega_{\delta}\setminus\overline{\omega})$ regular, as guaranteed by classical elliptic regularity results. 
We only needed this high regularity assumption on the domain and the data since we are applying the chain rule approach to obtain the shape derivative of the functional. 
These regularity assumptions can be relaxed by applying a different method to obtain the shape derivative, for example, by employing the rearrangement method \cite{IKP2006}.

We point out that equation~\eqref{derivative:lagrangian} is difficult to handle since we cannot find explicitly the direction $\VV$.
In fact, the computed expression with the shape derivative ${\un}^{\prime}$ is not useful for practical applications, especially in the numerical realization of the proposed shape optimization problem via an iterative procedure.
This is because the implementation requires the solution of~\eqref{eq:shape_derivative_of_the_state} for each velocity field $\VV$, at every iteration.
To get around this difficulty, we apply the adjoint method and introduce the variable $p$ -- in order to eliminate from the gradient expression the shape derivative $\un^{\prime}$ -- which solves the following adjoint problem
\begin{equation}\label{eq:adjoint_problem}
     		-\Delta p = \beta \big(\uo-v^{k}\big)+\lambda^{k} \ \text{in} \ \Omega\setminus\overline{\omega},\quad
	 	\dn{p} = \uo-f \ \text{on} \ \partial \Omega,\quad
		p = 0\ \text{on} \ \partial \omega.	
\end{equation}
This leads us to the following expression for the shape derivative of $\mathcal{G}^{k}$:
\begin{equation}\label{eq:shape_derivative_of_the_Lagrangian}
	d\mathcal{G}^{k}(\Omega\setminus\overline{\omega})[\VV] 
		= \int_{\partial \omega} D_{\omega} \mathcal{G}^{k}(\Omega\setminus\overline{\omega}) \cdot \VV \, ds
		= \int_{\partial \omega} \left( \partial_n p \partial_n u + \lambda^{k} v^{k} - \frac{\beta}{2} \big( v^{k} \big)^{2} \right) \vect{n} \cdot \VV \, ds,
\end{equation}
%
%
%
%
%

In practice, the computed shape derivative $D_{\omega} \mathcal{G}(\Omega\setminus\overline{\omega})$ is not used directly in a numerical procedure since it may cause some unwanted oscillations on the boundary during the approximation process, causing some instabilities in the algorithm. 
To address this issue, we have to apply the so-called Sobolev-gradient method \cite{Neuberger1997} which we will discuss next.
\subsubsection{Classical computation of the extended-regularized deformation fields}
The shape gradient of $J$ is only supported on $\partial \omega$ and may lack enough smoothness necessary in numerical realization (particularly, when employing finite element methods or FEMs).
To improve the regularity of the descent direction $D_{\omega} \mathcal{G}^{k}(\Omega\setminus\overline{\omega})$ and extend its definition to the entirety of $\Omega\setminus\overline{\omega}$, we make use of its Riesz representation which we obtained here by solving the following system of partial differential equations:
\begin{equation}\label{eq:deformation_field}
 		-\Delta \VV+\VV = 0 \ \text{in} \ \Omega\setminus\overline{\omega},\quad
	 	 \VV = 0 \ \text{on} \ \partial \Omega,\quad
		\partial_n\VV = -D_{\omega} \mathcal{G}^{k}(\Omega\setminus\overline{\omega}) \ \text{on} \ \partial \omega.	
\end{equation}
Accordingly, we can formulate a Sobolev gradient-based descent (SGD) algorithm laid out in Algorithm \ref{algo:SGD_algorithm} to solve our problem. 
\begin{algorithm} 
\caption{SGD algorithm for $\omega$-subproblem \eqref{eq:controle}}\label{algo:SGD_algorithm} 
\begin{enumerate}
	\item \textit{Input} Fix $\beta$, $a$,  $b$, and {$\varepsilon$} and set $\lambda^{k}$, $\mu_{m}$, $\omega^{k}_{m}=\omega^{k}$,  $u^{k}_{m}=u^{k}$, $v^{k}_{m}=v^{k}$. Also, set $m=0$. 	
	\item \textit{Iteration} For $m = 1, 2, \ldots$,
	\begin{enumerate}\itemsep0.3em 
		\item[2.1] solve the state problem \eqref{eq:Neumann_problem};
		\item[2.2] solve the  adjoint problem \eqref{eq:adjoint_problem};
		\item[2.2] compute the descent direction $\VV_{m}^{k}$ via equation \eqref{eq:deformation_field};
		\item[2.3] update the current boundary $\partial \omega_{m}^{k}$ by $\VV_{m}^{k}$ to obtain ${\partial \omega_{m+1}^k}$; i.e., for some small scalar $t^{k}>0$, set 
		\[
			\partial \omega_{m+1}^{k}:=\left\{x+t^{k} \VV_{m}^{k}(x) \mid x\in \partial \omega^{k}_{m}\right\}.
		\]
	\end{enumerate}
	\item \textit{Stop test} Repeat \textit{Iteration} until convergence; that is, \textbf{while} $\Vert  d\mathcal{G}^{k}(\omega_{m}^{k})[\VV^{k}_{m}] \Vert \geqslant \varepsilon$ \textbf{do} \textit{Iteration} 
	\item \textit{Output} $\omega^{k+1}=\omega_{m+1}^{k}$.
\end{enumerate}
\end{algorithm}
%
%
\begin{remark}\label{rem:step_size_conditions}
In Step 2.3 of Algorithm \ref{algo:SGD_algorithm}, the step size $t^{k}$ is initialized using the formula $t^{0} = \mu J^{0}/\|\VV^{0}\|_{H^{1}(\Omega\setminus\overline{\omega}^{0})^{d}}$ with $\mu=0.5$. 
We continue to use this step size in the succeeding iterations but further adjust it to avoid inverted triangles (or tetrahedrons) within the mesh after each update.
Note that a backtracking procedure (with the initial value for the step size $t^{k} = \mu J^{k}/\|\VV^{k}\|_{H^{1}(\Omega\setminus\overline{\omega}^{k})^{d}}$, where $\mu > 0$ is sufficiently small, at each iteration) based on a line search method for shape optimization similar to \cite[p. 281]{RabagoAzegami2020}, could also be employed. 
However, the previously mentioned choice of descent step size is more effective in providing a reconstruction of the unknown cavity. 
Furthermore, based on our experience, the cost function $J$ is insensitive to large deformations.
We leave the improvement of the choice of the step size for further research.
\end{remark}
\begin{remark}
The algorithm above is tailored for the finite element scheme. 
If one chooses to use the boundary element method instead, the descent direction choice of $-D_{\omega} \mathcal{G}^{k}(\Omega\setminus\overline{\omega})$ can be directly incorporated into a similar algorithm.
\end{remark}
%
%
\subsection{Solution of the $v$-subproblem \eqref{eq:etat}}
Now we turn our attention to the resolution of $v$-subproblem \eqref{eq:etat} by minimizing the augmented Lagrangian functional $\mathcal{L}_{\beta}$ with respect to $v$.
That is, we solve the $v$-subproblem $\eqref{eq:etat}$ given by 
\[
\begin{aligned}
  v^{k+1}&=\argmin_{v \in L^{2}(\Omega\setminus\overline{\omega})}\Big\{  J(\omega^{k+1})+U_{\mathcal{K}}(v)+\frac{\beta}{2}  \int_{\Omega\setminus\overline{\omega}} \vert u_{\omega^{k+1}}-v \vert^2 \, dx +  \int_{\Omega\setminus\overline{\omega}} \lambda^{k}  ( u_{\omega^{k+1}}-v ) \, dx \Big\}.
  \end{aligned}
\]
Applying the projection method, we obtain the equation $v^{k+1}=P_{\mathcal{K}}\big(u_{\omega^{k+1}}+ \lambda^{k}/\beta \big)$,
where $P_{\mathcal{K}}(w) := \max(a, \min(b, w))$, for all $w \in L^{2}(\Omega \setminus \overline{\omega})$ is the projection operator onto the admissible set $\mathcal{K}$.
\subsection{ADMM-SGD algorithm}
Finally, based on the discussions above, we can now propose a modification of Algorithm \ref{algo:ADMM_algorithm} for the numerical solution of the constrained shape optimal control problem \eqref{eq:control} with an equality or inequality constraints subject to \eqref{eq:Neumann_problem}.
More precisely, Algorithm \ref{algo:ADMM_algorithm} can be specified as a nested iterative ADMM-SGD scheme for the optimal control problem \eqref{eq:control_Uad} following the instructions given in Algorithm \ref{algo:ADMM-SGD}. 

\begin{algorithm}[!htp] \caption{ADMM-SGD}\label{algo:ADMM-SGD}
\begin{enumerate} \itemsep0.1em 
    \item \textit{Initialization} Specify the input data $g$, and choose $\omega^{0}$, $\lambda^{0}$, $\beta$, $a$, $b$, $v^{0}$, and $\varepsilon$.
    \item \textit{Iteration} For $k=0,\ldots,N$,
    \begin{enumerate}\itemsep0.1em
   	 \item[2.1] compute $u_{\omega^{k}}$  solution of the state \eqref{eq:Neumann_problem} associated to $\omega^{k}$;
   	 \item[2.2] compute $p^{k}$ solution of the adjoint state \eqref{eq:adjoint_problem};
   	 \item[2.3] update $\omega^{k+1}$ by the gradient-descent method in Algorithm \ref{algo:SGD_algorithm};
   	 \item[2.4] update $v^{k+1}$ as $\displaystyle v^{k+1}=\max \left( a, \min \left(u_{\omega^{k+1}}+\lambda^{k}/\beta, b \right) \right)$;
   	 \item[2.5] set $\lambda^{k+1}=\lambda^{k}+\beta (u_{\omega^{k+1}}-v^{k+1} )$.
	\end{enumerate}
	\item \textit{Stop test} Repeat \textit{Iteration} until convergence.
\end{enumerate}
\end{algorithm}

\begin{remark}
The methods and algorithms presented above can be easily modified for the case of noisy data. 
If one intends to add a regularization term, whether or not the data is contaminated by noise (for instance, via perimeter or volume regularization), these terms will be incorporated into the Lagrangian functional. 
This addition results in the inclusion of their respective shape derivatives in \eqref{derivative:lagrangian}.
\end{remark}
%
%
%
\section{Numerical Implementation and Examples}\label{sec:numerics}
We now illustrate the feasibility of the proposed scheme and its advantages over the classical optimization approach `$J(\Omega \setminus \overline{\omega}) \to \inf$.' 
For this purpose, we carry out experiments not only with the case of exact data but also with noisy data.  
In the case of noise-contaminated data, we will employ a regularization method using the area (or volume) functional for the 2D cases and the perimeter (or more appropriately, the surface area) functional for the 3D cases in the detection process. 
In all test cases, the specimen's shape is that of the unit ball $B(\vect{0},1)$ — with a unit radius centered at the origin.
Moreover, the prescribed flux is set to $ g = 1 $, and the additional data $ f $ on $ {\partial}{\omega} $ are obtained by numerically solving the forward problem \eqref{eq:Neumann_problem} using very fine meshes and $ P_{2} $ finite element basis functions. 
To avoid `inverse crimes' (see \cite[p. 179]{ColtonKress2019}), we use coarser meshes and $ P_{1} $ finite elements in the inversion process. 
Meanwhile, we stop the algorithm as soon as it reaches a maximum number of iterations $N$ (which means precisely the maximum number of successful mesh deformations).
This, of course, pertains to the concept of `convergence' as mentioned in the \textit{Stop test} step of Algorithm \ref{algo:ADMM-SGD}.
Obviously, this criterion can be modified and even improved, but this simple one already permits us to obtain effective results.
The computations are performed on a MacBook Pro with an Apple M1 chip and 16GB RAM main memory, via {\sc FreeFem++} \cite{Hecht2012}.

\medskip
\textbf{2D case.} 
For test cases in two spatial dimensions, we consider a flower-like shape, a peanut-shaped, and an \textsf{L}-shaped exact cavity. 
The specifications in the algorithms are as follows: $N = 300$, $ \lambda^{0} = 0.001 $, $ a = 0.5\min u(\Omega \setminus \overline{\omega}^{\ast}) $, $ b = 1.5\max u(\Omega \setminus \overline{\omega}^{\ast}) $, $ v^{0} = 1 $, $ \varepsilon = 10^{-6} $, and $ \omega^{0} = B(\mathbf{0}, 0.8) $.

\begin{figure}[htp!]
\centering
\resizebox{0.325\linewidth}{!}{\includegraphics{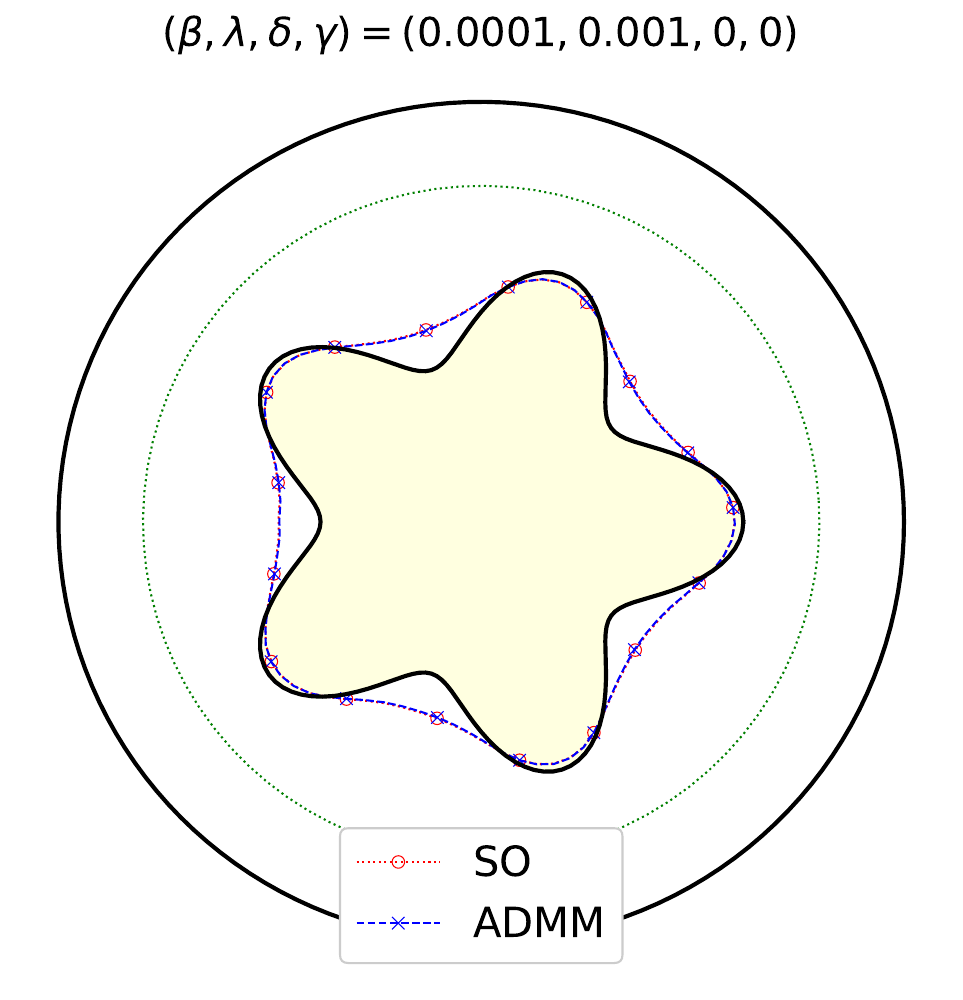}} 
\resizebox{0.325\linewidth}{!}{\includegraphics{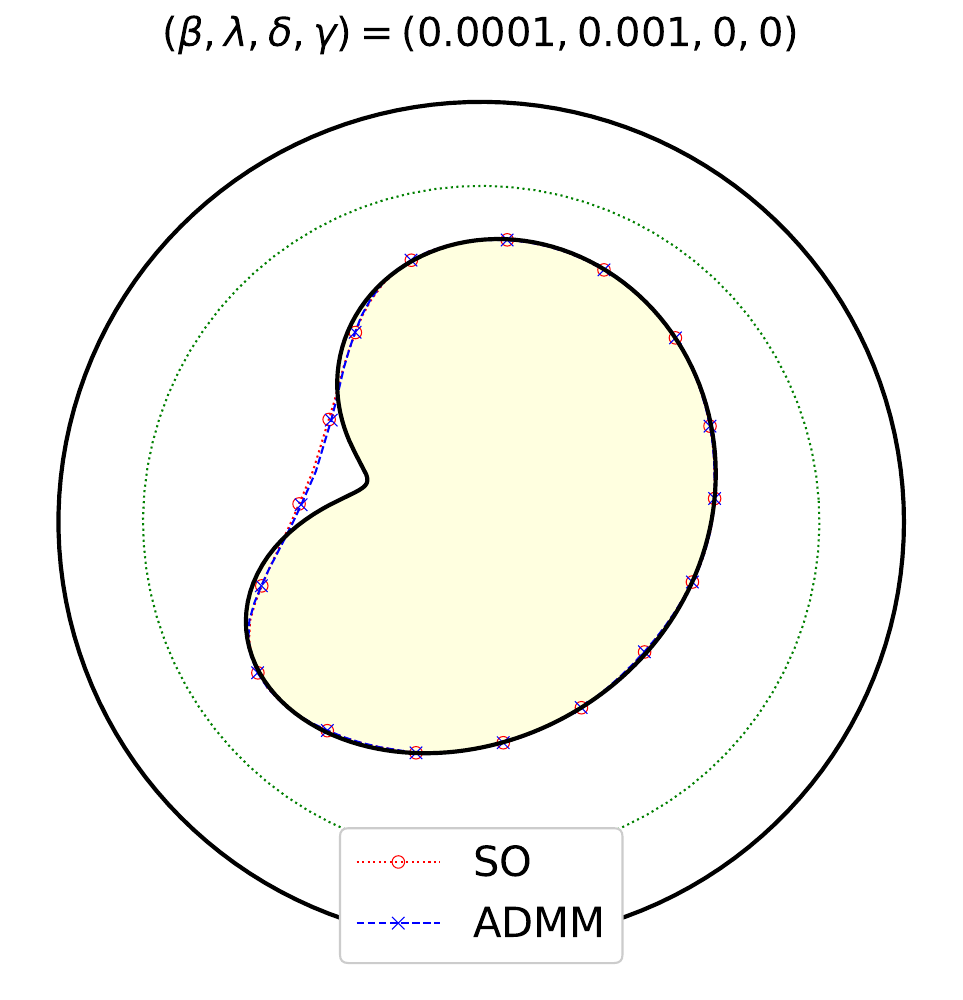}}  
\resizebox{0.325\linewidth}{!}{\includegraphics{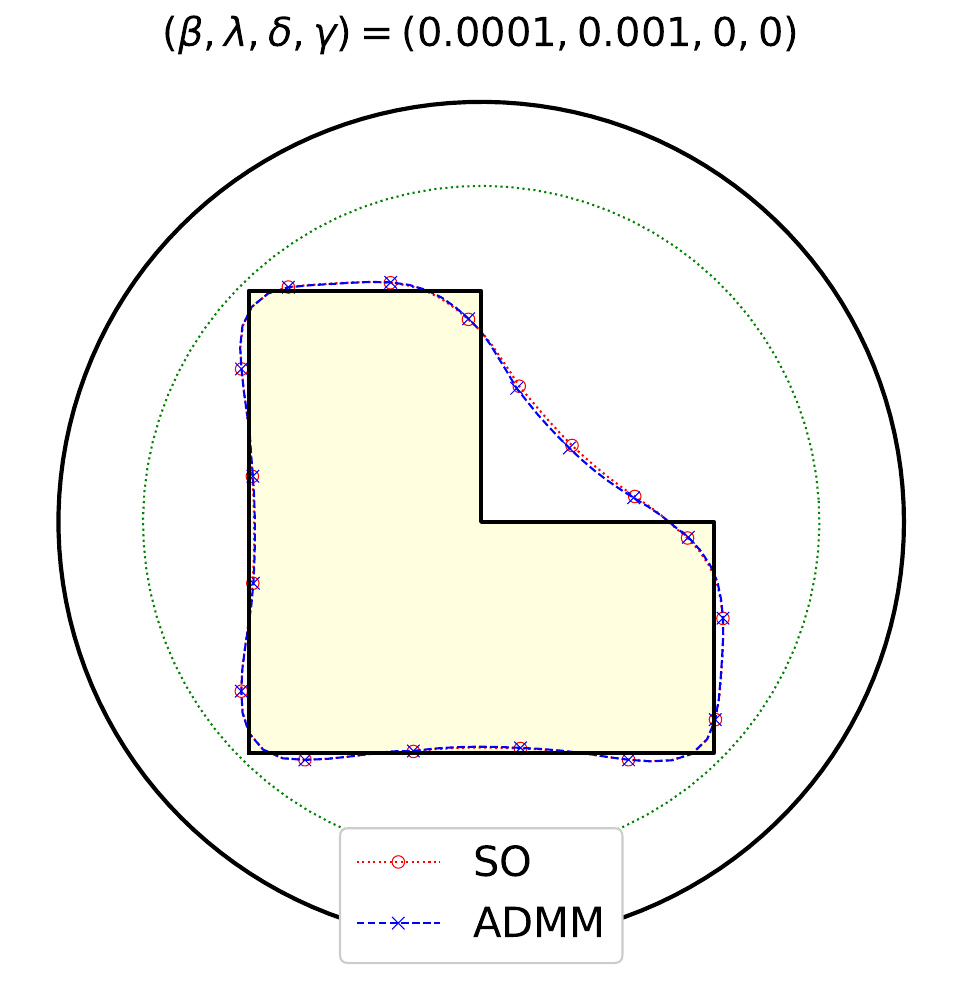}} 
\\[0.7em] 
\resizebox{0.325\linewidth}{!}{\includegraphics{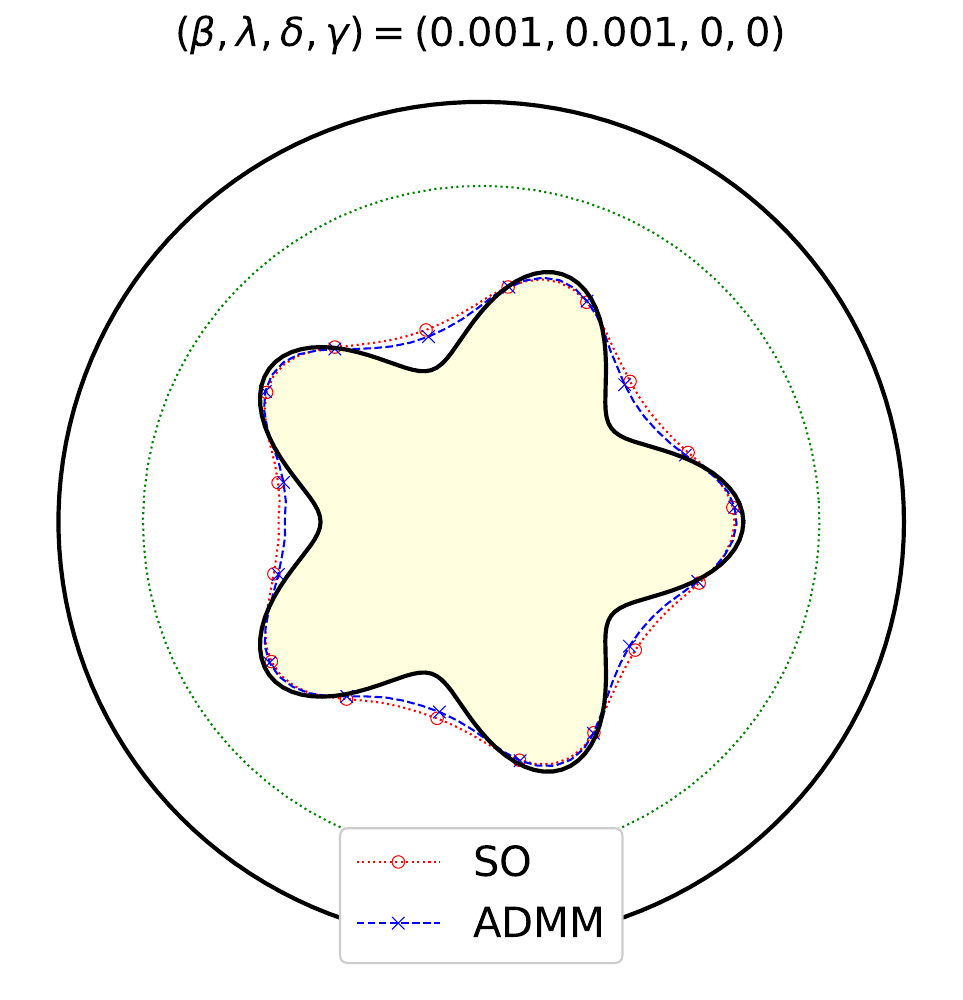}} 
\resizebox{0.325\linewidth}{!}{\includegraphics{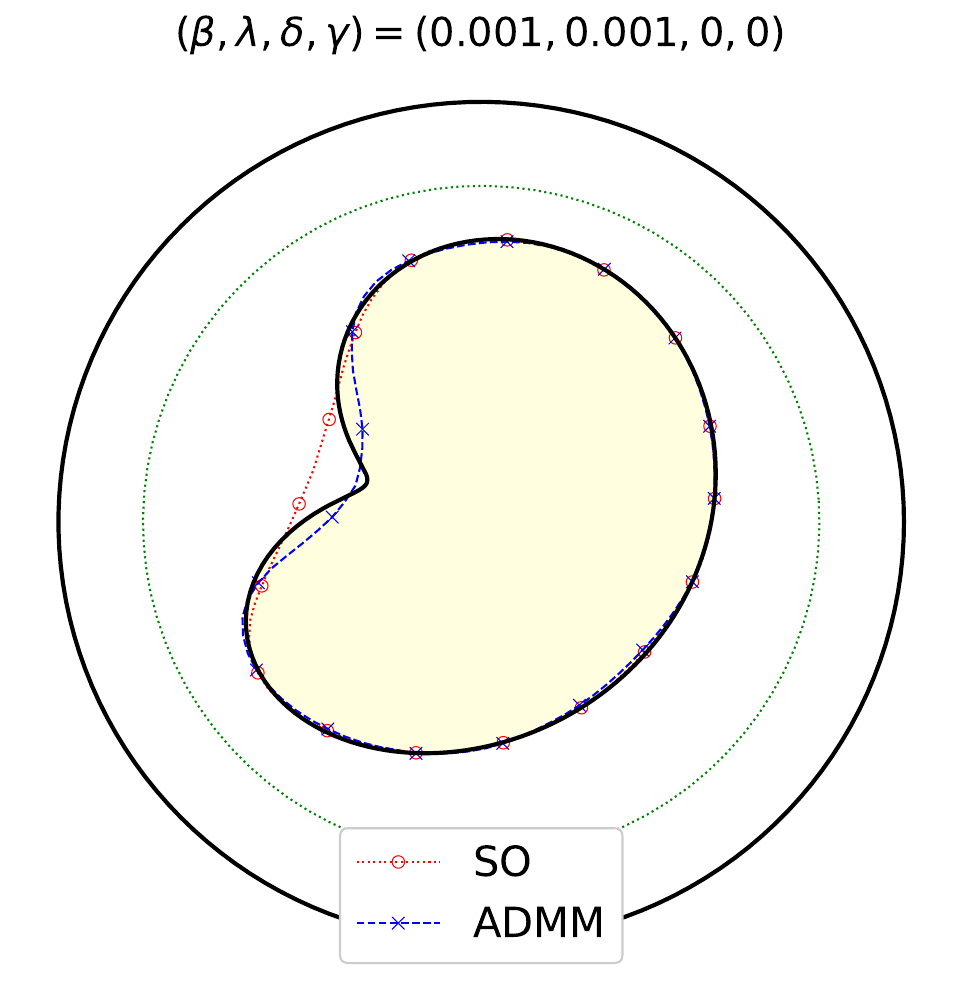}}  
\resizebox{0.325\linewidth}{!}{\includegraphics{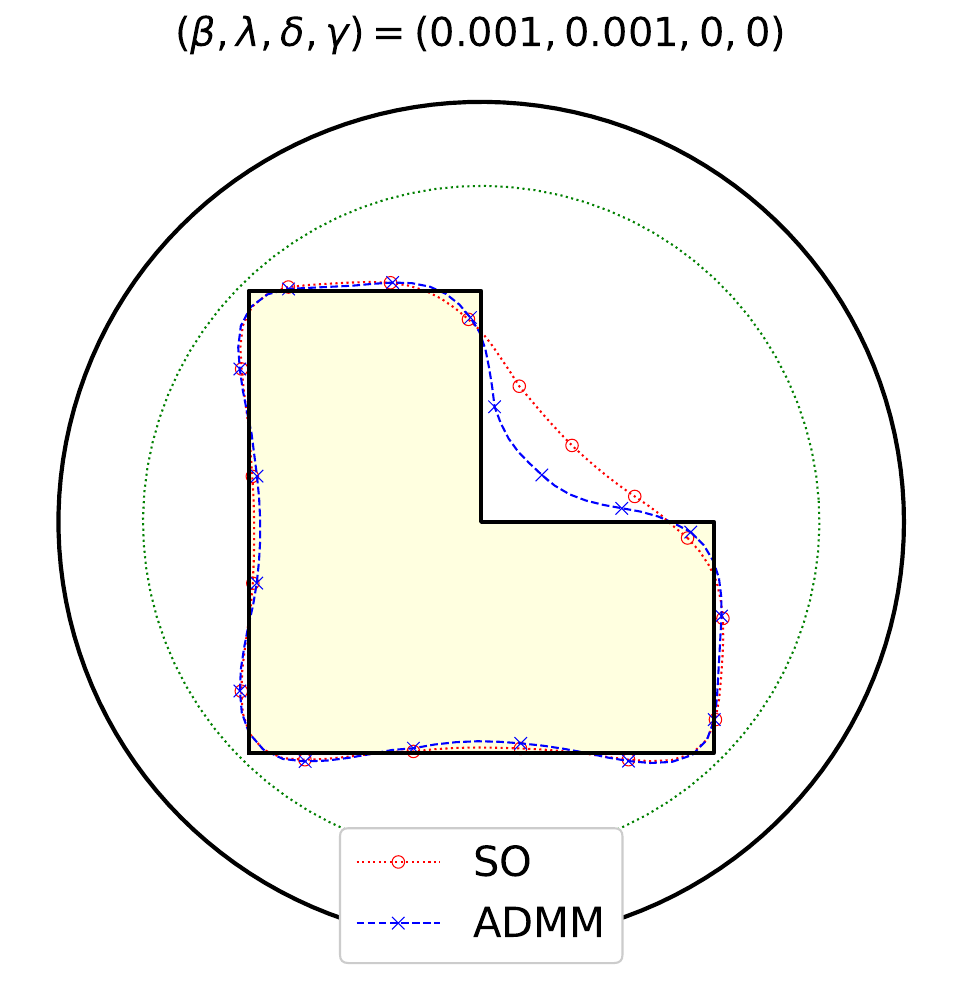}} 
\\[0.7em] 
\resizebox{0.325\linewidth}{!}{\includegraphics{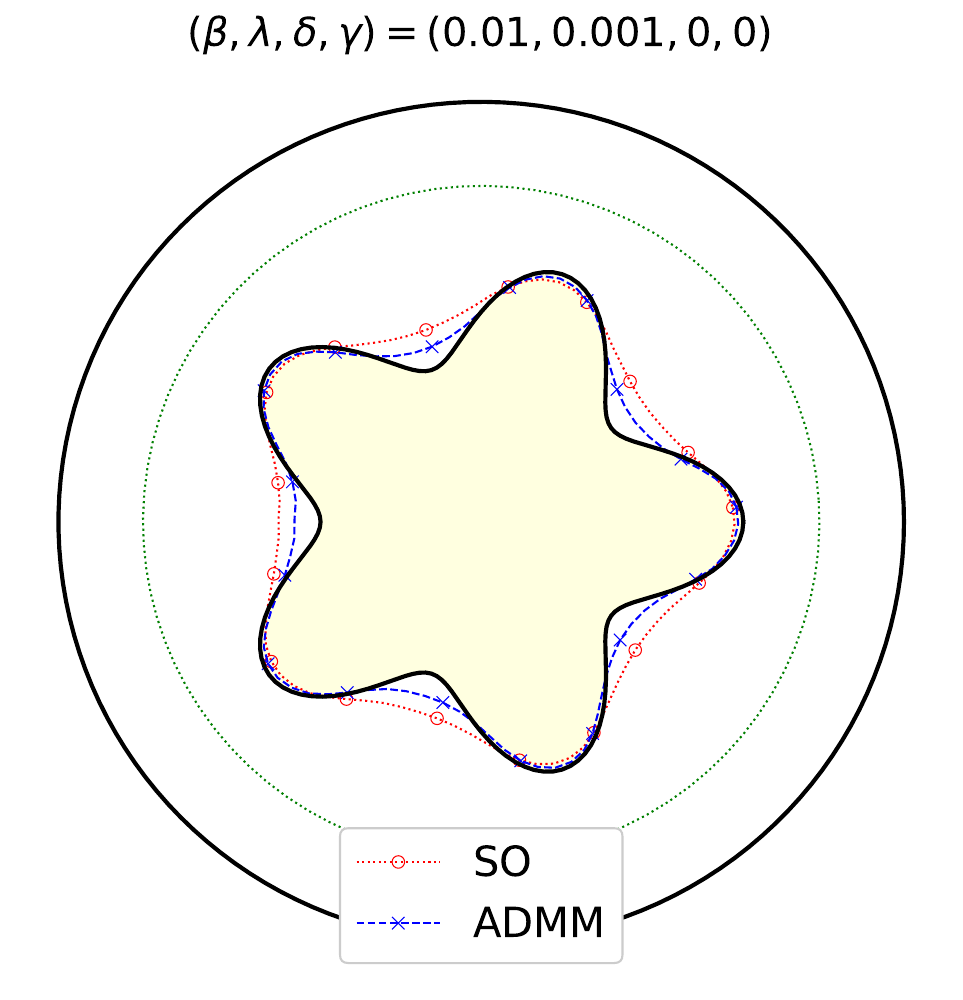}} 
\resizebox{0.325\linewidth}{!}{\includegraphics{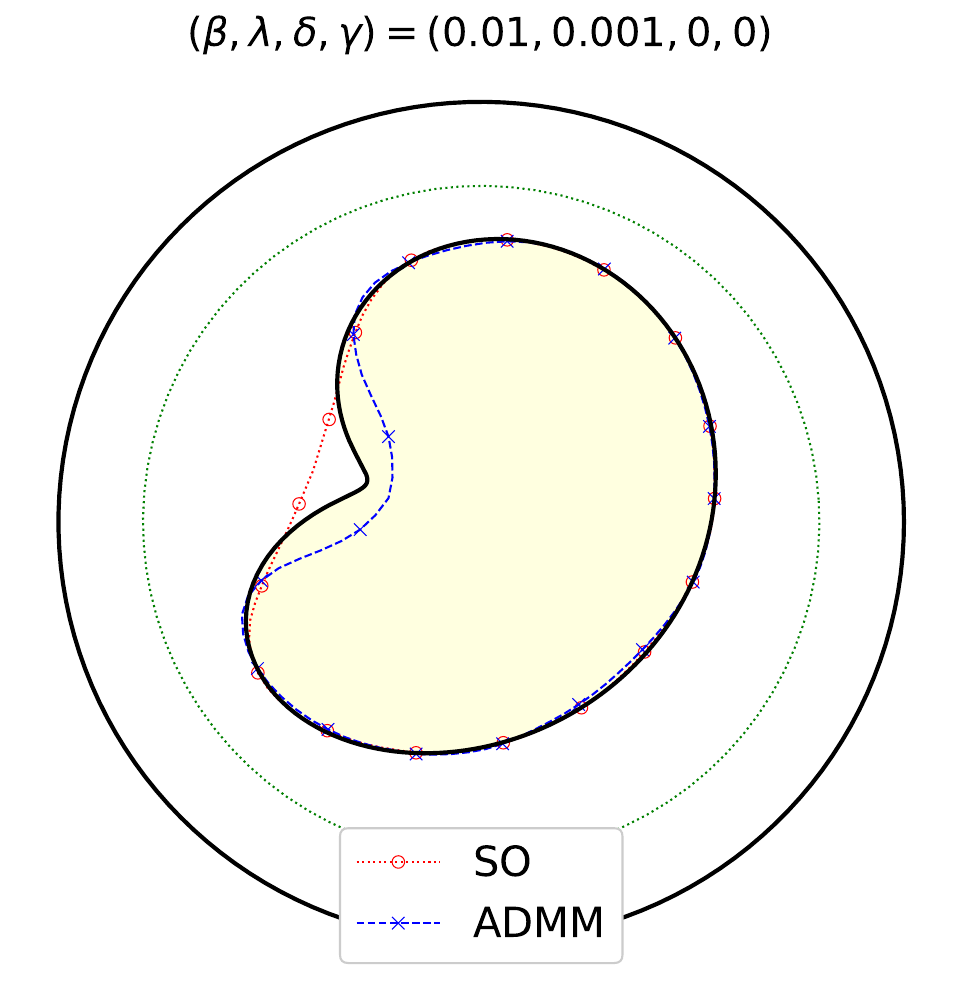}}  
\resizebox{0.325\linewidth}{!}{\includegraphics{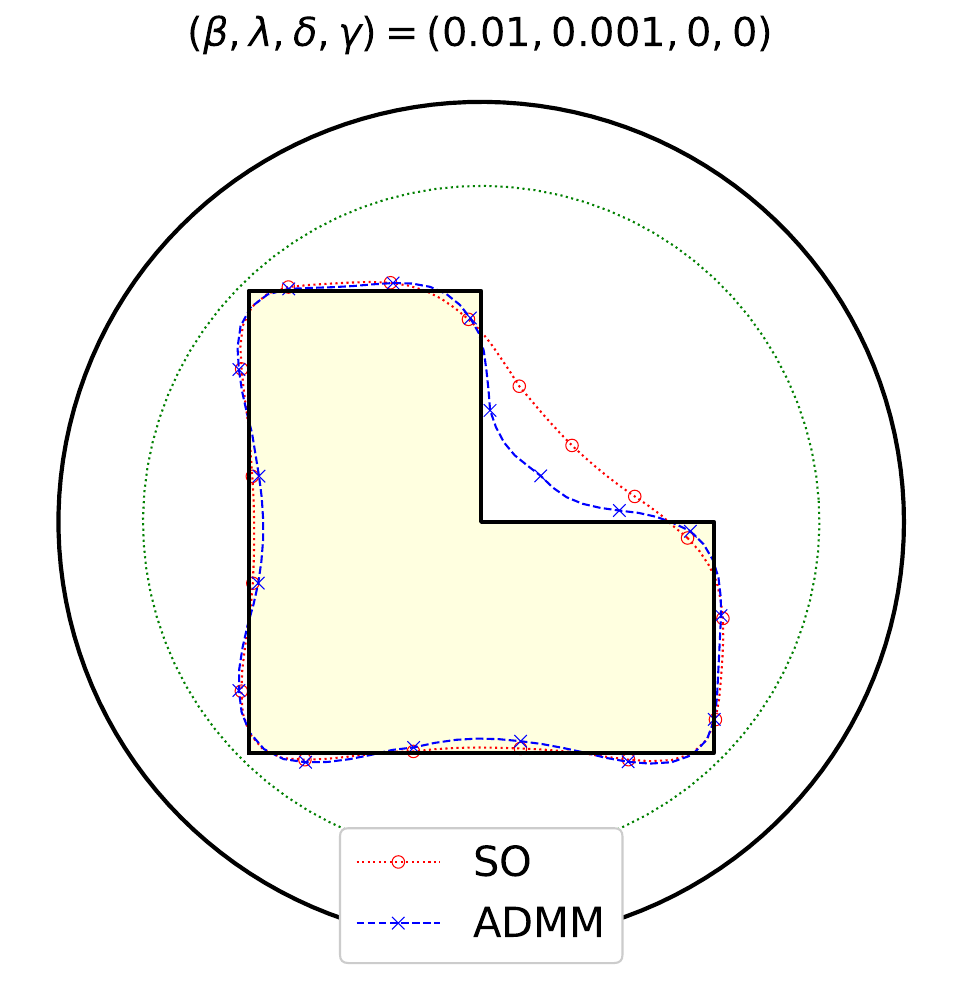}} 
\caption{Reconstructions with exact data at varying values of $\beta$}
\label{fig:figure1a}
\end{figure}

\begin{figure}[htp!]
\centering 
\resizebox{0.24\linewidth}{!}{\includegraphics{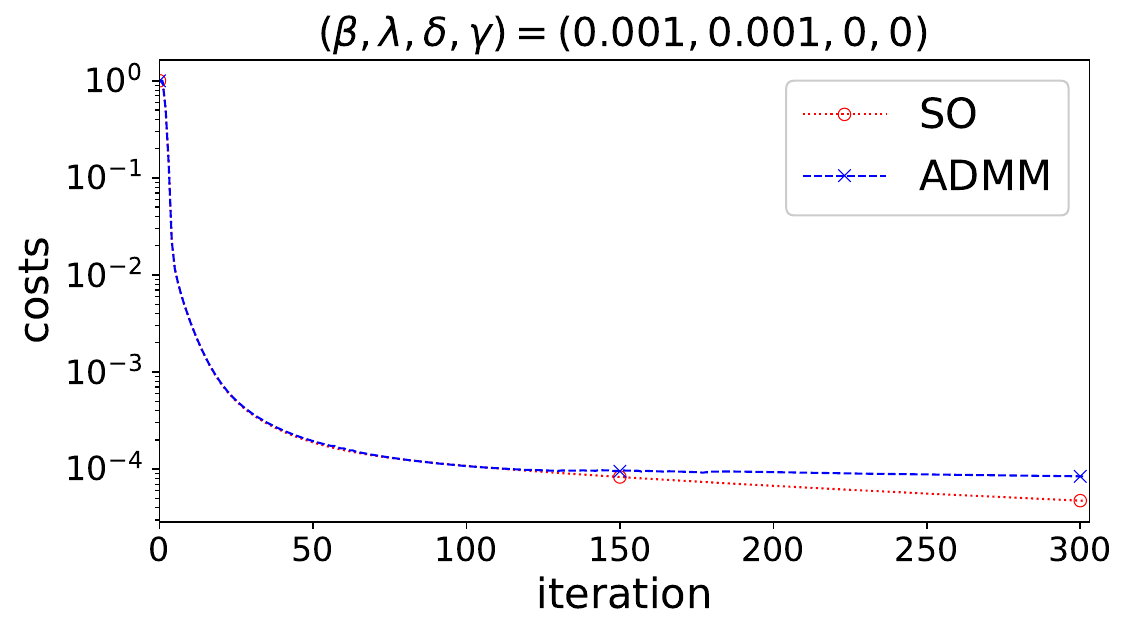}}  
\resizebox{0.24\linewidth}{!}{\includegraphics{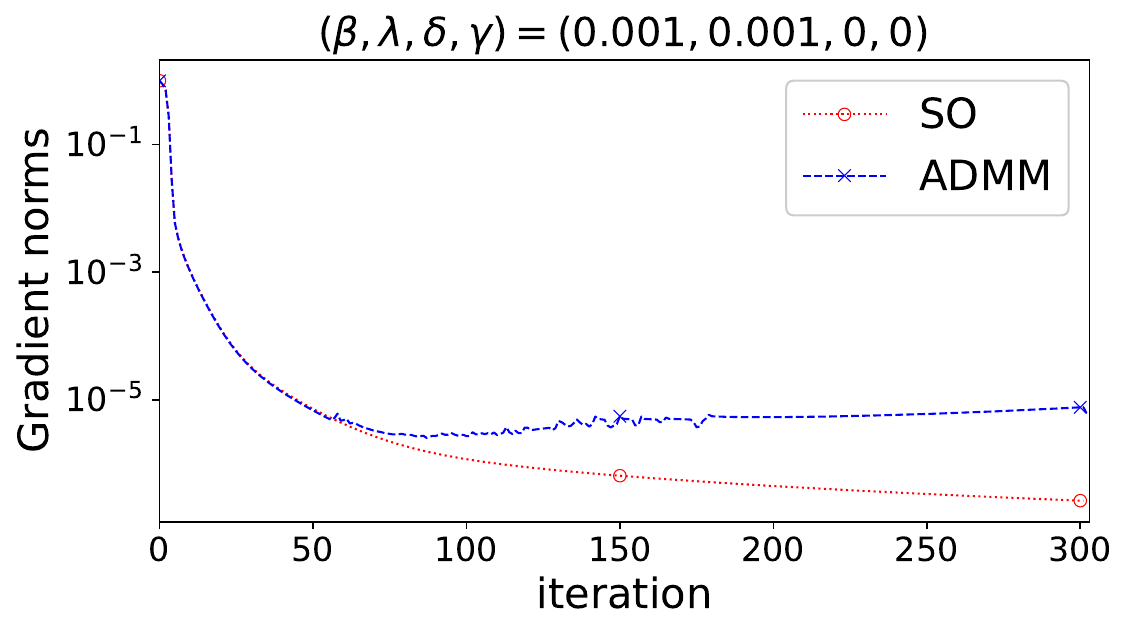}}  
\resizebox{0.24\linewidth}{!}{\includegraphics{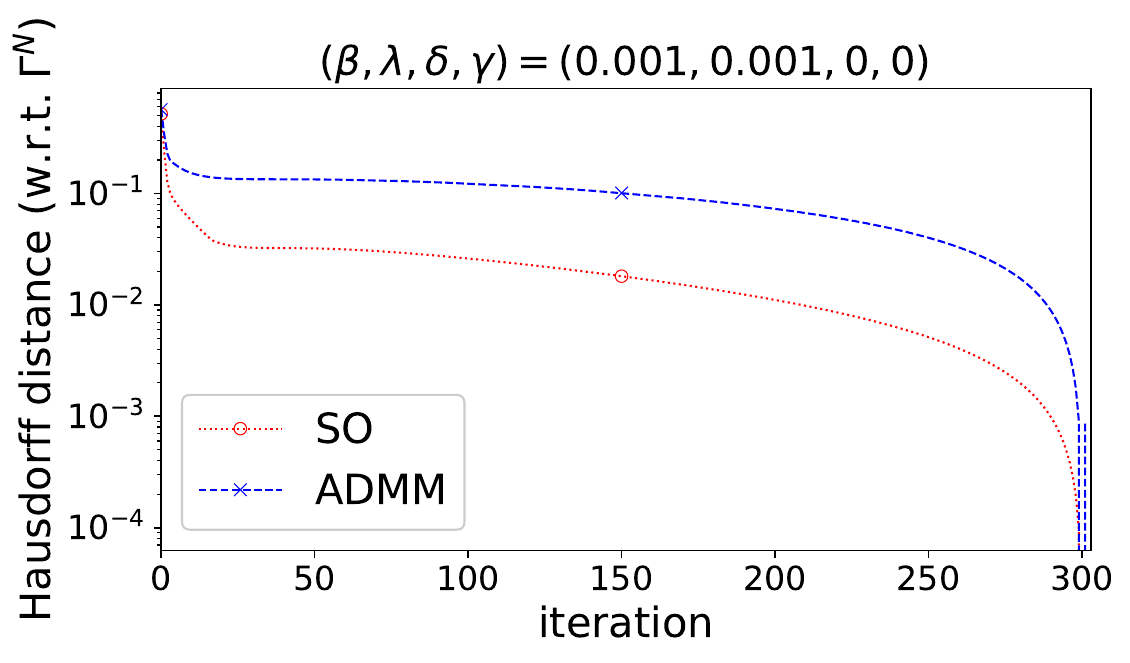}}  
\resizebox{0.255\linewidth}{!}{\includegraphics{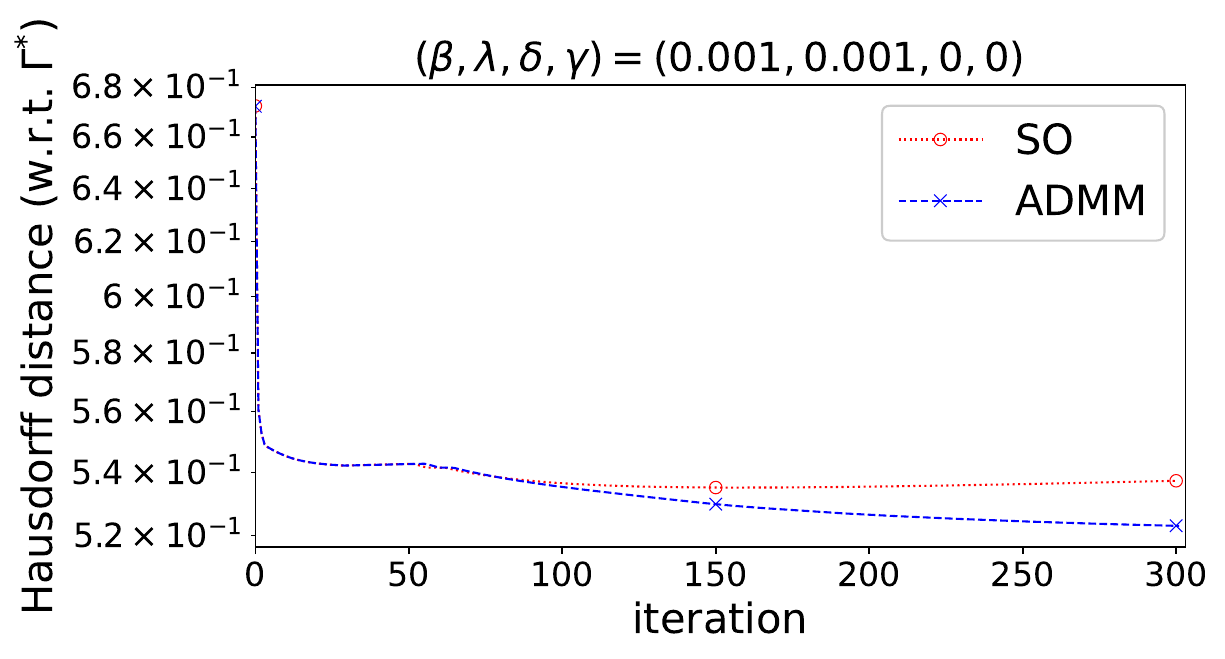}} 
\caption{Histories of costs, gradient norms, and Hausdorff distances with respect to the final computed shape ${\partial}{\omega}^{N}$ (third row), $N=300$, and the exact cavity shape ${\partial}{\omega}^{\ast}$ (last row) corresponding to the case of the \textsf{L}-shape cavity shown in Figure \ref{fig:figure1a} with $\beta = 0.001$} 
\label{fig:figure1b}
\end{figure}

\begin{figure}[htp!]
\centering
\resizebox{0.325\linewidth}{!}{\includegraphics{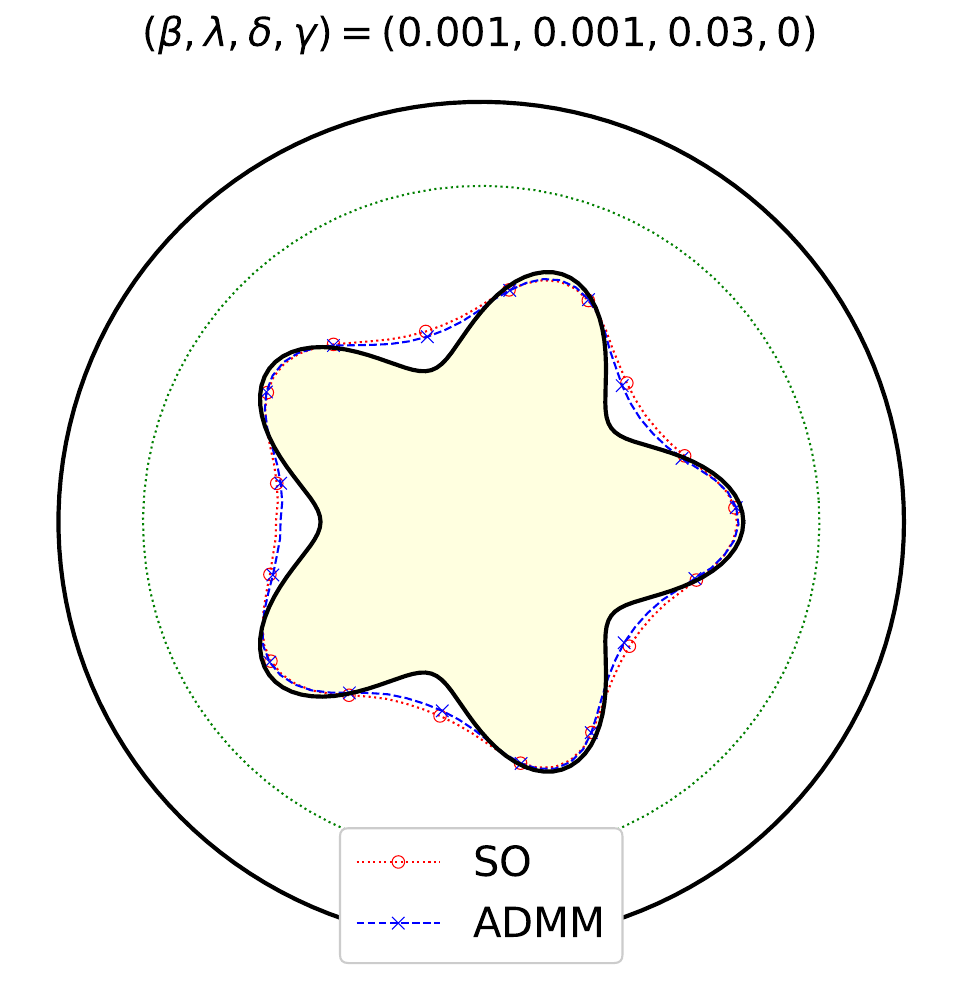}} 
\resizebox{0.325\linewidth}{!}{\includegraphics{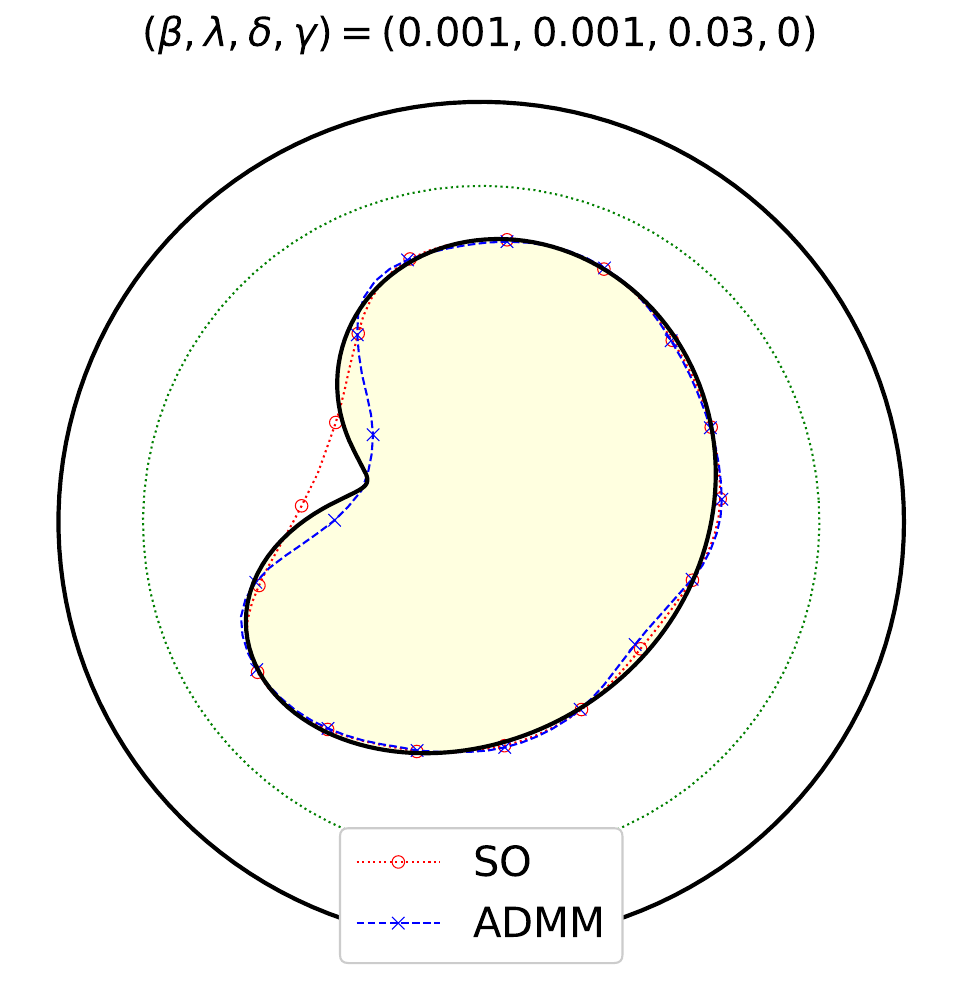}}  
\resizebox{0.325\linewidth}{!}{\includegraphics{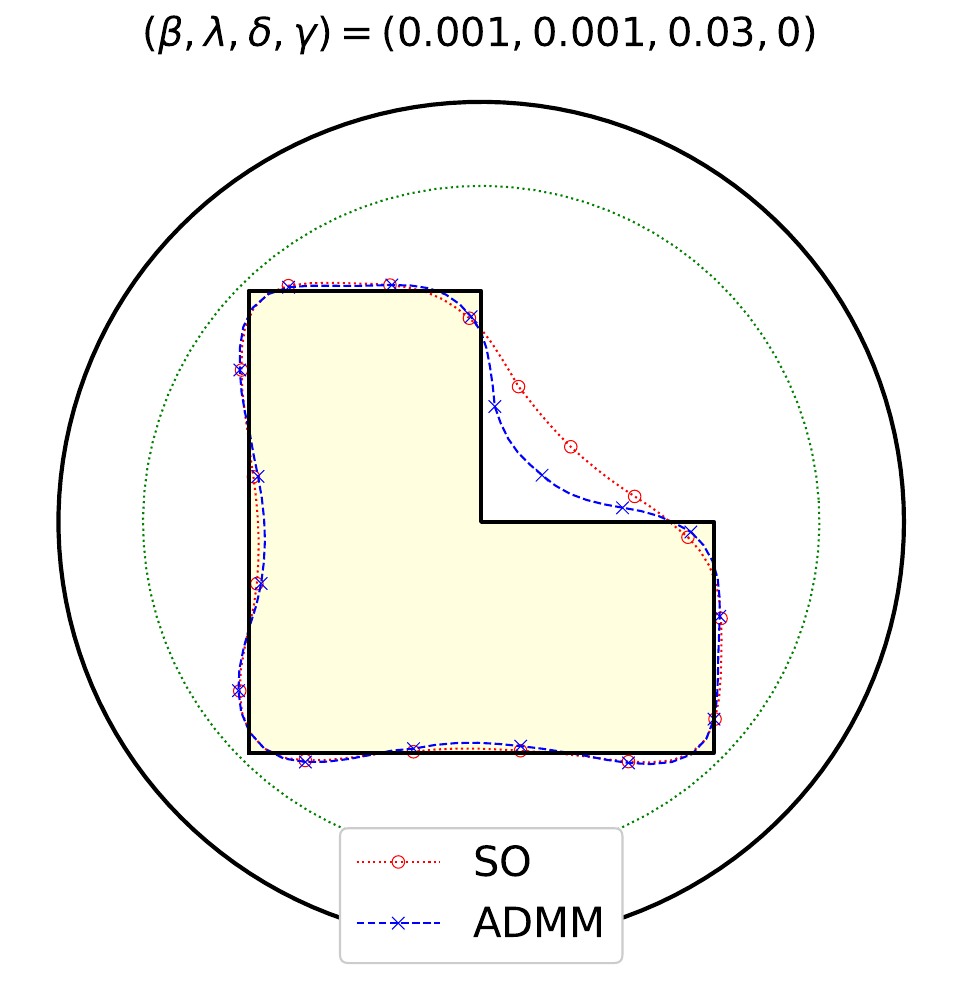}} 
\\[0.7em] 
\resizebox{0.325\linewidth}{!}{\includegraphics{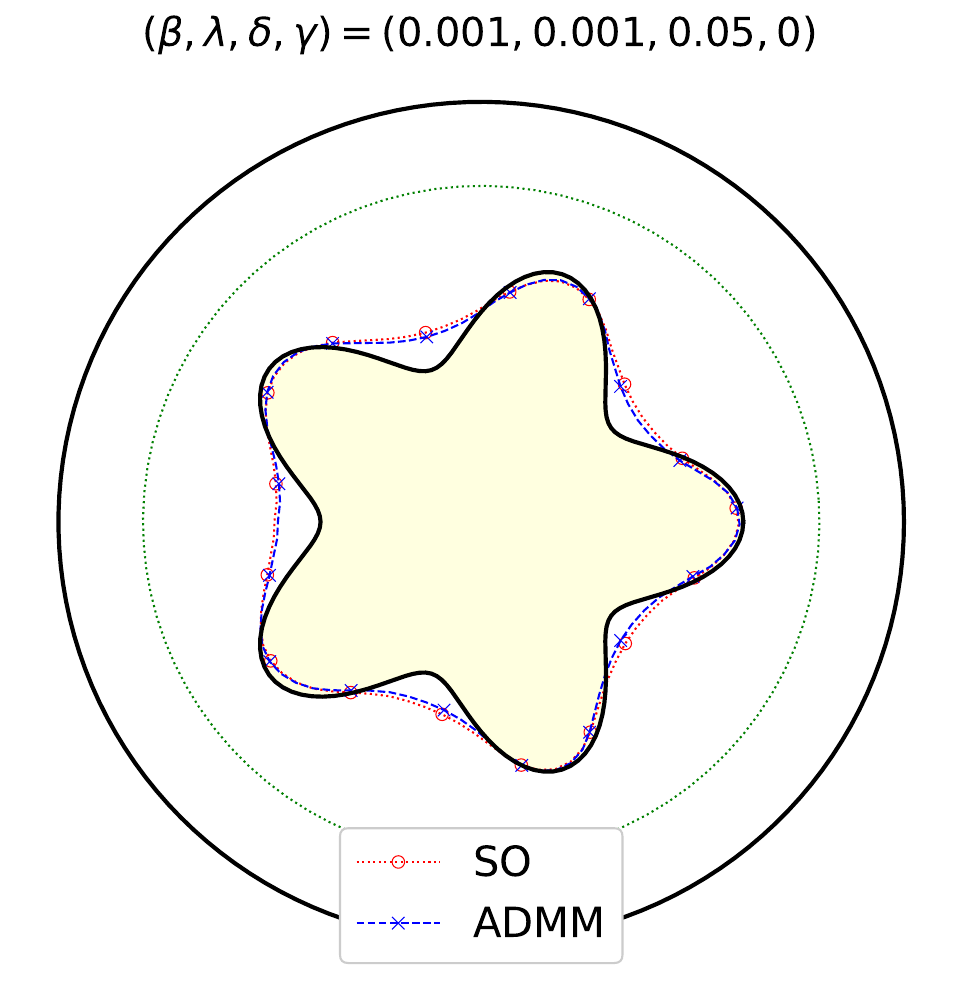}} 
\resizebox{0.325\linewidth}{!}{\includegraphics{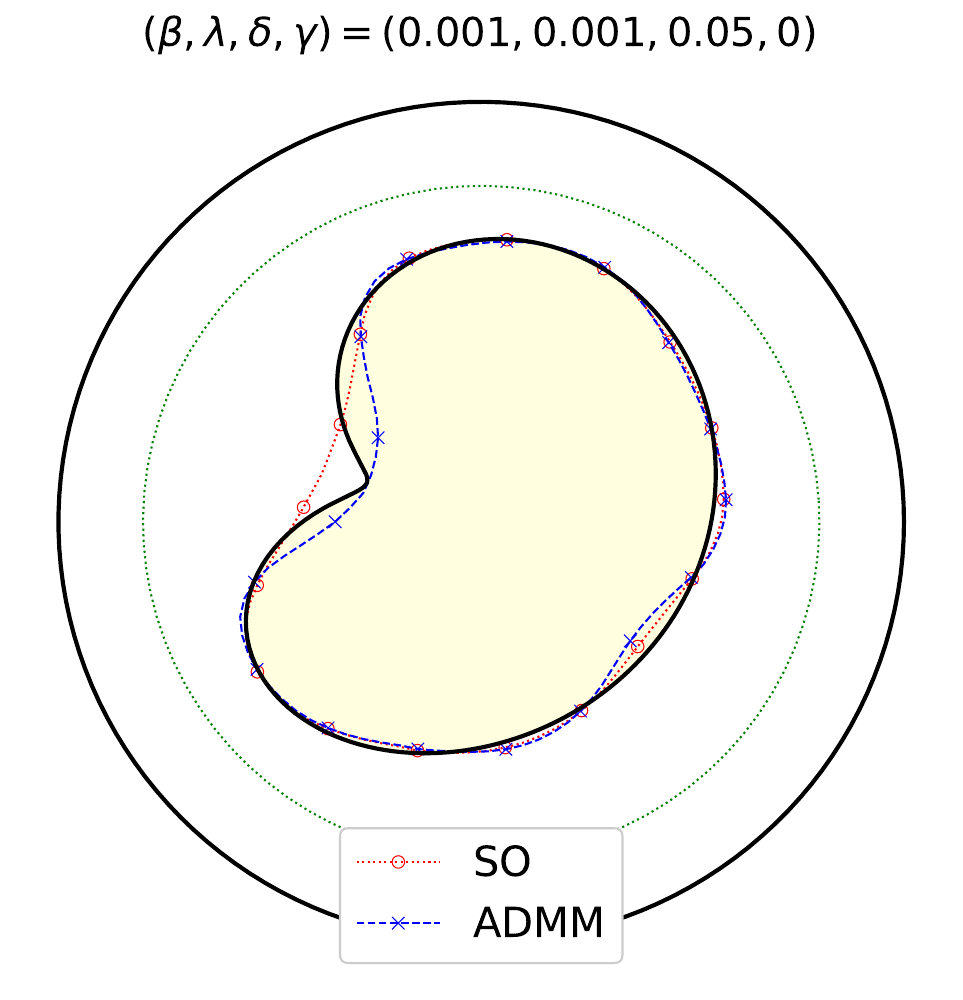}}  
\resizebox{0.325\linewidth}{!}{\includegraphics{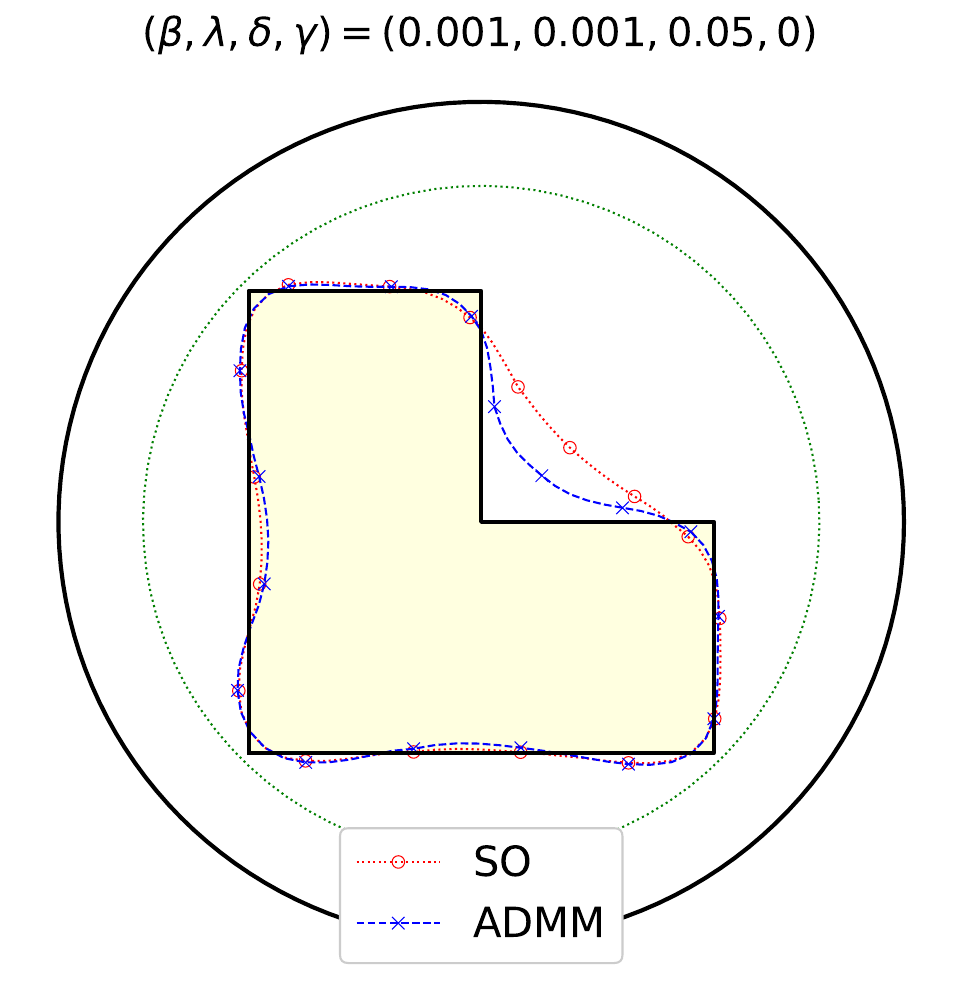}} 
\\[0.7em] 
\resizebox{0.325\linewidth}{!}{\includegraphics{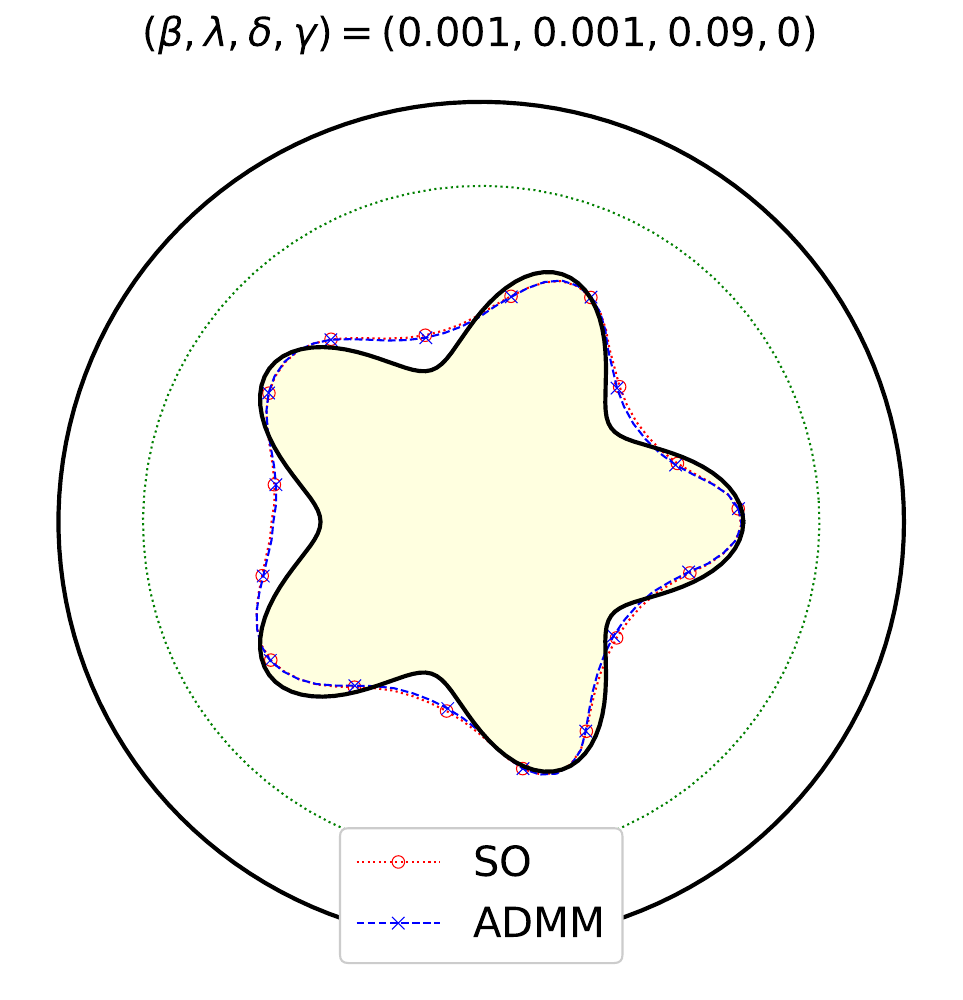}} 
\resizebox{0.325\linewidth}{!}{\includegraphics{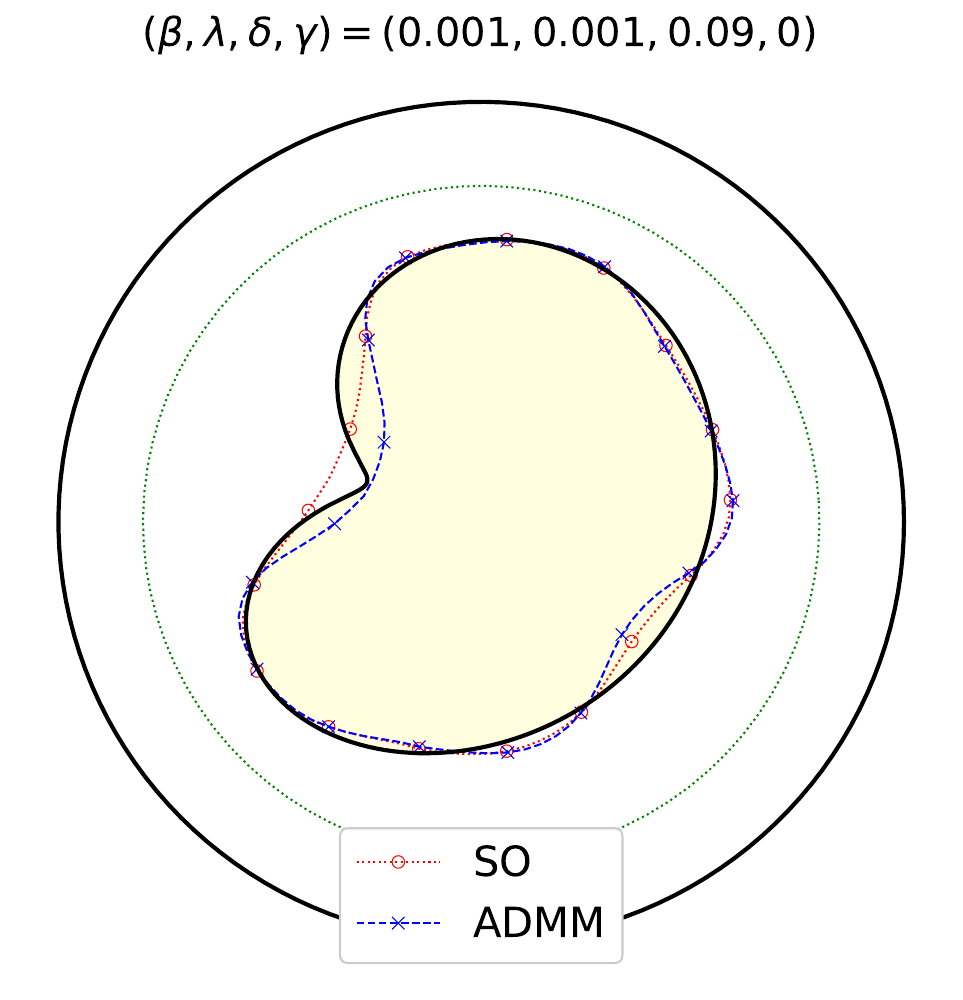}}  
\resizebox{0.325\linewidth}{!}{\includegraphics{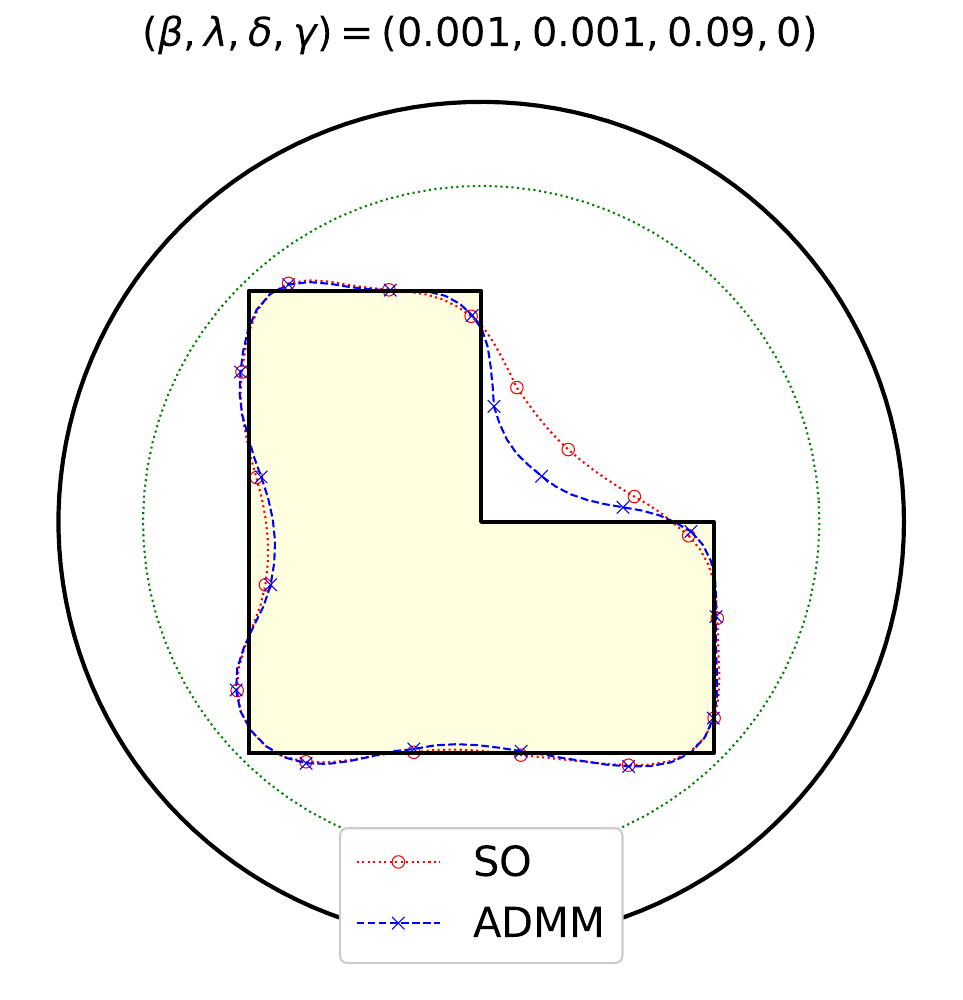}} 
\caption{Reconstructions with noisy data at varying noise levels ($\delta = 3\%, 5\%, 9\%$) with $\beta=0.001$ and without regularization (i.e., $\gamma = 0$)}
\label{fig:figure1c}
\end{figure}

\begin{figure}[htp!]
\centering 
\resizebox{0.24\linewidth}{!}{\includegraphics{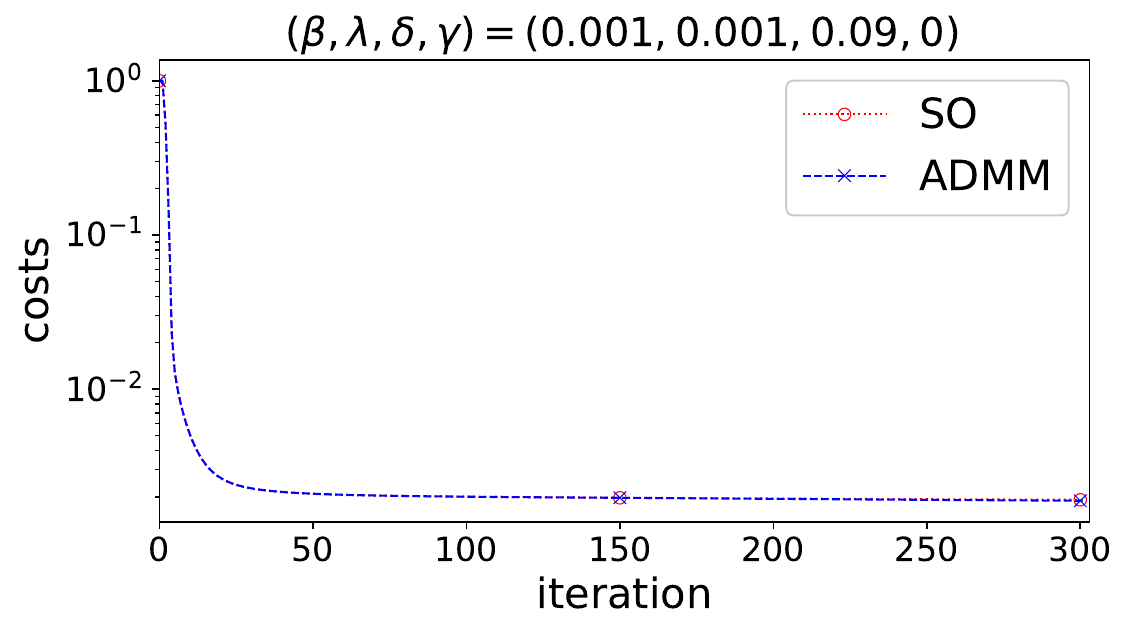}}  
\resizebox{0.24\linewidth}{!}{\includegraphics{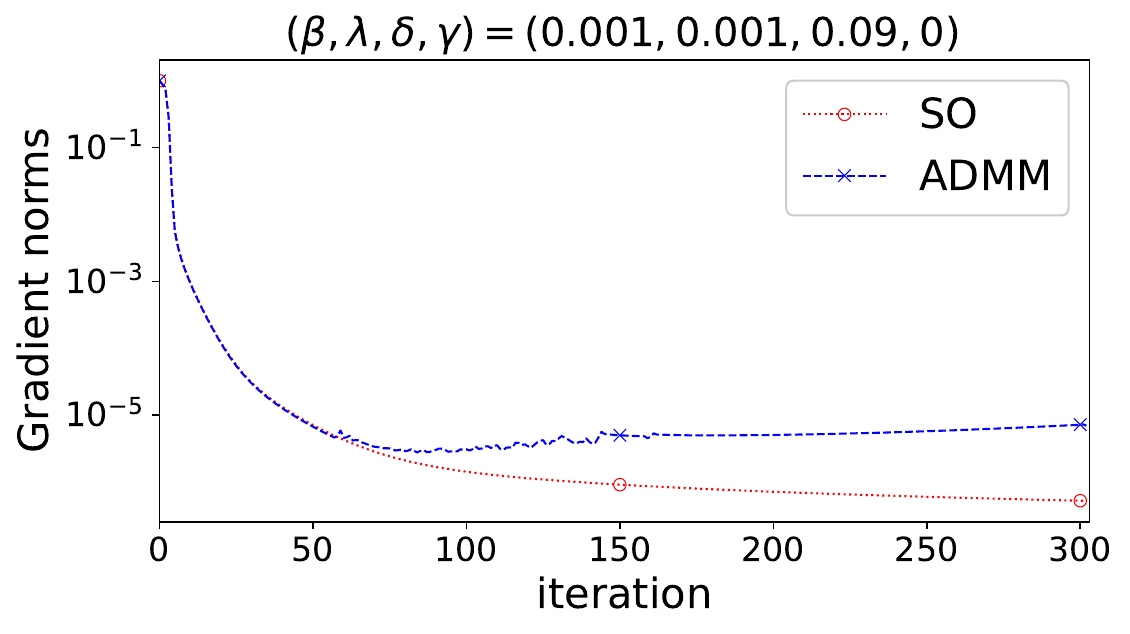}}  
\resizebox{0.24\linewidth}{!}{\includegraphics{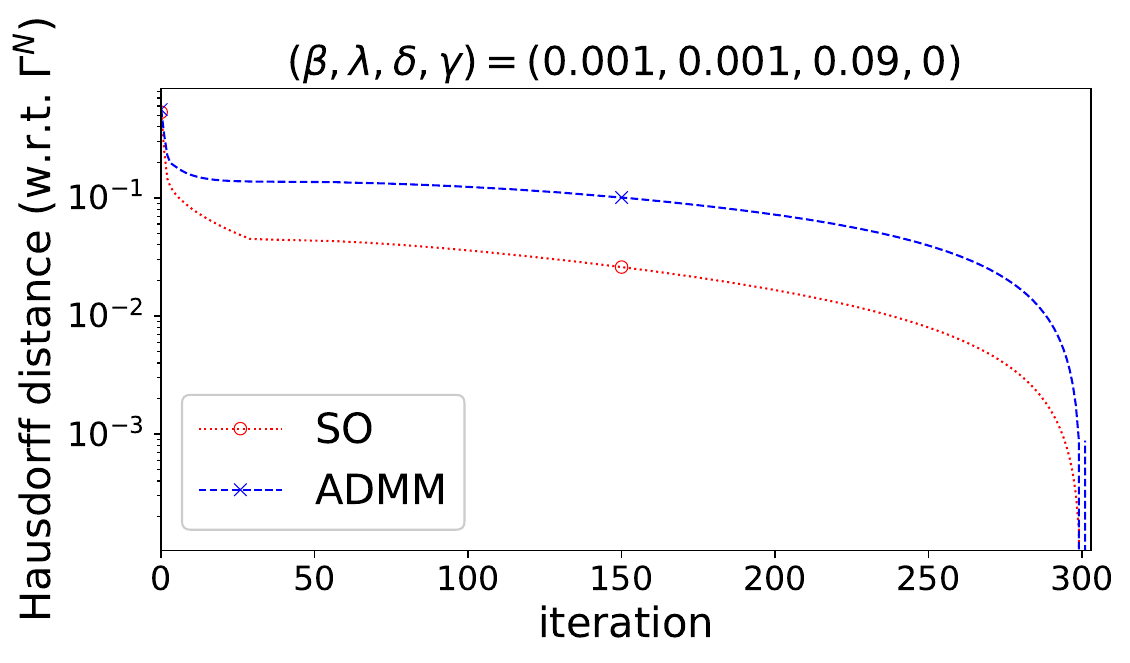}}  
\resizebox{0.255\linewidth}{!}{\includegraphics{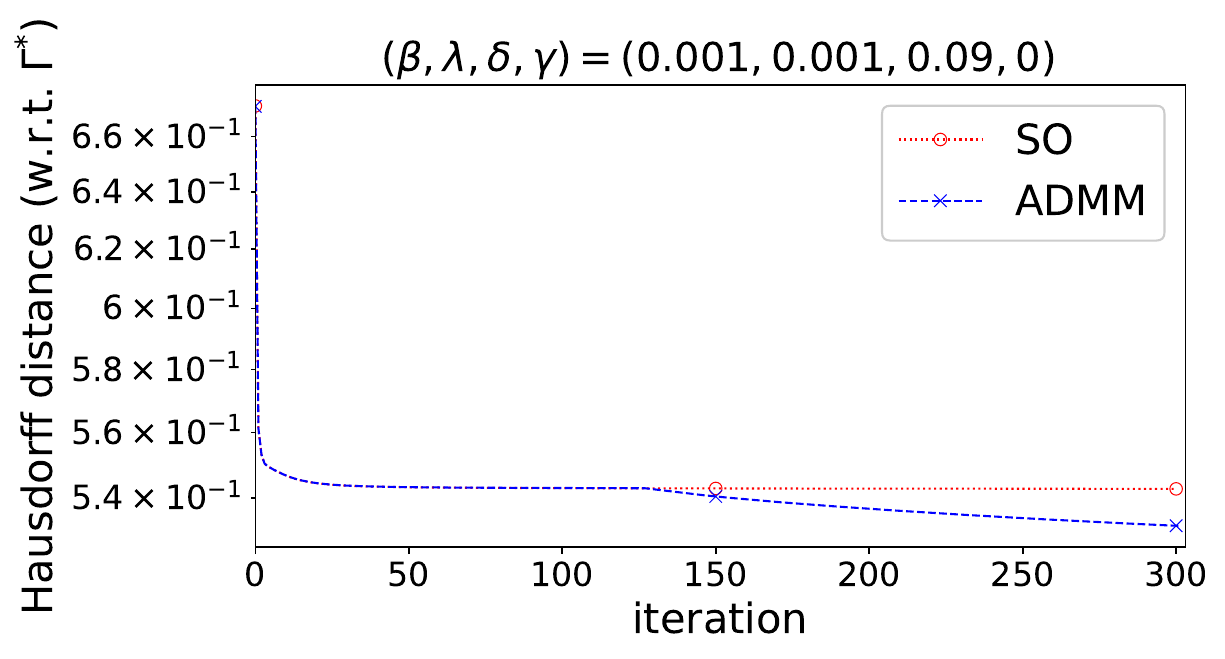}} 
\caption{Histories of costs, gradient norms, and Hausdorff distances with respect to the final computed shape ${\partial}{\omega}^{N}$ (third row), $N=300$, and the exact cavity shape ${\partial}{\omega}^{\ast}$ (last row) corresponding to the case of the \textsf{L}-shape cavity shown in Figure \ref{fig:figure1c} when $\delta = 9\%$} 
\label{fig:figure1d}
\end{figure}

\begin{figure}[htp!]
\centering
\resizebox{0.325\linewidth}{!}{\includegraphics{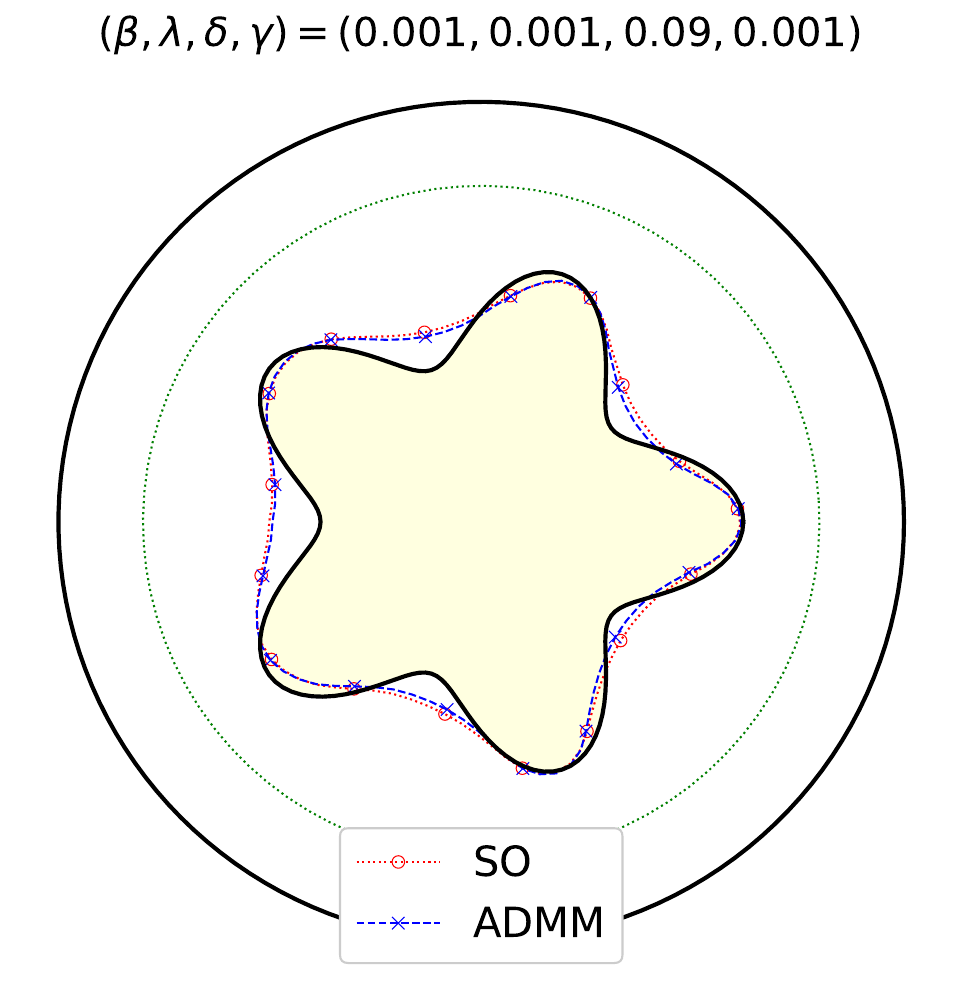}} 
\resizebox{0.325\linewidth}{!}{\includegraphics{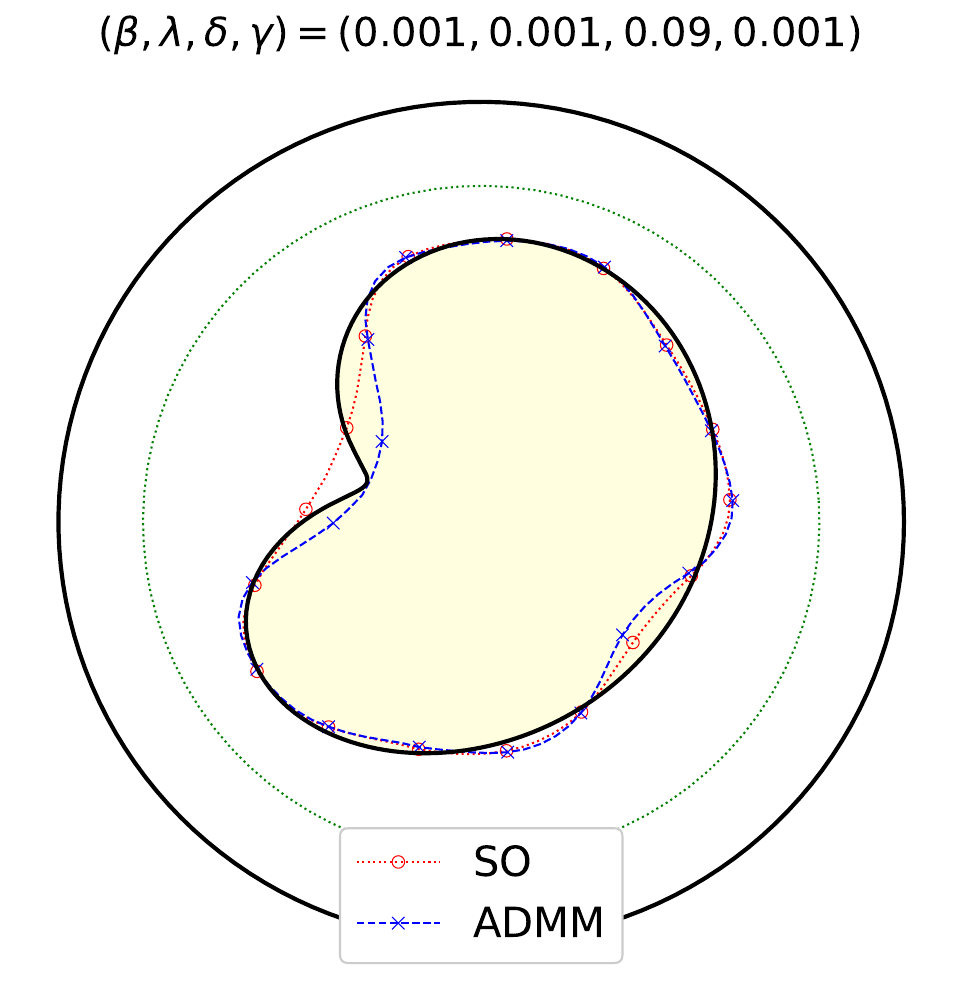}}  
\resizebox{0.325\linewidth}{!}{\includegraphics{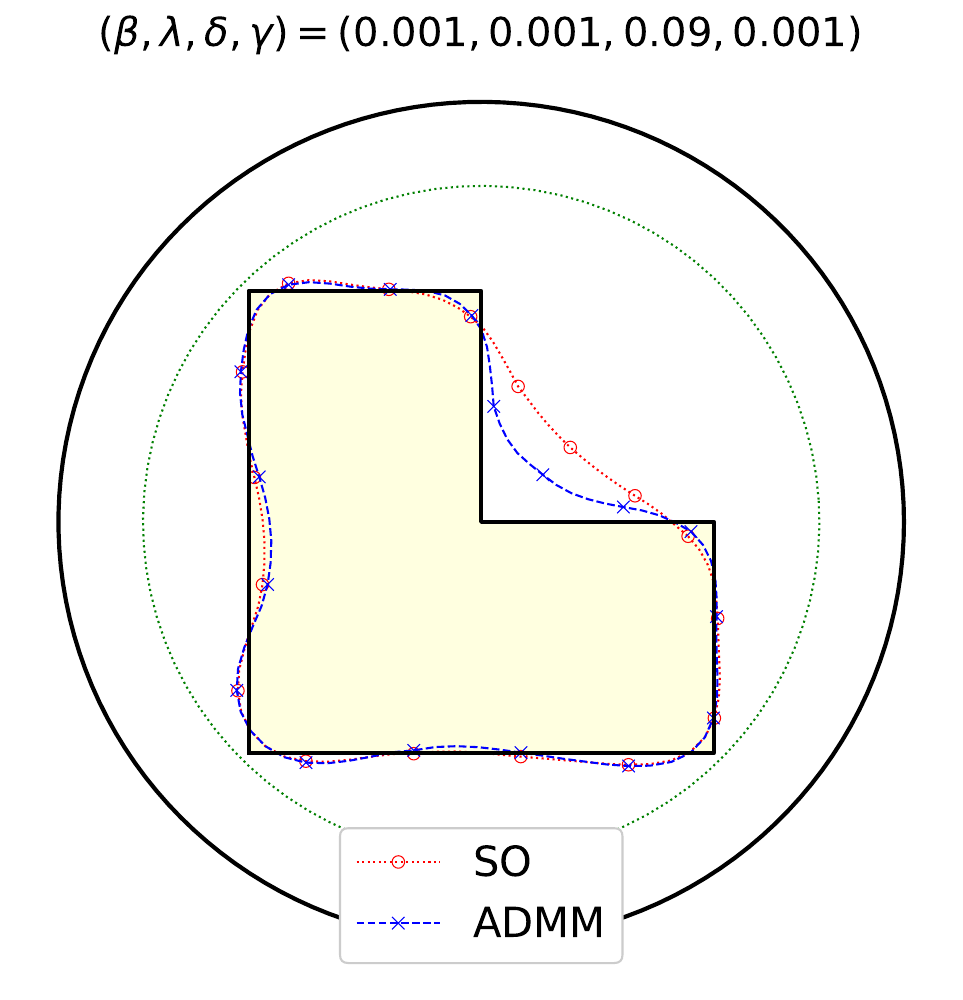}} 
\\[0.7em] 
\resizebox{0.325\linewidth}{!}{\includegraphics{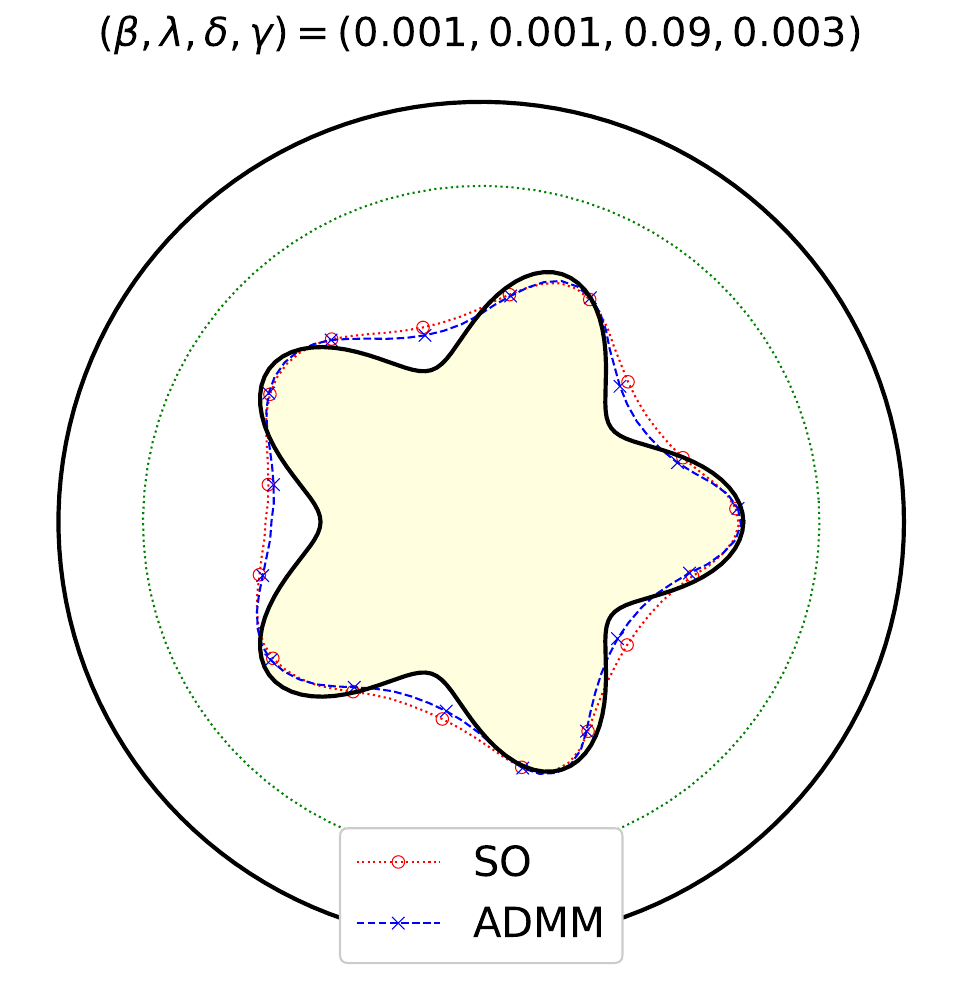}} 
\resizebox{0.325\linewidth}{!}{\includegraphics{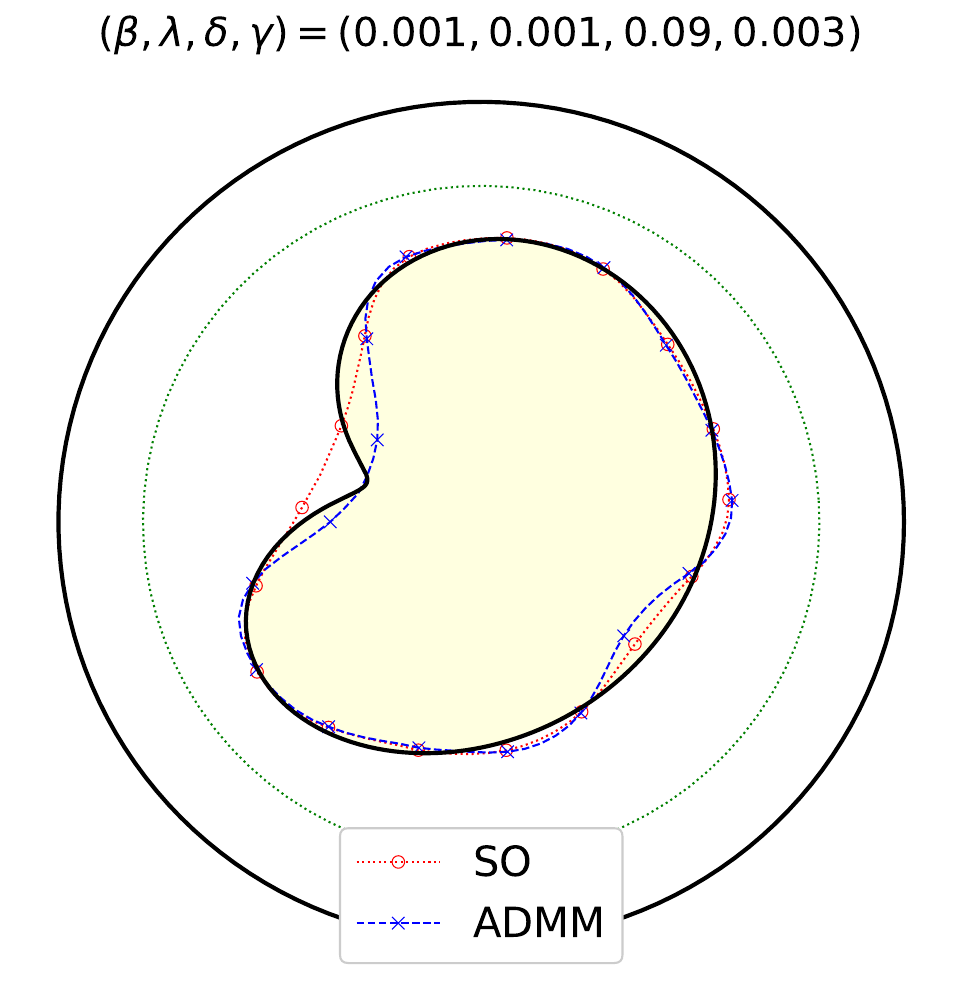}}  
\resizebox{0.325\linewidth}{!}{\includegraphics{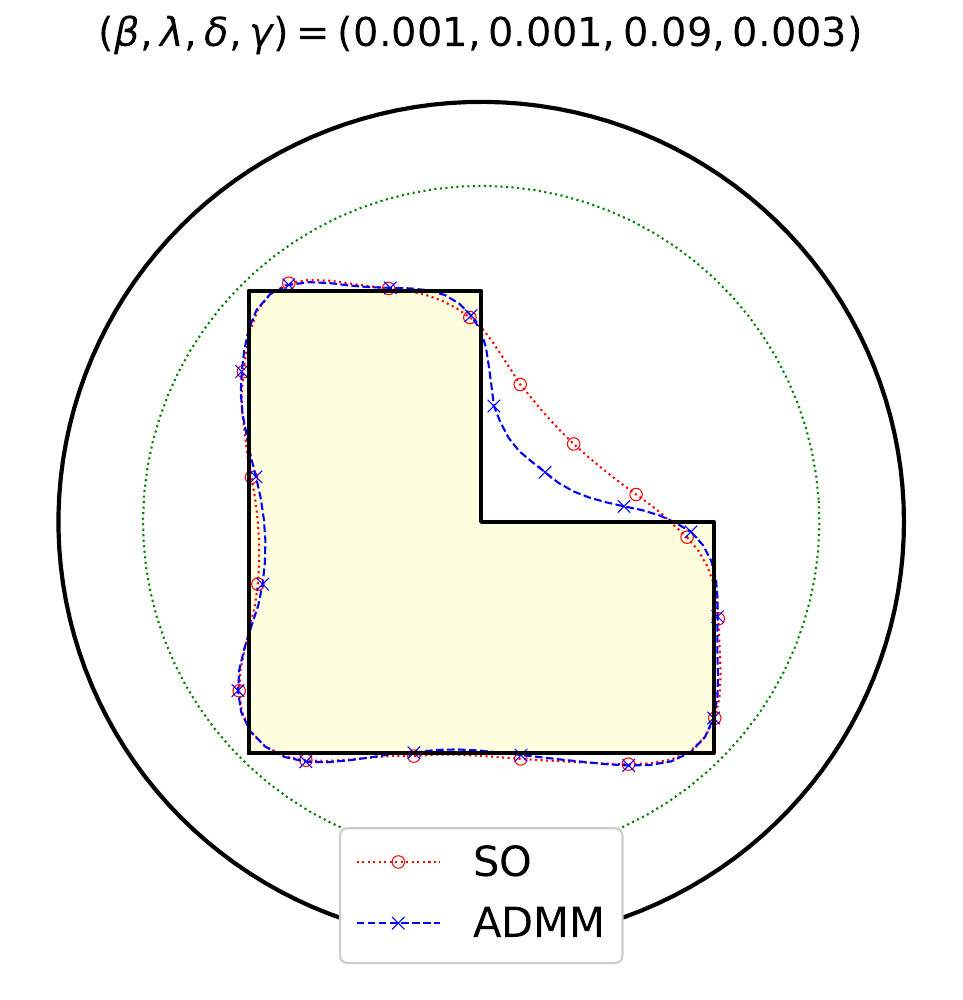}} 
\\[0.7em] 
\resizebox{0.325\linewidth}{!}{\includegraphics{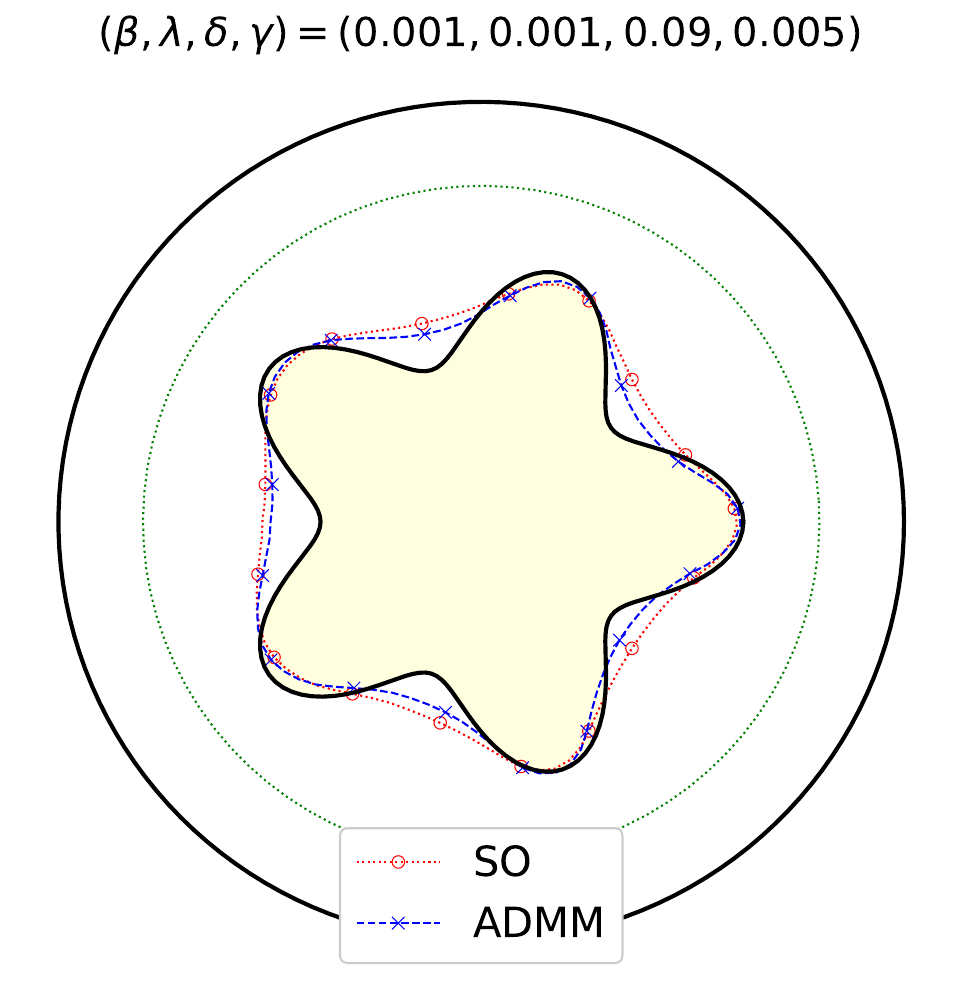}} 
\resizebox{0.325\linewidth}{!}{\includegraphics{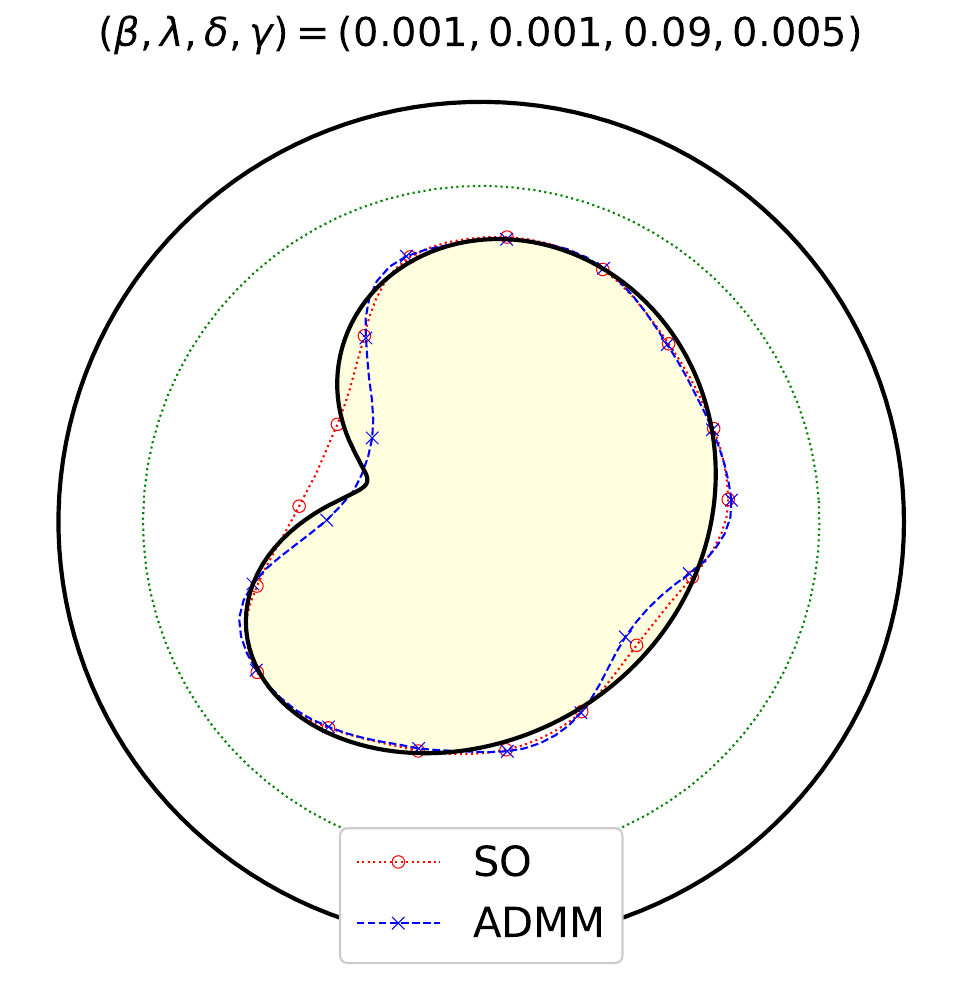}}  
\resizebox{0.325\linewidth}{!}{\includegraphics{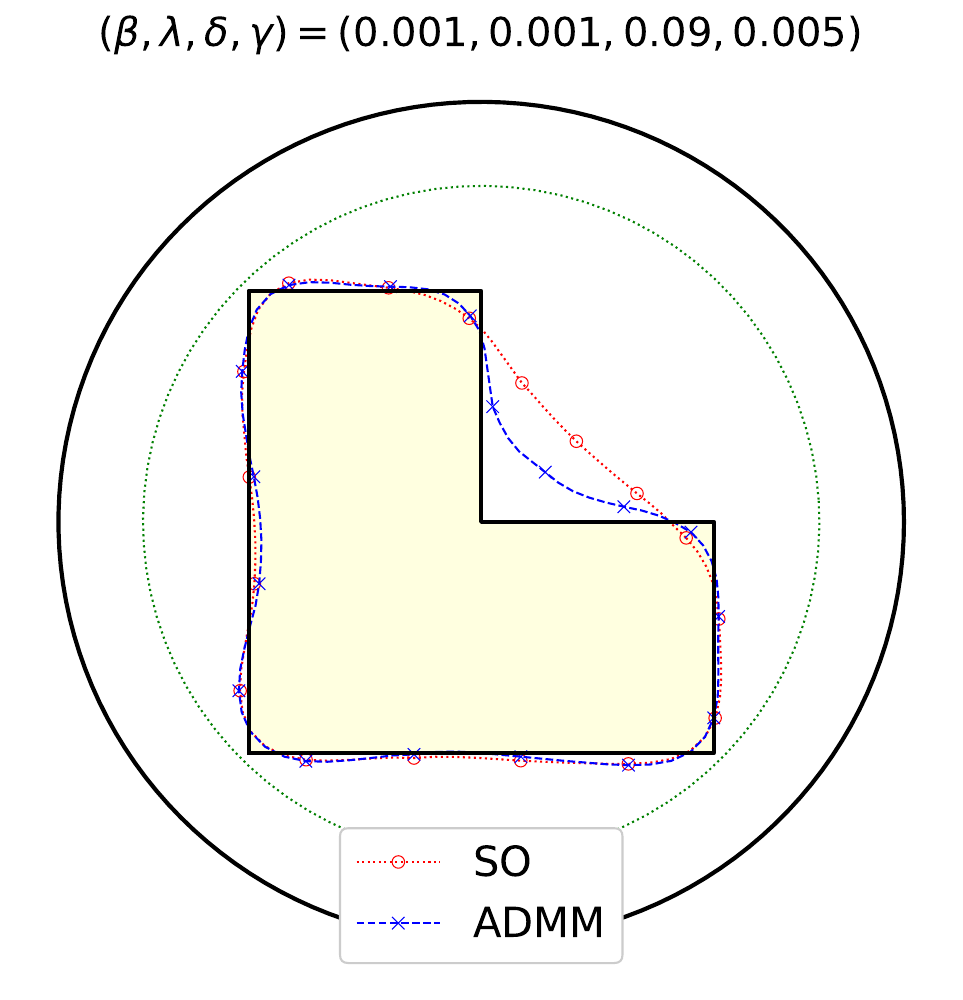}} 
\caption{Reconstructions with noisy data at noise level $\delta = 9\%$, with $\beta=0.001$ at varying levels of regularization ($\gamma = 0.001, 0.003, 0.005$)}
\label{fig:figure1e}
\end{figure}

\begin{figure}[htp!]
\centering 
\resizebox{0.24\linewidth}{!}{\includegraphics{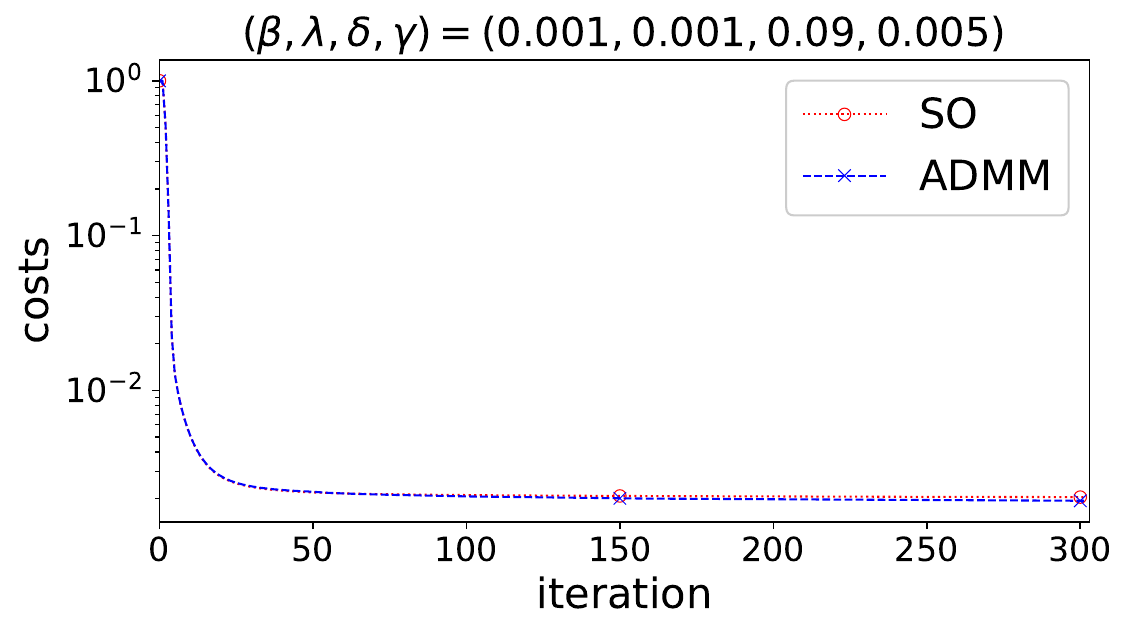}}  
\resizebox{0.24\linewidth}{!}{\includegraphics{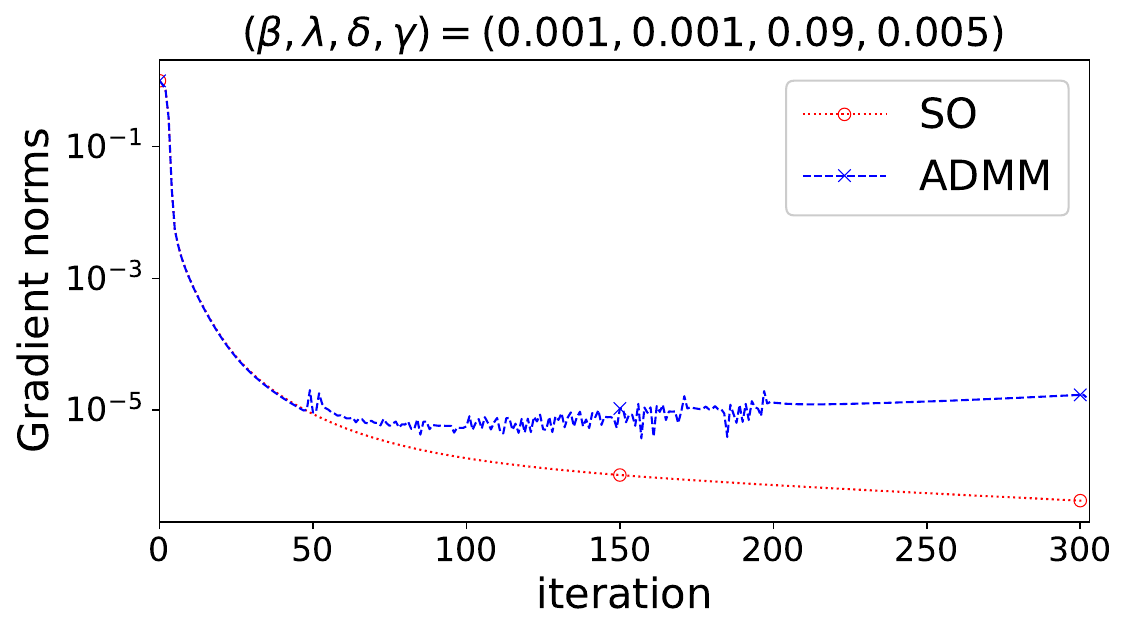}}  
\resizebox{0.24\linewidth}{!}{\includegraphics{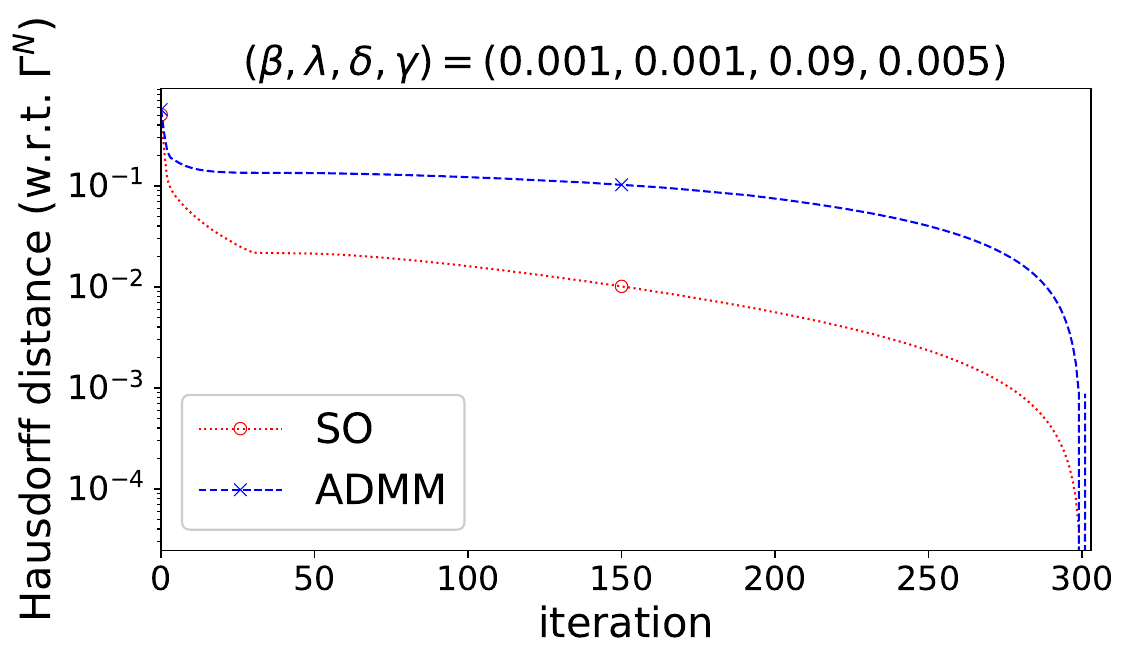}}  
\resizebox{0.255\linewidth}{!}{\includegraphics{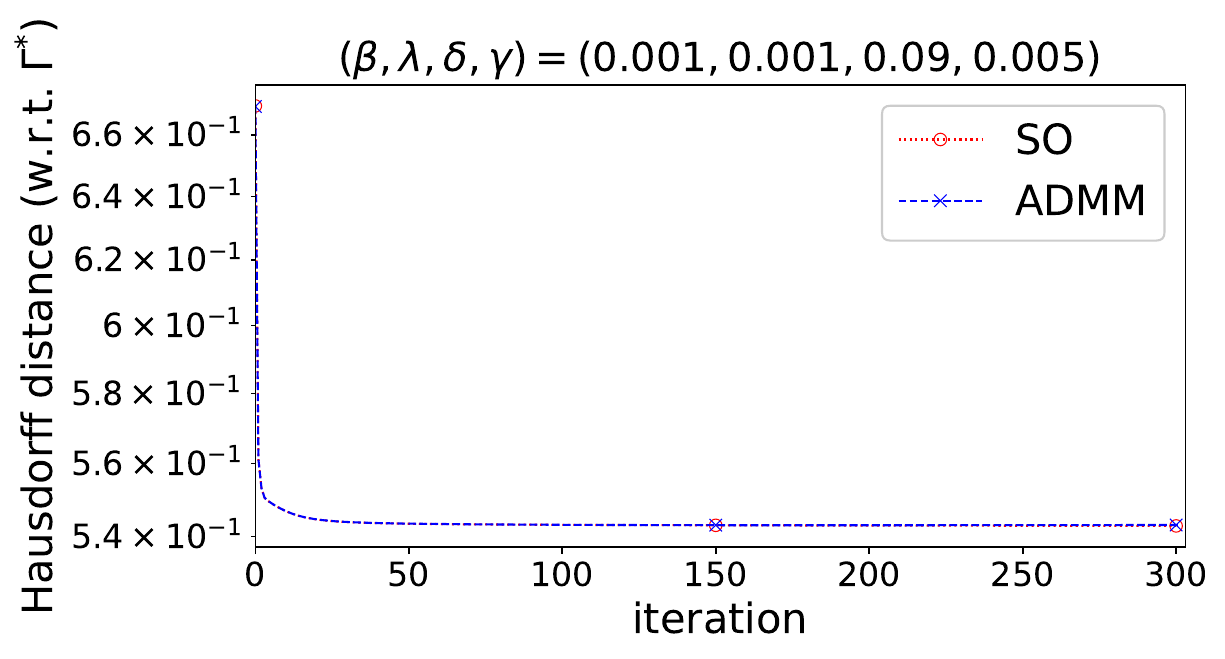}} 
\caption{Histories of costs, gradient norms, and Hausdorff distances with respect to the final computed shape ${\partial}{\omega}^{N}$ (third row), $N=300$, and the exact cavity shape ${\partial}{\omega}^{\ast}$ (last row) corresponding to the case of the \textsf{L}-shape cavity shown in Figure \ref{fig:figure1e} when $\gamma = 0.005$} 
\label{fig:figure1f}
\end{figure}

\begin{figure}[htp!]
\centering
\resizebox{0.325\linewidth}{!}{\includegraphics{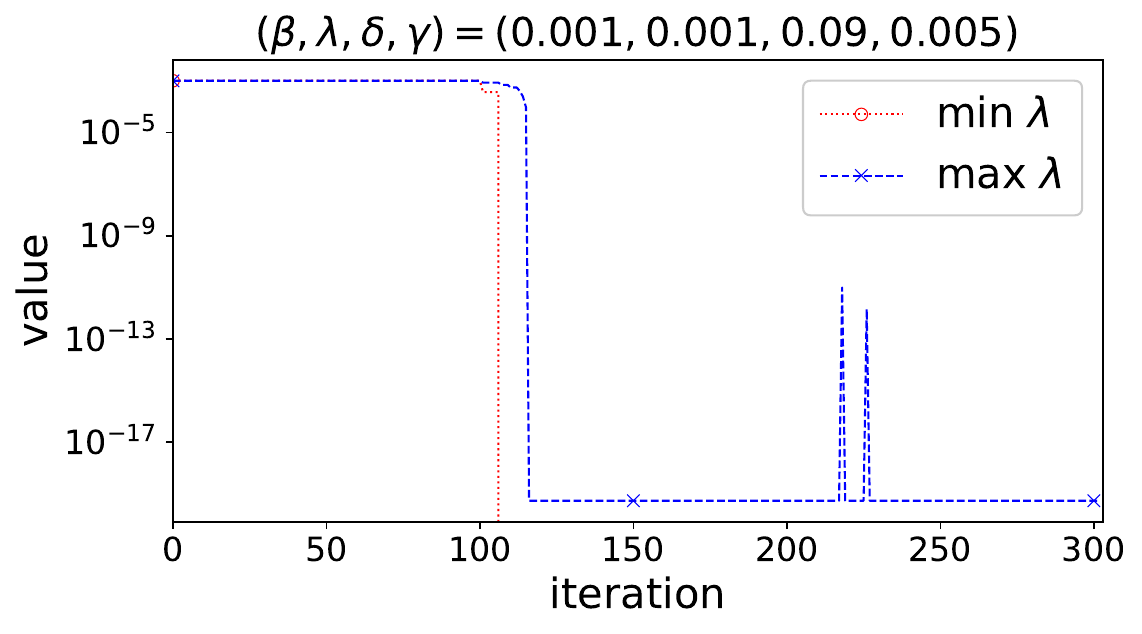}} 
\resizebox{0.325\linewidth}{!}{\includegraphics{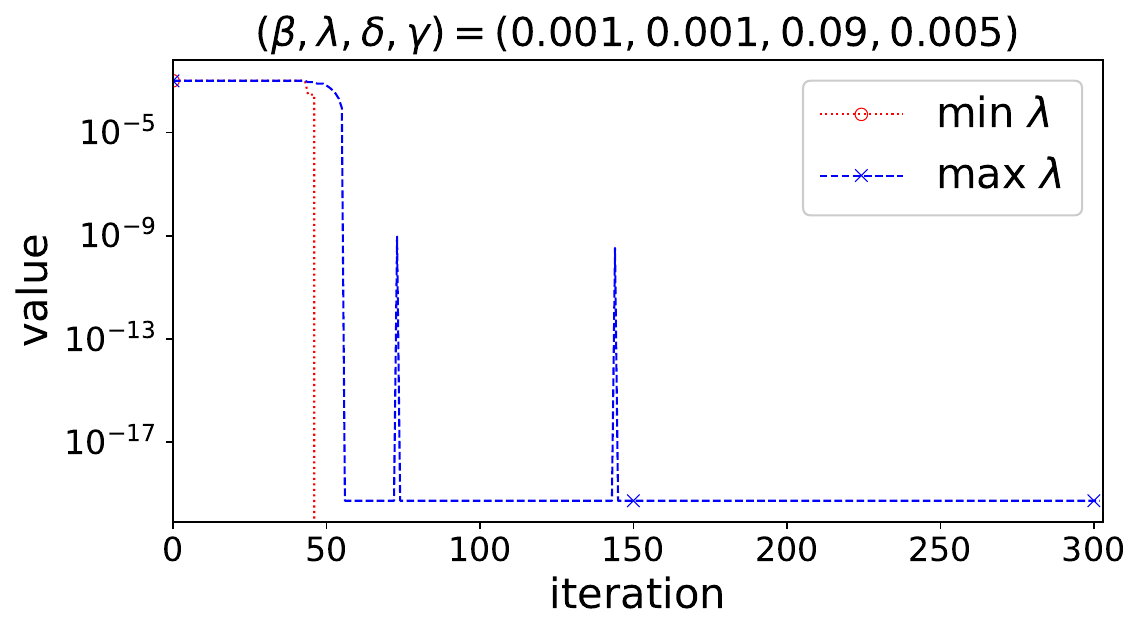}}  
\resizebox{0.325\linewidth}{!}{\includegraphics{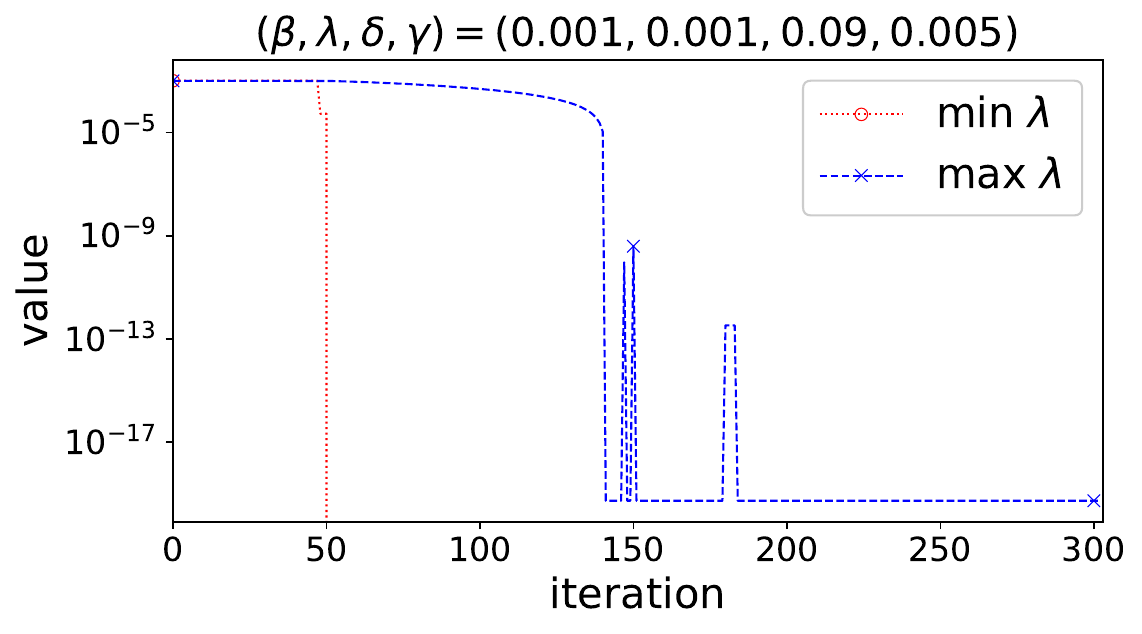}}
\caption{Histories of minimum and maximum values of $\lambda$ corresponding to the plots shown Figure \ref{fig:figure1f}} 
\label{fig:figure1g}
\end{figure}

\begin{figure}[htp!]
\centering
\resizebox{0.325\linewidth}{!}{\includegraphics{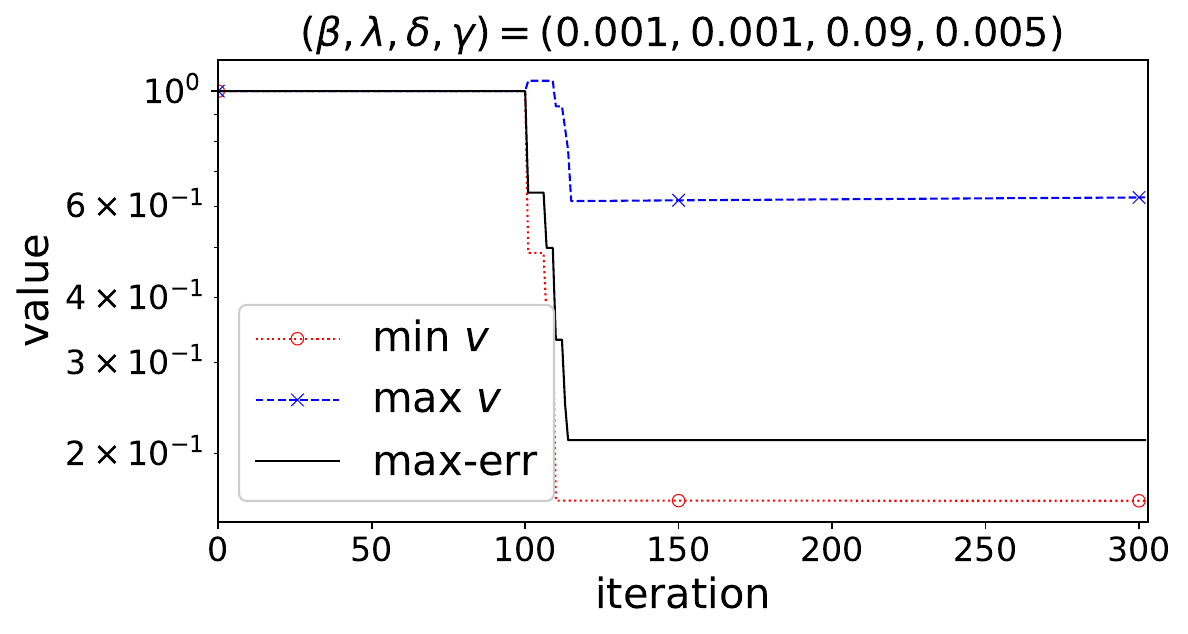}} 
\resizebox{0.325\linewidth}{!}{\includegraphics{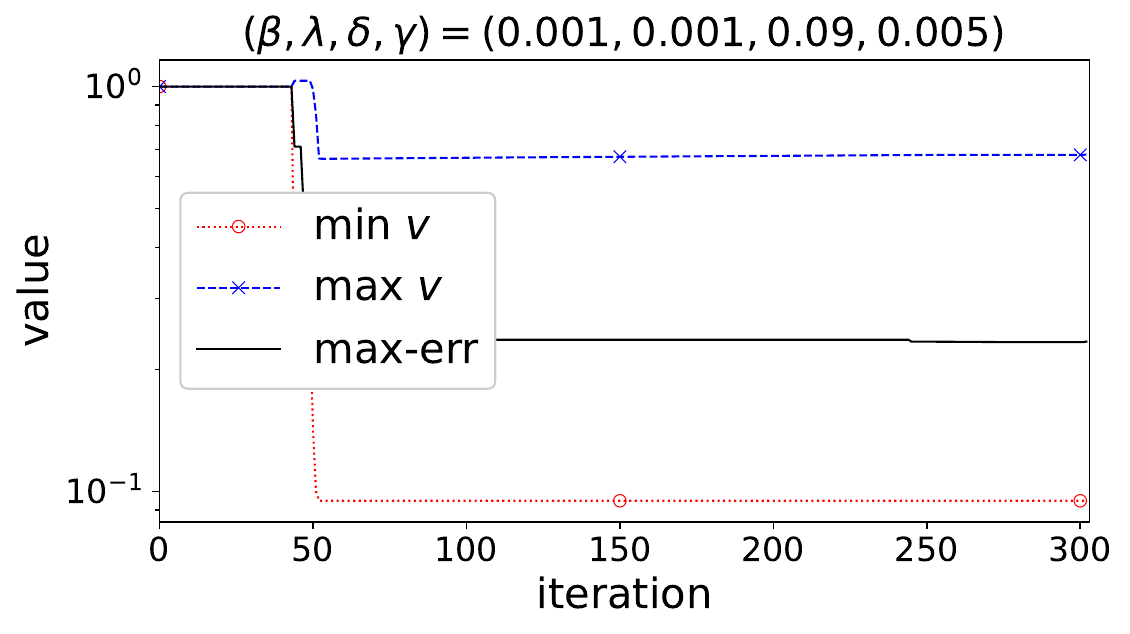}}  
\resizebox{0.325\linewidth}{!}{\includegraphics{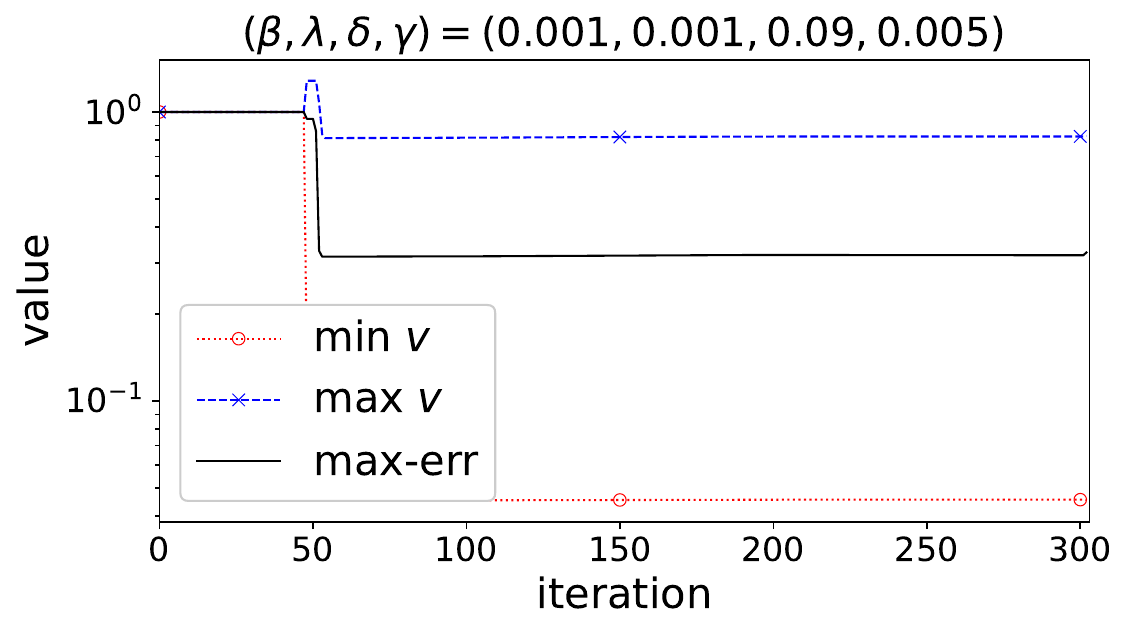}}
\caption{Histories minimum and maximum values of $v$ and the maximum error $\text{max-err}:=\|\uo - v\|_{L^{\infty}(\Omega\setminus\overline{\omega})}$ corresponding to the plots shown Figure \ref{fig:figure1f}} 
\label{fig:figure1h}
\end{figure}

\begin{remark}
	The choice of $a$ (and possibly $b$) can be determined by maximum principle.
	Since $u=0$ on ${\partial}{\omega}$, we can take $a = 0$, while $b$ can chosen depending on the input data $g$.
	For the present case, we can in fact take $a=0$ and overestimate $b$ by taking $b=2$.
	In our experience, these choices are also effective for ADMM.
\end{remark}
The results of the detections for the considered test cases with exact data are displayed in Figures \ref{fig:figure1a} through Figure \ref{fig:figure1h}.
In the plotted shapes (Figures \ref{fig:figure1a}, \ref{fig:figure1c}, and \ref{fig:figure1e}), the black solid lines represent the exact medium, while the green dotted lines depict the initial guess. 
The red dotted lines with circle markers represent the shape obtained by the classical shape optimization method (hereinafter abbreviated as SO), while the blue dashed lines with cross markers represent the shape obtained using ADMM. 
The main findings from our numerical experiments are as follows.
\begin{itemize}
\item As evident in Figure \ref{fig:figure1a}, when using exact data for the inversion, the proposed ADMM approach provides more accurate cavity detections -- as expected -- compared to classical shape optimization (SO). 
This observation remains consistent even in the presence of noisy data, as shown in Figure \ref{fig:figure1c} (see also Figure \ref{fig:figure1e}).
\item Indeed, with ADMM, we can achieve a more pronounced detection of the concave part(s) of the exact cavities, regardless of whether the data is contaminated with noise. Furthermore, reconstructions are faster with ADMM, as evidenced by the histories of Hausdorff distances $d_{H}({{\partial}{\omega}}^{\ast},{{\partial}{\omega}}^{k})$ between the exact ${{\partial}{\omega}}^{\ast}$ and the $k$th computed shape ${\partial}{\omega}^{k}$ of the cavity. 
For a visual representation, refer to Figure \ref{fig:figure1b} and \ref{fig:figure1d}.
\item We notice, however, that, in some situations, the cost computed due to ADMM is larger compared to SO, as seen in Figure \ref{fig:figure1b}. 
Nevertheless, in general, the cost values converge after some iterations.
\item We also observed, as one would expect (see \cite{Boydetal2011}), that the efficiency of ADMM depends highly on the magnitude of the free parameter $\beta$, at least when $\varepsilon$ is not set too small (refer again to Figure \ref{fig:figure1a}).
Note that, in the case of a very small $\varepsilon$, Steps 2.4 and 2.5 of Algorithm \ref{algo:ADMM-SGD} might never be reached or activated.
For very small values of $\beta$, we noticed from our experience that the proposed scheme is ineffective and exhibits almost the same convergence behavior and accuracy as that of SO (classical shape optimization). 
For larger values of $\beta$, on the other hand, ADMM has the tendency to overshoot the exact shape (this is also true when $\lambda$ is initially taken very large).
These are not surprising, as a balance between the primal and dual residuals in the approximation is necessary, similar to many ADMM-type algorithms (see, e.g., \cite{GlowinskiSongYuan2020}).
From the values we have tested, we found that ADMM performs well within the range $\beta \in (0.001, 0.01)$ but proves ineffective when taken far from this interval -- at least in the case of the present test scenarios.
\end{itemize}
For illustration purposes, the histories of gradient norms and Hausdorff distances, denoted as $d_{H}({{\partial}{\omega}}^{N},{{\partial}{\omega}}^{k})$, between the final computed shape ${{\partial}{\omega}}^{N}$ and the $k$th approximation ${\partial}{\omega}^{k}$ of the exact cavity shape are plotted in Figures \ref{fig:figure1b}, \ref{fig:figure1d}, and \ref{fig:figure1f}.
Based on these histories, it appears that ADMM and SO almost have the same rate of convergence in most situations, especially in the case of noisy data. 
However, from the histories of Hausdorff distances, it is evident that ADMM provides a faster and more accurate approximation of the exact cavity than SO.
The impact of the regularization parameter $\gamma$ when considering noisy data, for specific noise levels, is summarized in Figures \ref{fig:figure1e} and \ref{fig:figure1f}.
Additionally, for further illustration, the histories of minimum and maximum values of $\lambda$, the histories of minimum and maximum values of $v$, and the maximum error $\text{max-err}:=\|\uo - v\|_{L^{\infty}(\Omega\setminus\overline{\omega})}$ corresponding to the plots in Figure \ref{fig:figure1f} are respectively shown in Figures \ref{fig:figure1g} and \ref{fig:figure1h}.
%
%
%
%
%
%
%
%
%
%
%
%
%
%
%
%
%
%
%
%
%
%
%
%
%
%
%
%
%
%
%
%
%
%
%
%
%
%
%
%
%

\medskip
\textbf{3D case.}Let us now consider test cases in the three spatial dimensions to further evaluate our algorithm. 
The computational setup remains unchanged, with only a few modifications. Specifically, we set $N = 600$, $\lambda^{0} = 0.001$, $a = 0.5\min u(\Omega \setminus \overline{\omega}^{\ast})$, $b = 1.5\max u(\Omega \setminus \overline{\omega}^{\ast})$, $v^{0} = 1$, $\varepsilon = 10^{-6}$, and $\omega^{0} = B(\mathbf{0}, 0.8)$.

\begin{figure}[htp!]
\centering
\resizebox{0.16\linewidth}{!}{\includegraphics{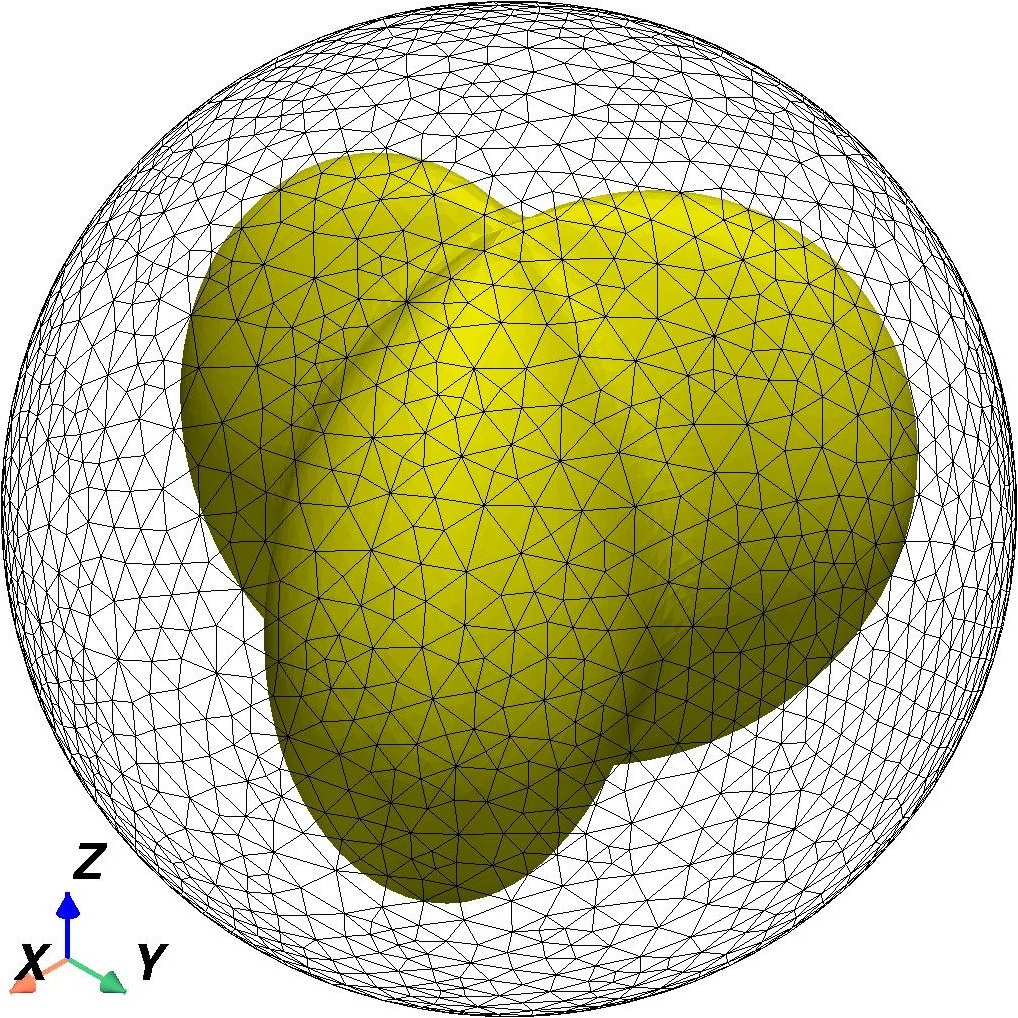}} \hfill
\resizebox{0.16\linewidth}{!}{\includegraphics{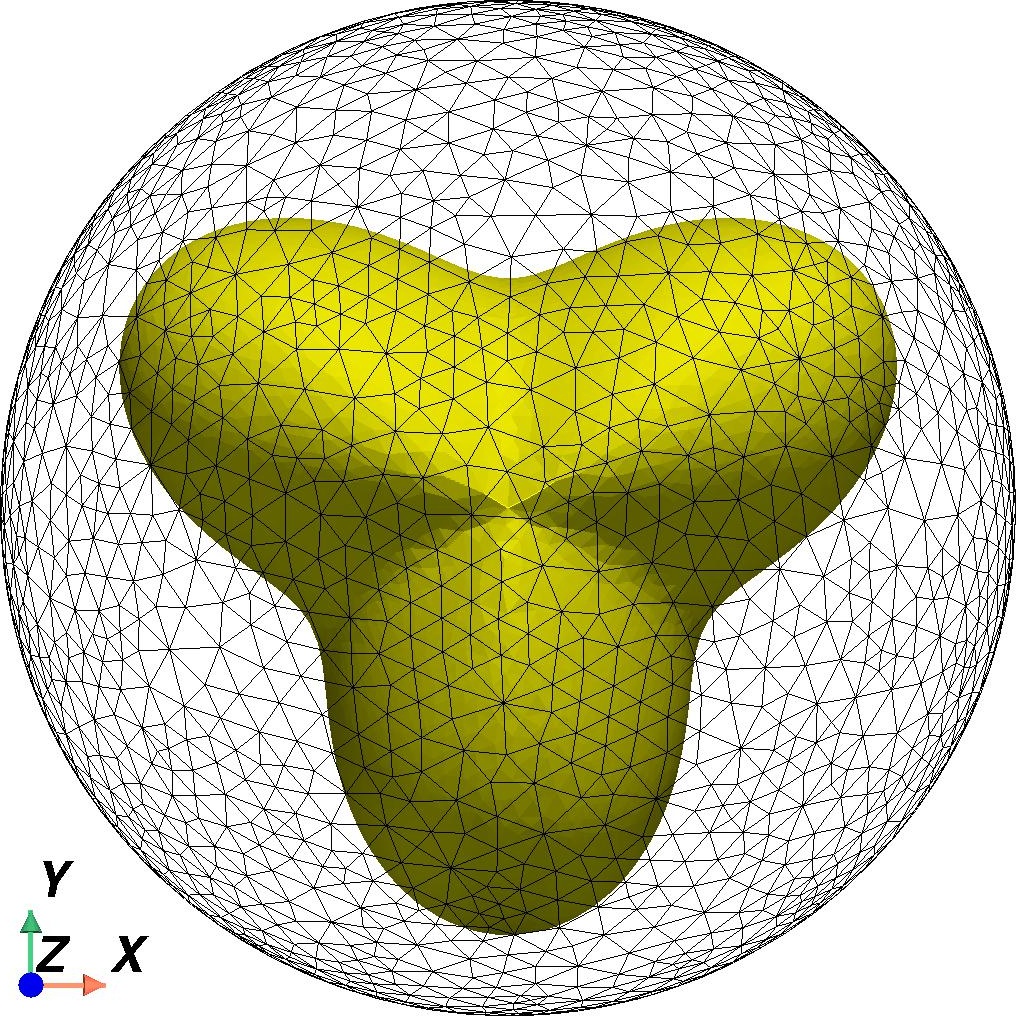}} \hfill
\resizebox{0.16\linewidth}{!}{\includegraphics{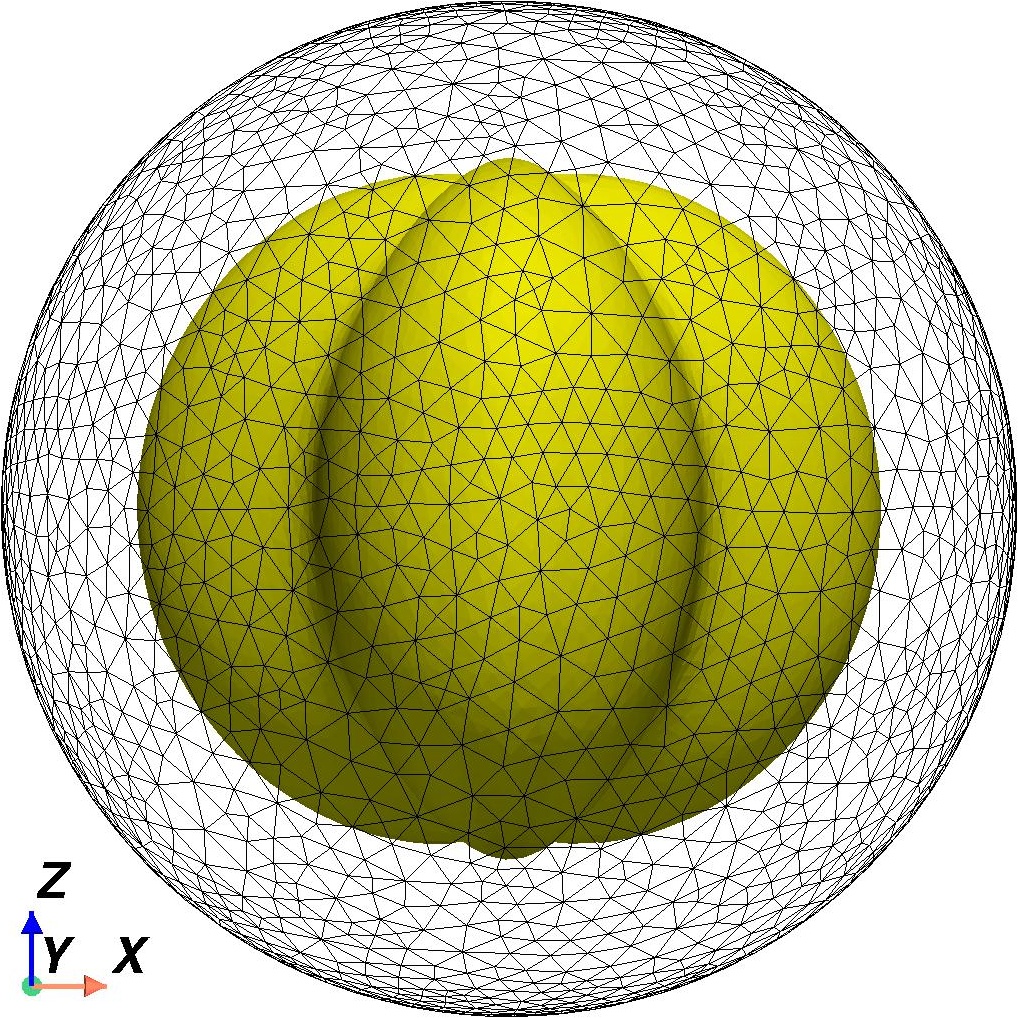}} \hfill
\resizebox{0.16\linewidth}{!}{\includegraphics{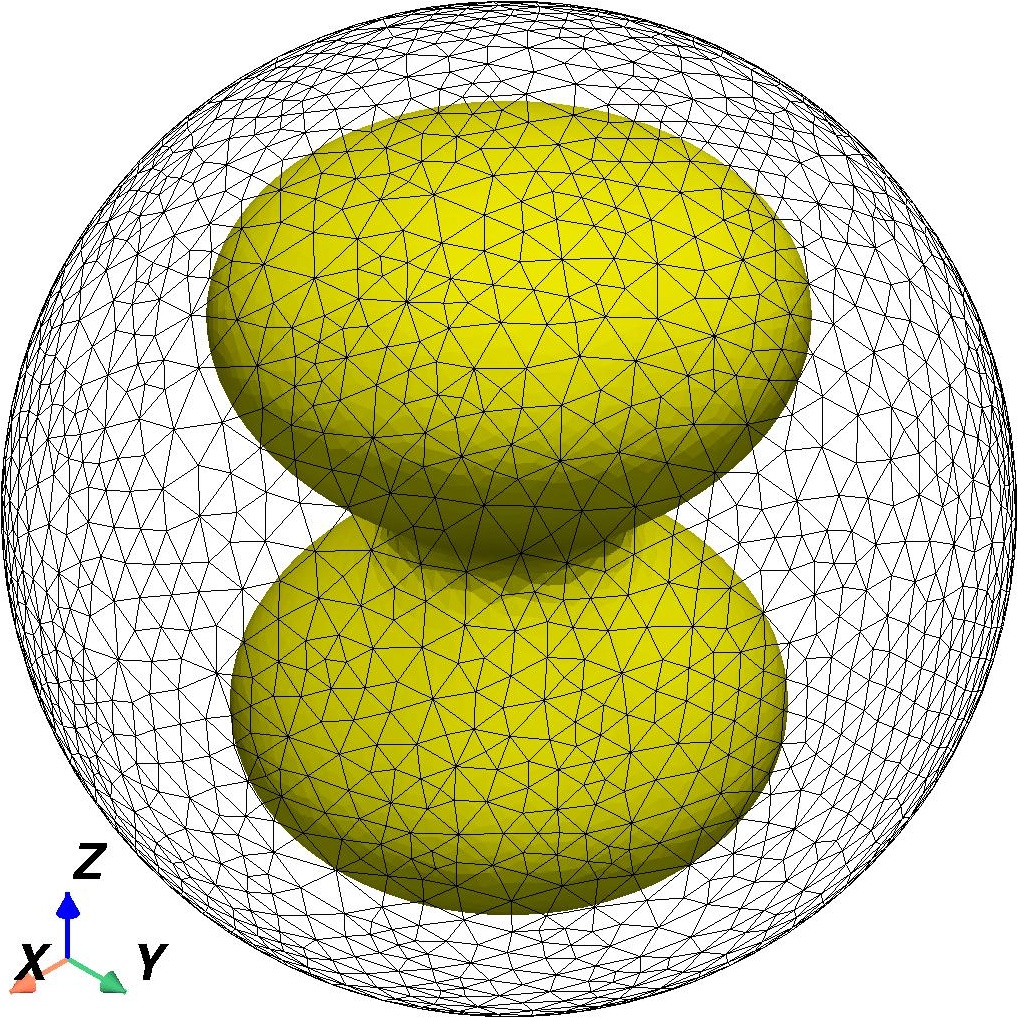}} \hfill
\resizebox{0.16\linewidth}{!}{\includegraphics{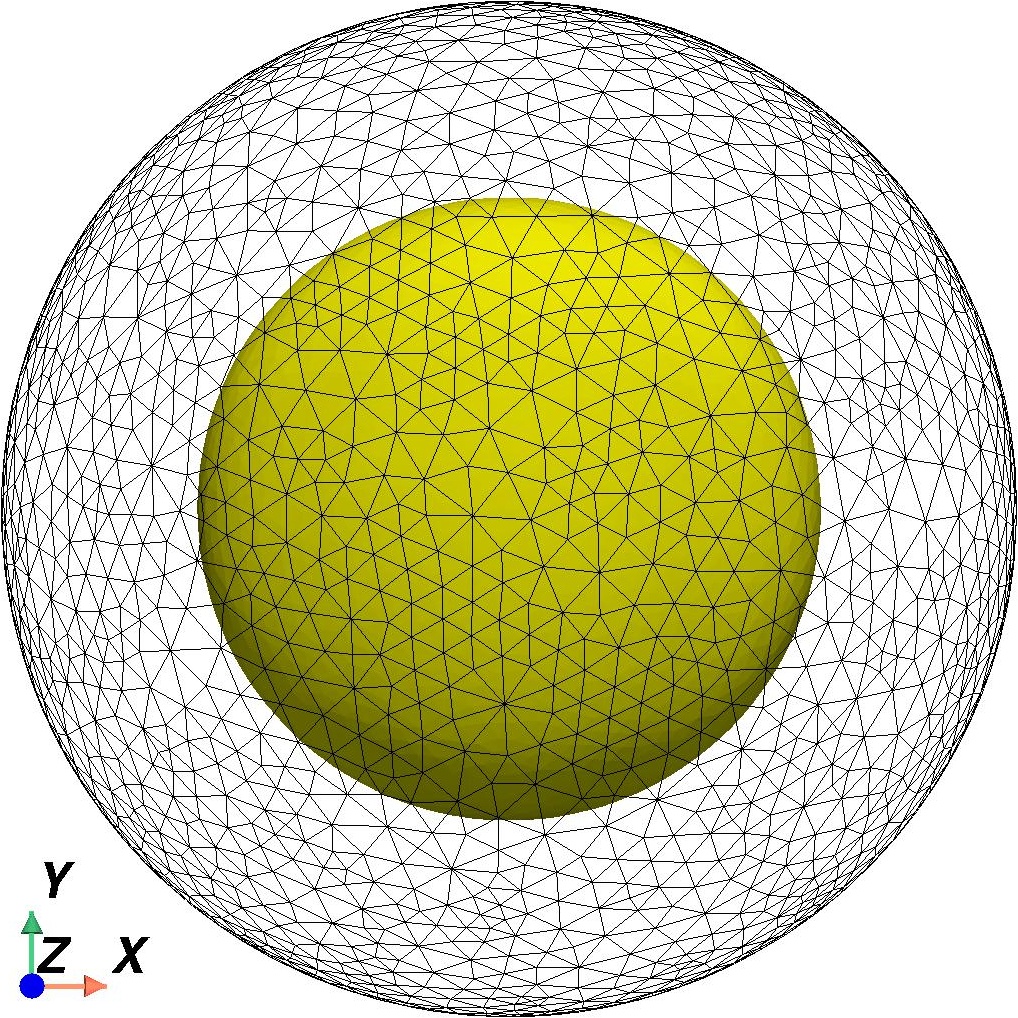}} \hfill
\resizebox{0.16\linewidth}{!}{\includegraphics{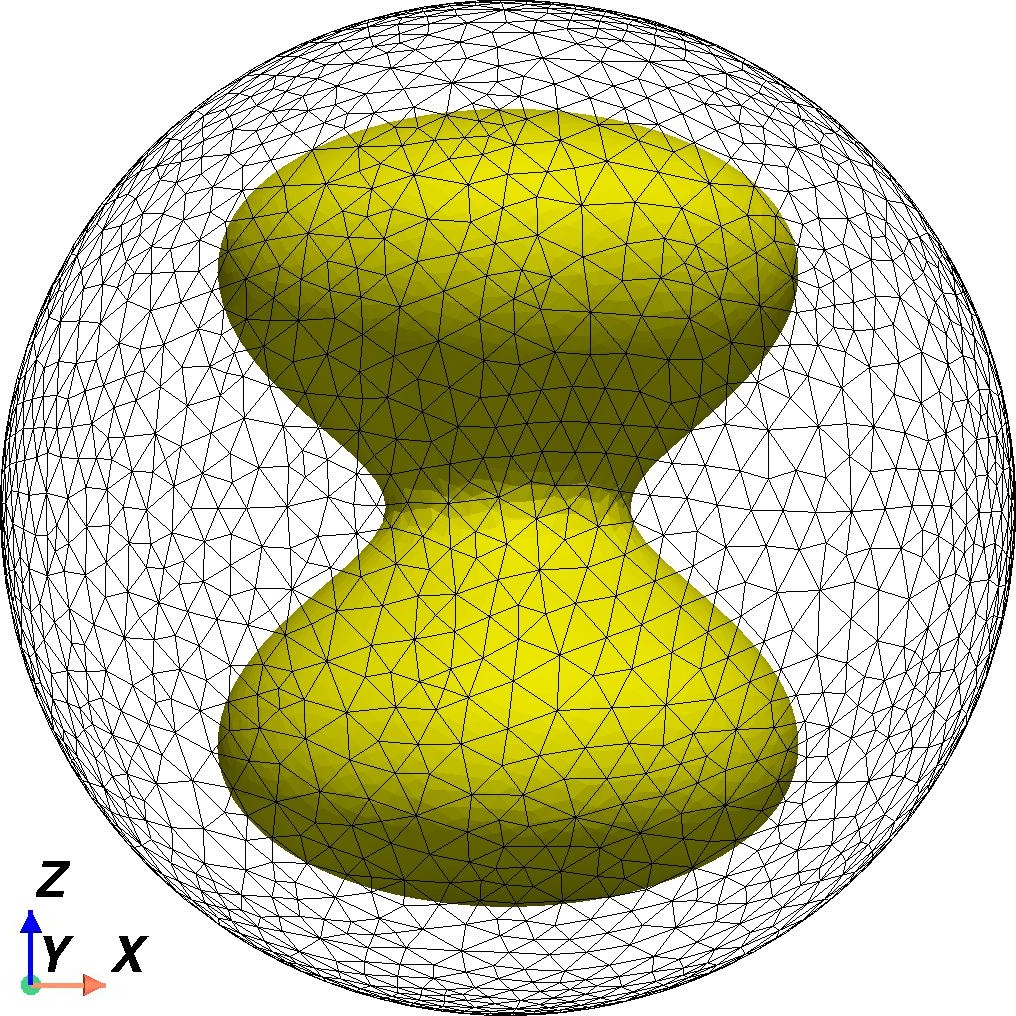}} 
 \\[1em] 
\resizebox{0.16\linewidth}{!}{\includegraphics{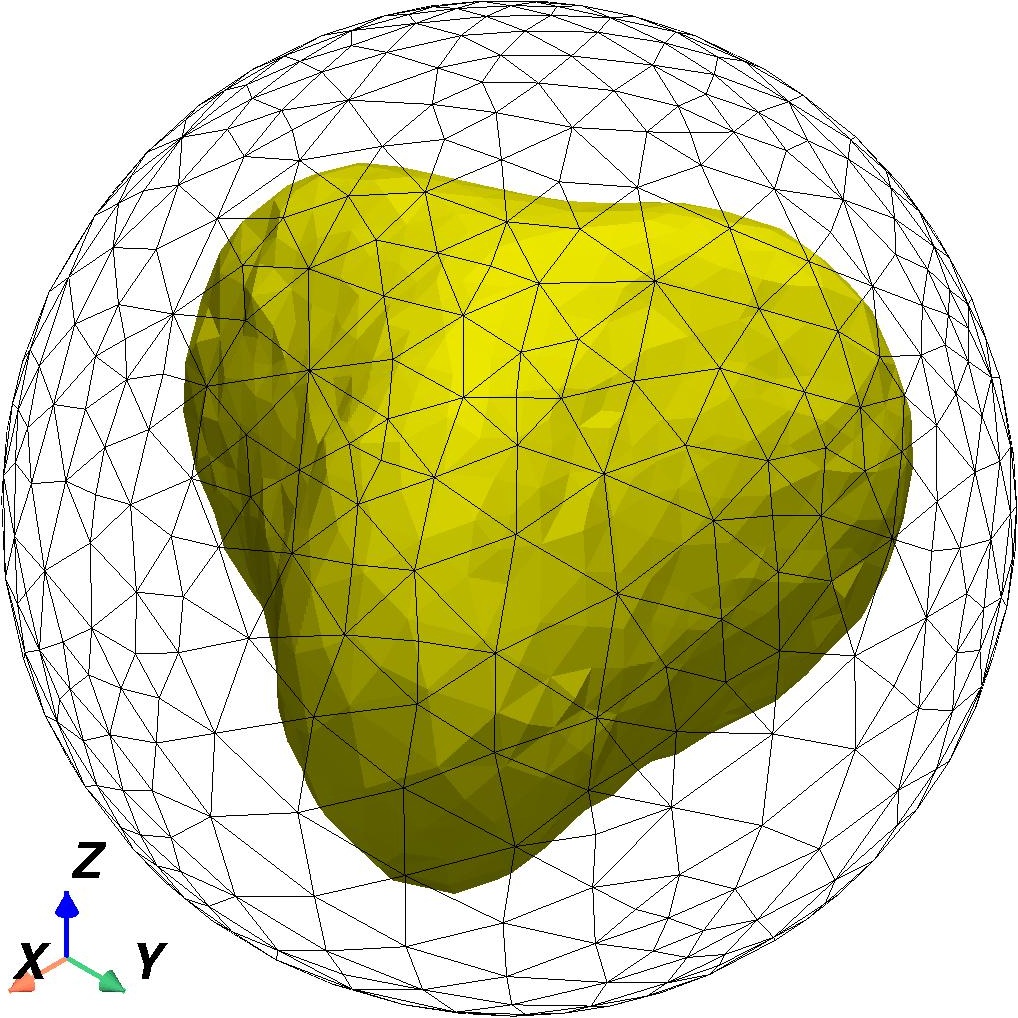}} \hfill
\resizebox{0.16\linewidth}{!}{\includegraphics{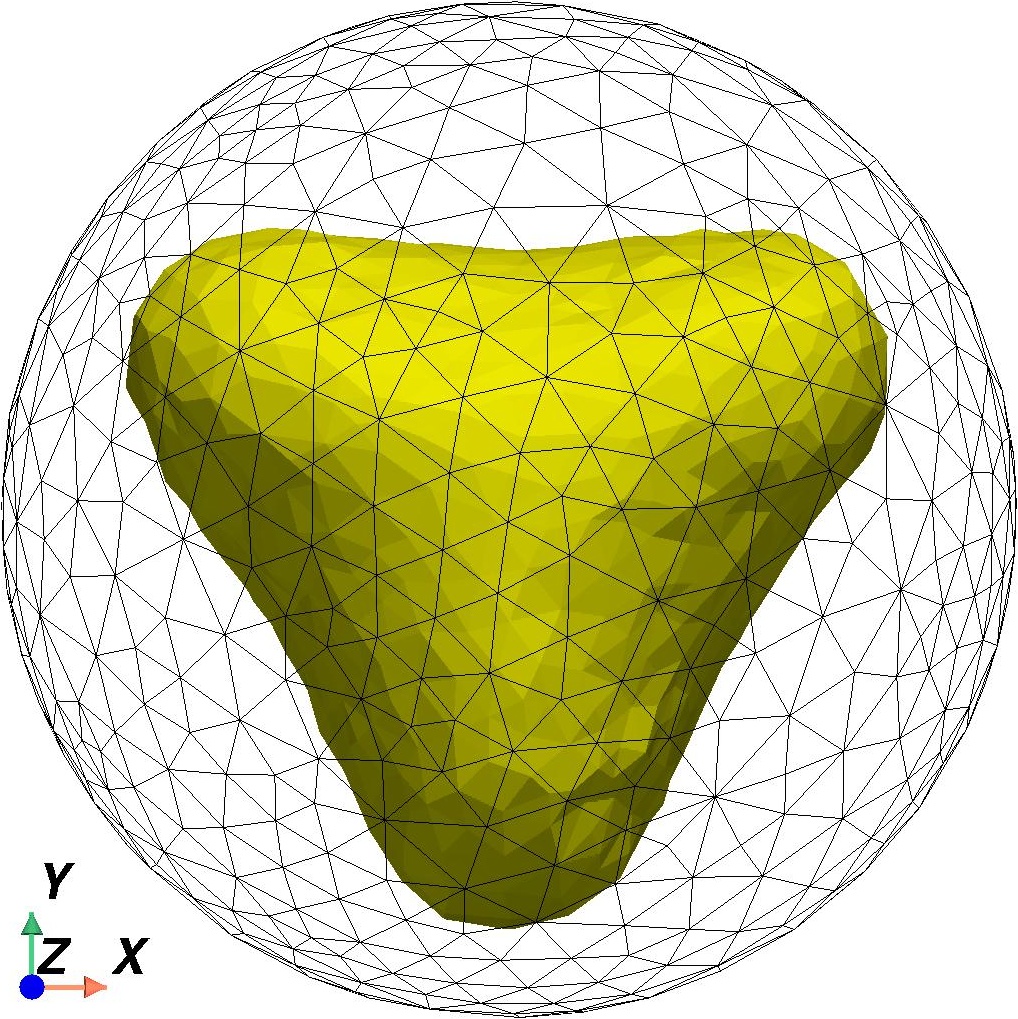}} \hfill
\resizebox{0.16\linewidth}{!}{\includegraphics{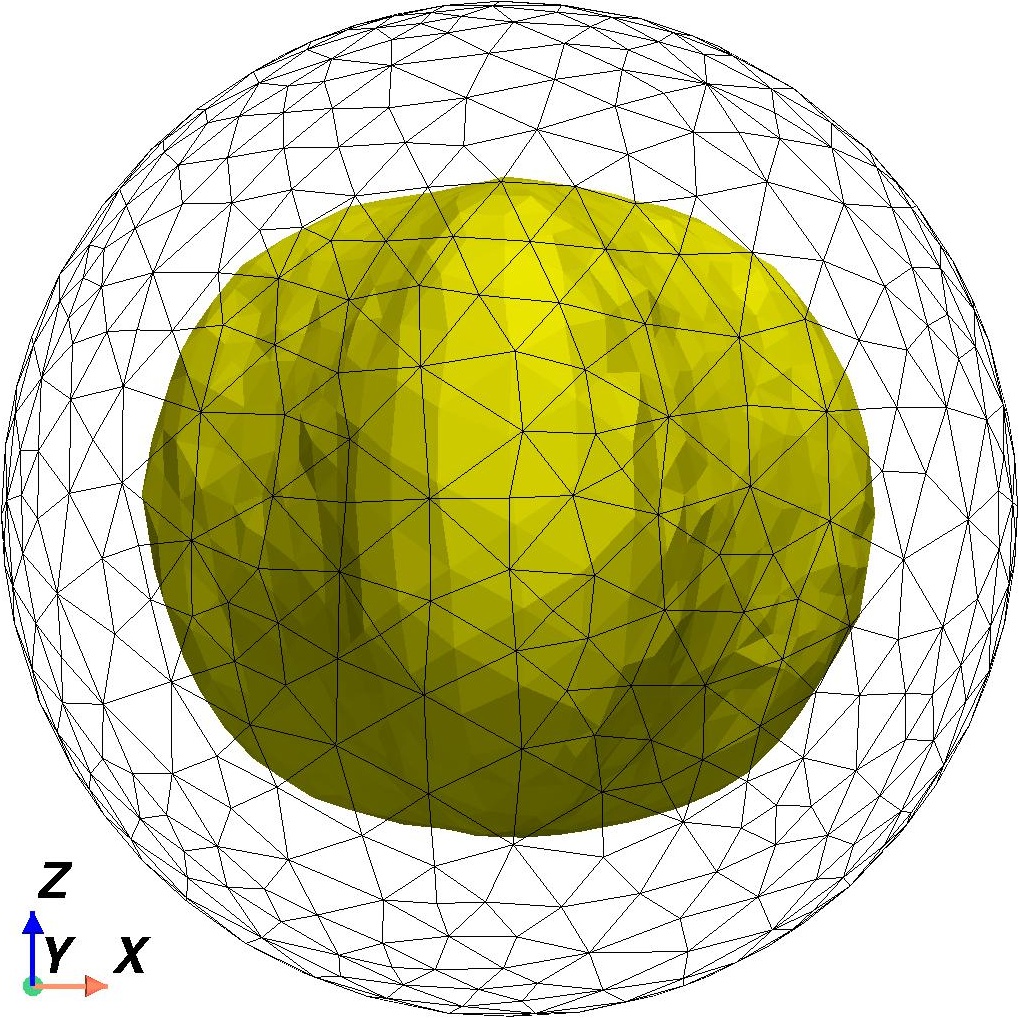}} \hfill
\resizebox{0.16\linewidth}{!}{\includegraphics{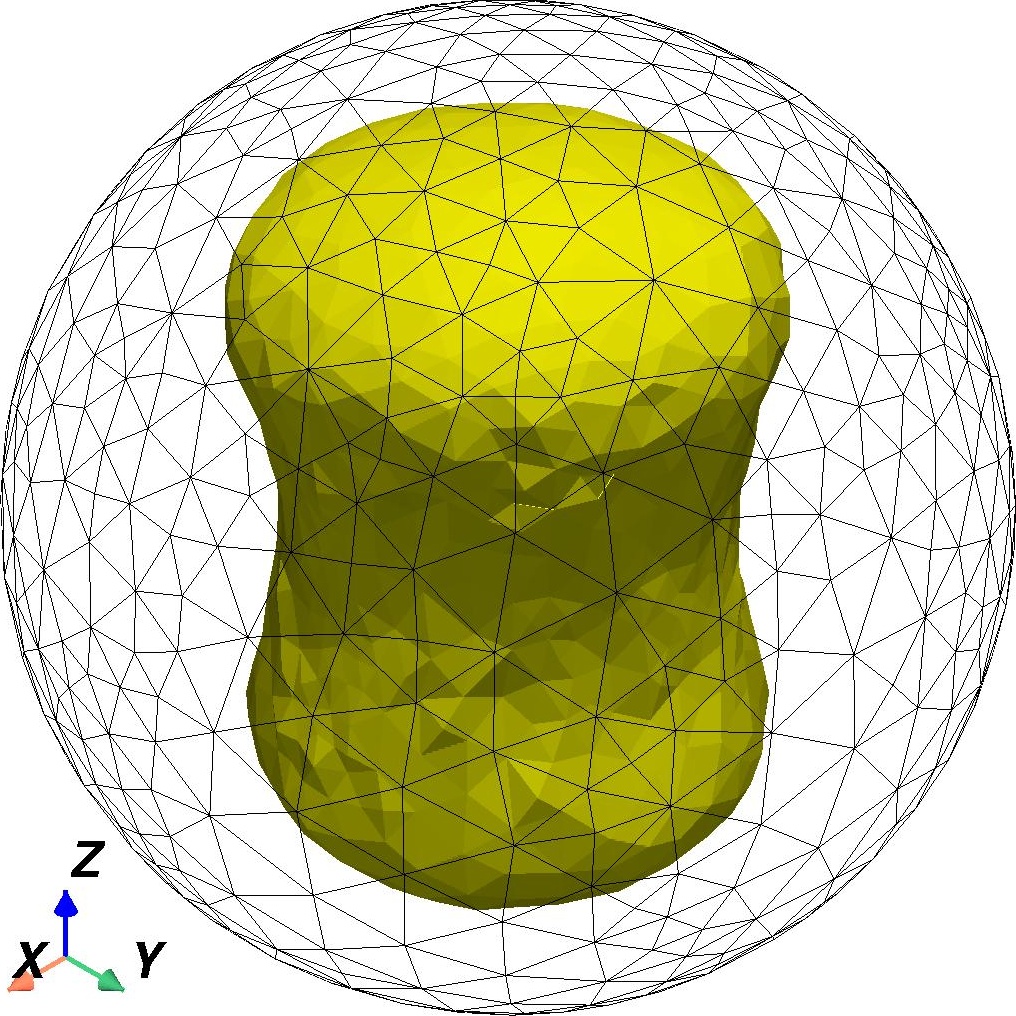}} \hfill
\resizebox{0.16\linewidth}{!}{\includegraphics{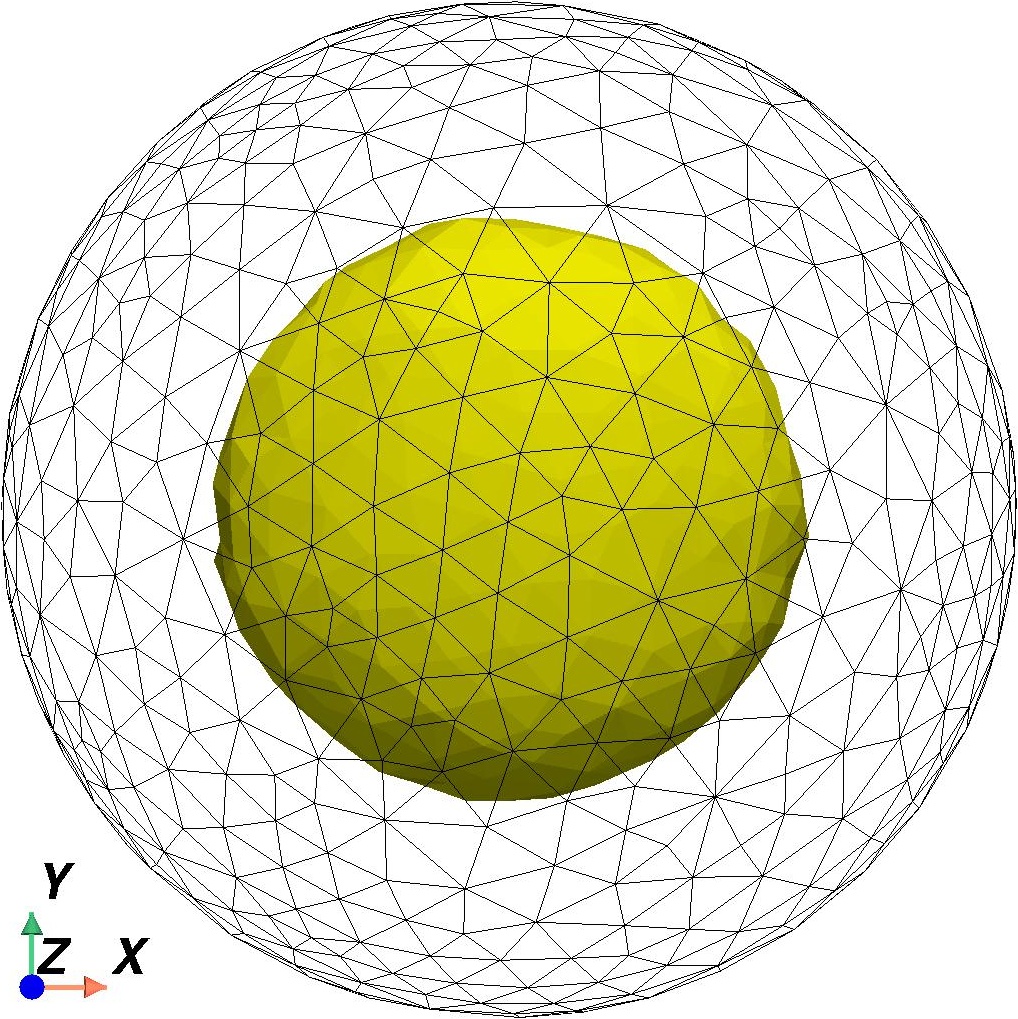}} \hfill
\resizebox{0.16\linewidth}{!}{\includegraphics{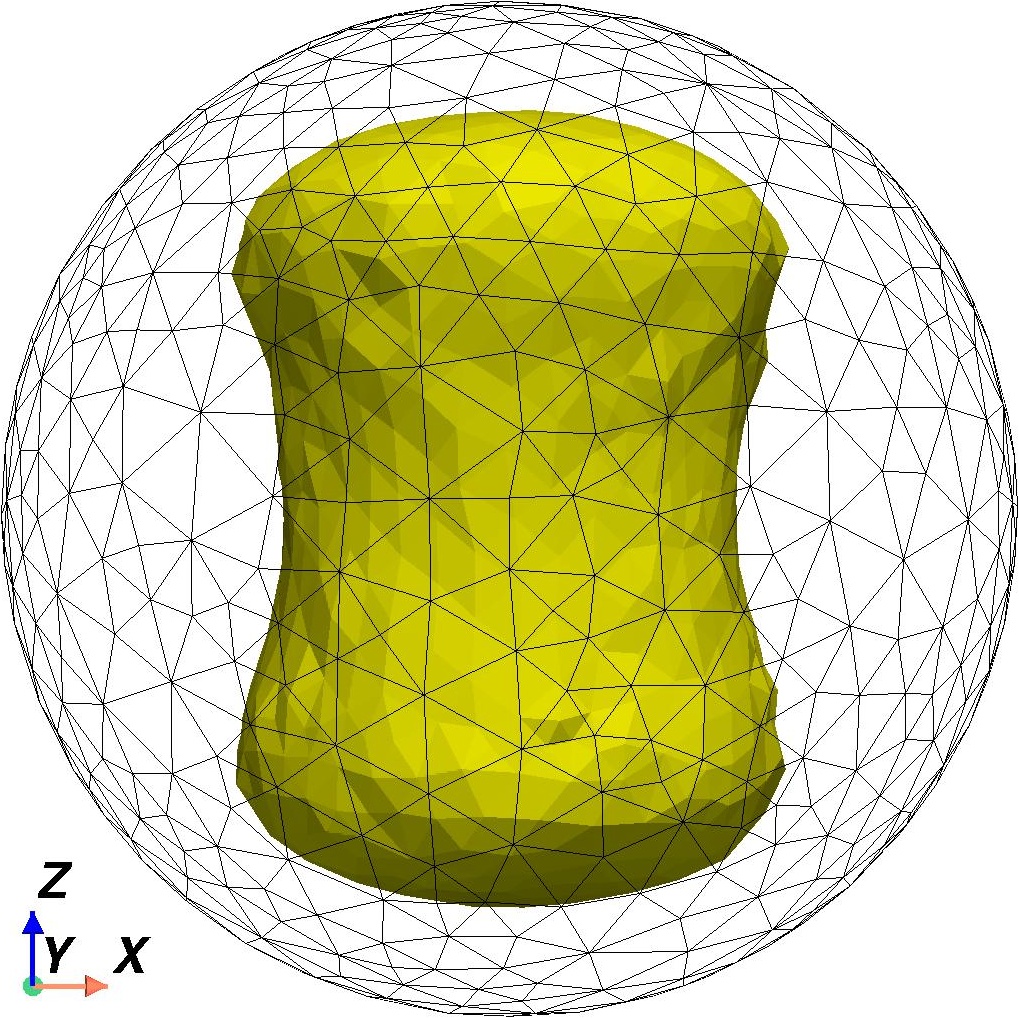}} 
 \\[1em] 
\resizebox{0.16\linewidth}{!}{\includegraphics{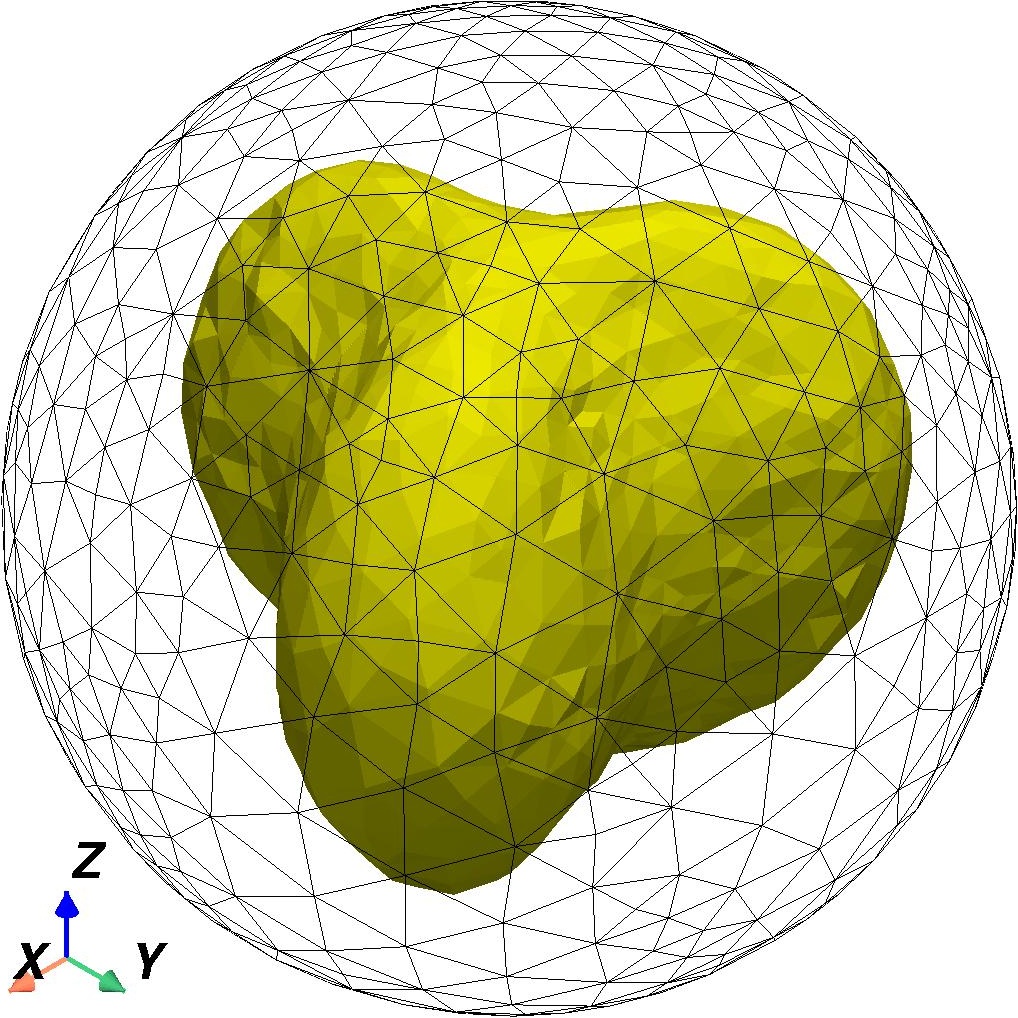}} \hfill
\resizebox{0.16\linewidth}{!}{\includegraphics{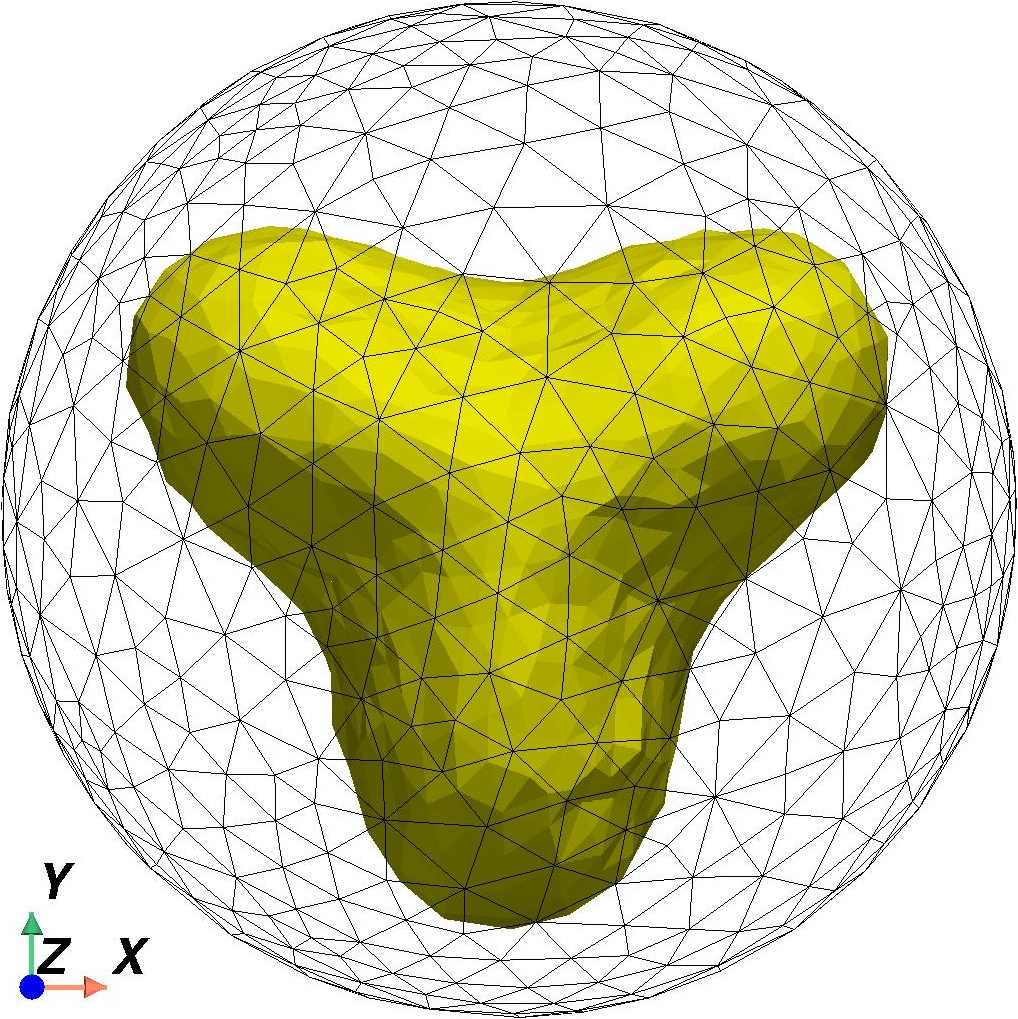}} \hfill
\resizebox{0.16\linewidth}{!}{\includegraphics{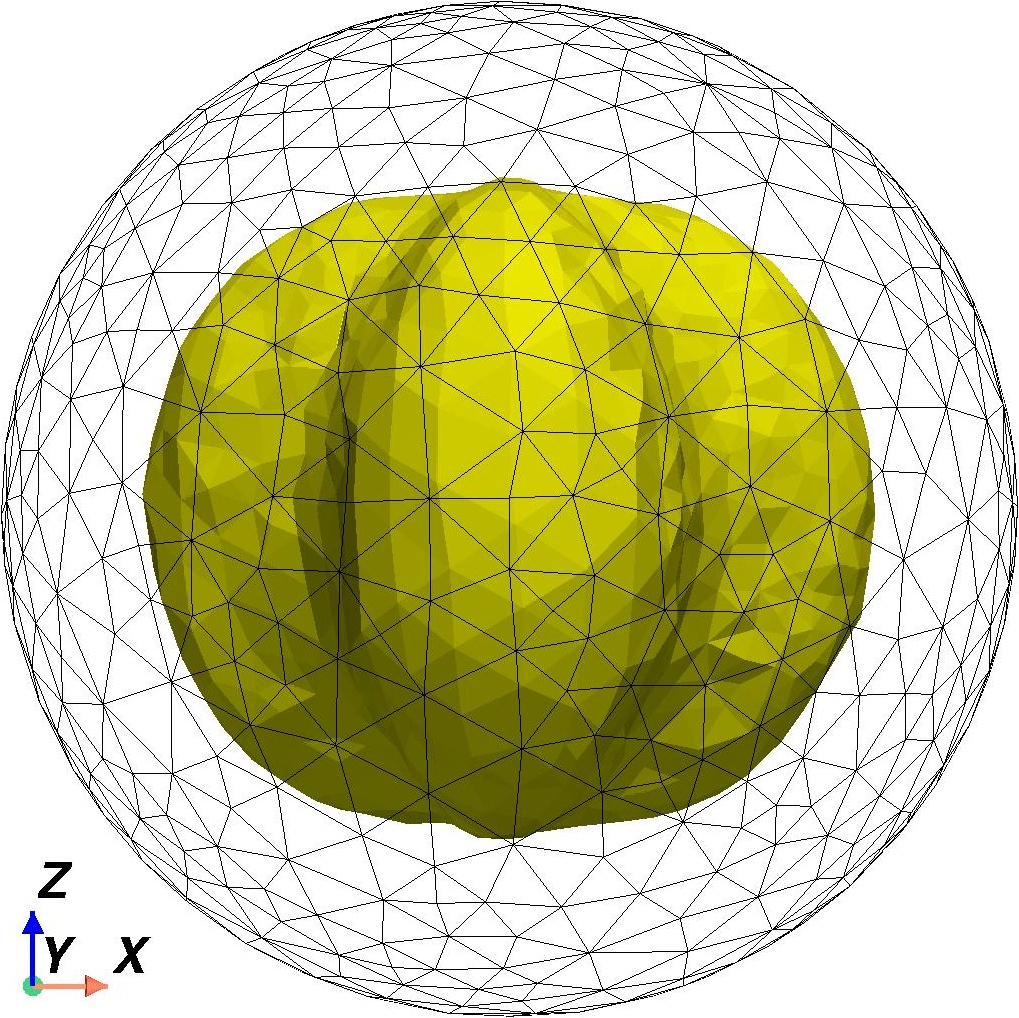}} \hfill
\resizebox{0.16\linewidth}{!}{\includegraphics{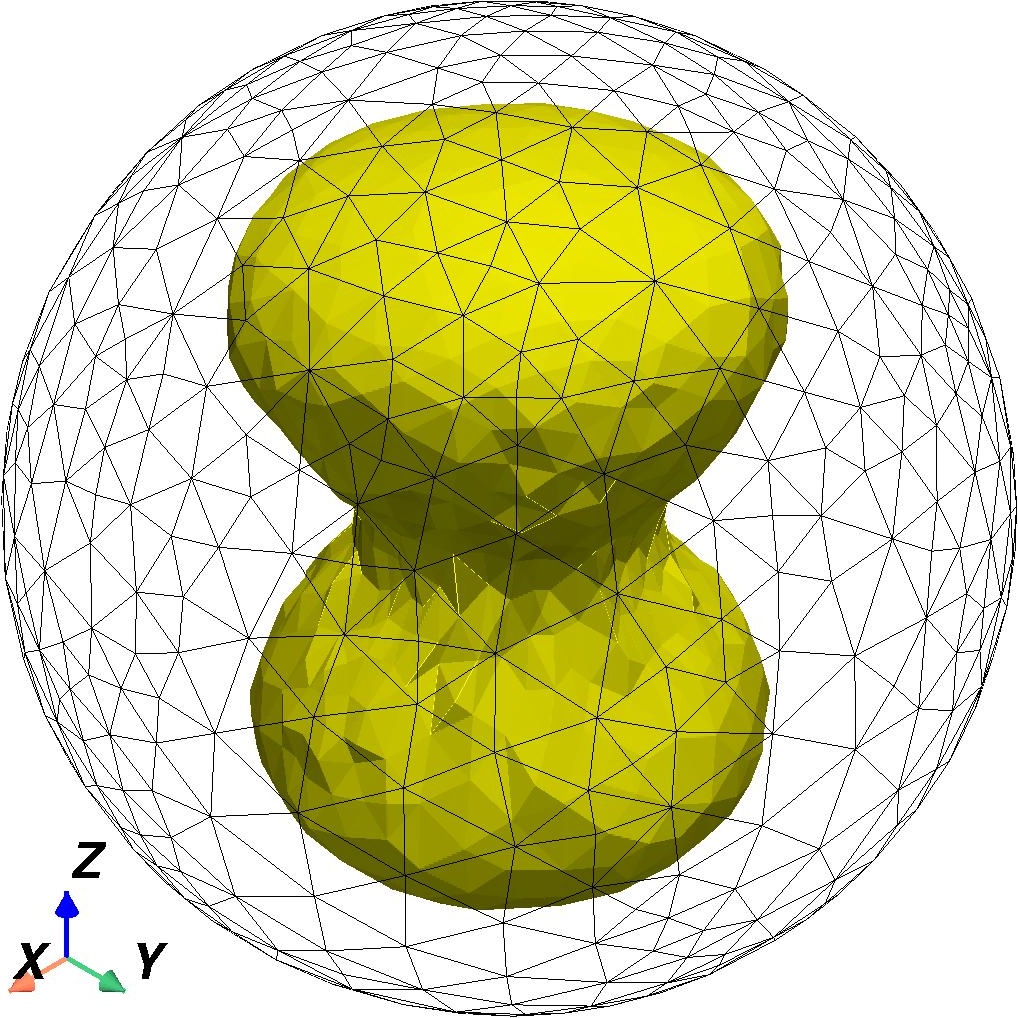}} \hfill
\resizebox{0.16\linewidth}{!}{\includegraphics{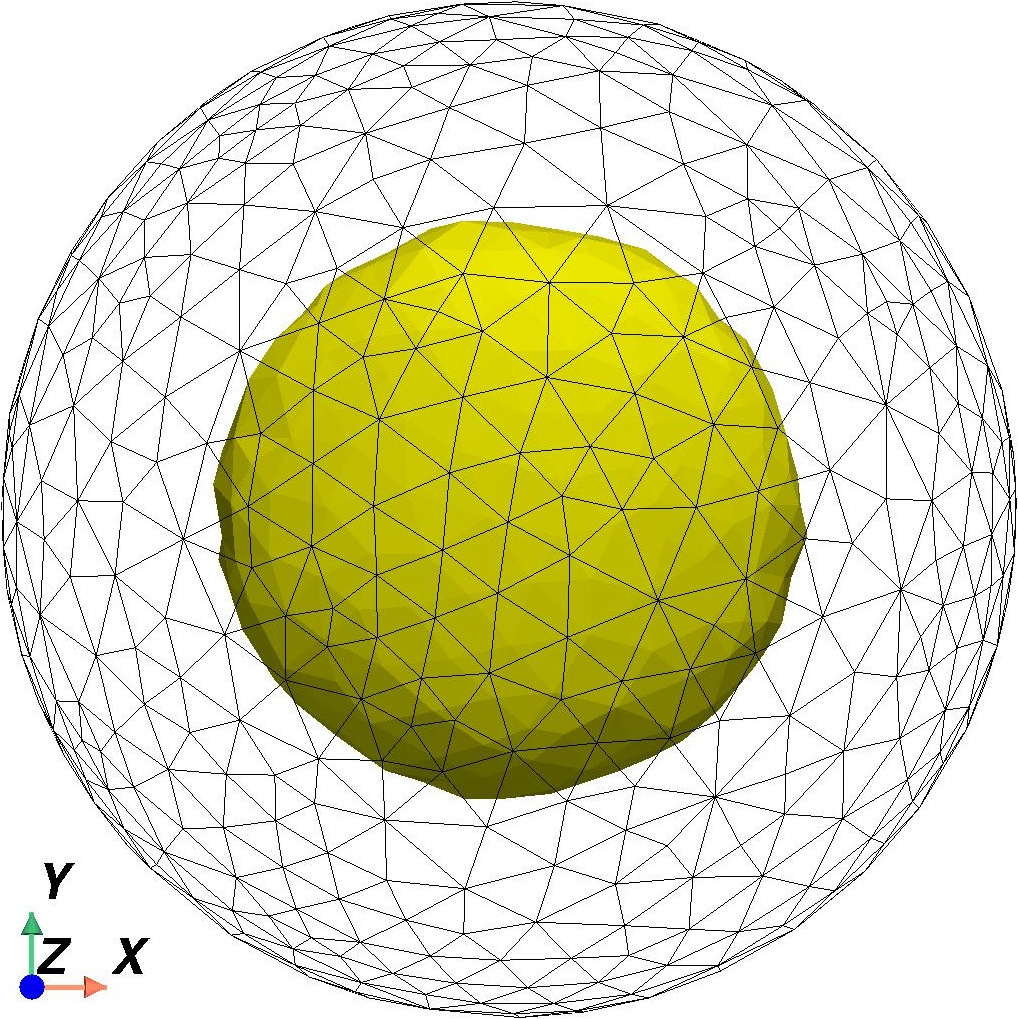}} \hfill
\resizebox{0.16\linewidth}{!}{\includegraphics{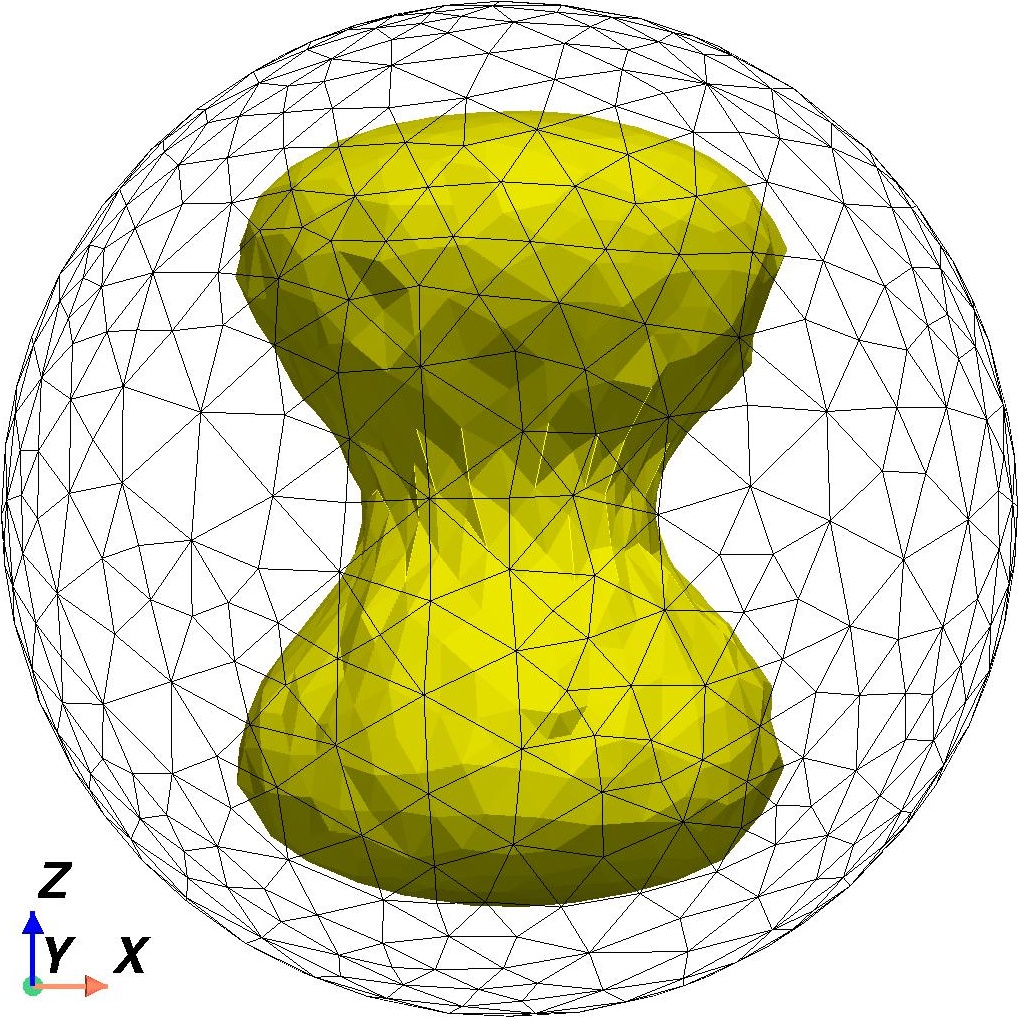}}
\caption{Exact geometries of the cavities (top/first row) and reconstructed shapes obtained via SO (middlesecond row) and ADMM (bottom/third row) with exact data}
\label{fig:figure2a}
\end{figure}

For the exact cavity, we examine two test cases featuring shapes that are strictly non-convex: a flower-like cavity and a dumbbell-like cavity (refer to the first row of plots in Figure \ref{fig:figure2a}). 
Additionally, we explore a case with a smaller cavity size and more pronounced concavities (see Figure \ref{fig:figure3a}). 
However, in this test case, we set $N = 500$ for the stopping condition.
For the flower-like cavity, we choose $\beta = 0.3$, while for the other cavity, we choose $\beta = 0.1$. 
These choices, based on our experience, provide the best reconstructions of the shapes.

For the forward problem, we discretize the exact domain with minimum and maximum mesh sizes $h_{\min}^{\ast} = 0.05$ and $h_{\max}^{\ast} = 0.1$ (see the first row of plots in Figure \ref{fig:figure2a}), respectively with tetrahedrons of (maximum) volume $0.001$.
For the inversion procedure, we discretize $(\Omega\setminus\overline{\omega})^{0}$ using a coarse mesh with $h_{\min} = 0.15$ and $h_{\max} = 0.2$, with tetrahedrons of volume $0.005$. 

\begin{figure}[htp!]
\centering
\resizebox{0.16\linewidth}{!}{\includegraphics{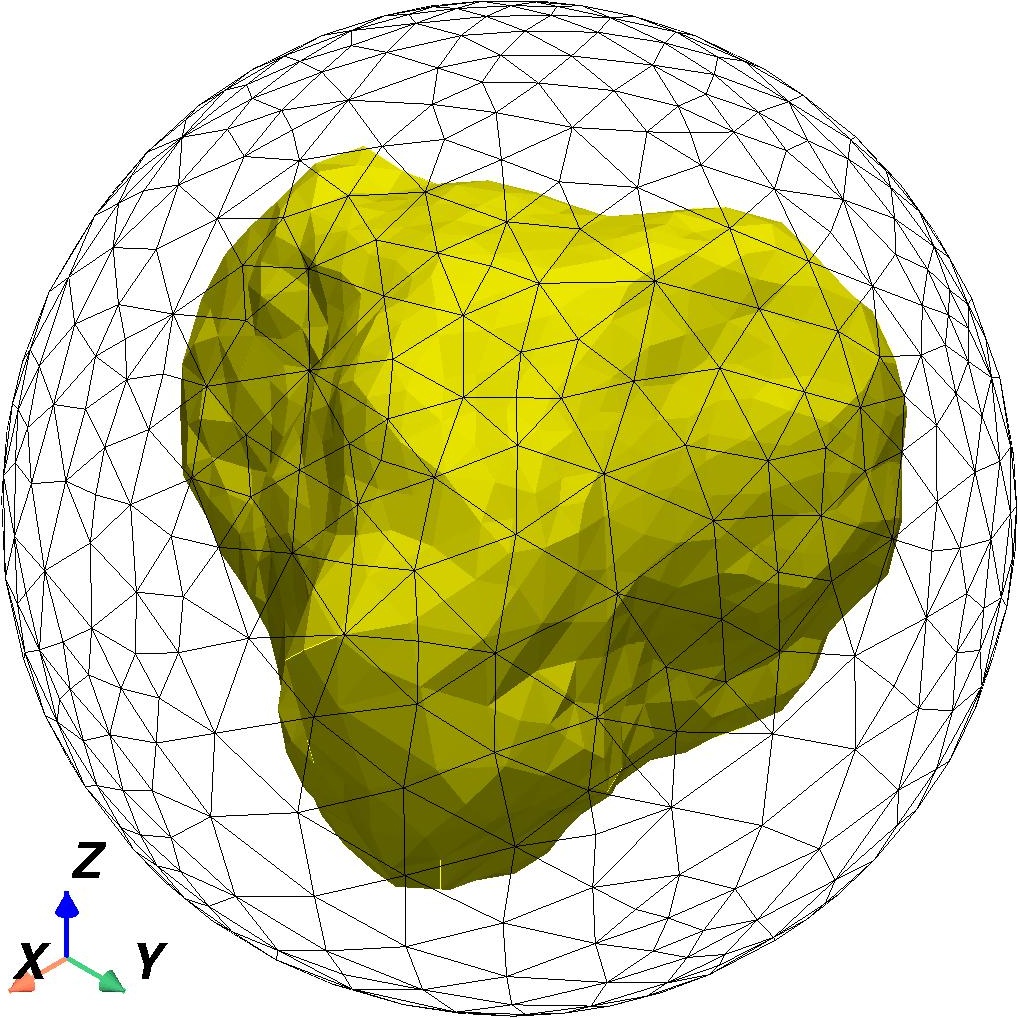}} \hfill
\resizebox{0.16\linewidth}{!}{\includegraphics{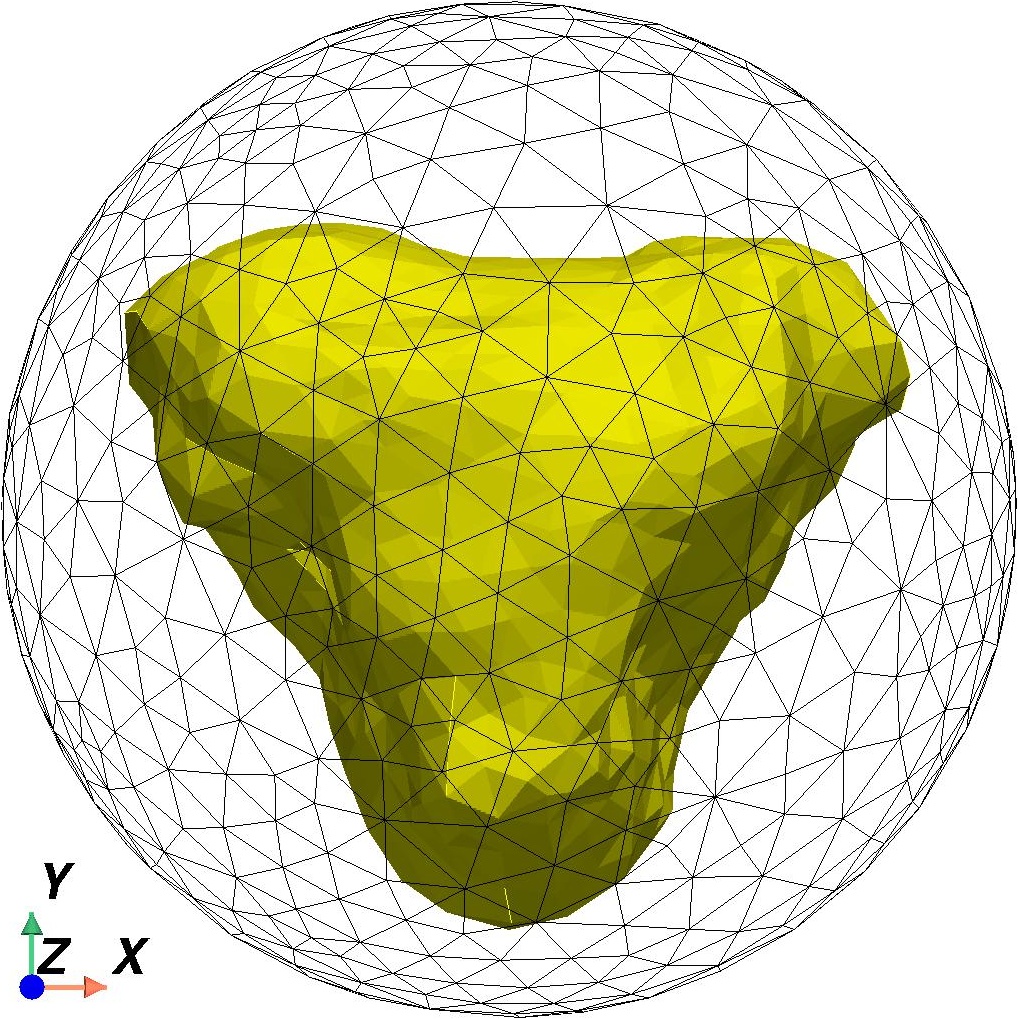}} \hfill
\resizebox{0.16\linewidth}{!}{\includegraphics{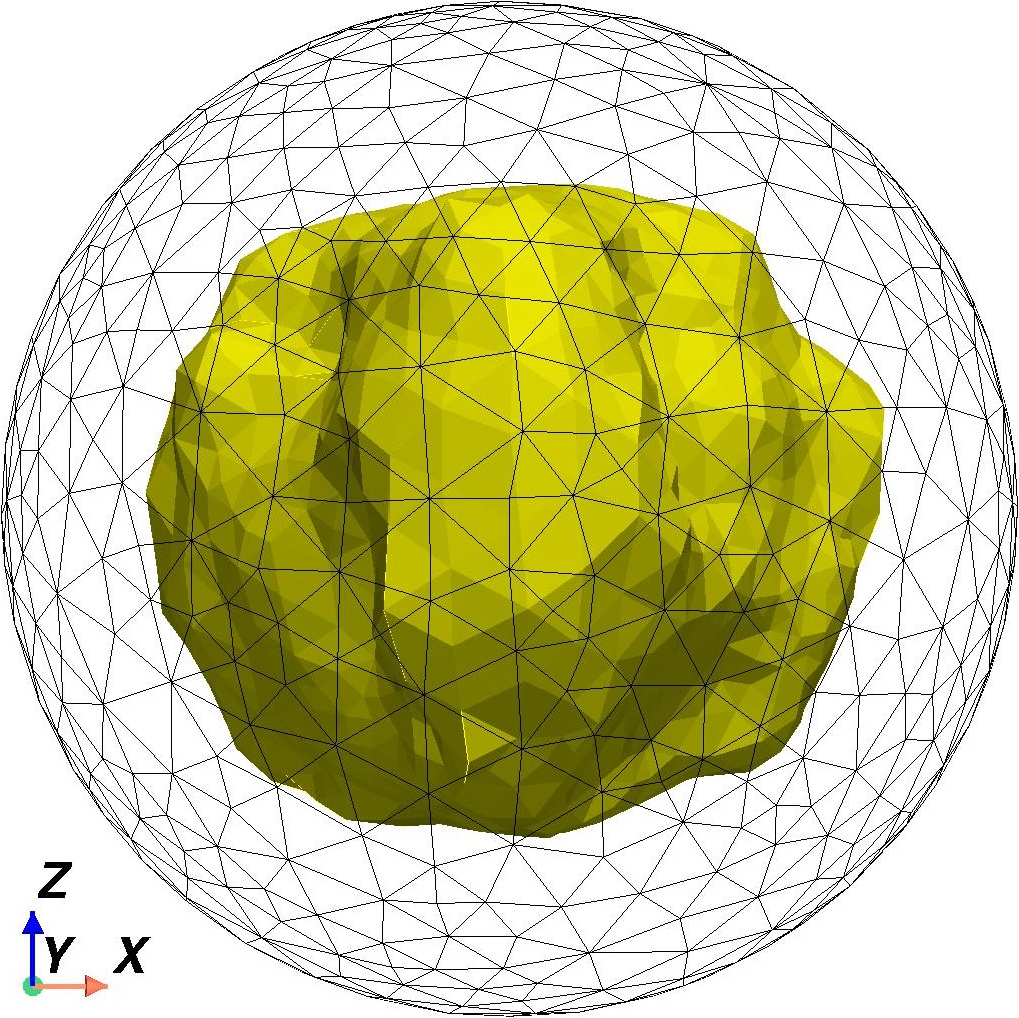}} \hfill
\resizebox{0.16\linewidth}{!}{\includegraphics{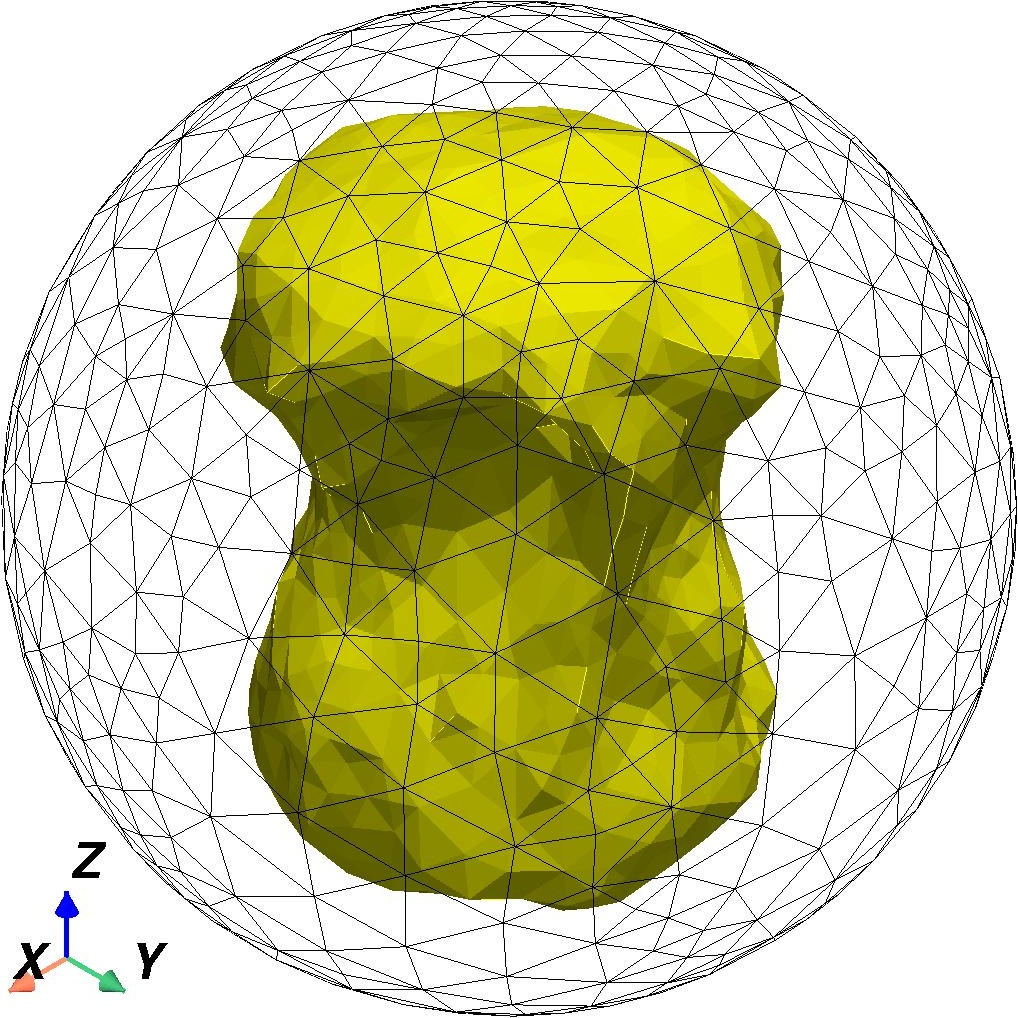}} \hfill
\resizebox{0.16\linewidth}{!}{\includegraphics{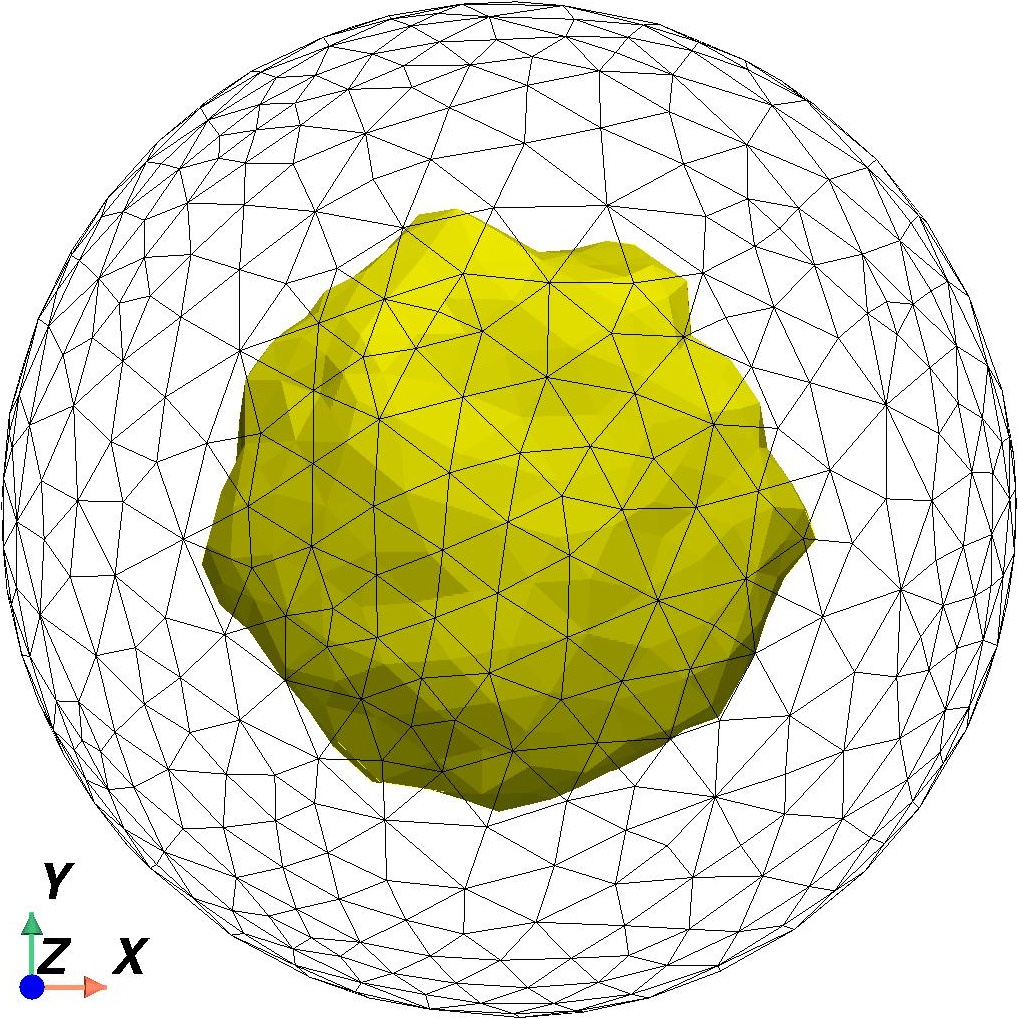}} \hfill
\resizebox{0.16\linewidth}{!}{\includegraphics{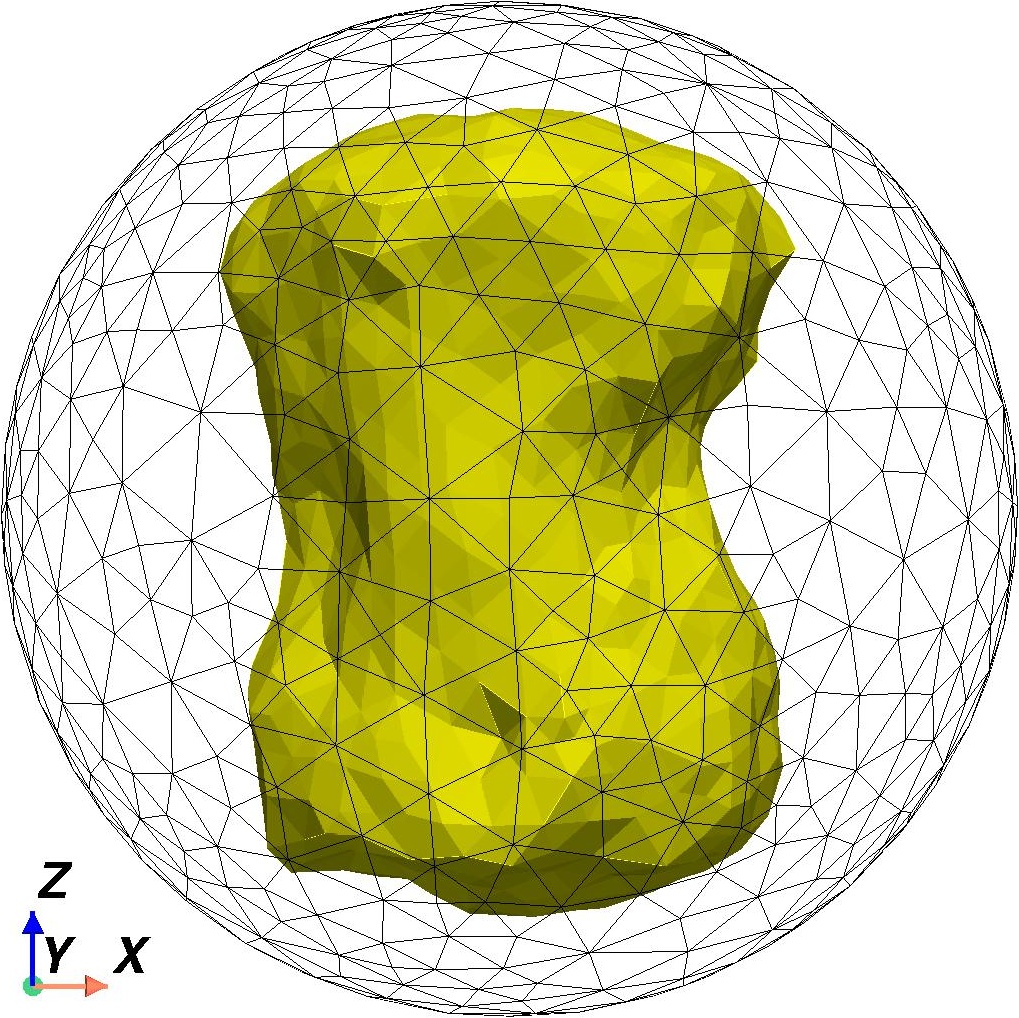}} 
 \\[1em] 
 \resizebox{0.16\linewidth}{!}{\includegraphics{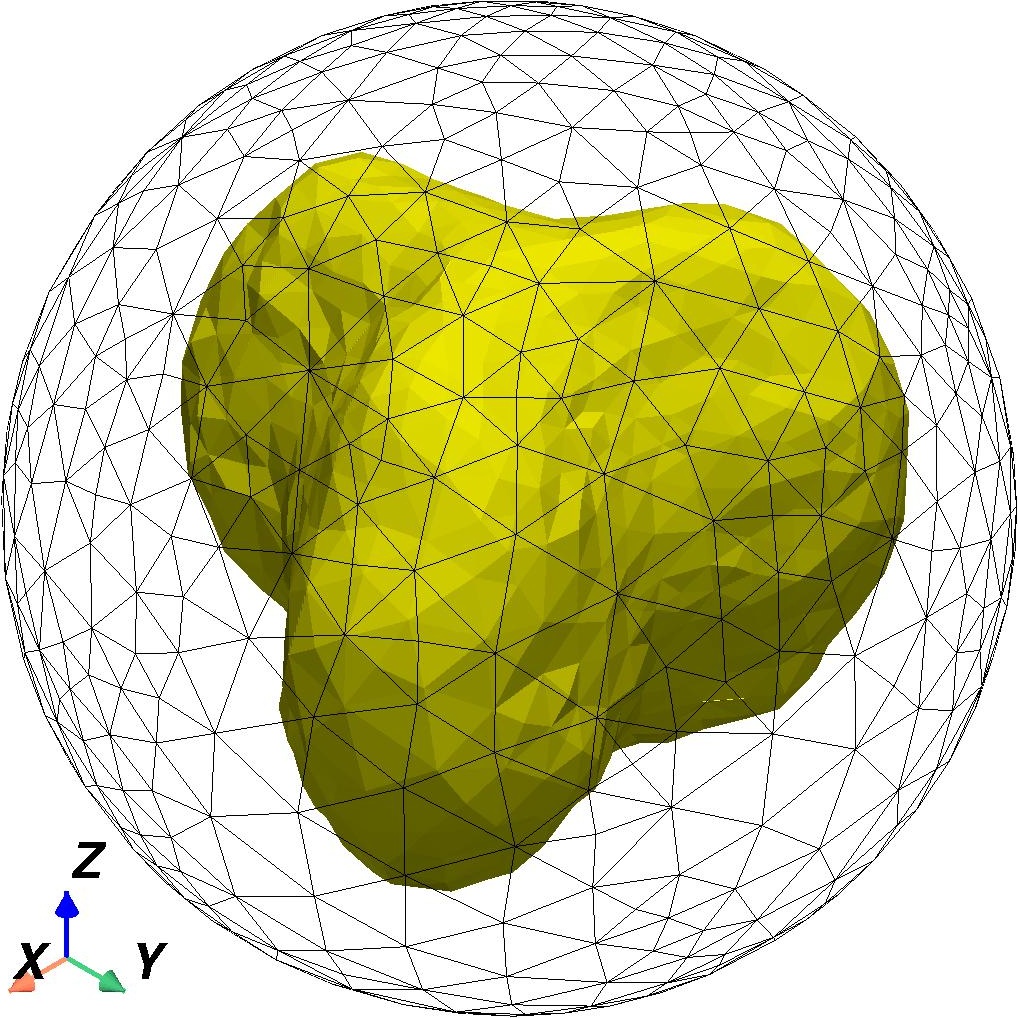}} \hfill
\resizebox{0.16\linewidth}{!}{\includegraphics{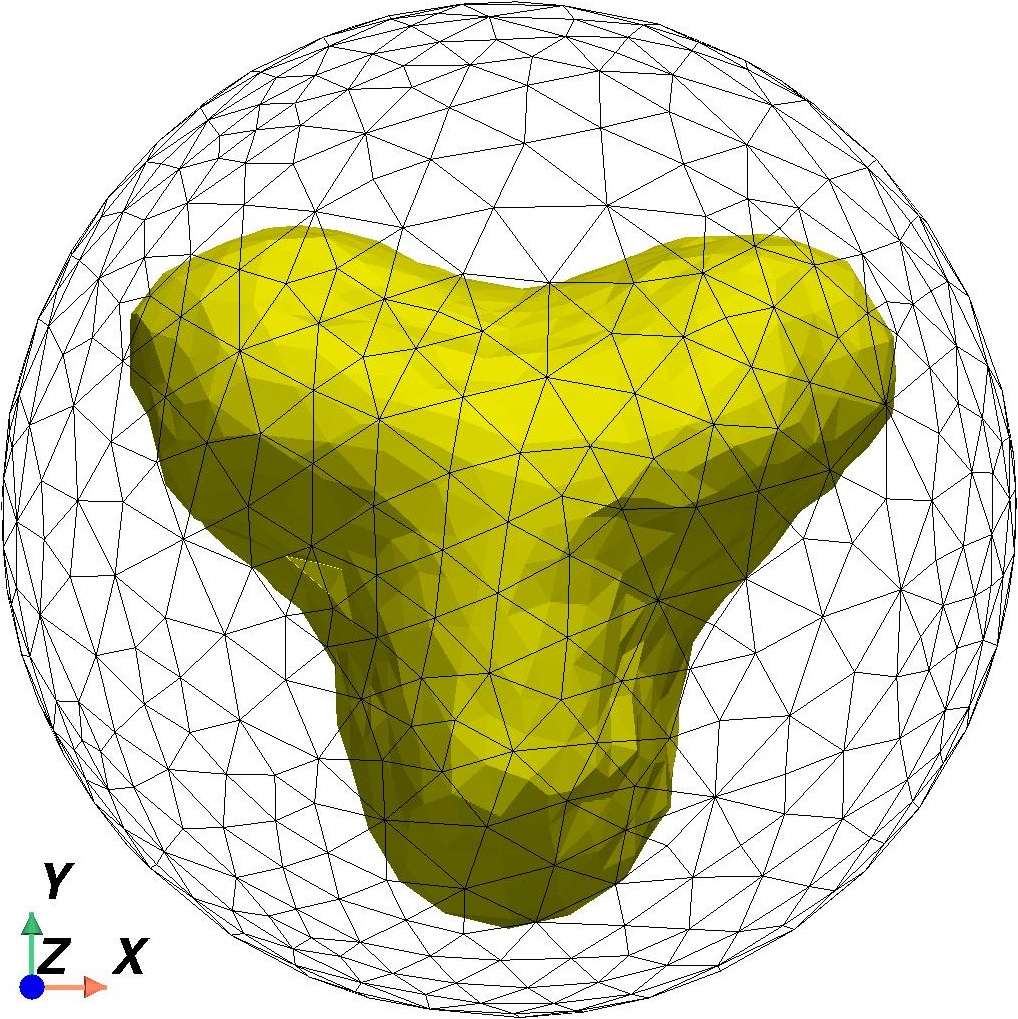}} \hfill
\resizebox{0.16\linewidth}{!}{\includegraphics{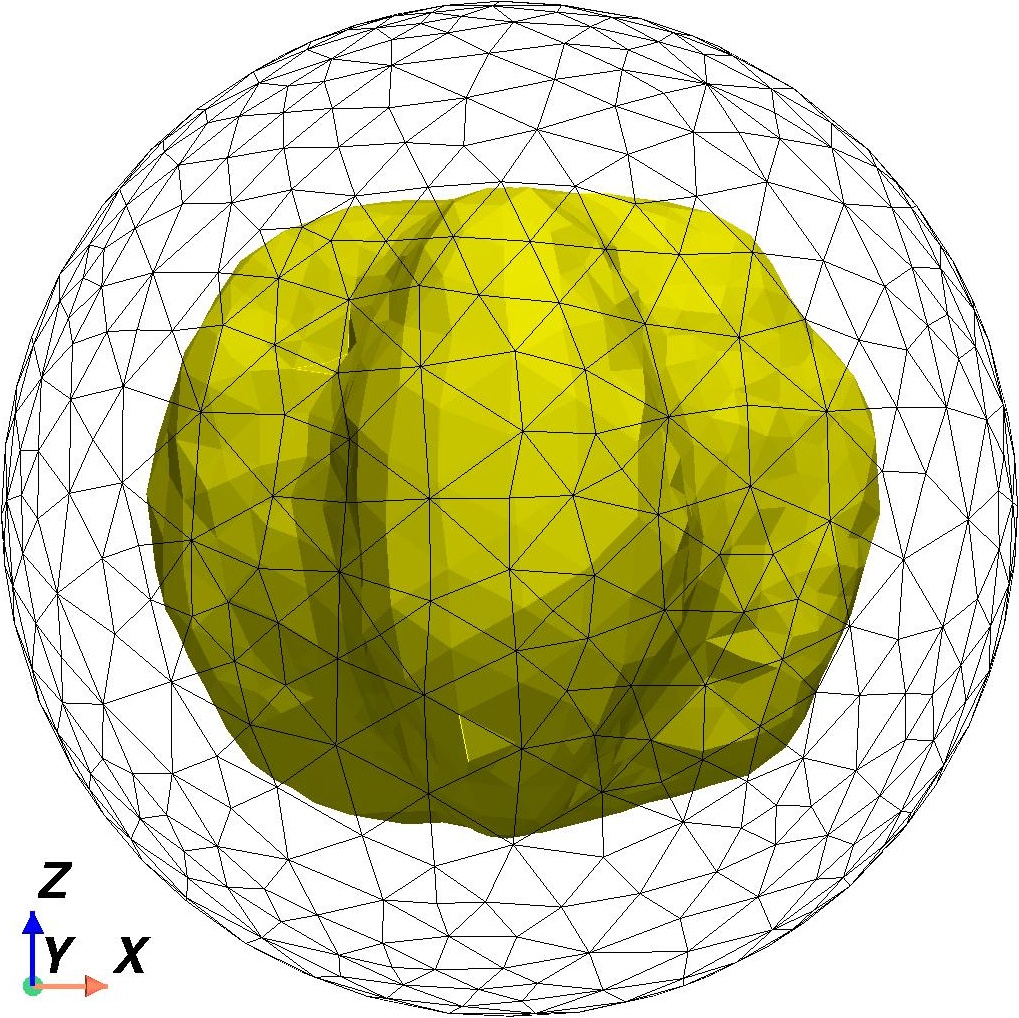}} \hfill
\resizebox{0.16\linewidth}{!}{\includegraphics{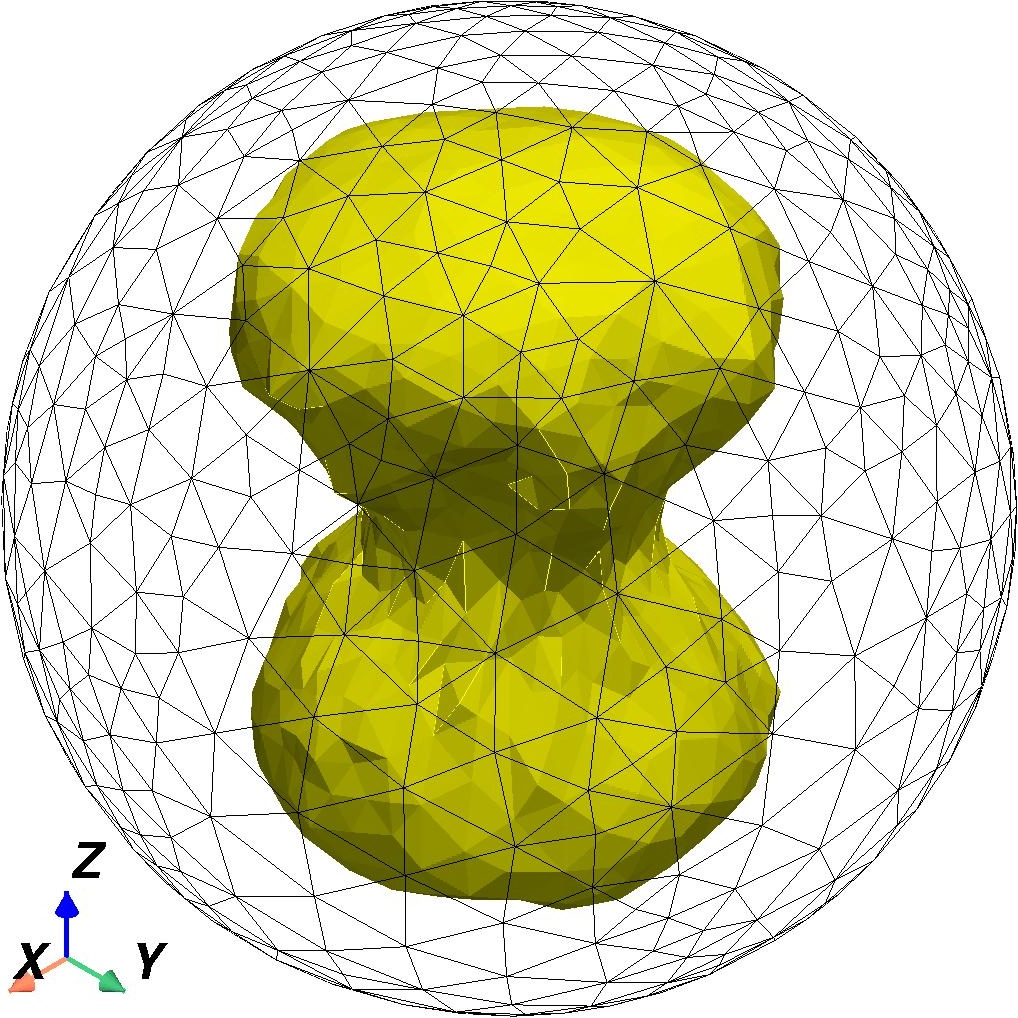}} \hfill
\resizebox{0.16\linewidth}{!}{\includegraphics{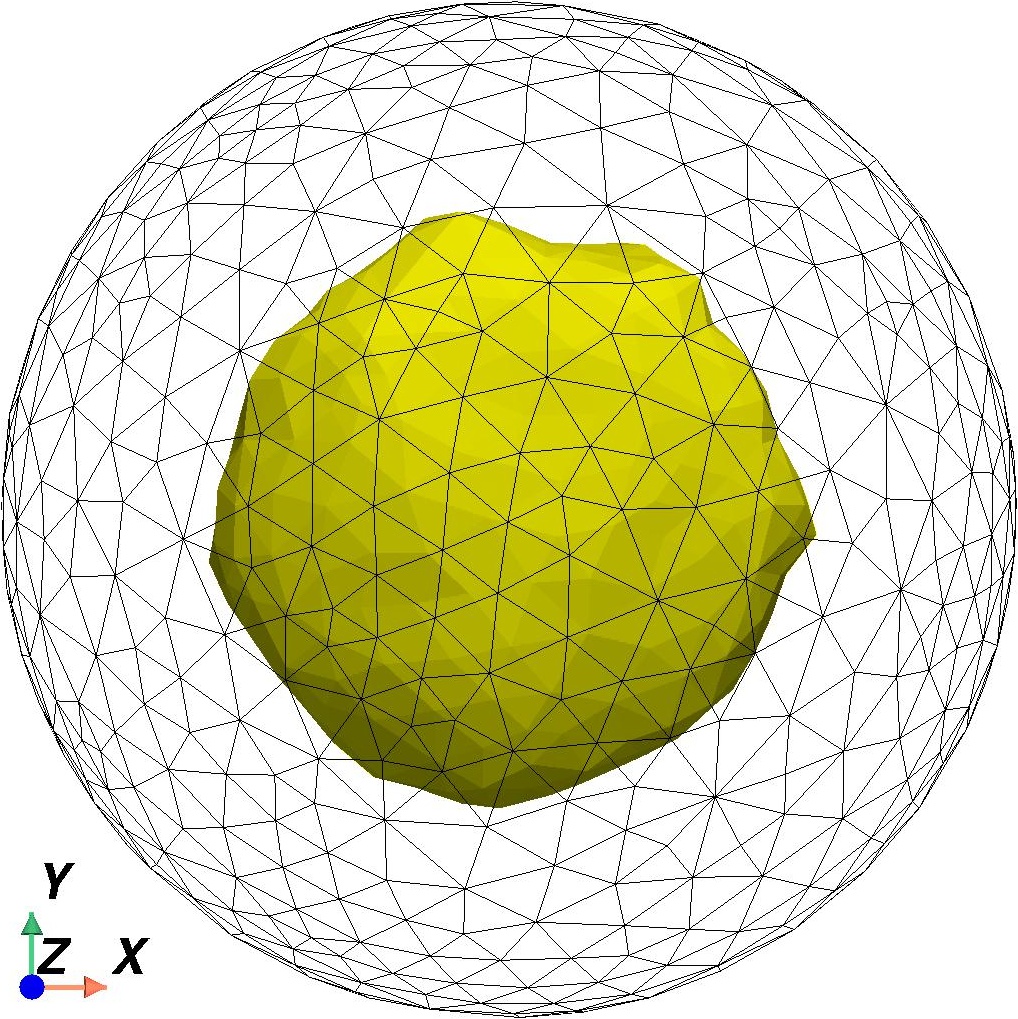}} \hfill
\resizebox{0.16\linewidth}{!}{\includegraphics{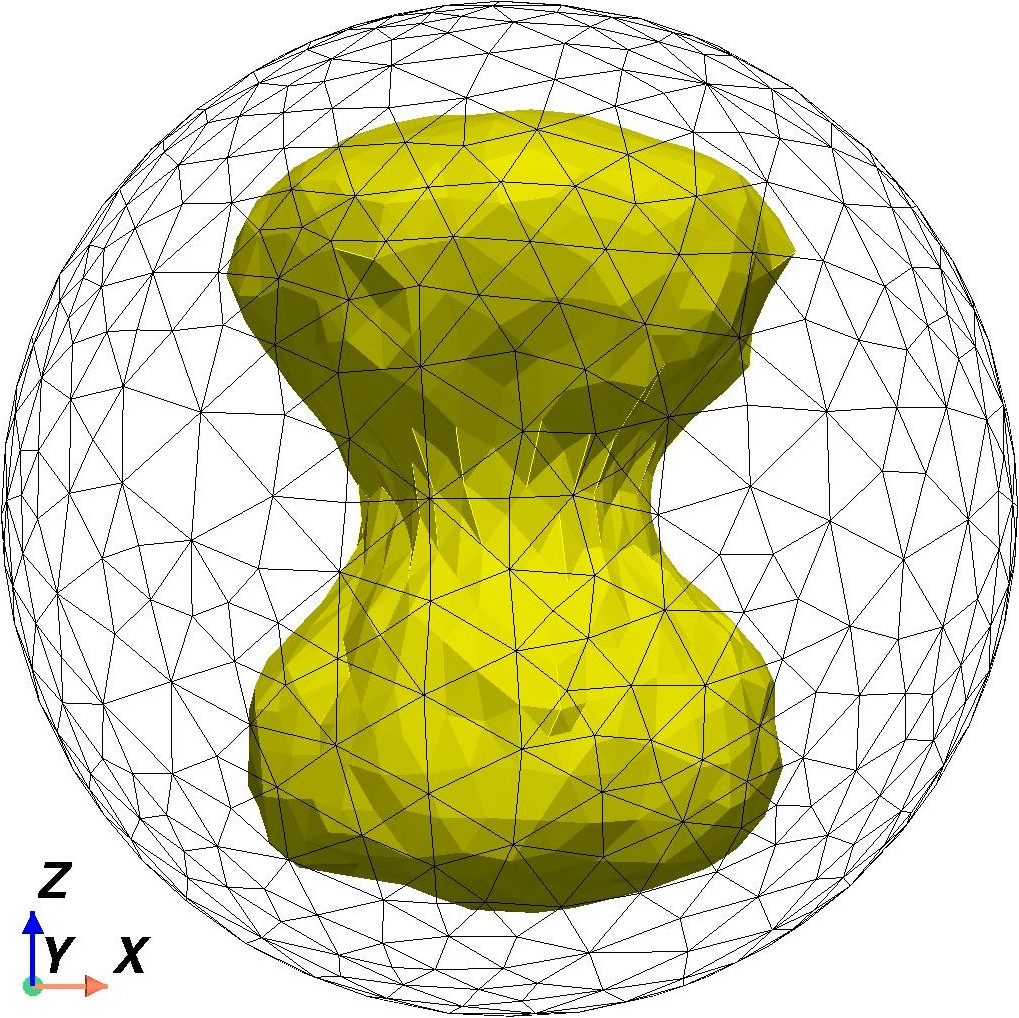}}
\caption{Reconstructions via SO (top/first row) and via ADMM (bottom/second row) with noised data at a noise level of $\delta = 15\%$ without regularization (i.e., $\gamma = 0$)}
\label{fig:figure2b}
\end{figure}

\begin{figure}[htp!]
\centering
 \resizebox{0.16\linewidth}{!}{\includegraphics{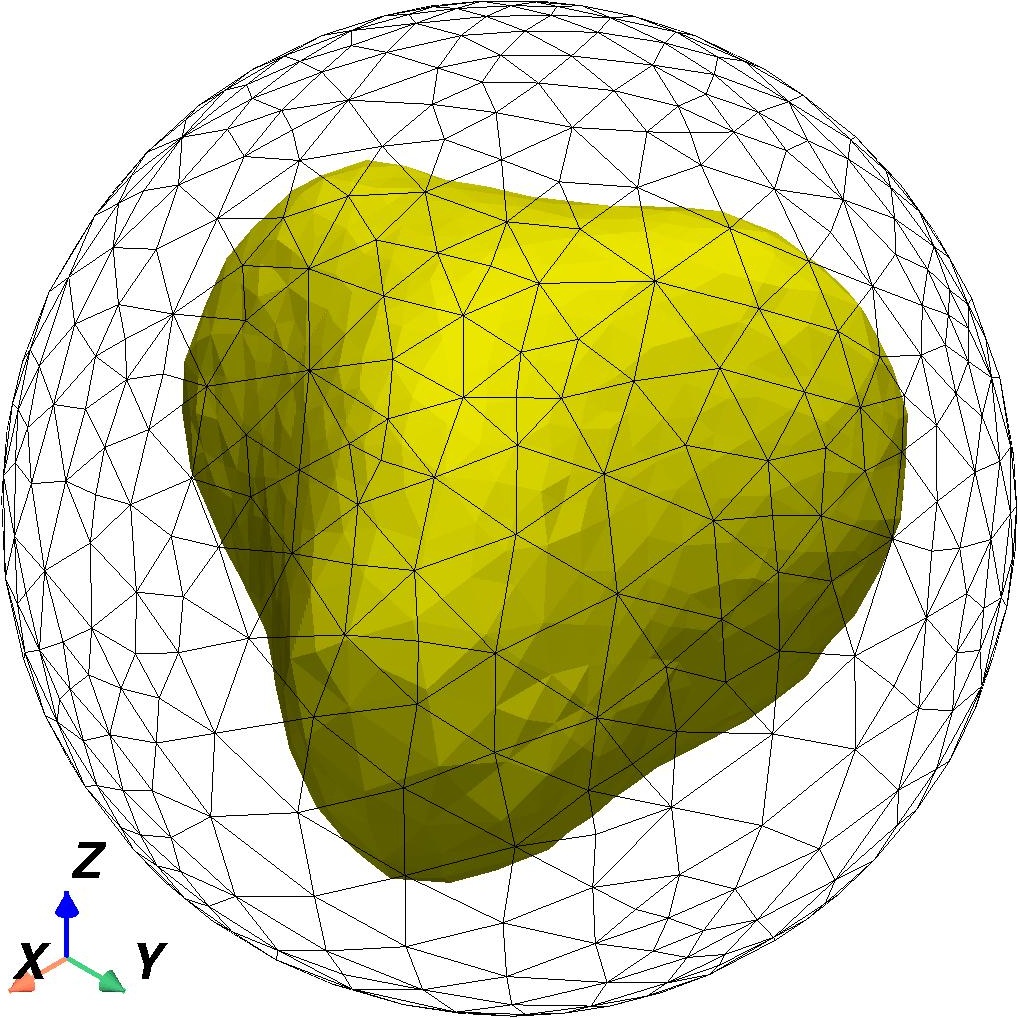}} \hfill
\resizebox{0.16\linewidth}{!}{\includegraphics{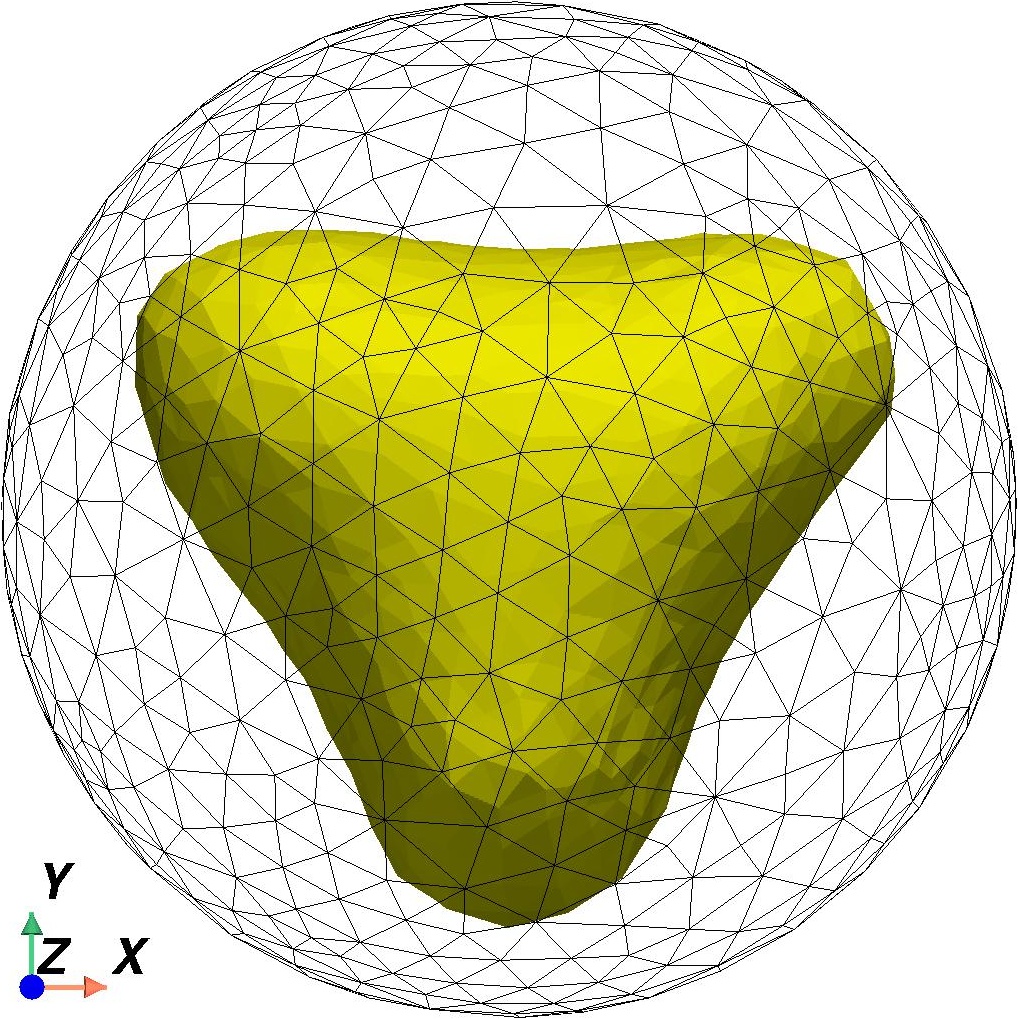}} \hfill
\resizebox{0.16\linewidth}{!}{\includegraphics{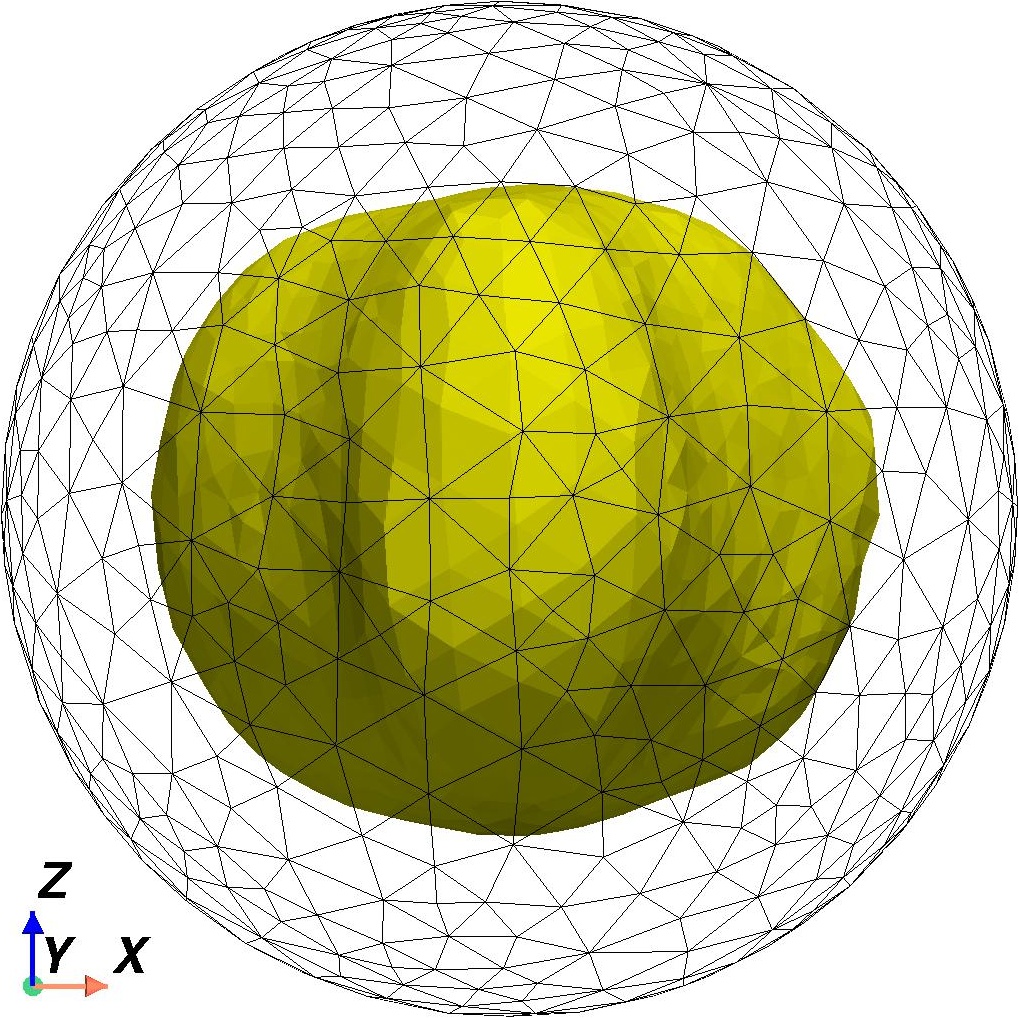}} \hfill
\resizebox{0.16\linewidth}{!}{\includegraphics{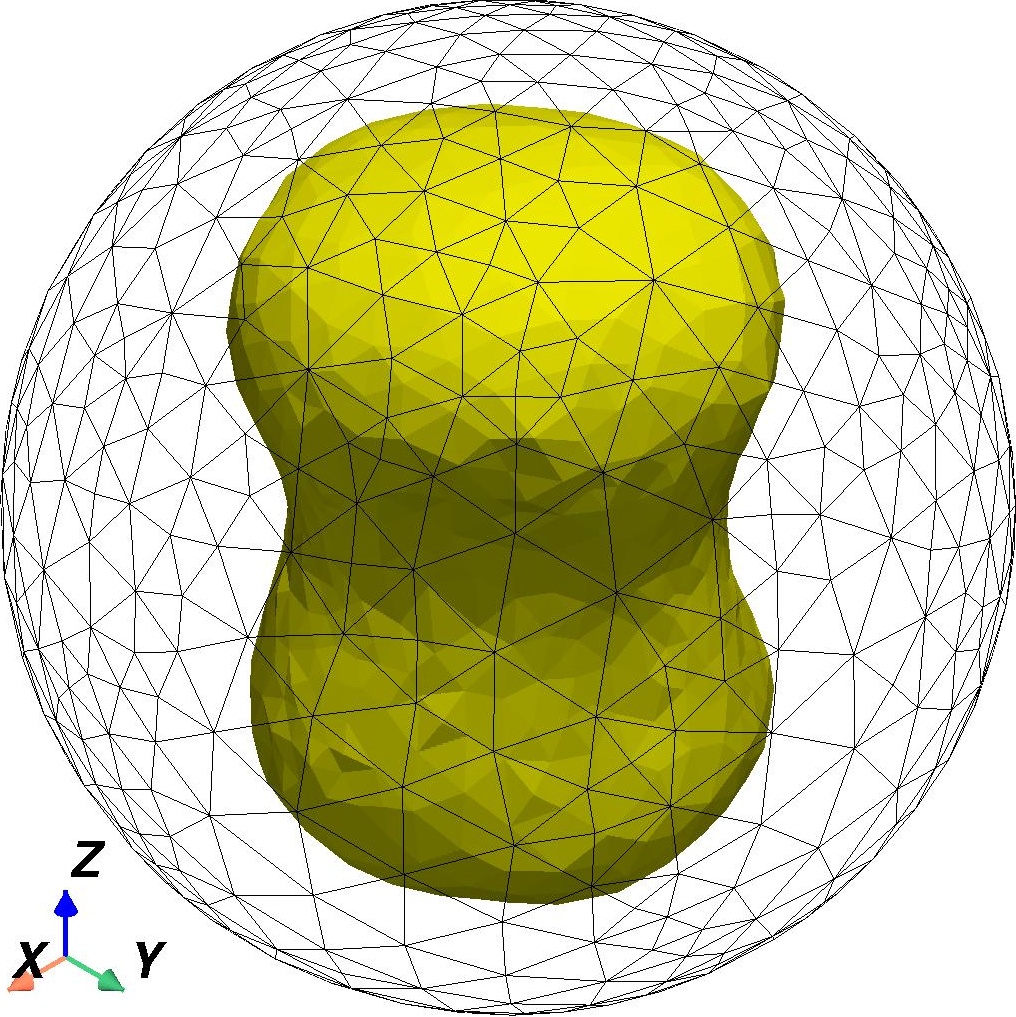}} \hfill
\resizebox{0.16\linewidth}{!}{\includegraphics{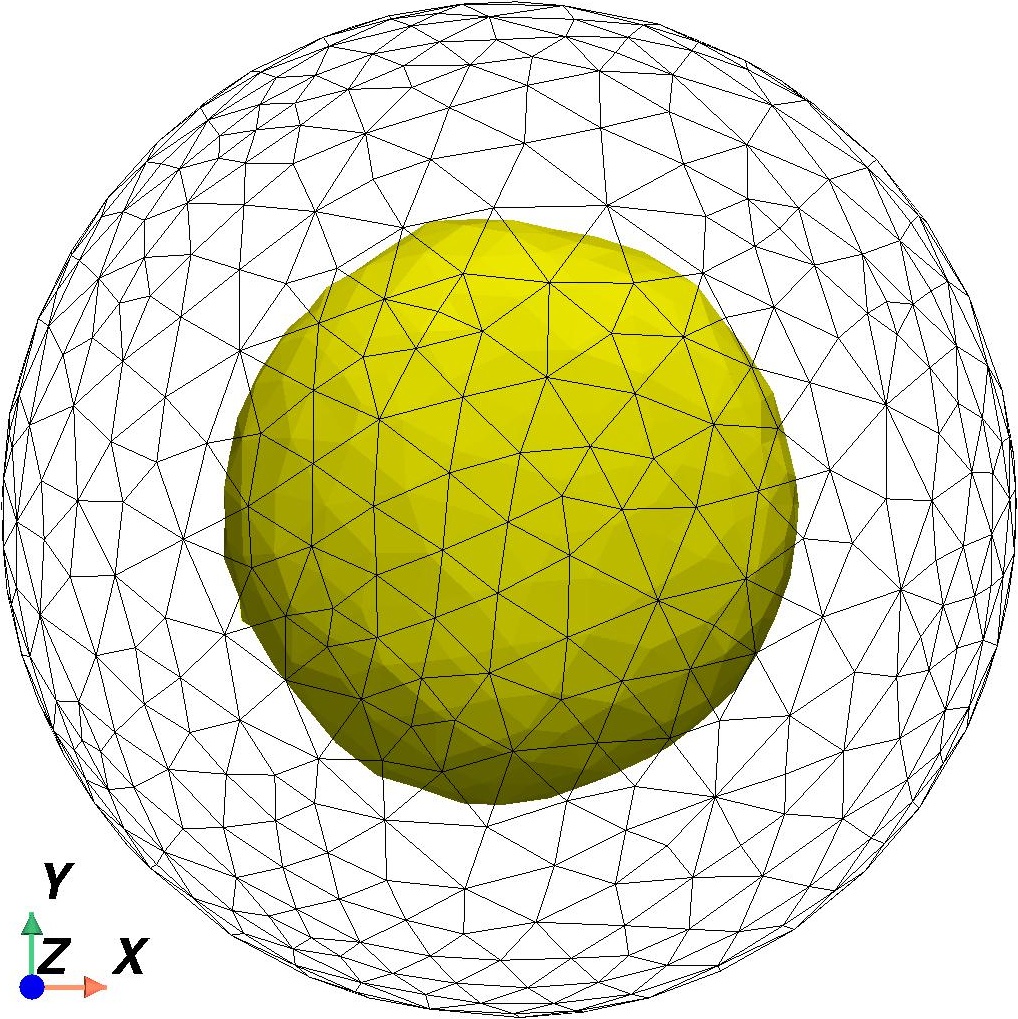}} \hfill
\resizebox{0.16\linewidth}{!}{\includegraphics{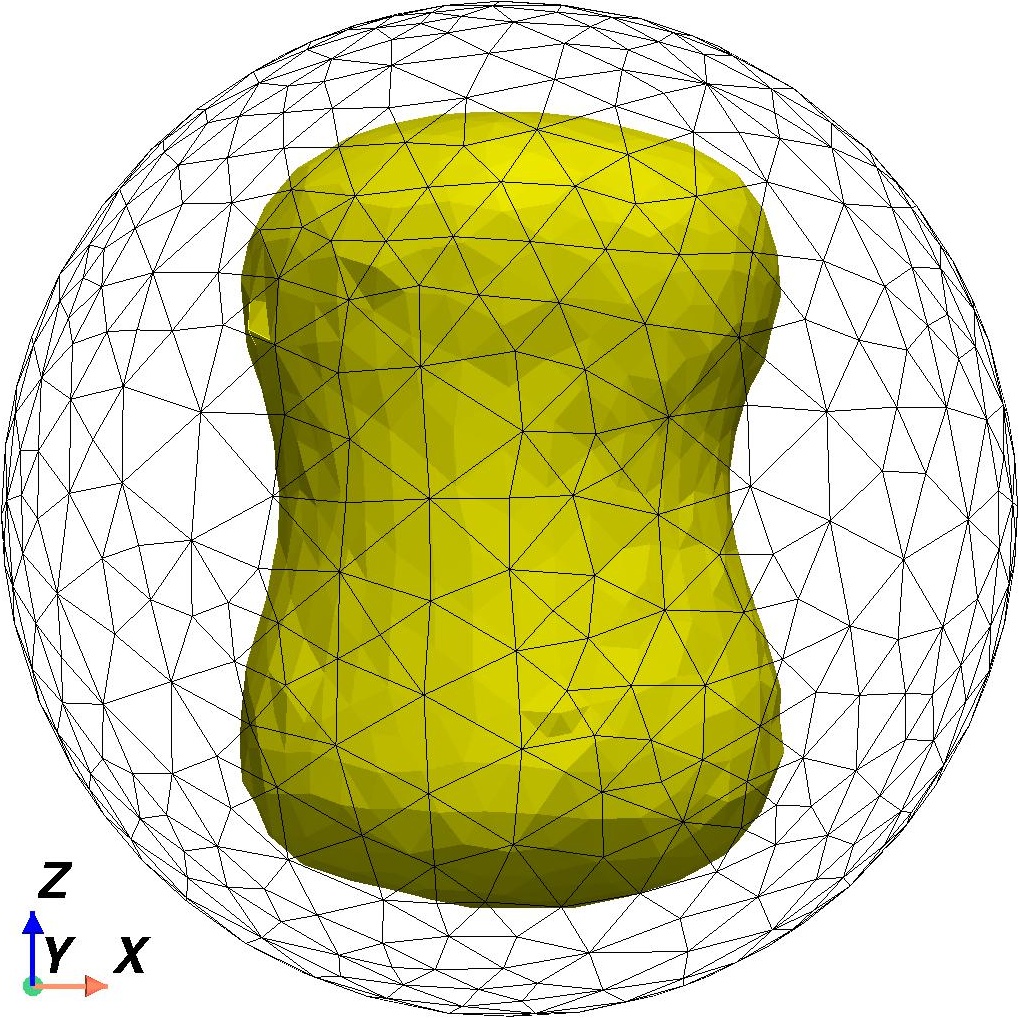}} 
 \\[1em] 
 \resizebox{0.16\linewidth}{!}{\includegraphics{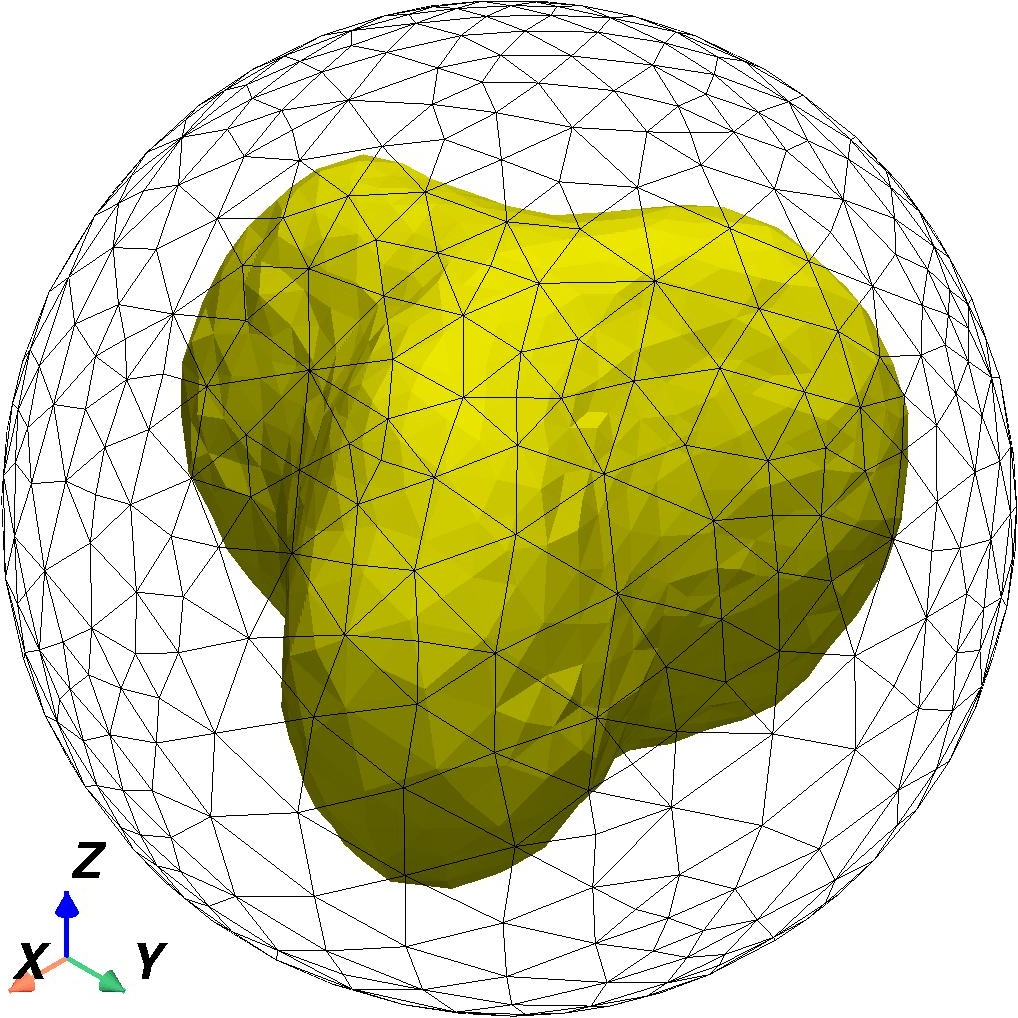}} \hfill
\resizebox{0.16\linewidth}{!}{\includegraphics{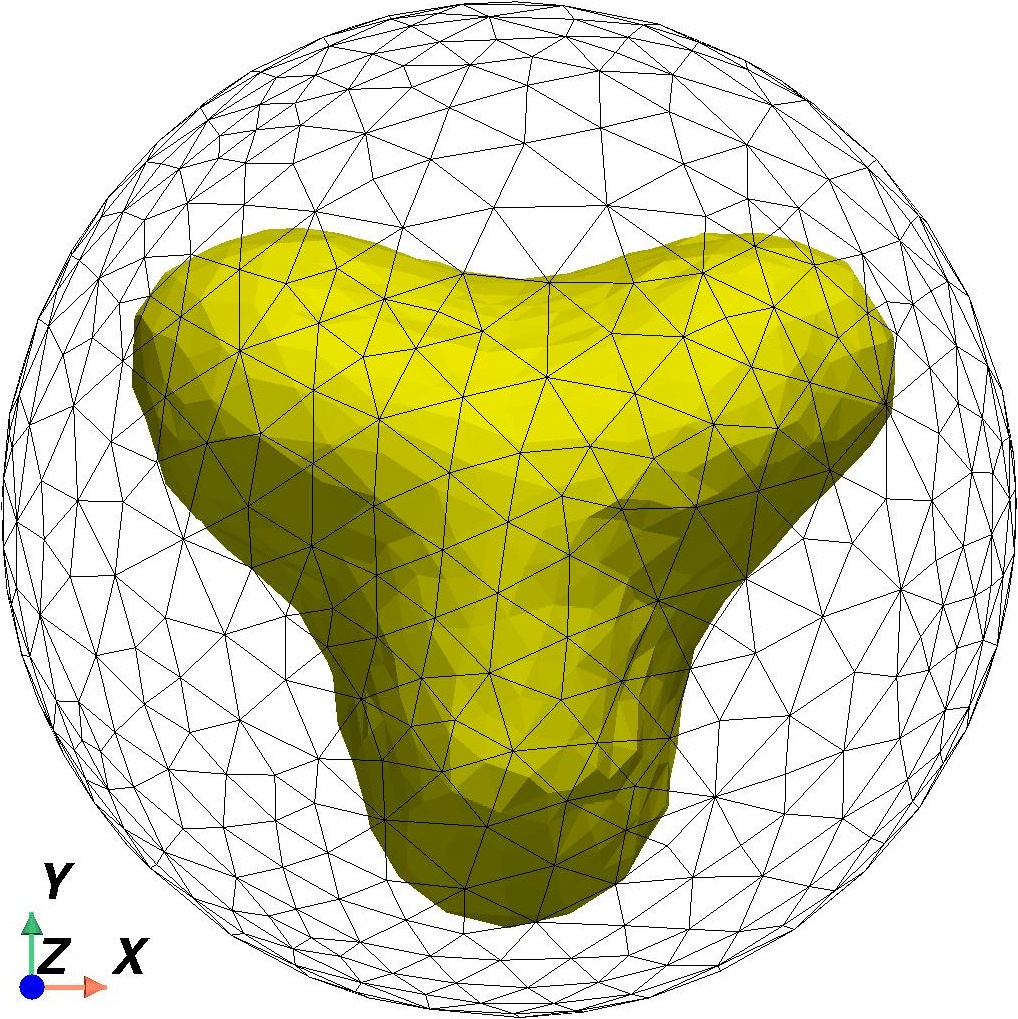}} \hfill
\resizebox{0.16\linewidth}{!}{\includegraphics{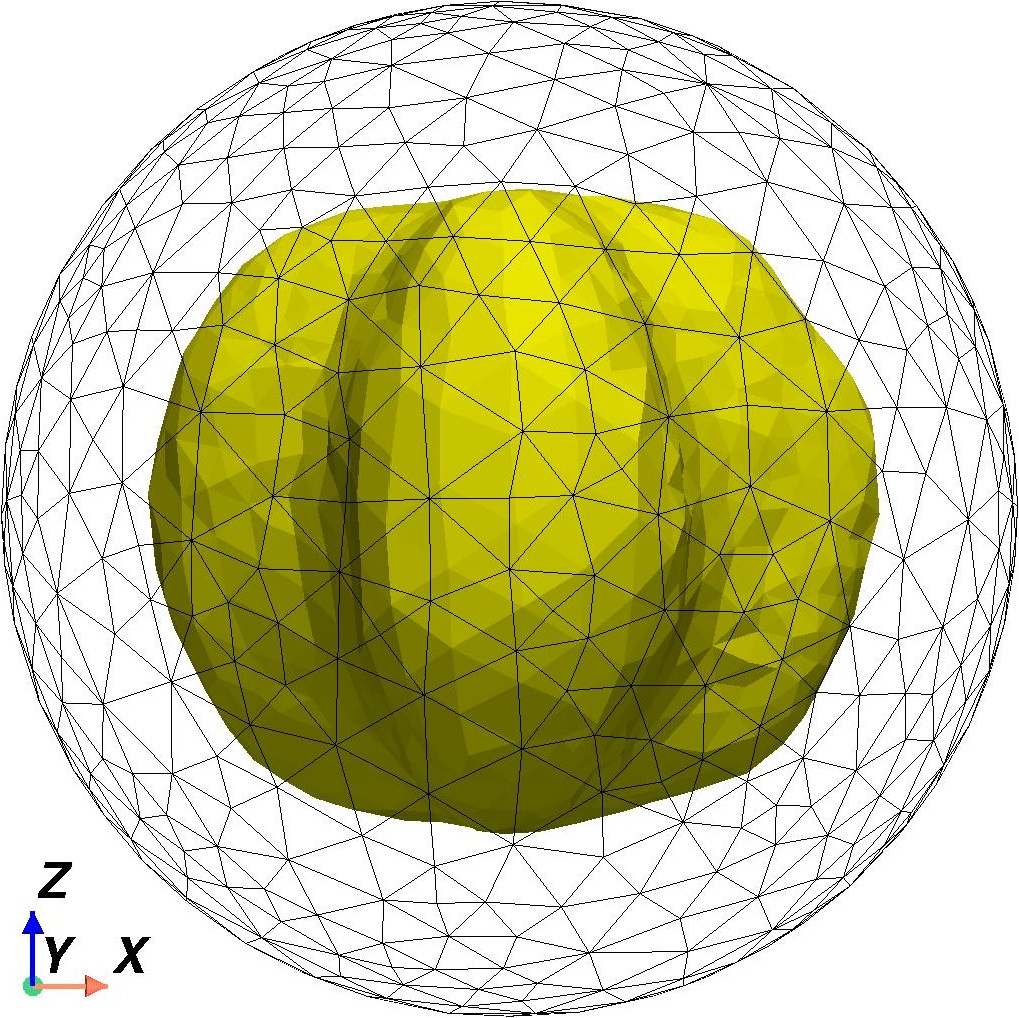}} \hfill
 \resizebox{0.16\linewidth}{!}{\includegraphics{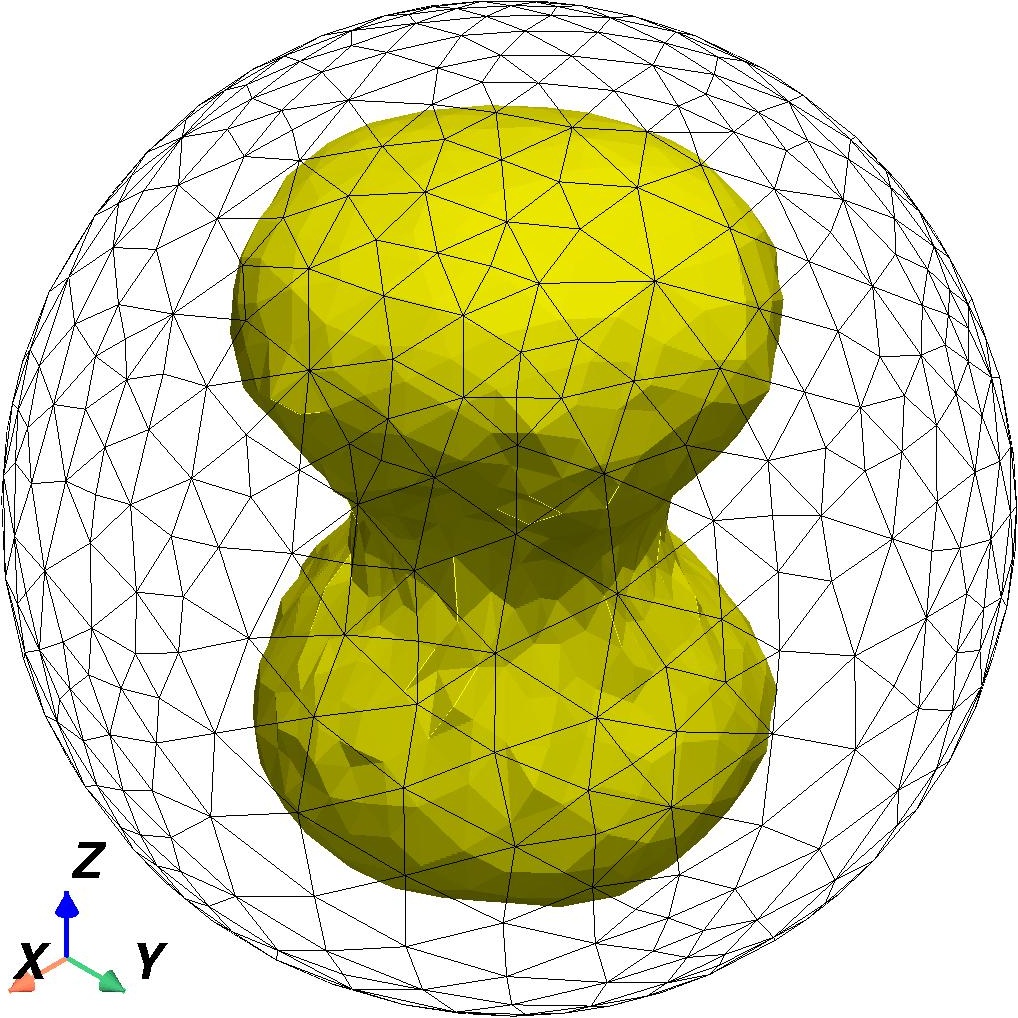}} \hfill
\resizebox{0.16\linewidth}{!}{\includegraphics{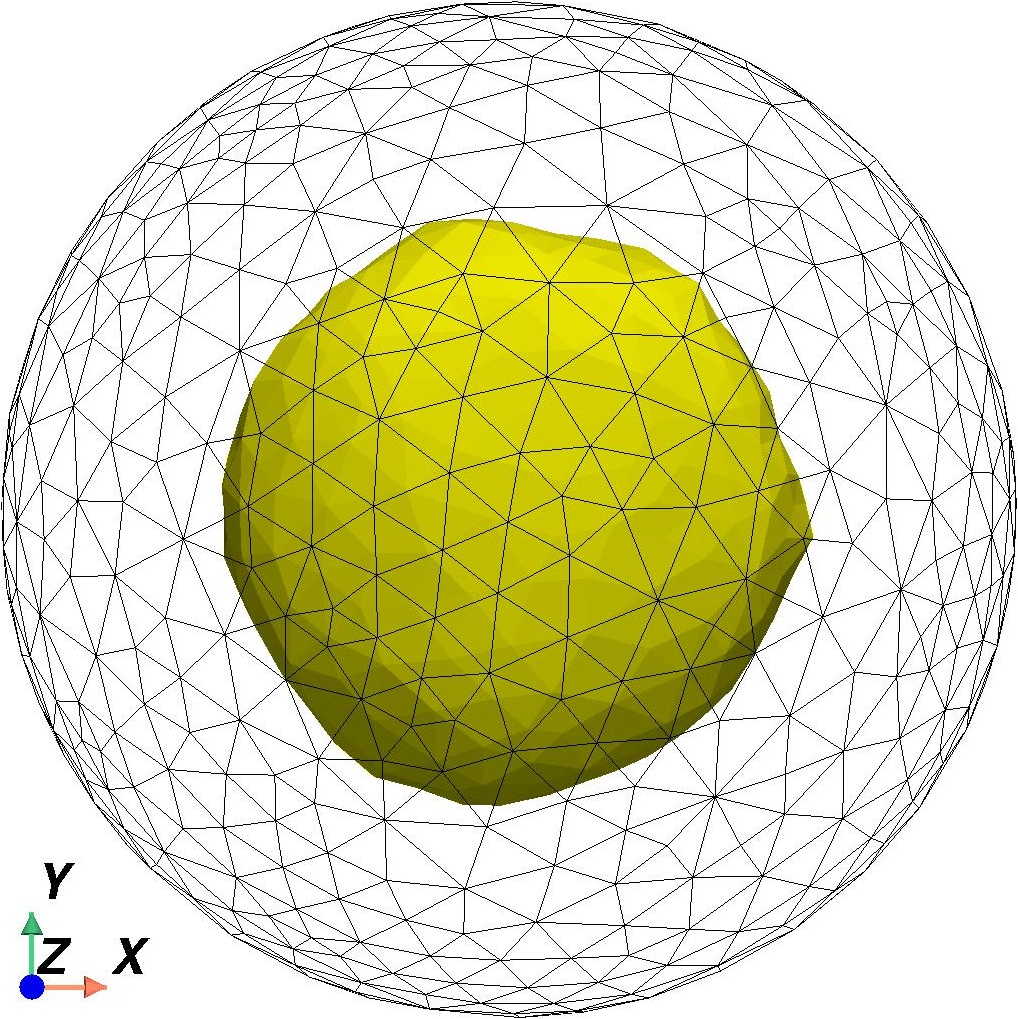}} \hfill
\resizebox{0.16\linewidth}{!}{\includegraphics{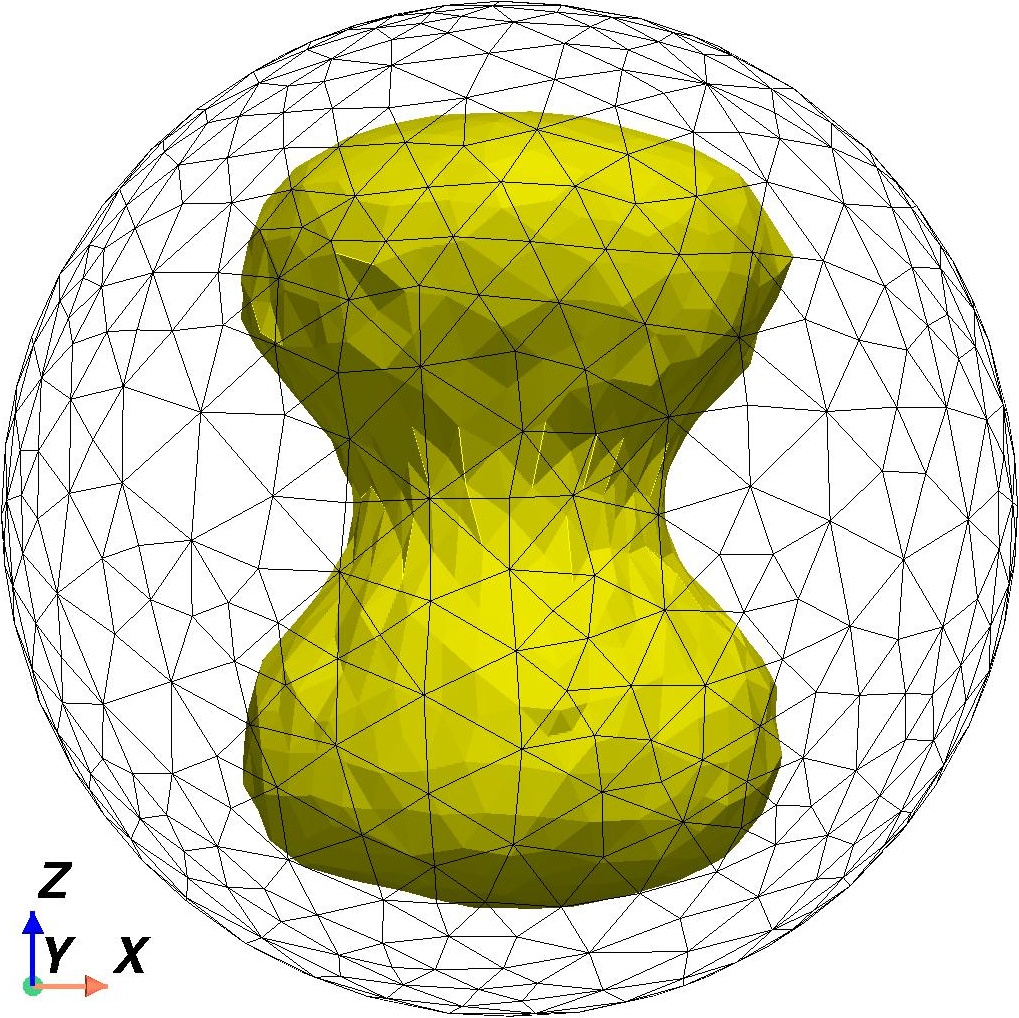}}
\caption{Reconstructions via SO (top/first row) and via ADMM (bottom/second row) with noised data at a noise level of $\delta = 15\%$ with regularization where $\gamma = 0.003$}
\label{fig:figure2c}
\end{figure}

\begin{figure}[htp!]
\centering
\resizebox{0.16\linewidth}{!}{\includegraphics{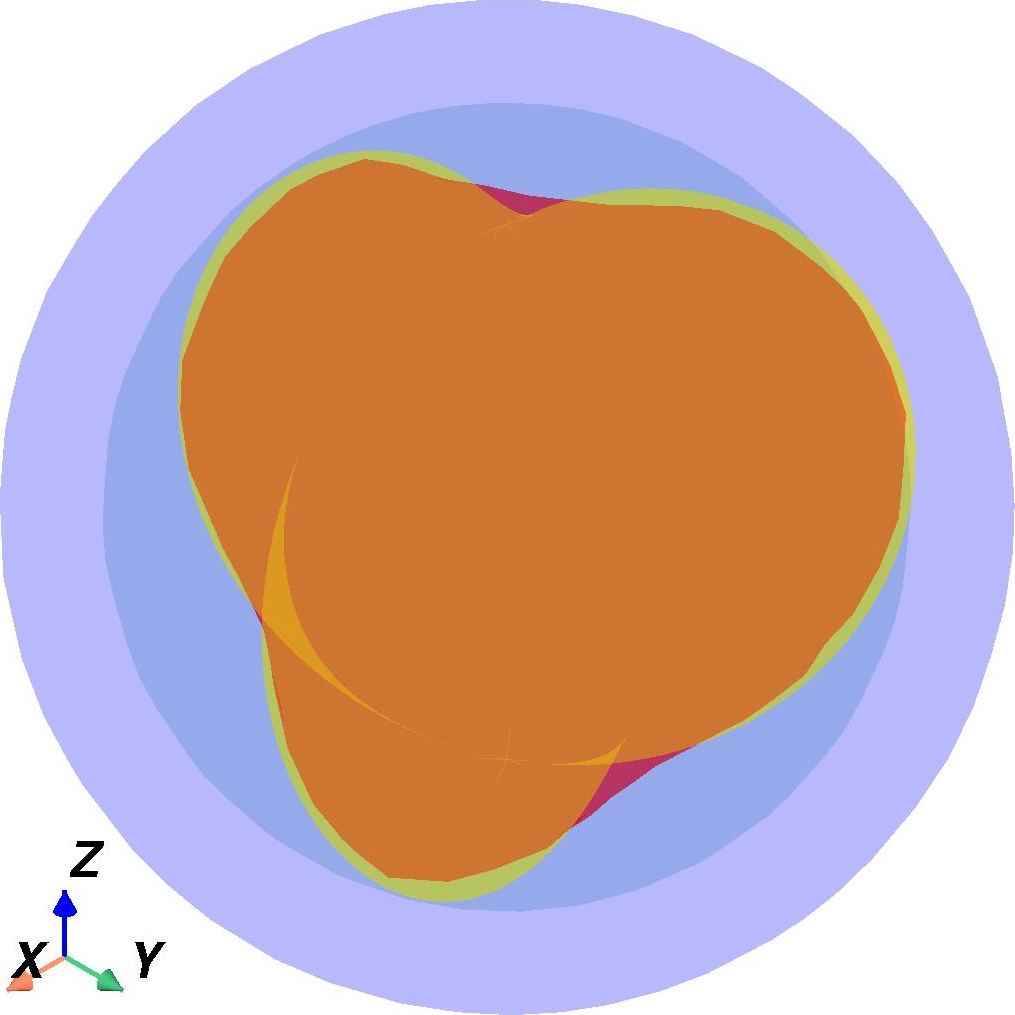}} \hfill
\resizebox{0.16\linewidth}{!}{\includegraphics{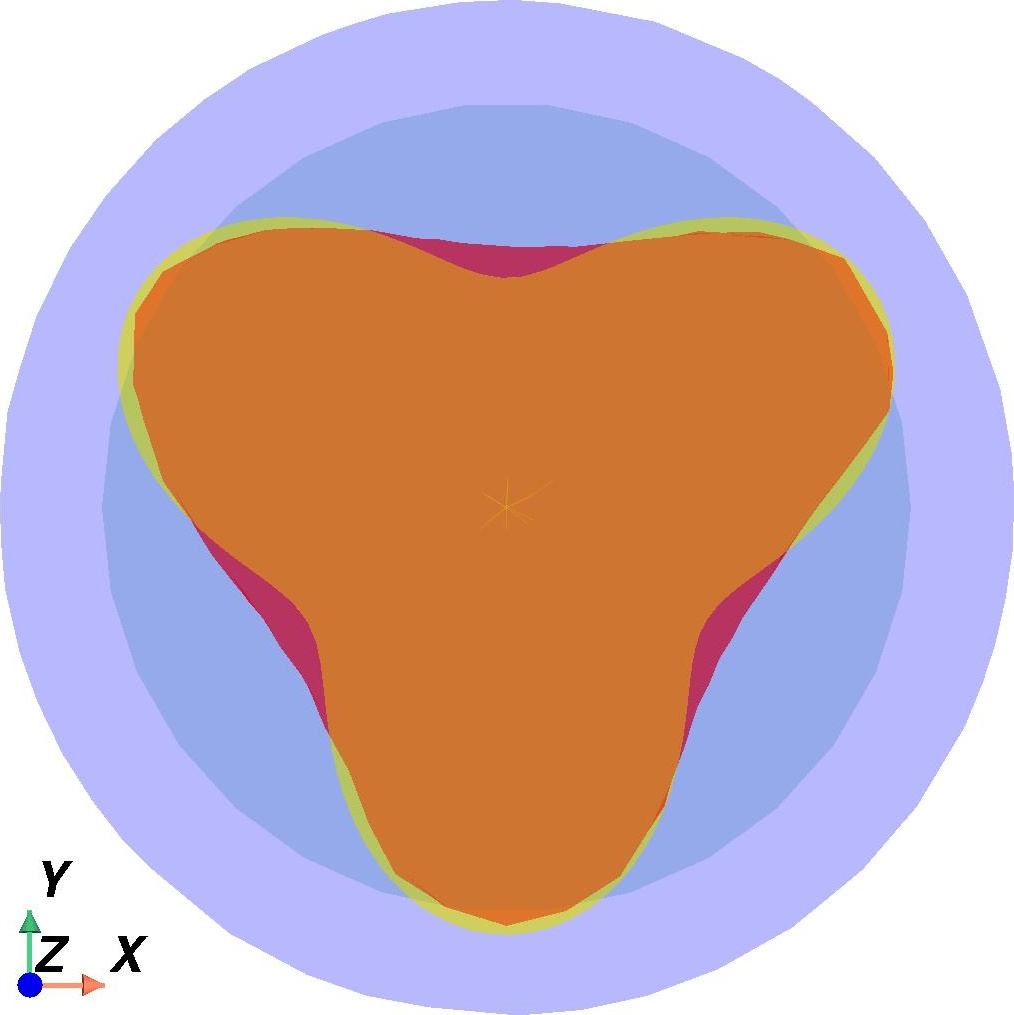}} \hfill
\resizebox{0.16\linewidth}{!}{\includegraphics{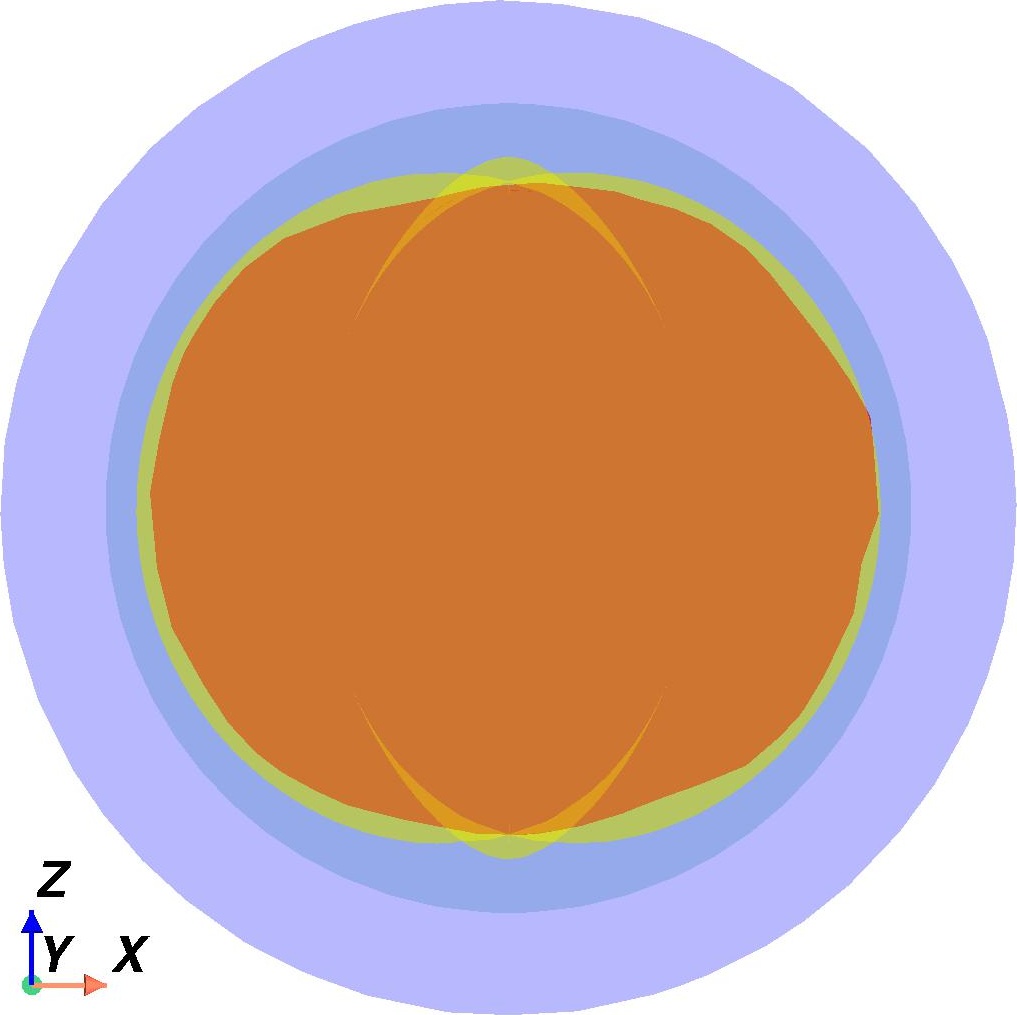}} \hfill
\resizebox{0.16\linewidth}{!}{\includegraphics{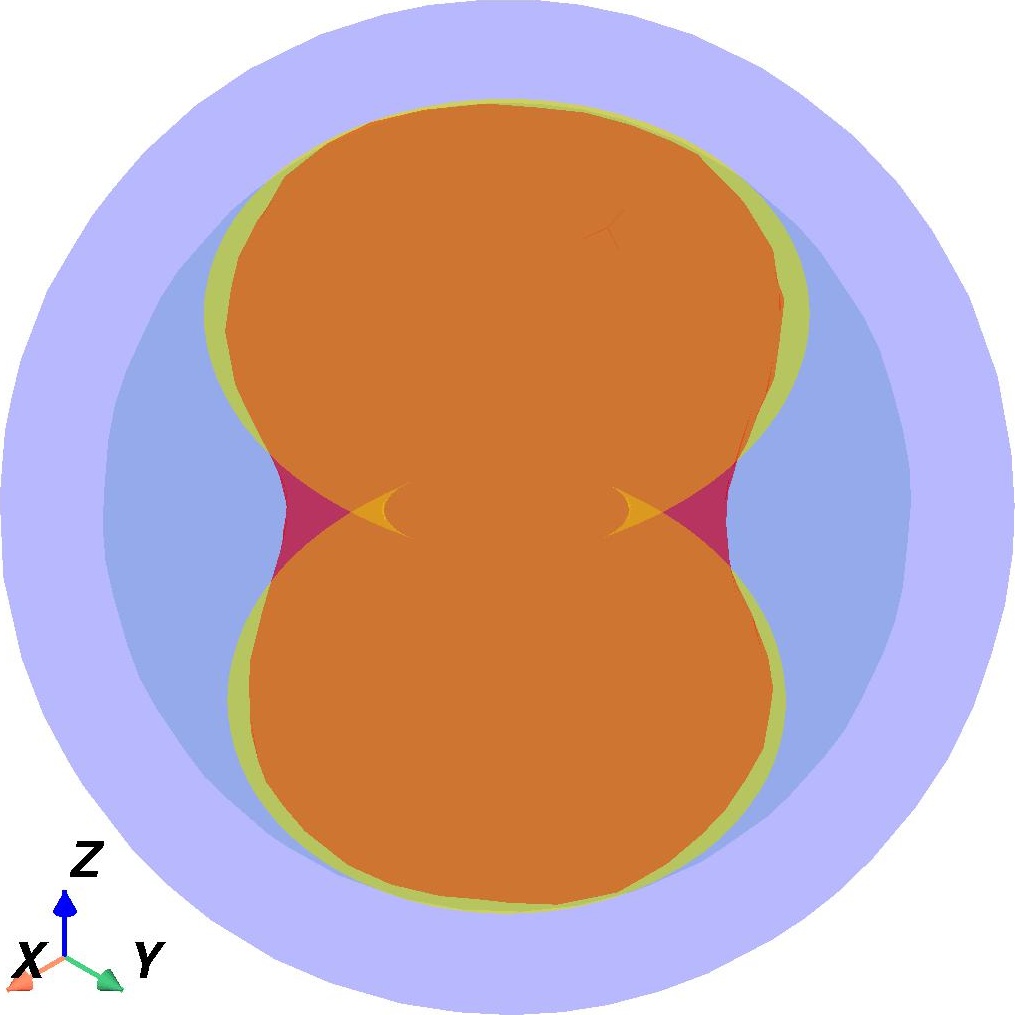}} \hfill
\resizebox{0.16\linewidth}{!}{\includegraphics{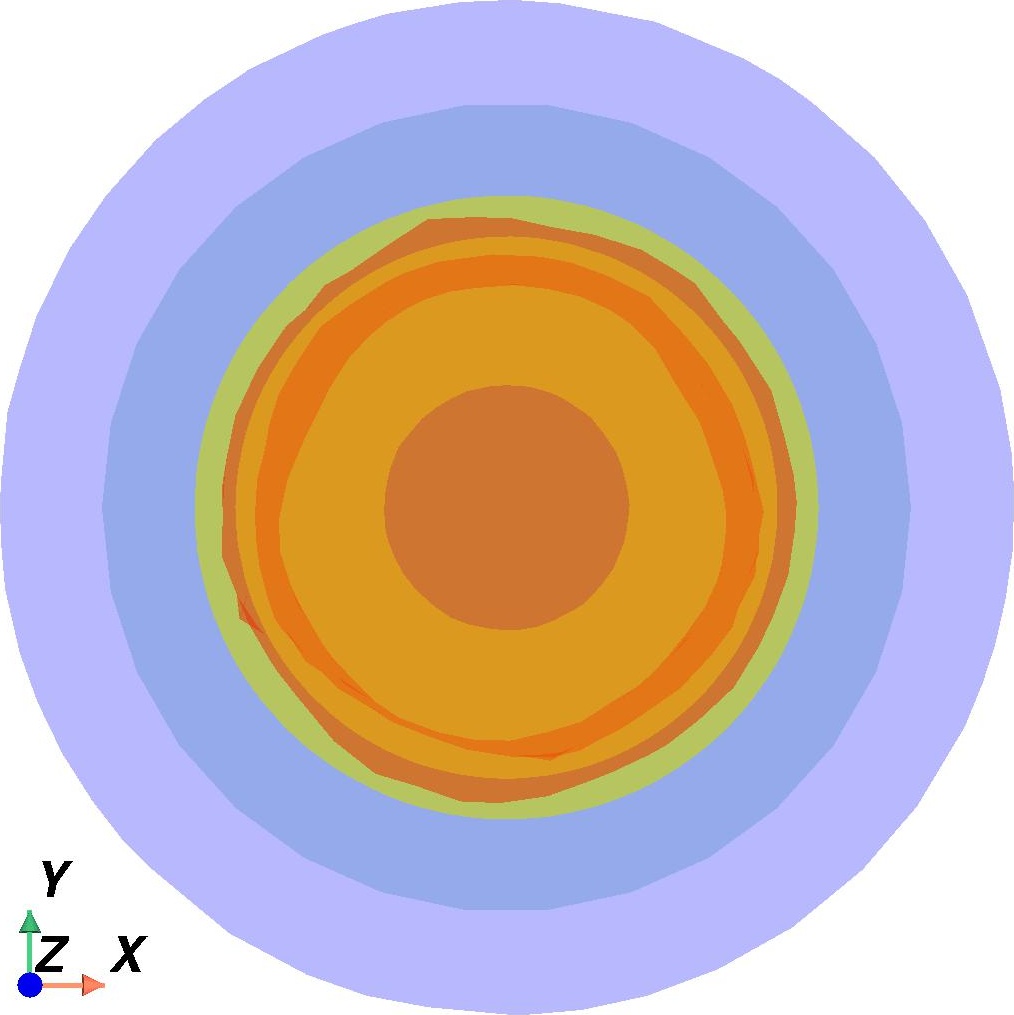}} \hfill
\resizebox{0.16\linewidth}{!}{\includegraphics{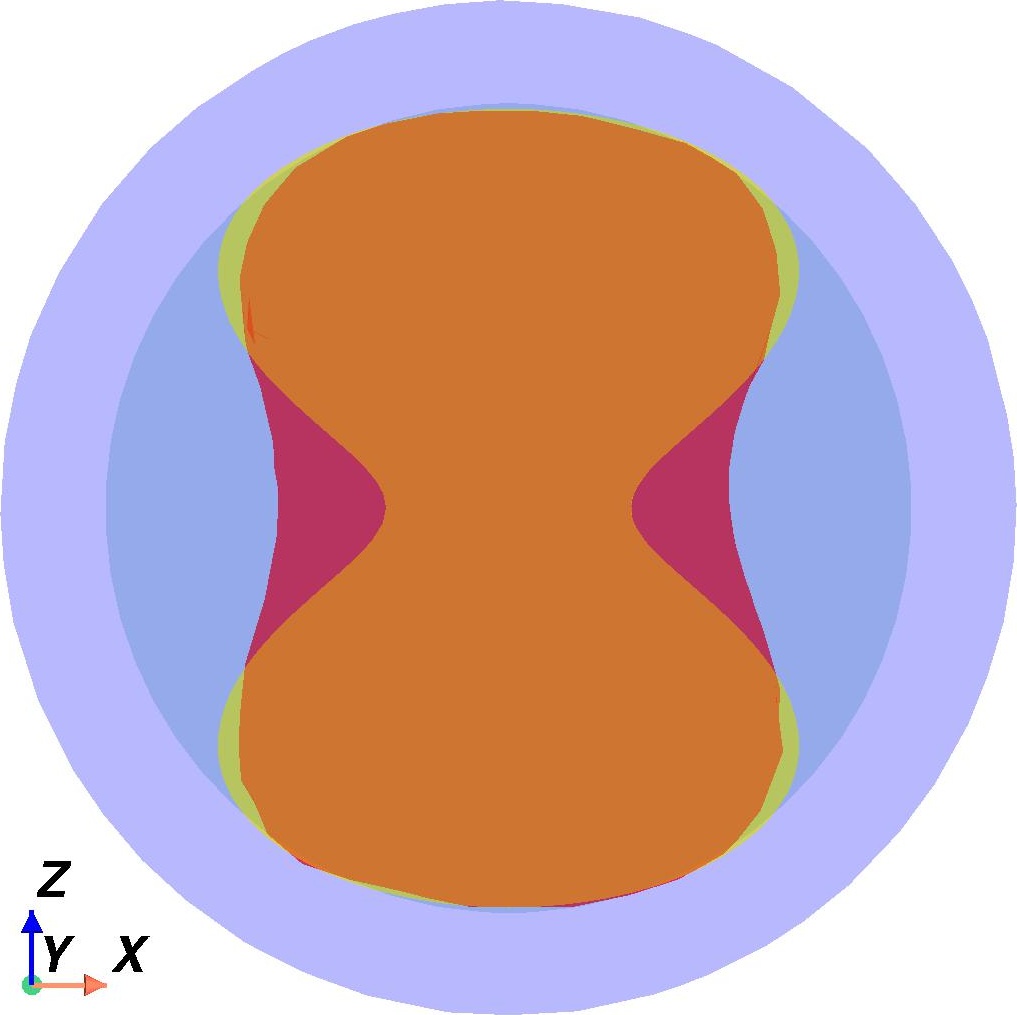}}
 \\[1em] 
\resizebox{0.16\linewidth}{!}{\includegraphics{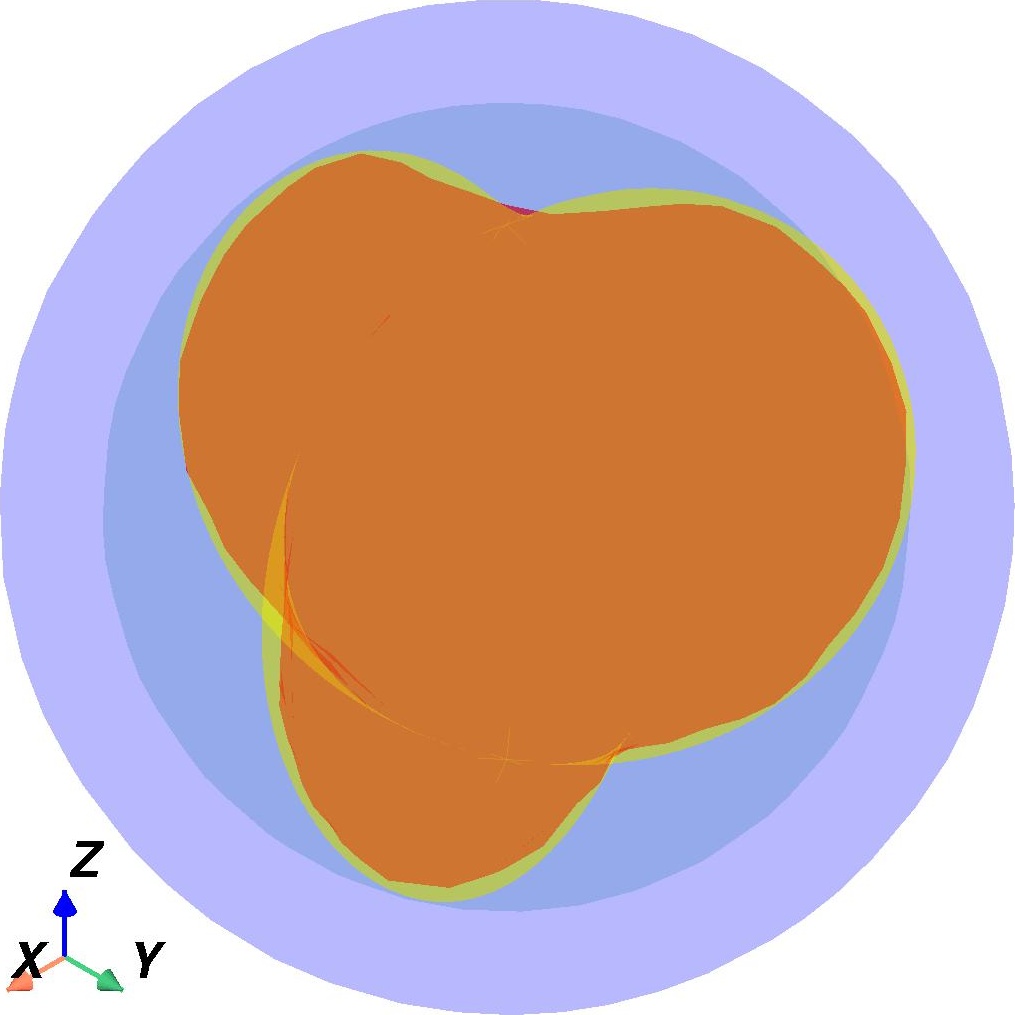}} \hfill
\resizebox{0.16\linewidth}{!}{\includegraphics{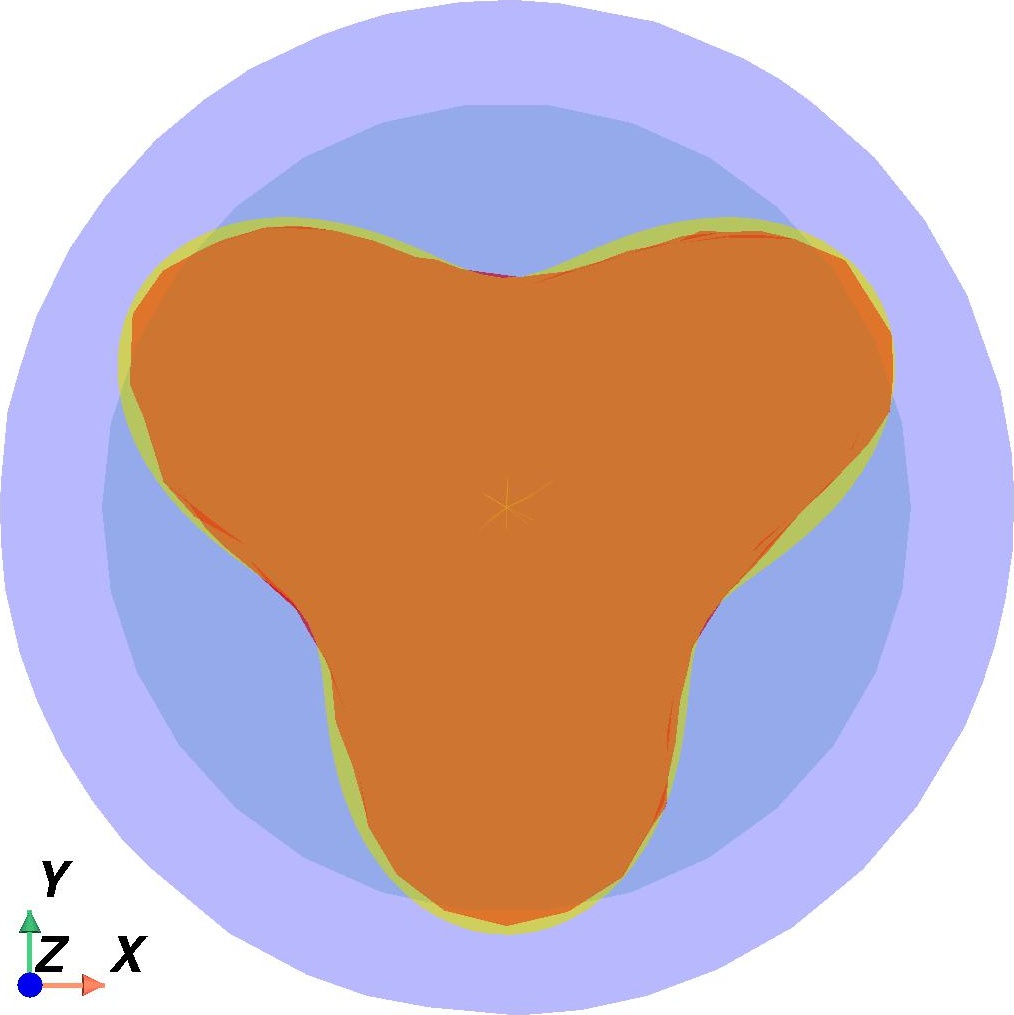}} \hfill
\resizebox{0.16\linewidth}{!}{\includegraphics{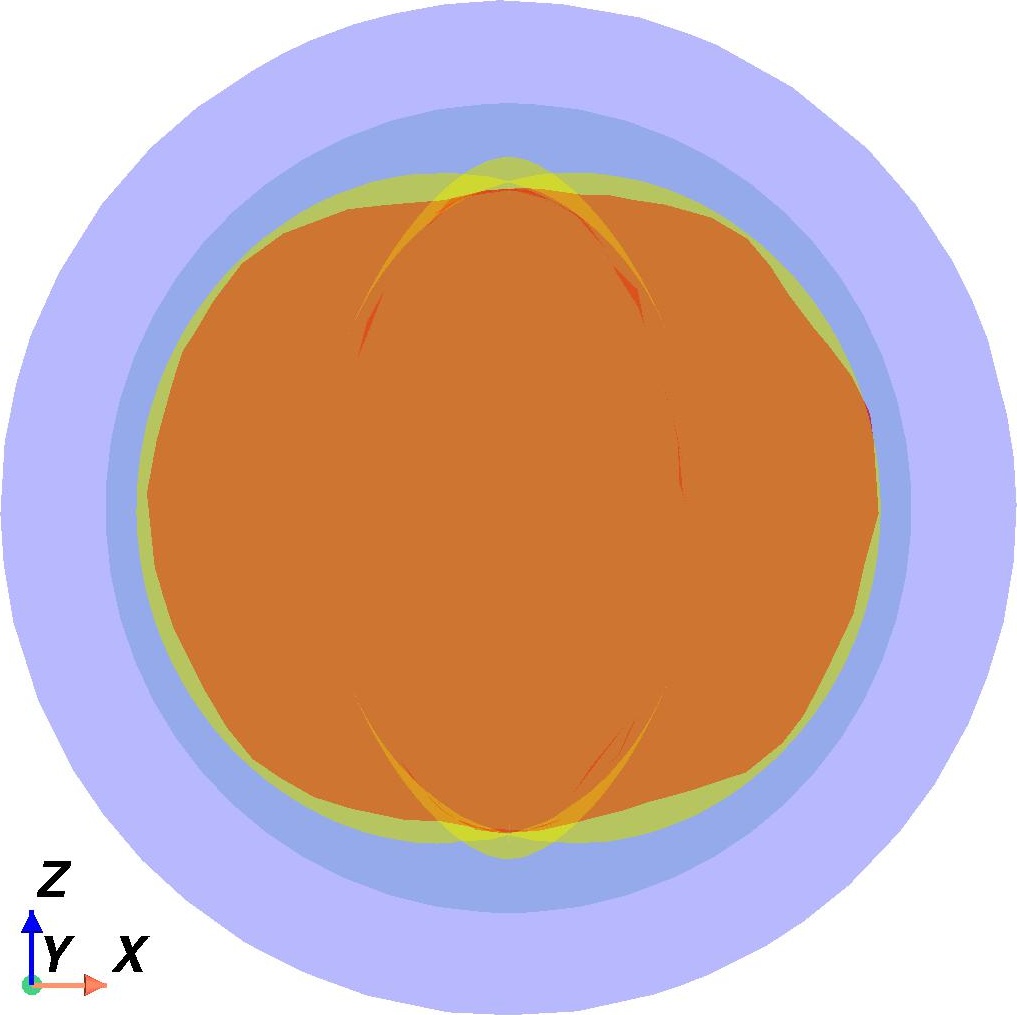}}\hfill
\resizebox{0.16\linewidth}{!}{\includegraphics{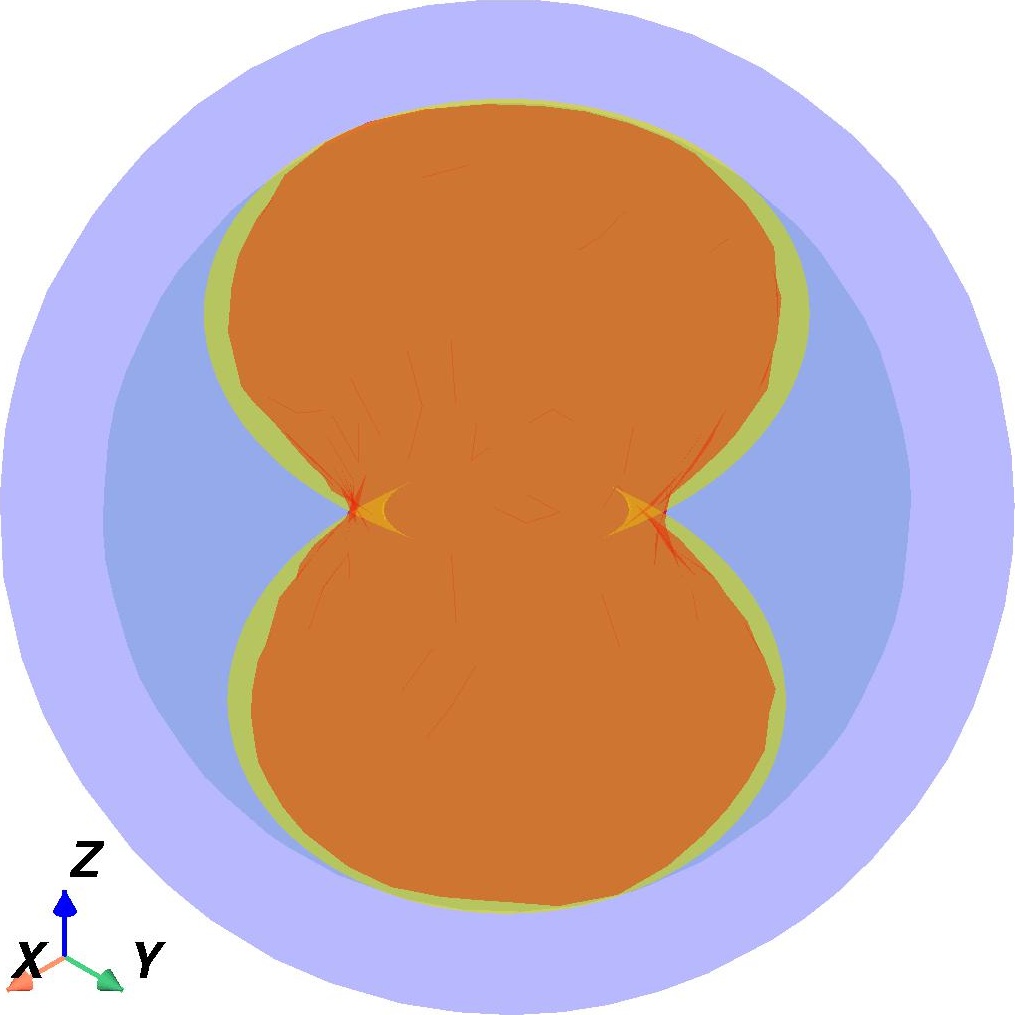}} \hfill
\resizebox{0.16\linewidth}{!}{\includegraphics{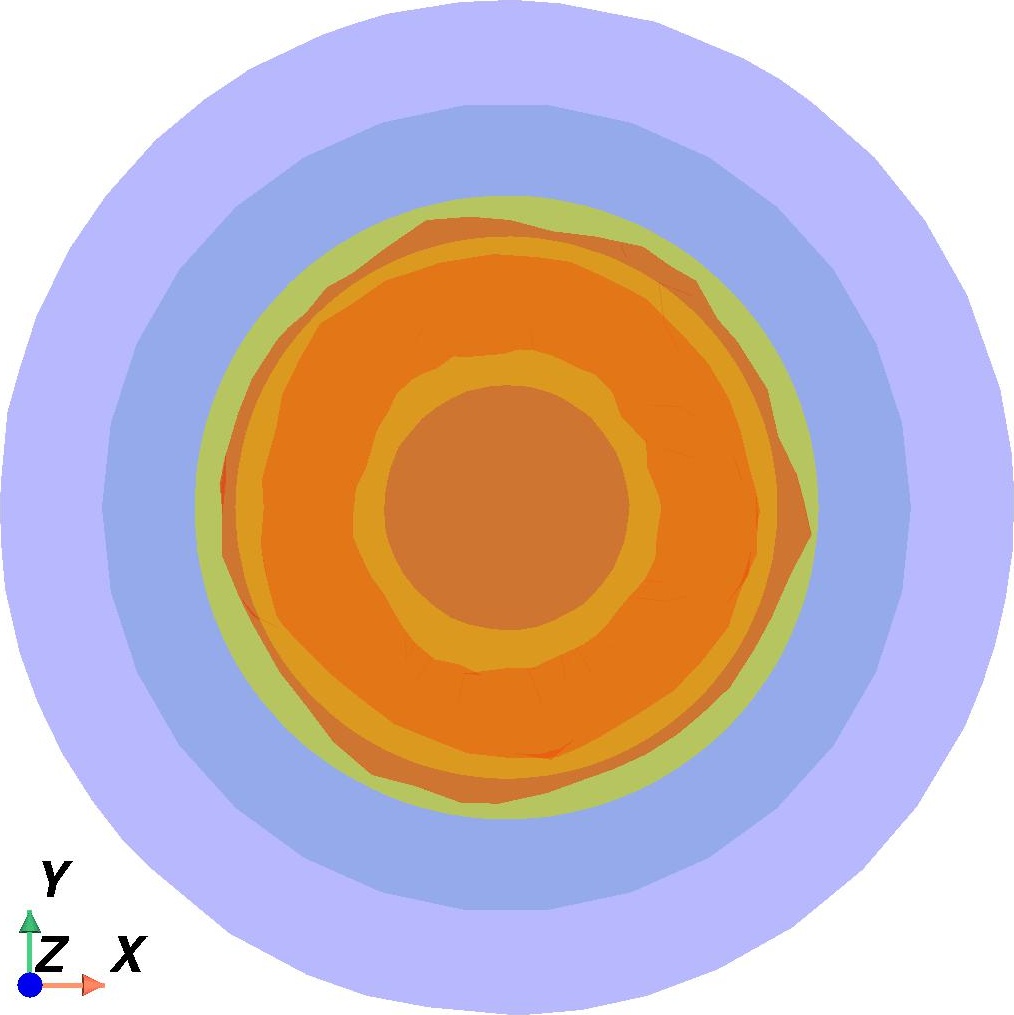}} \hfill
\resizebox{0.16\linewidth}{!}{\includegraphics{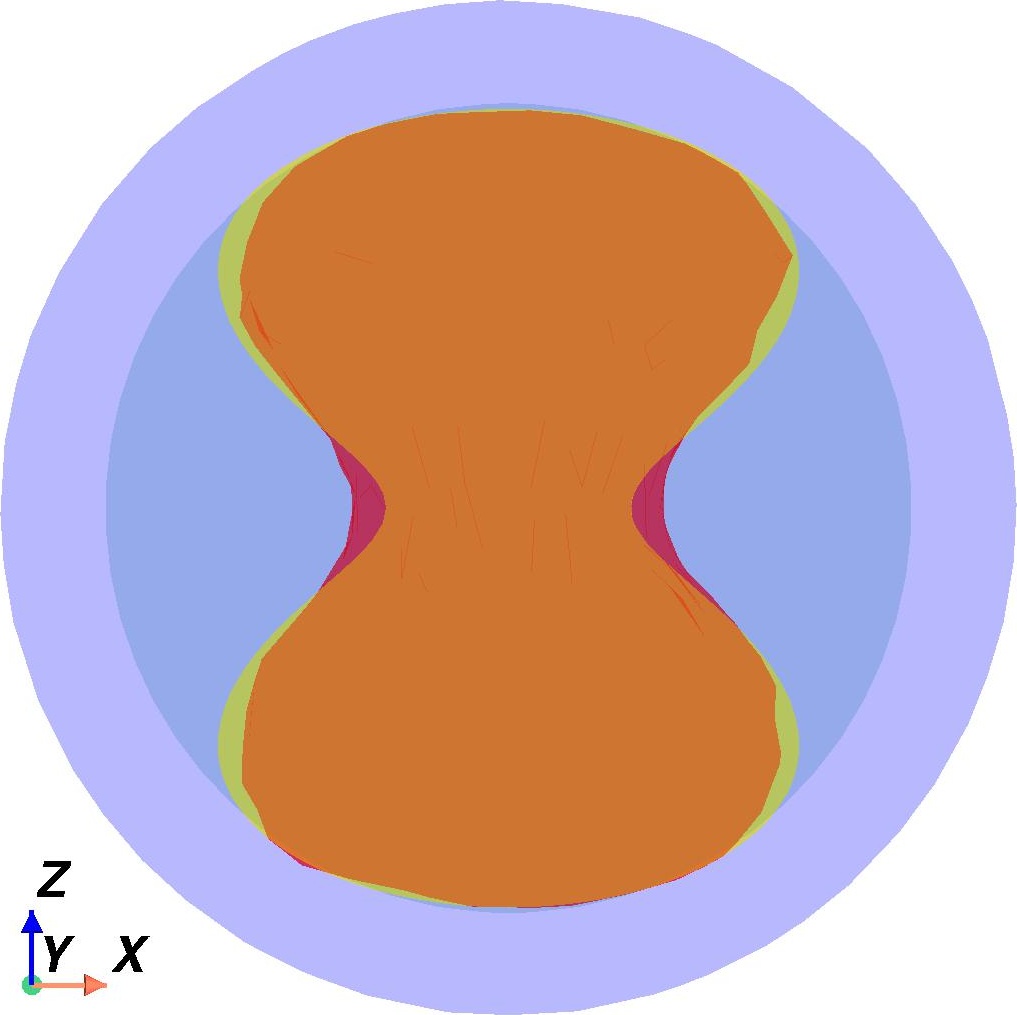}} 
\caption{Cross comparisons of exact and reconstructed shapes via SO (top/first row) and via ADMM (bottom/second row) with noised data at a noise level of $\delta = 15\%$ with regularization where $\gamma = 0.003$}
\label{fig:figure2d}
\end{figure}

\begin{figure}[htp!]
\centering
\resizebox{0.27\linewidth}{!}{\includegraphics{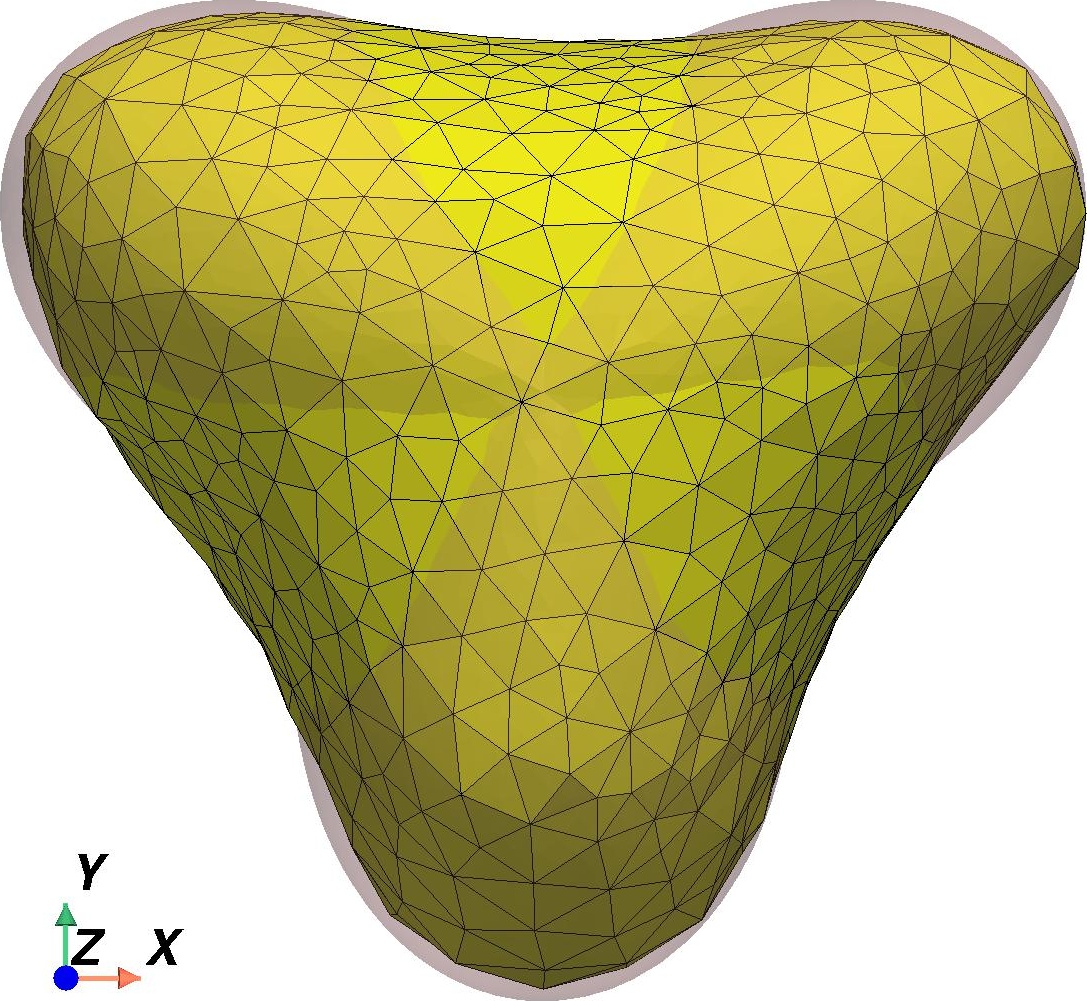}}\
\resizebox{0.27\linewidth}{!}{\includegraphics{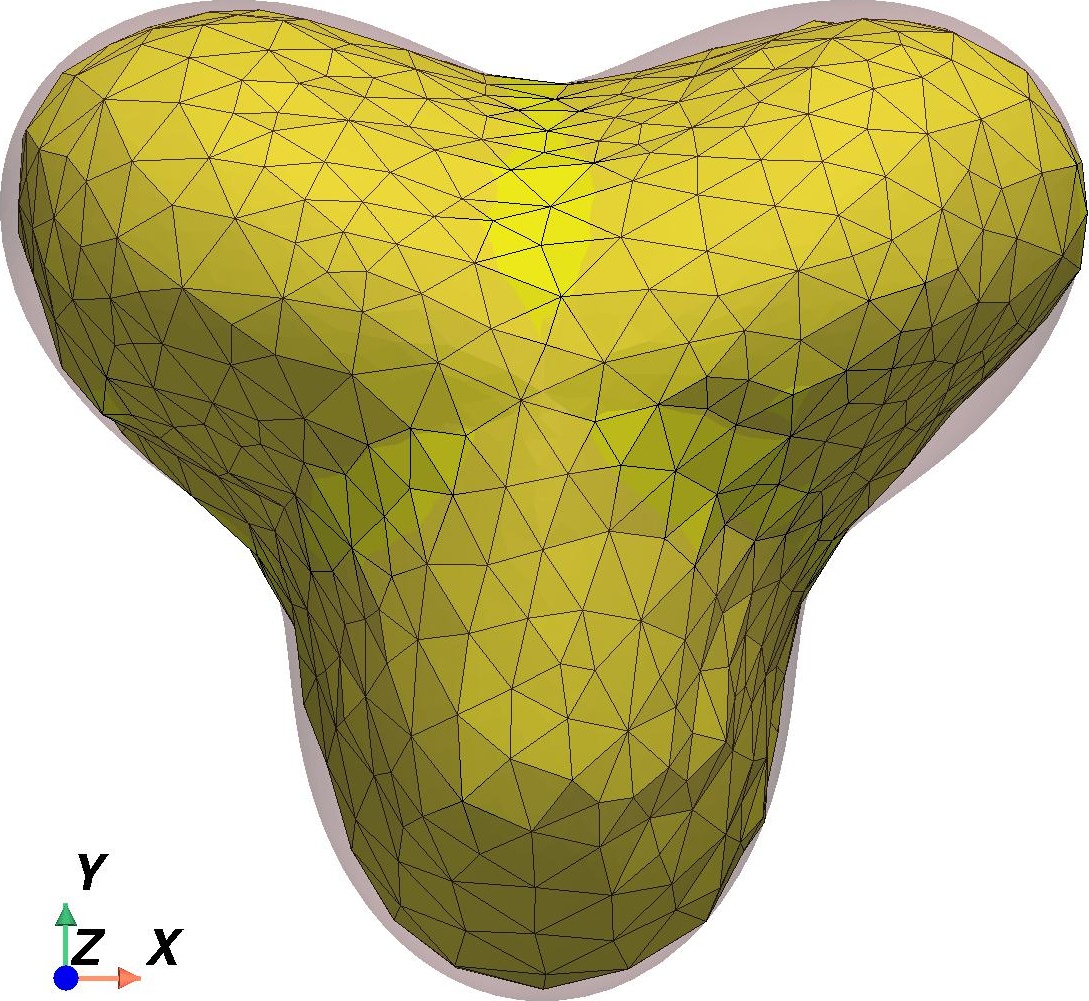}} 
\quad
\resizebox{0.2\linewidth}{!}{\includegraphics{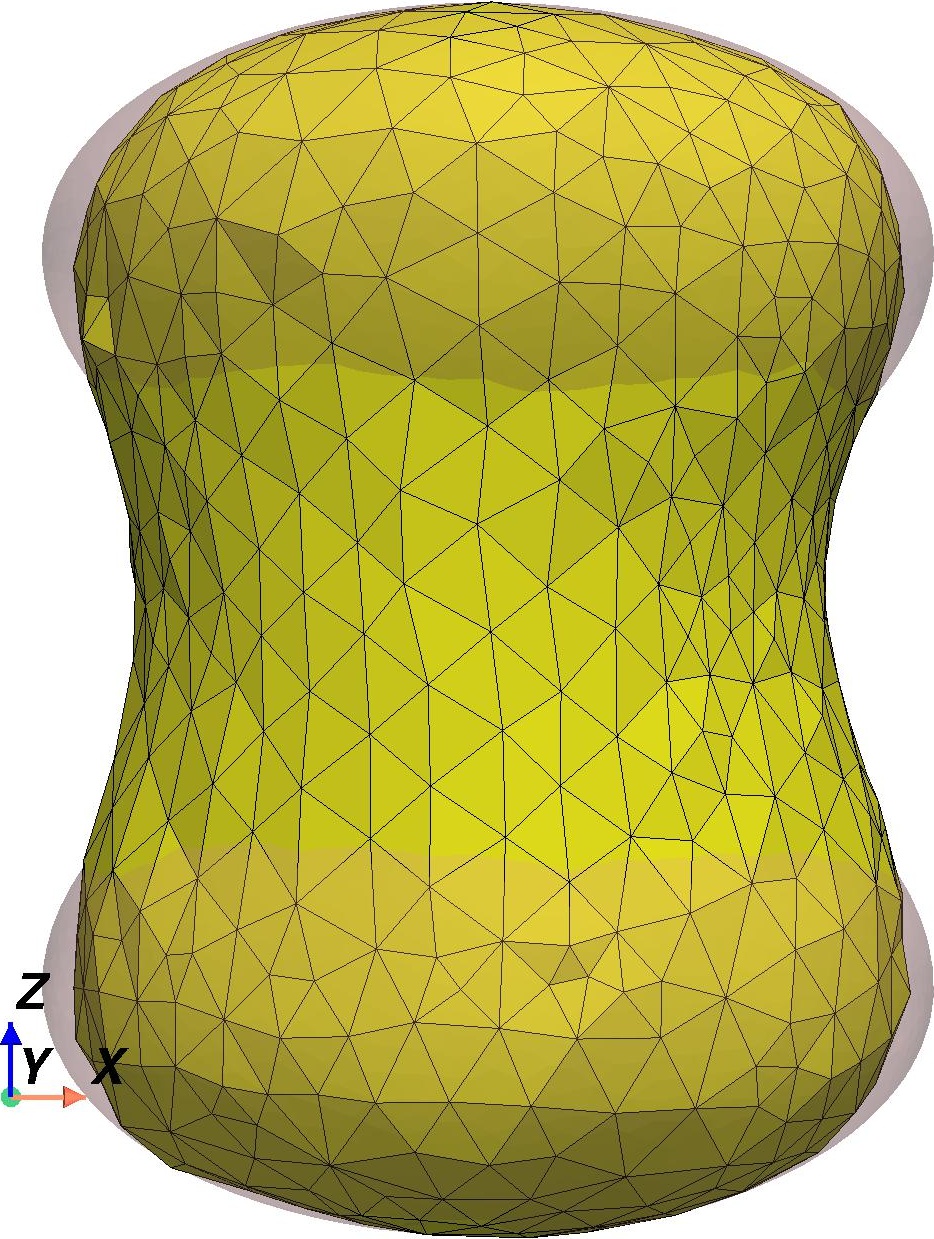}}\
\resizebox{0.2\linewidth}{!}{\includegraphics{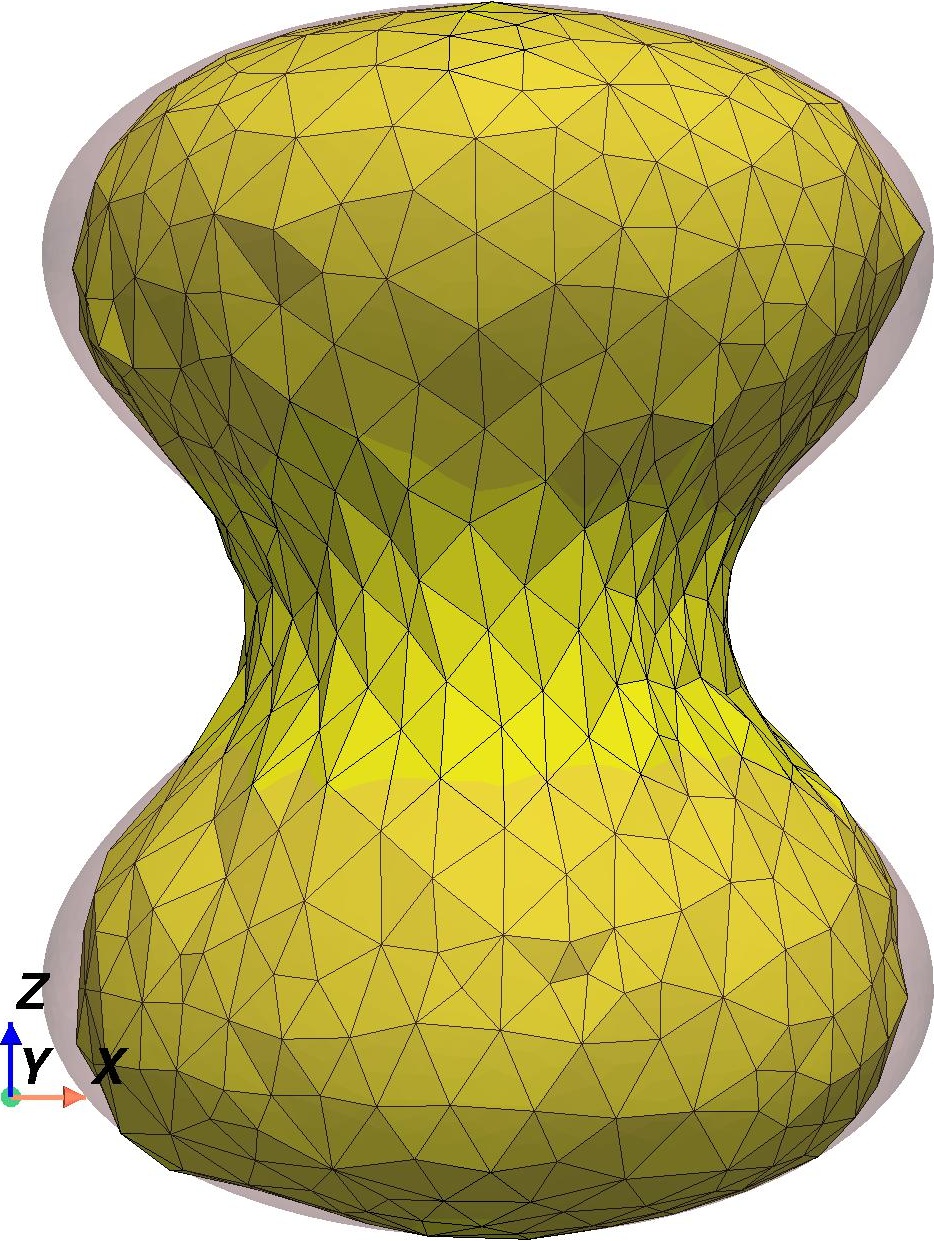}}
\caption{A closer examination of cross-comparisons between exact and computed cavities using SO (left column) and ADMM (right column) with noisy data at a noise level of $\delta = 15\%$, incorporating regularization with $\gamma = 0.003$}
\label{fig:figure2e}
\end{figure}

\begin{figure}[htp!]
\centering
\resizebox{0.24\linewidth}{!}{\includegraphics{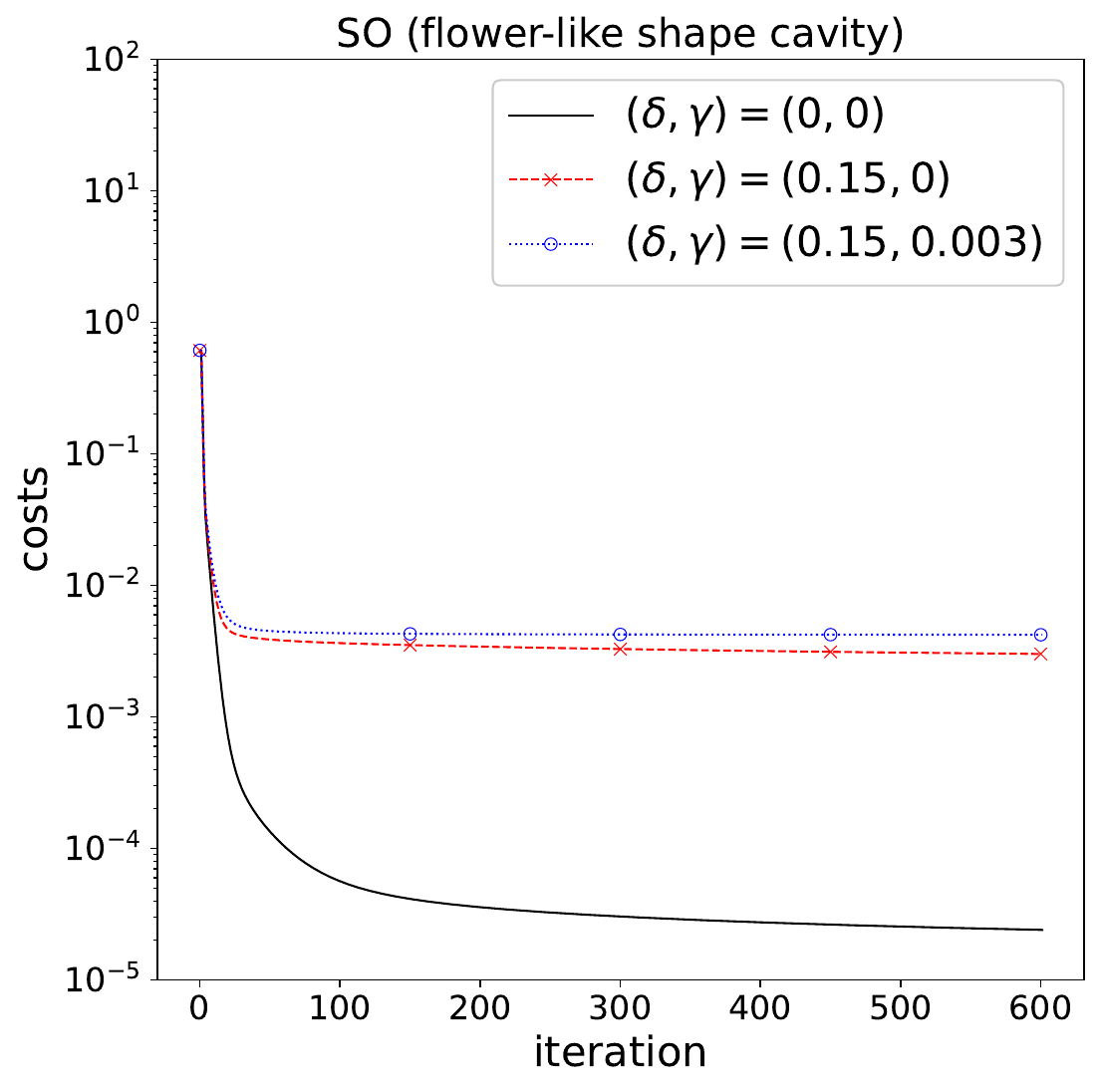}} \hfill
\resizebox{0.24\linewidth}{!}{\includegraphics{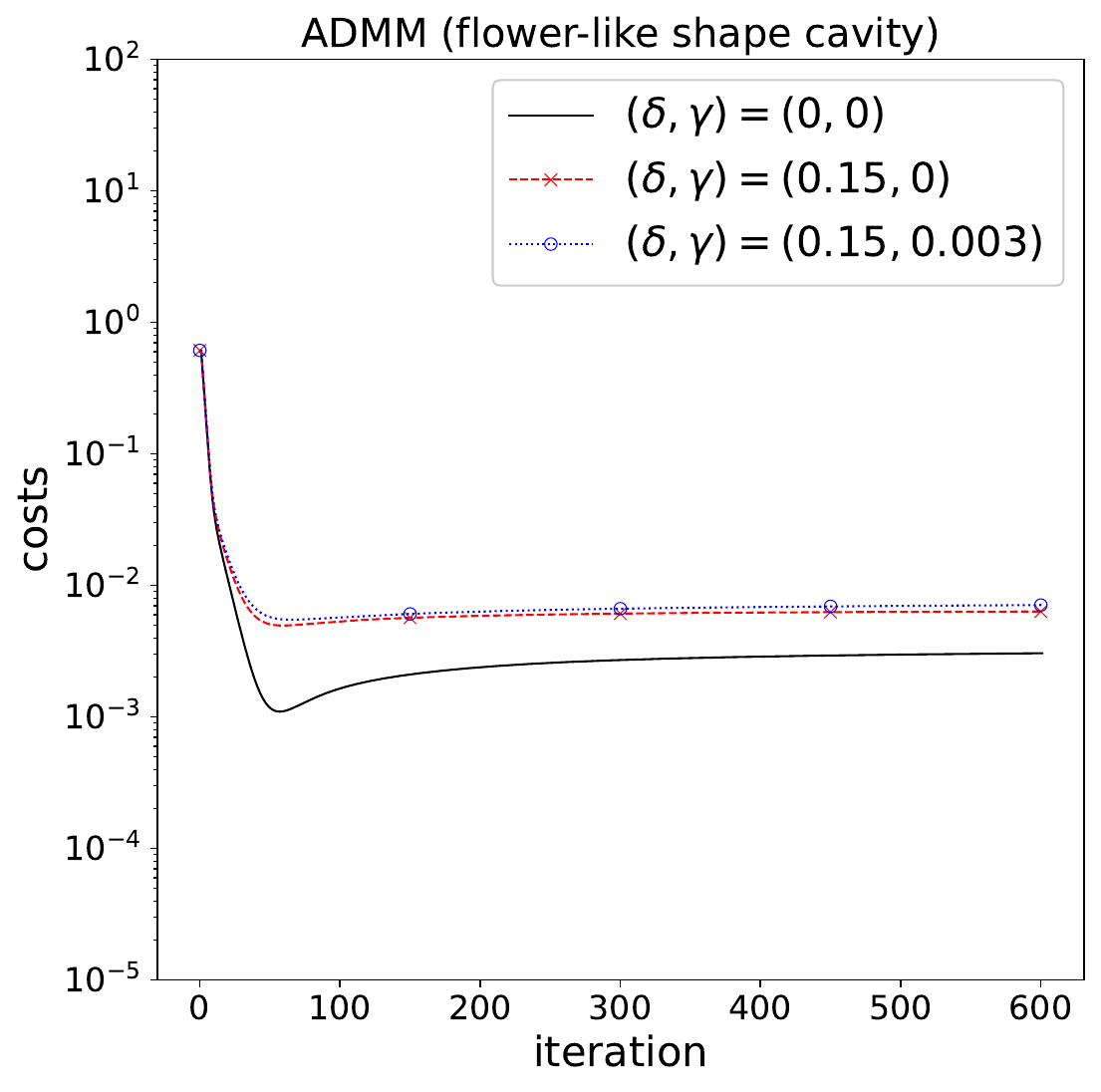}} \hfill
\resizebox{0.24\linewidth}{!}{\includegraphics{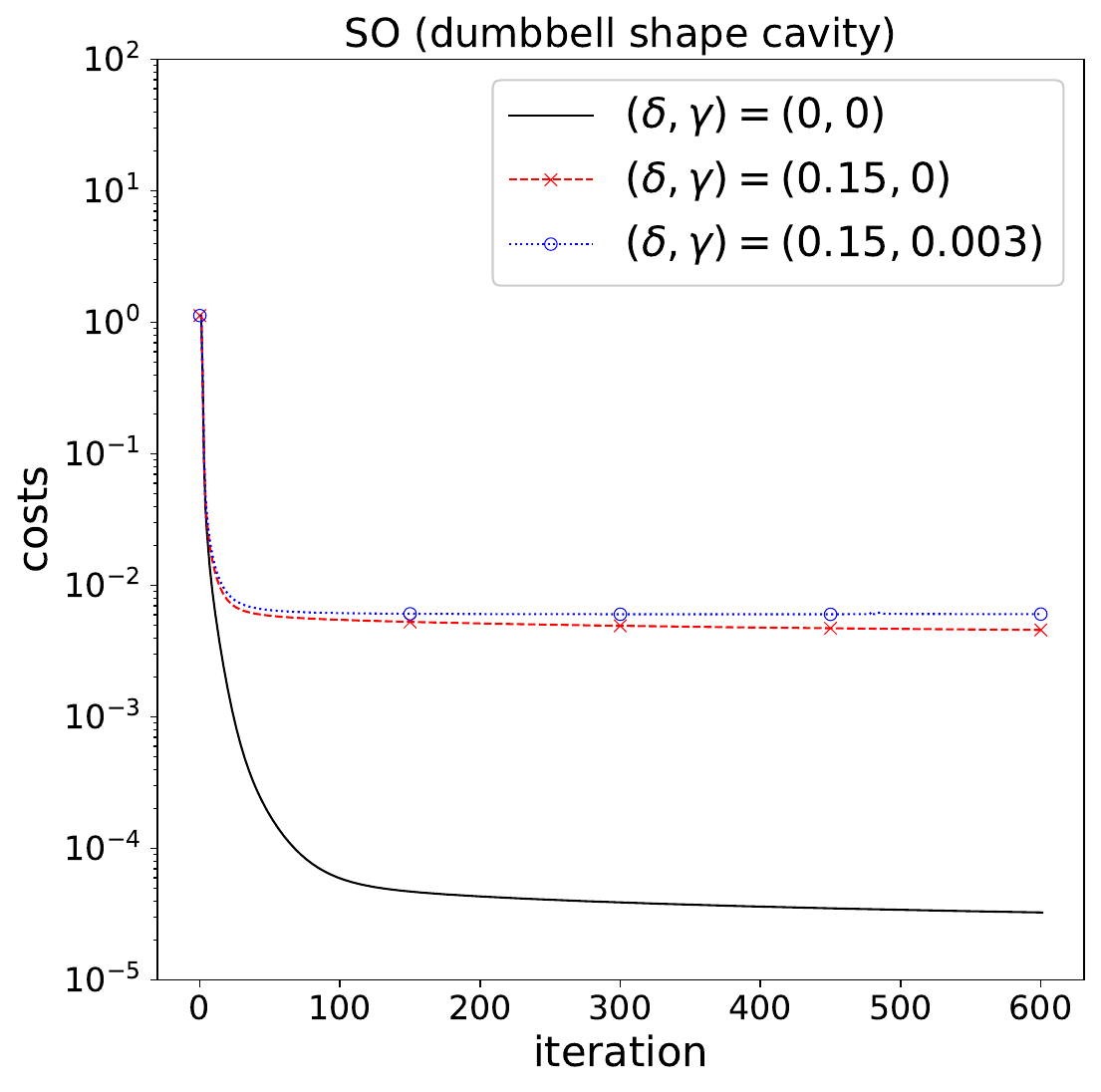}} \hfill
\resizebox{0.24\linewidth}{!}{\includegraphics{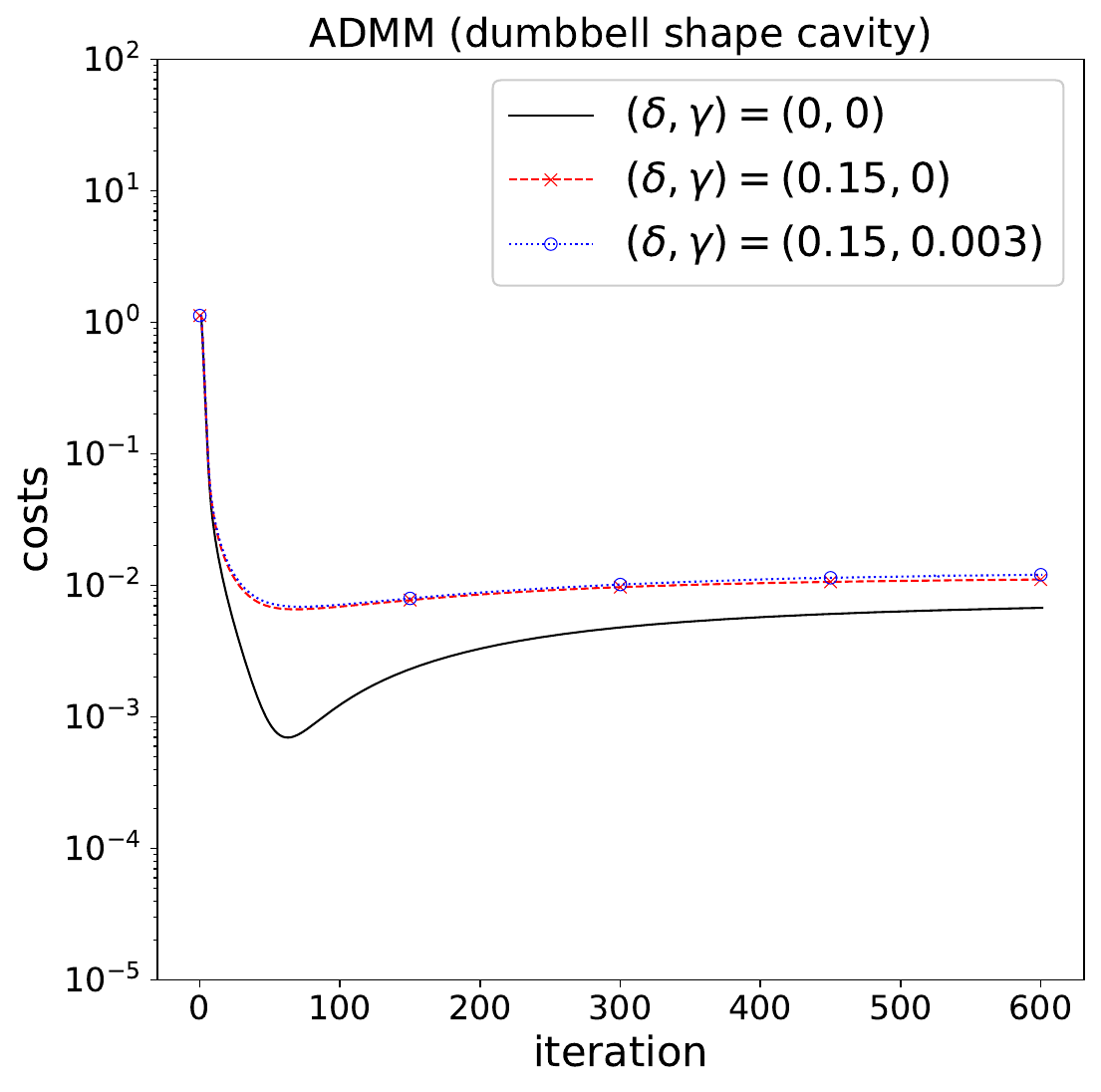}}
 \\[1em] 
\resizebox{0.24\linewidth}{!}{\includegraphics{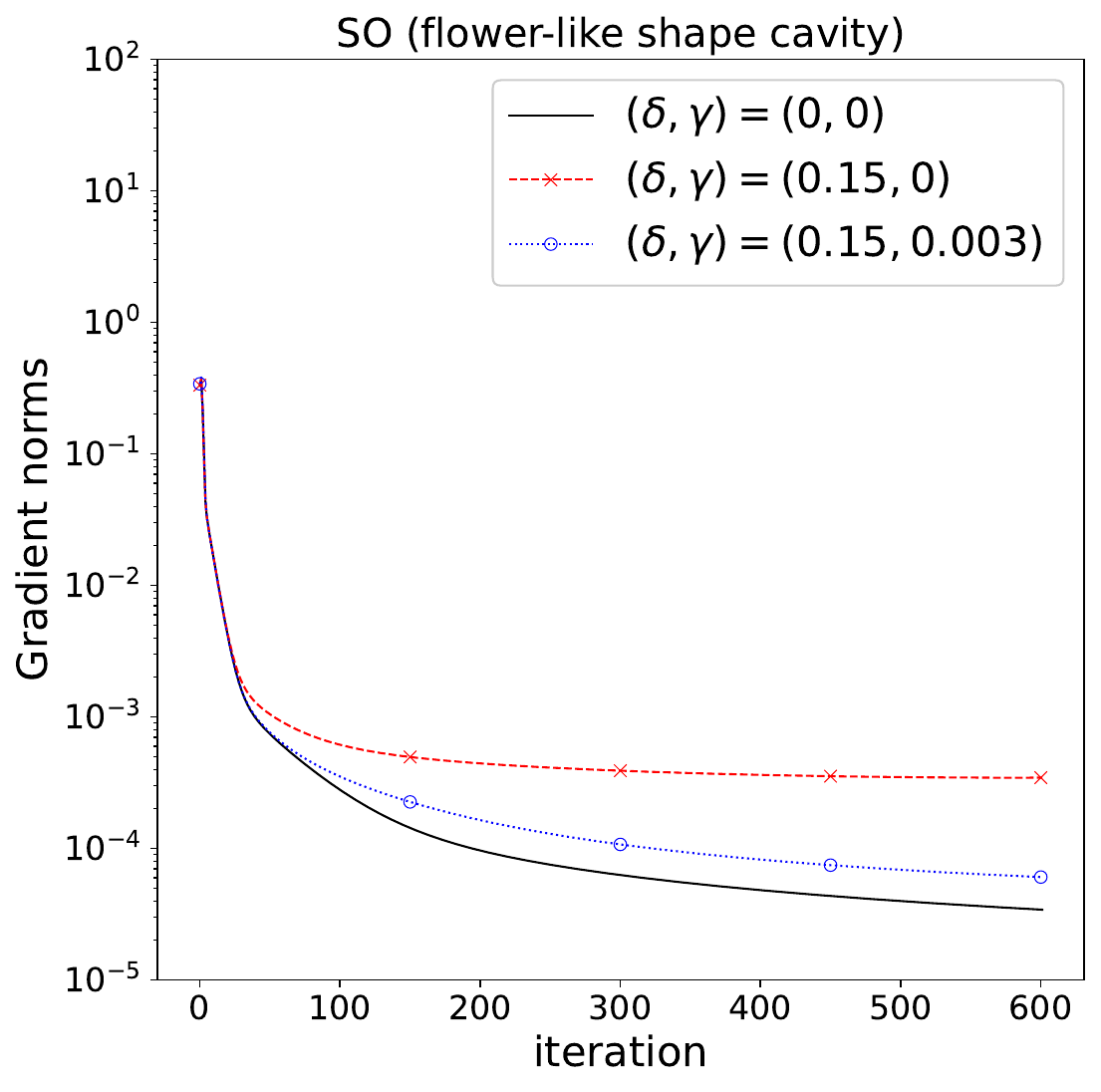}} \hfill
\resizebox{0.24\linewidth}{!}{\includegraphics{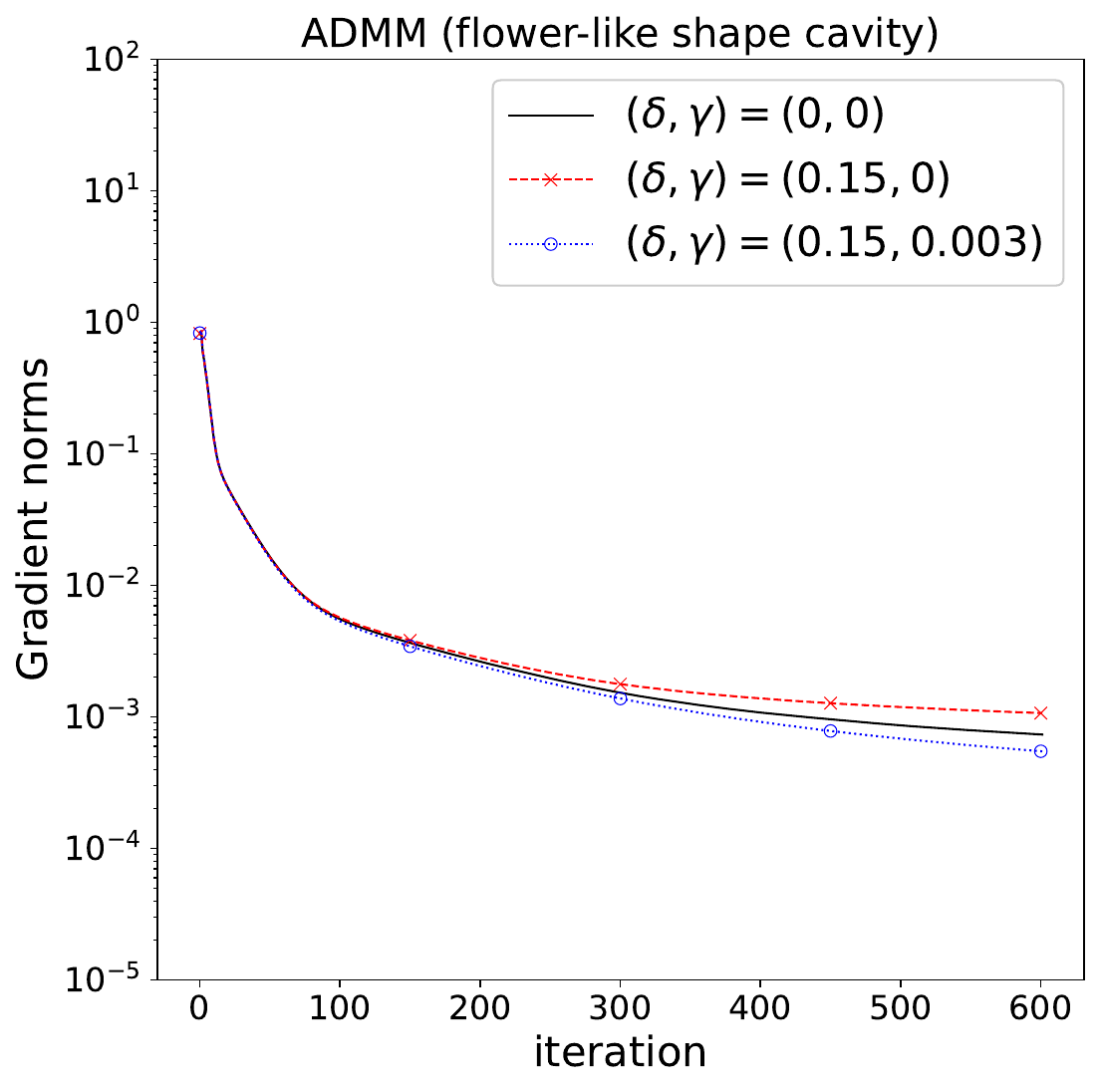}} \hfill
\resizebox{0.24\linewidth}{!}{\includegraphics{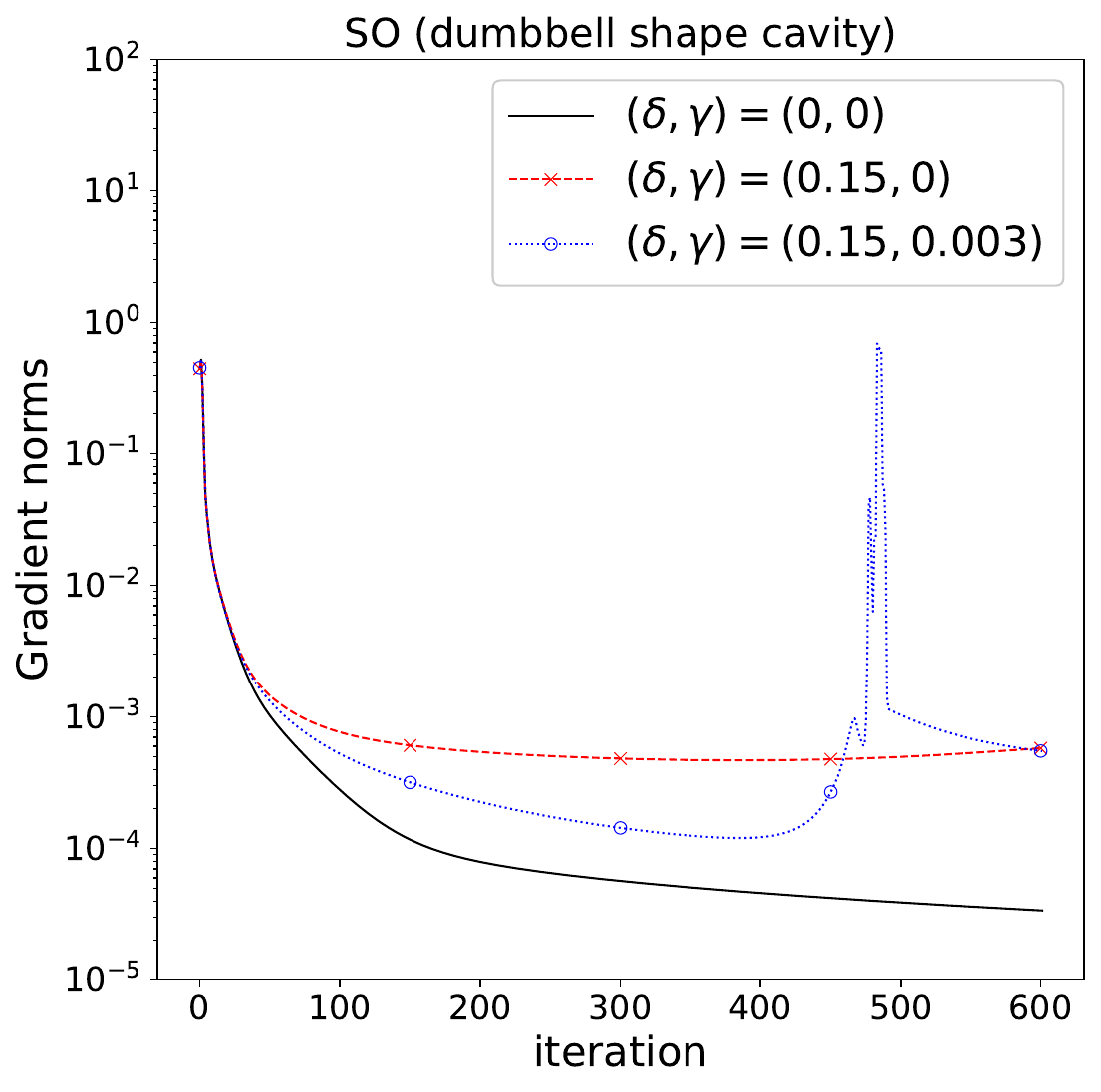}} \hfill
\resizebox{0.24\linewidth}{!}{\includegraphics{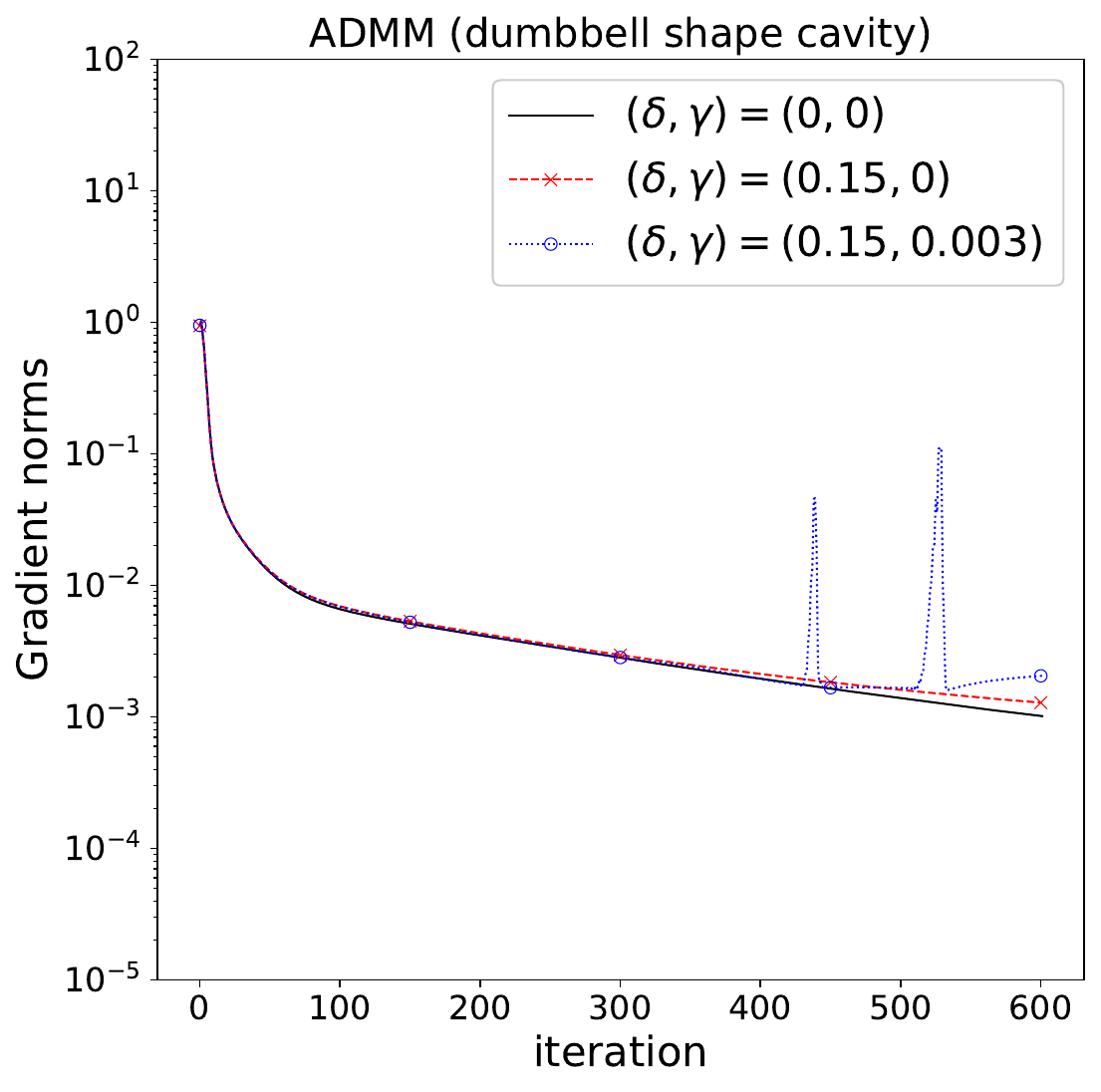}}
\caption{Histories of costs (top/first row) and gradient norms (bottom/second row) via SO and via ADMM}
\label{fig:figure2f}
\end{figure}

\begin{figure}[htp!]
\centering
\resizebox{0.25\linewidth}{!}{\includegraphics{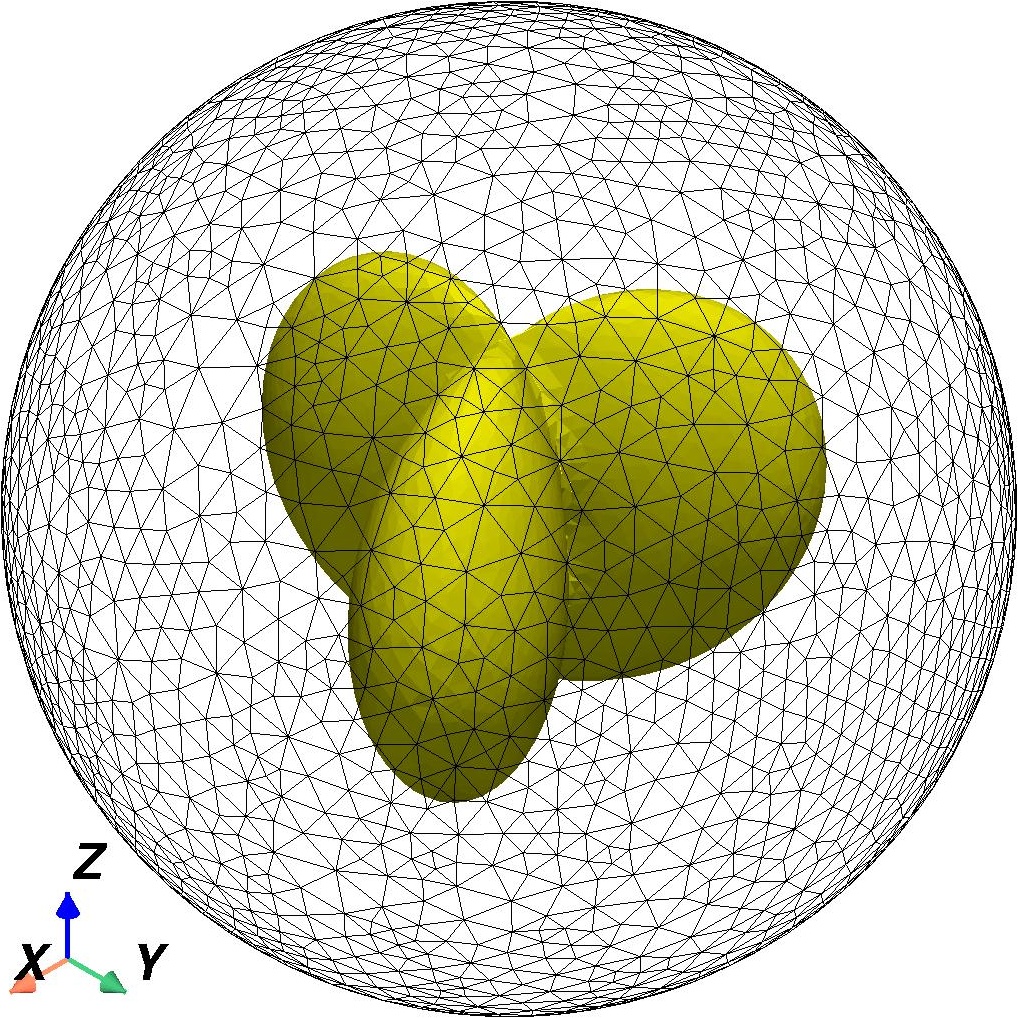}} \quad
\resizebox{0.25\linewidth}{!}{\includegraphics{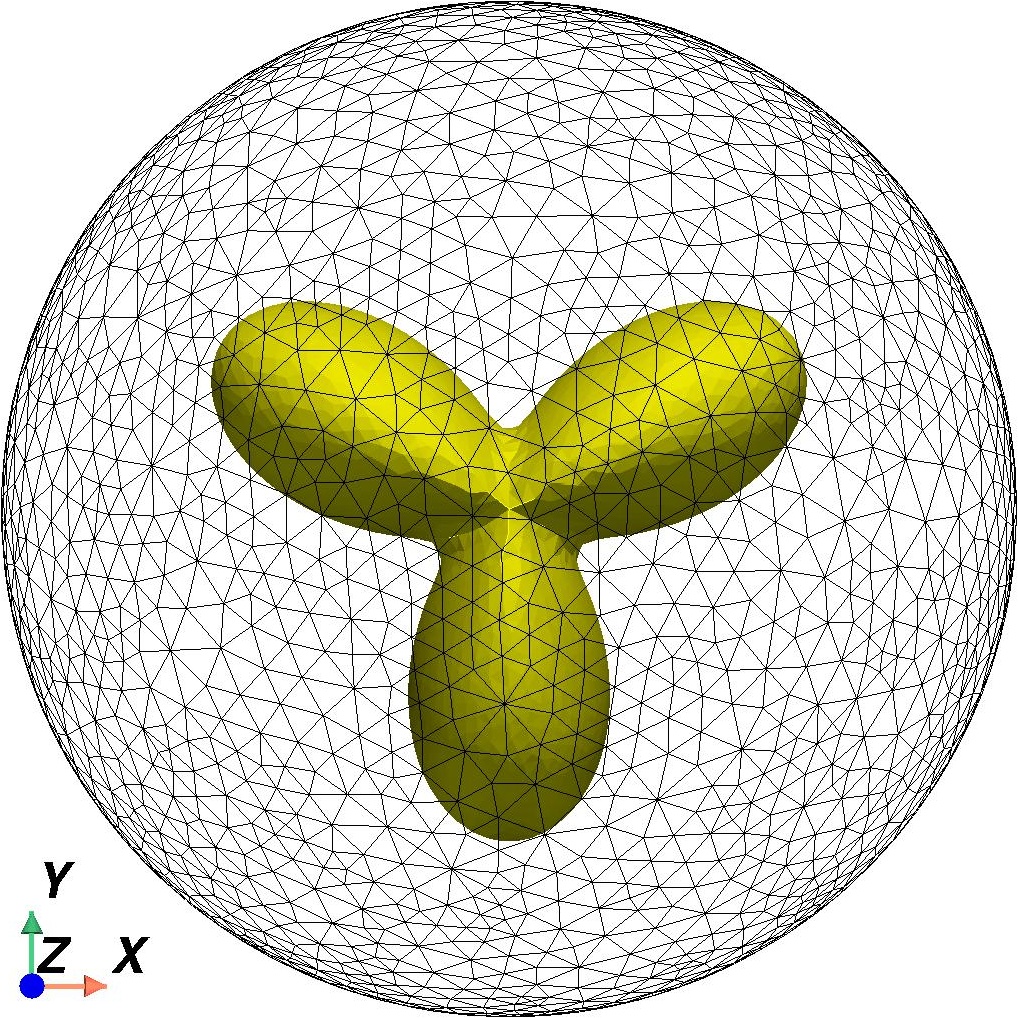}} \quad
\resizebox{0.25\linewidth}{!}{\includegraphics{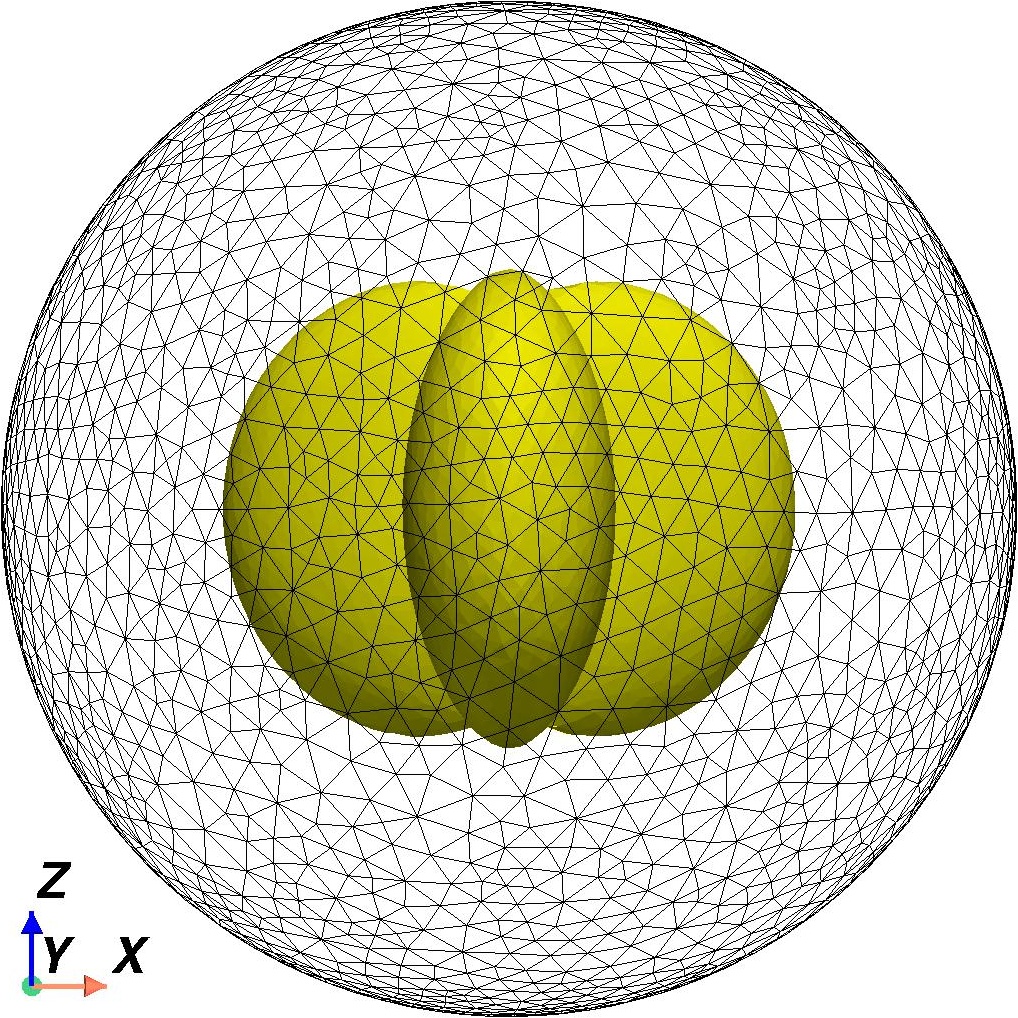}}
\caption{Exact geometries of the smaller size cavity}
\label{fig:figure3a}
\end{figure}

\begin{figure}[htp!]
\centering
\resizebox{0.16\linewidth}{!}{\includegraphics{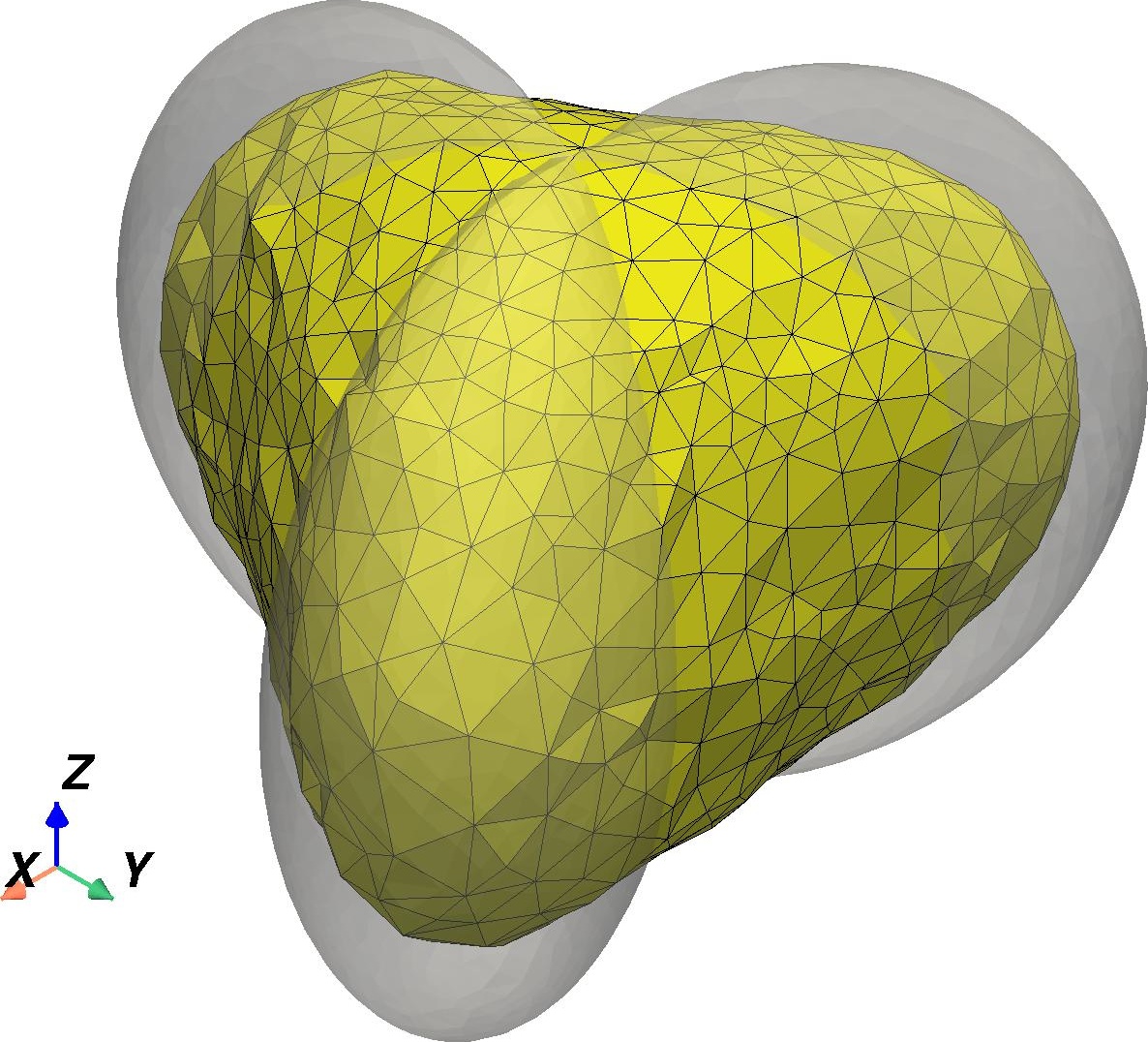}} 
\resizebox{0.16\linewidth}{!}{\includegraphics{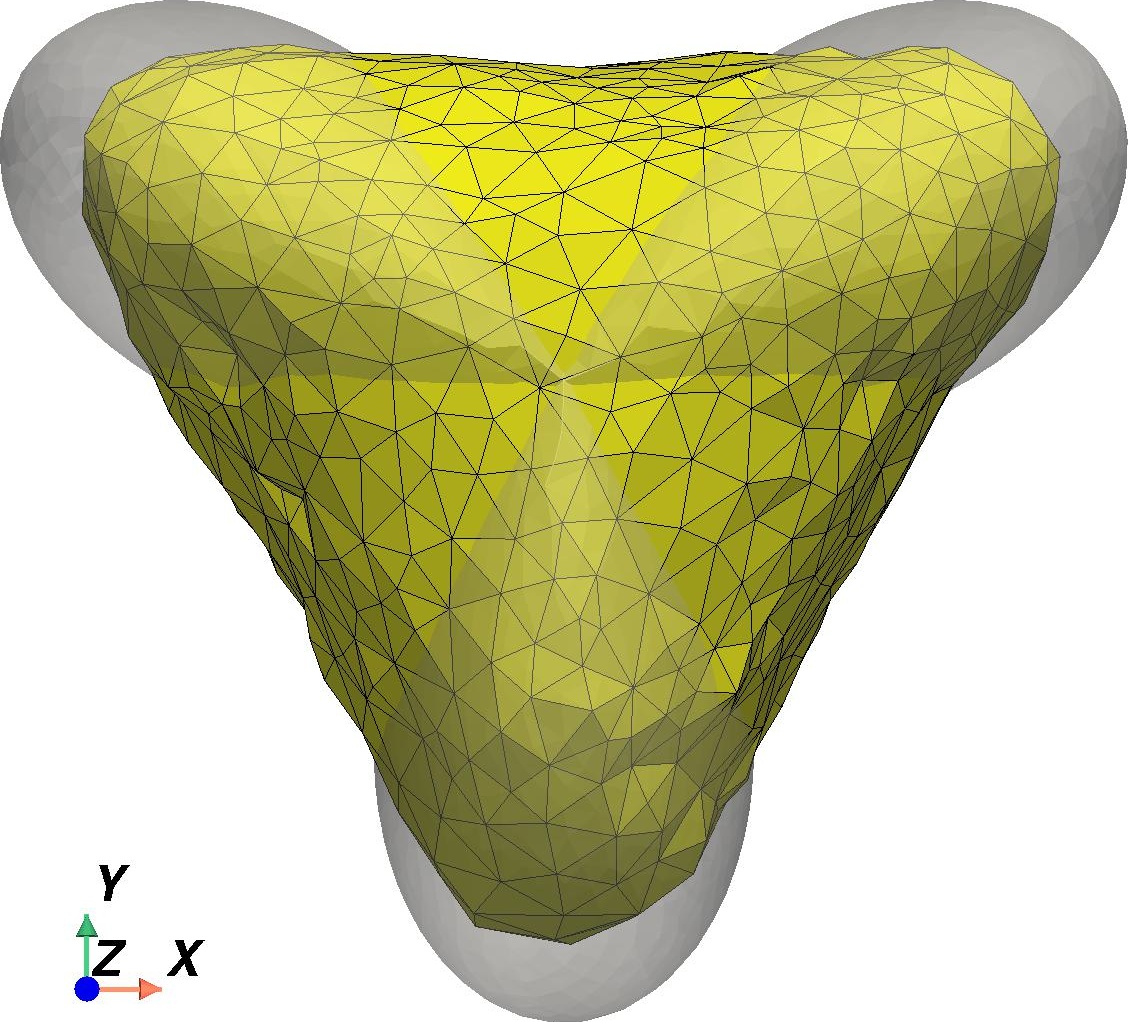}} 
\resizebox{0.16\linewidth}{!}{\includegraphics{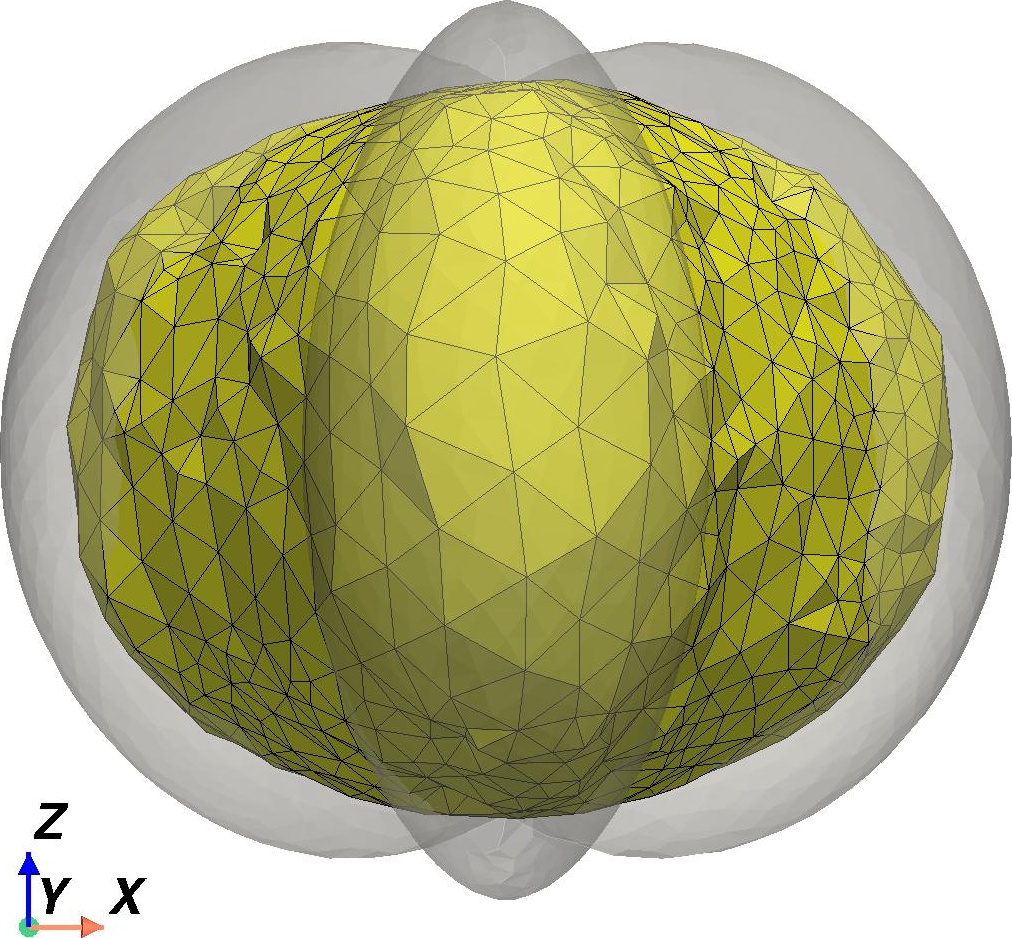}}\hfill
\resizebox{0.16\linewidth}{!}{\includegraphics{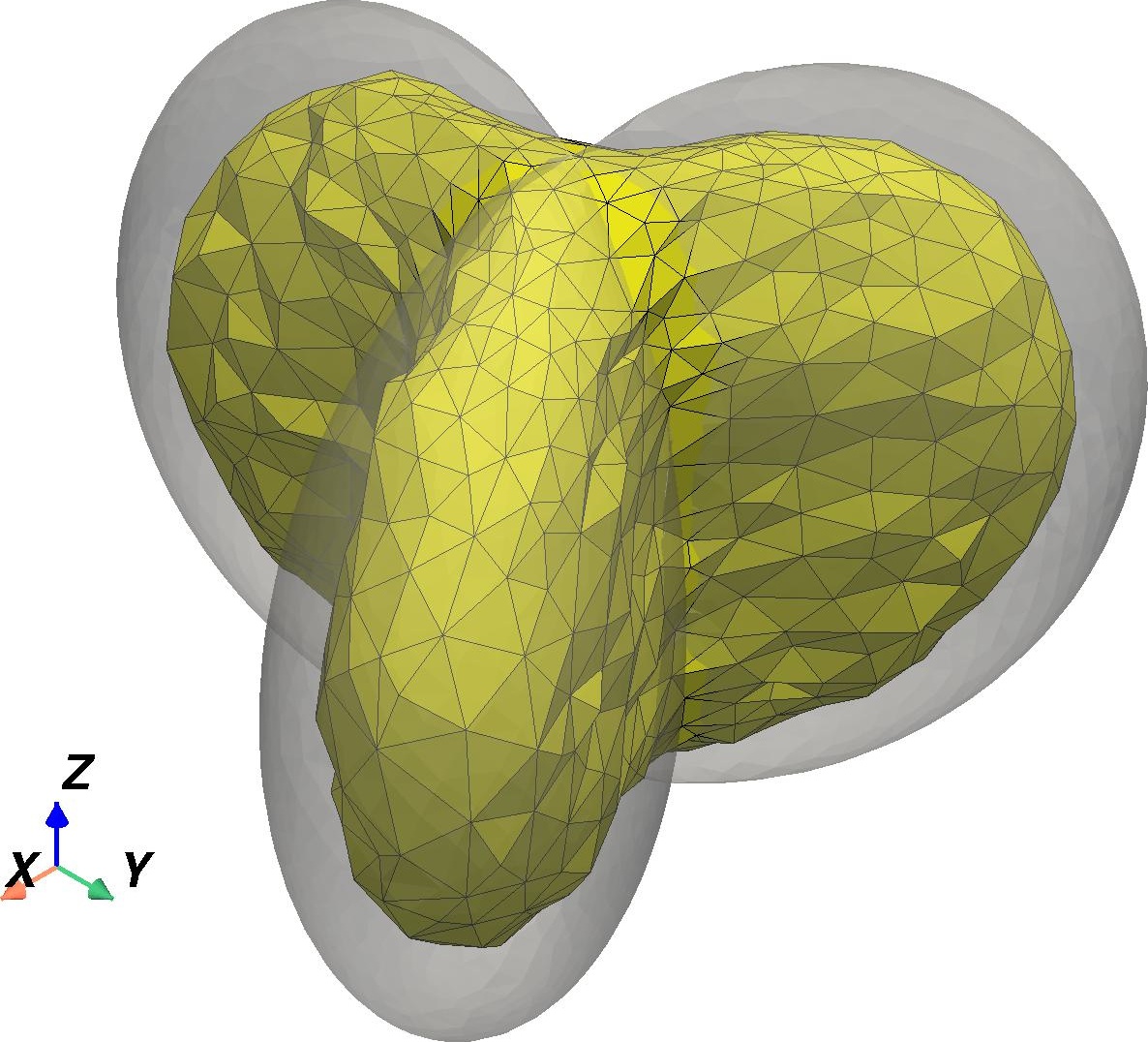}} 
\resizebox{0.16\linewidth}{!}{\includegraphics{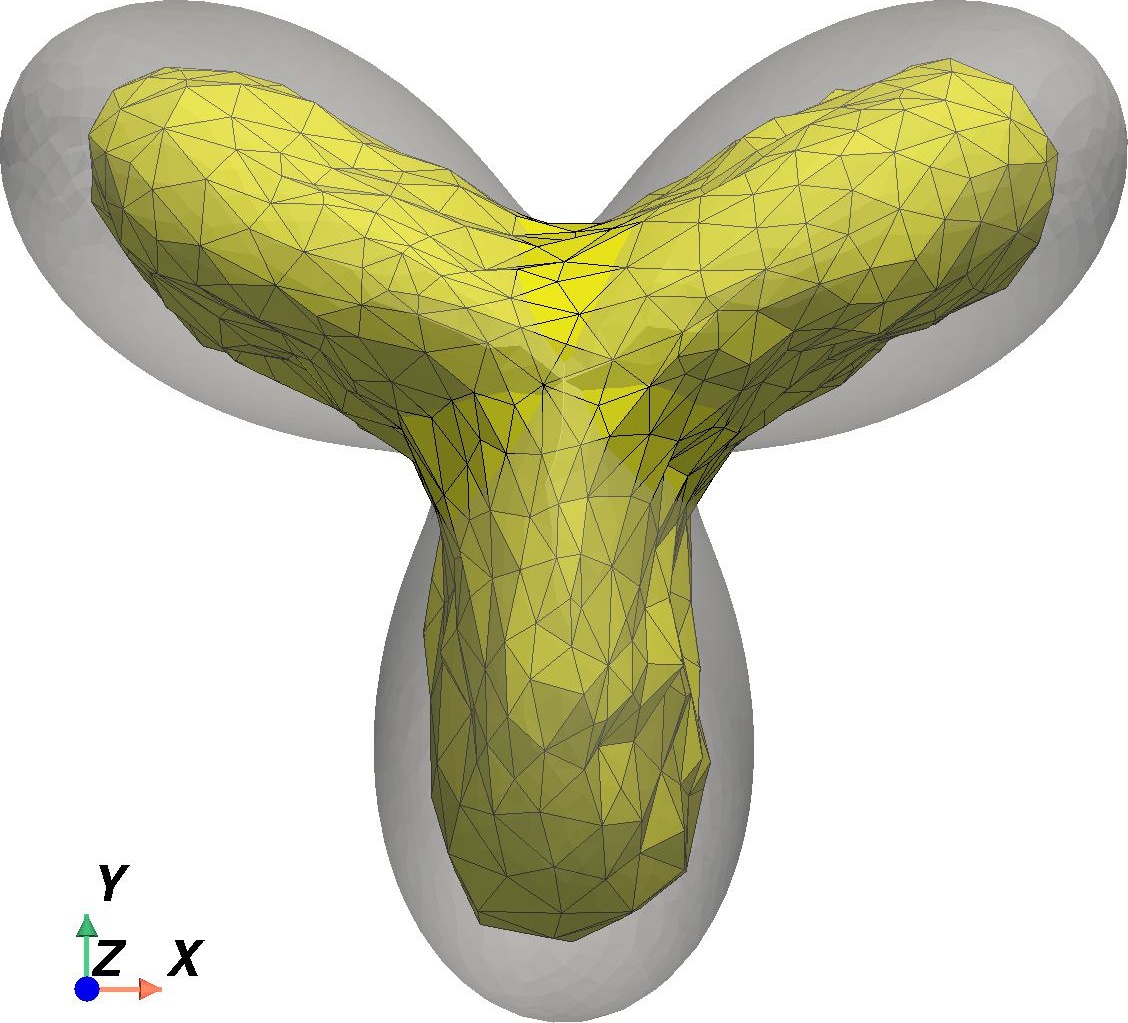}} 
\resizebox{0.16\linewidth}{!}{\includegraphics{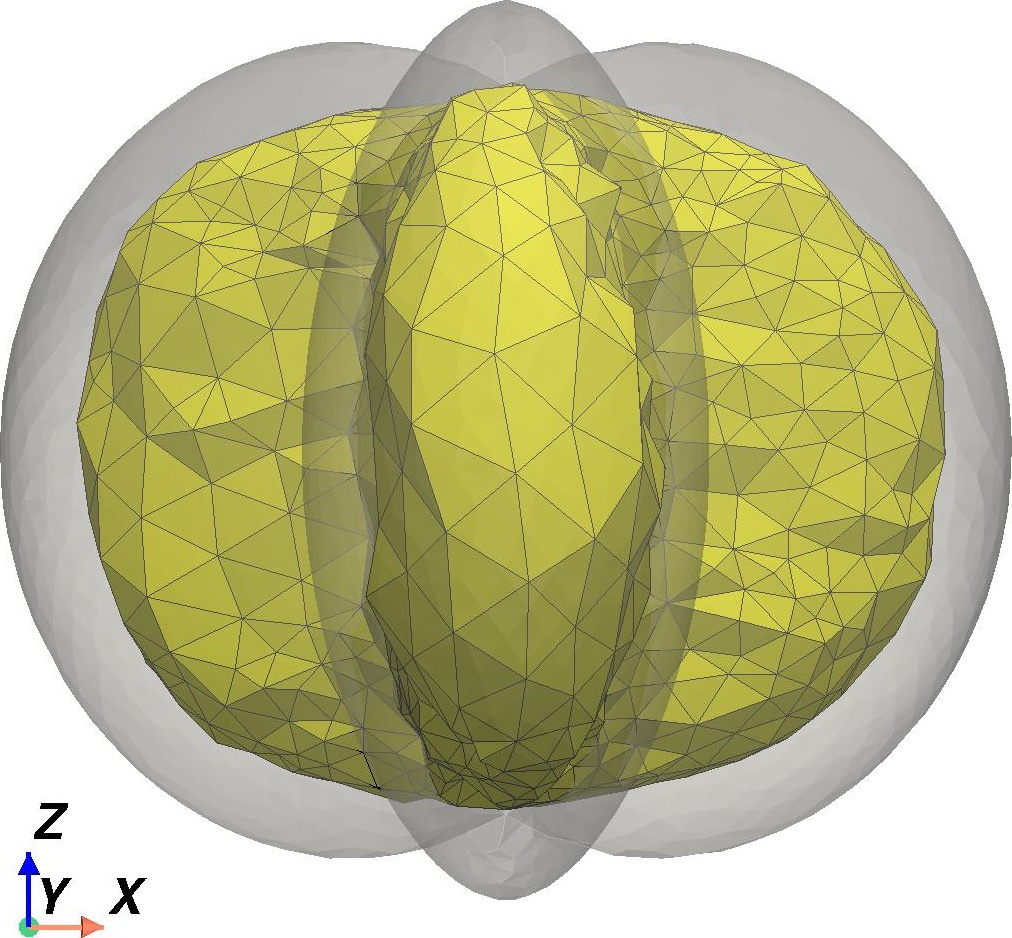}}
 \\[1em] 
\resizebox{0.16\linewidth}{!}{\includegraphics{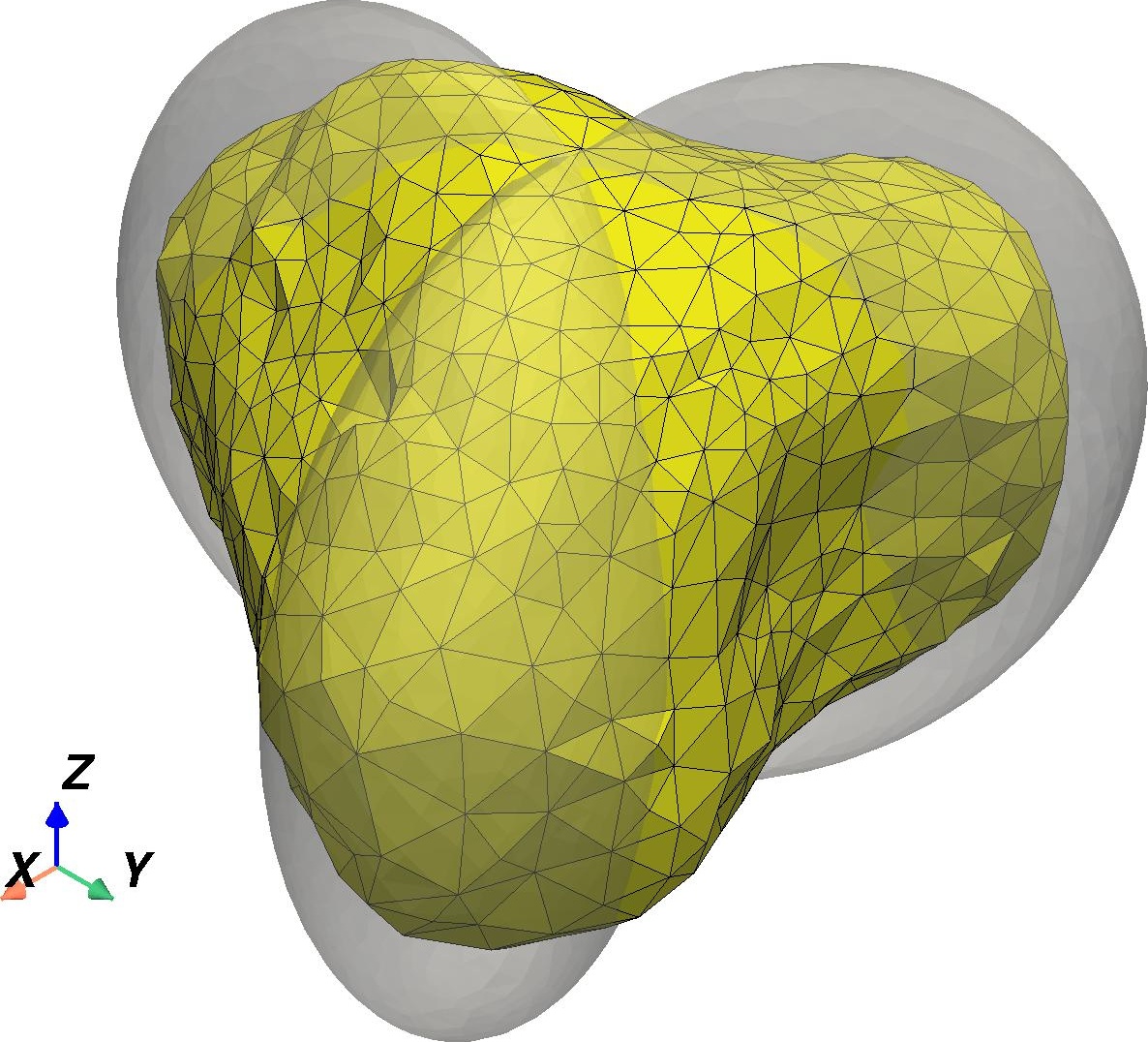}} 
\resizebox{0.16\linewidth}{!}{\includegraphics{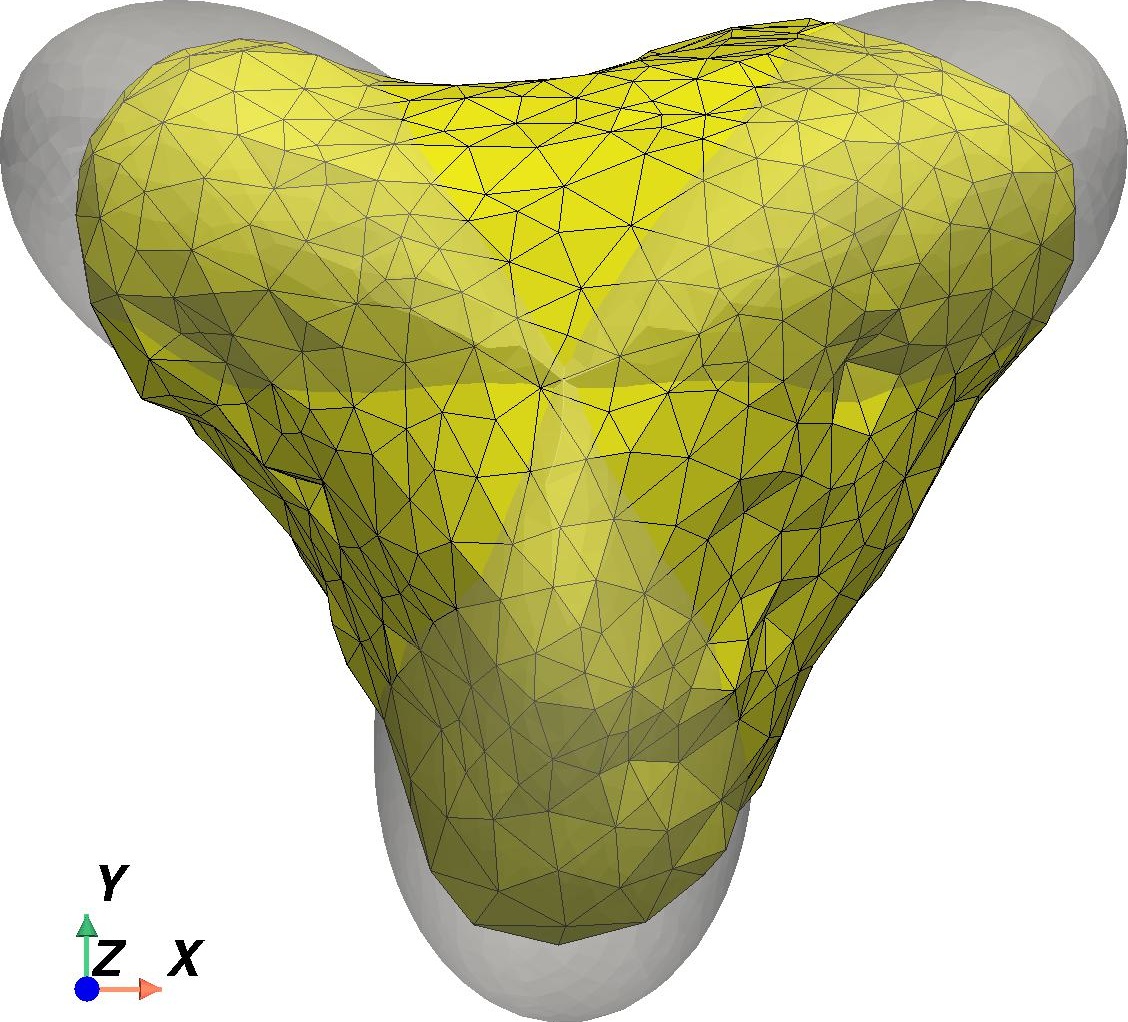}} 
\resizebox{0.16\linewidth}{!}{\includegraphics{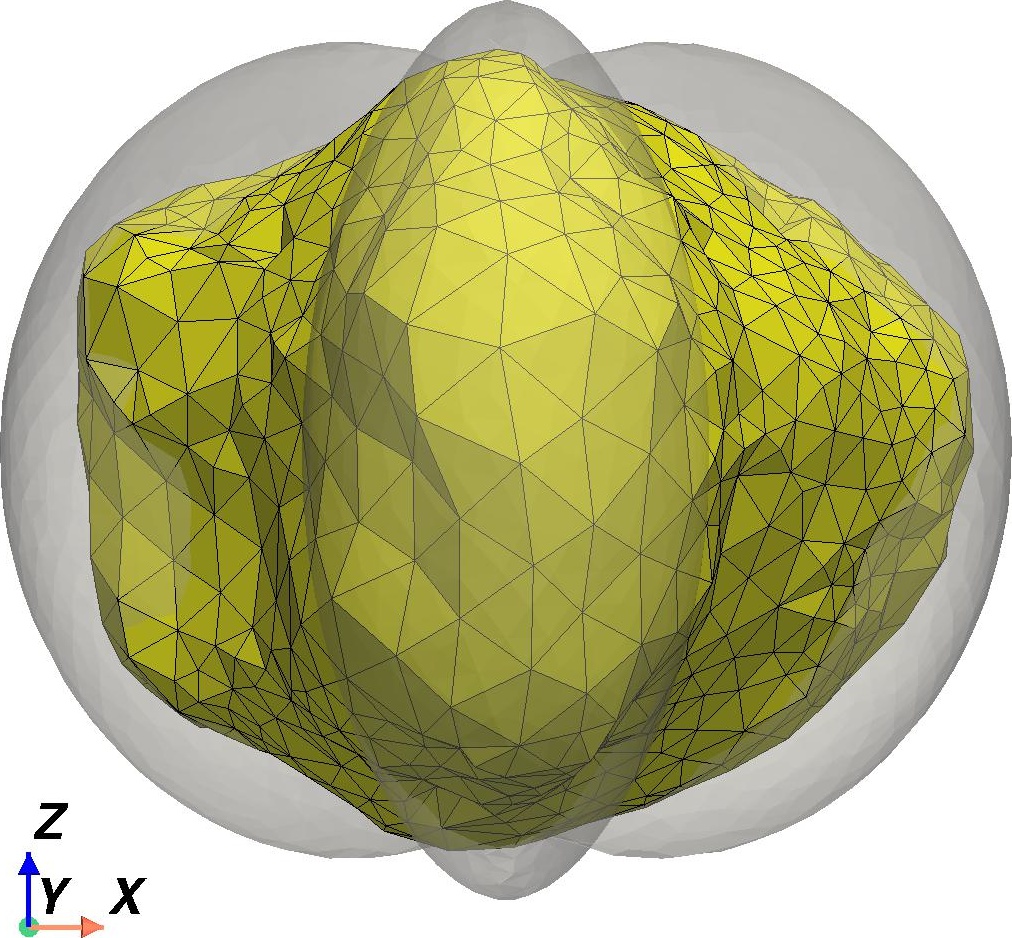}} \hfill
\resizebox{0.16\linewidth}{!}{\includegraphics{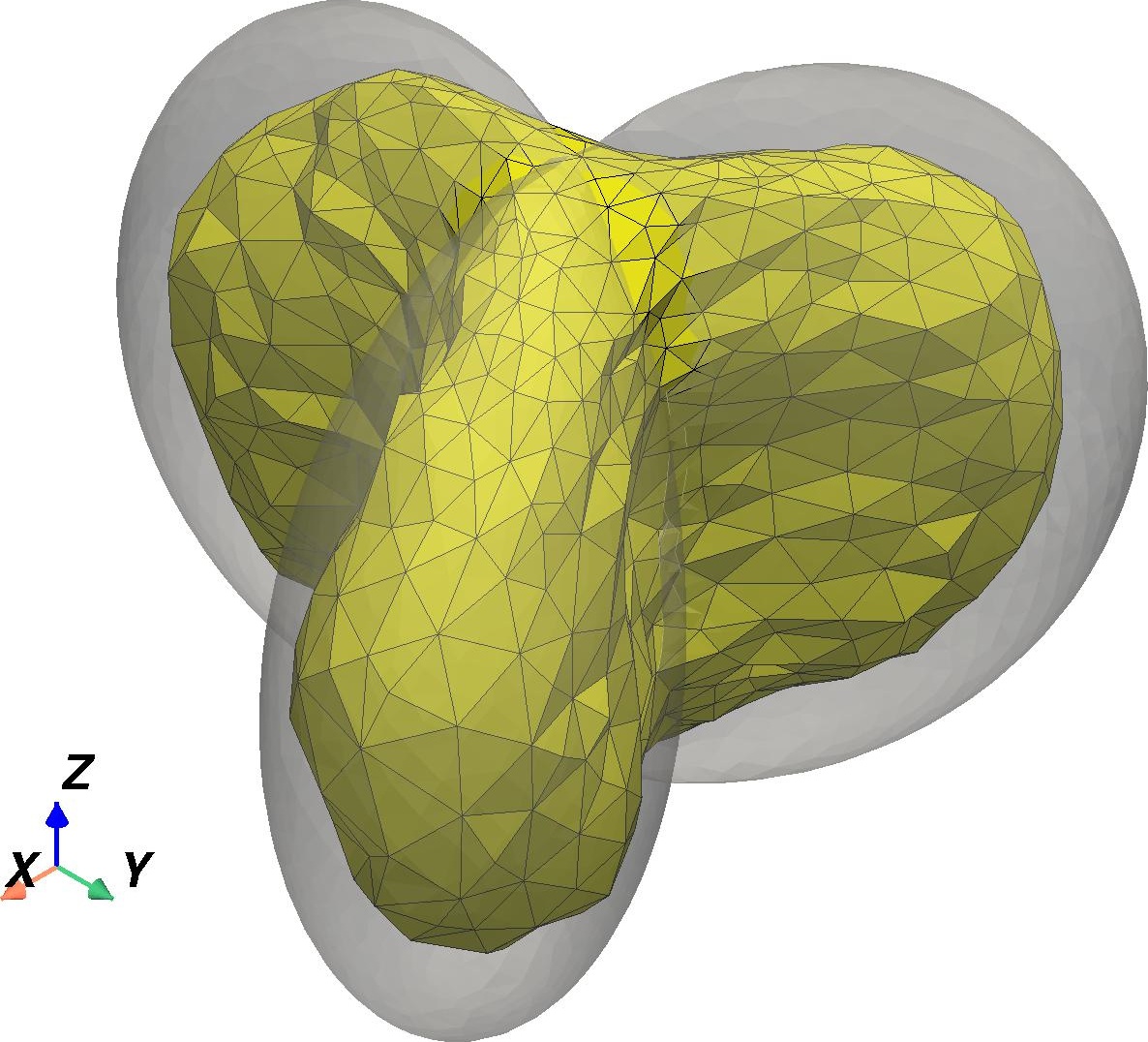}} 
\resizebox{0.16\linewidth}{!}{\includegraphics{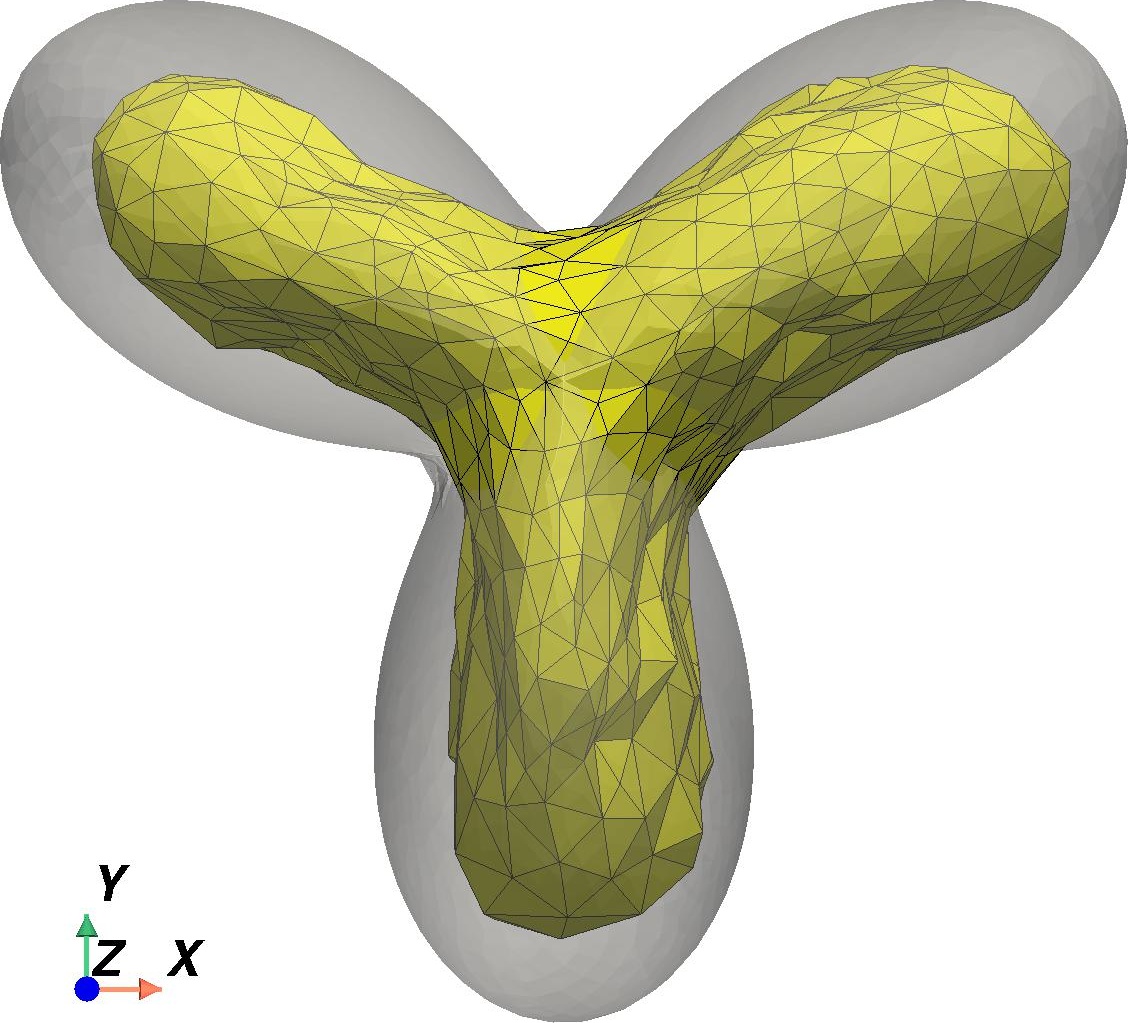}} 
\resizebox{0.16\linewidth}{!}{\includegraphics{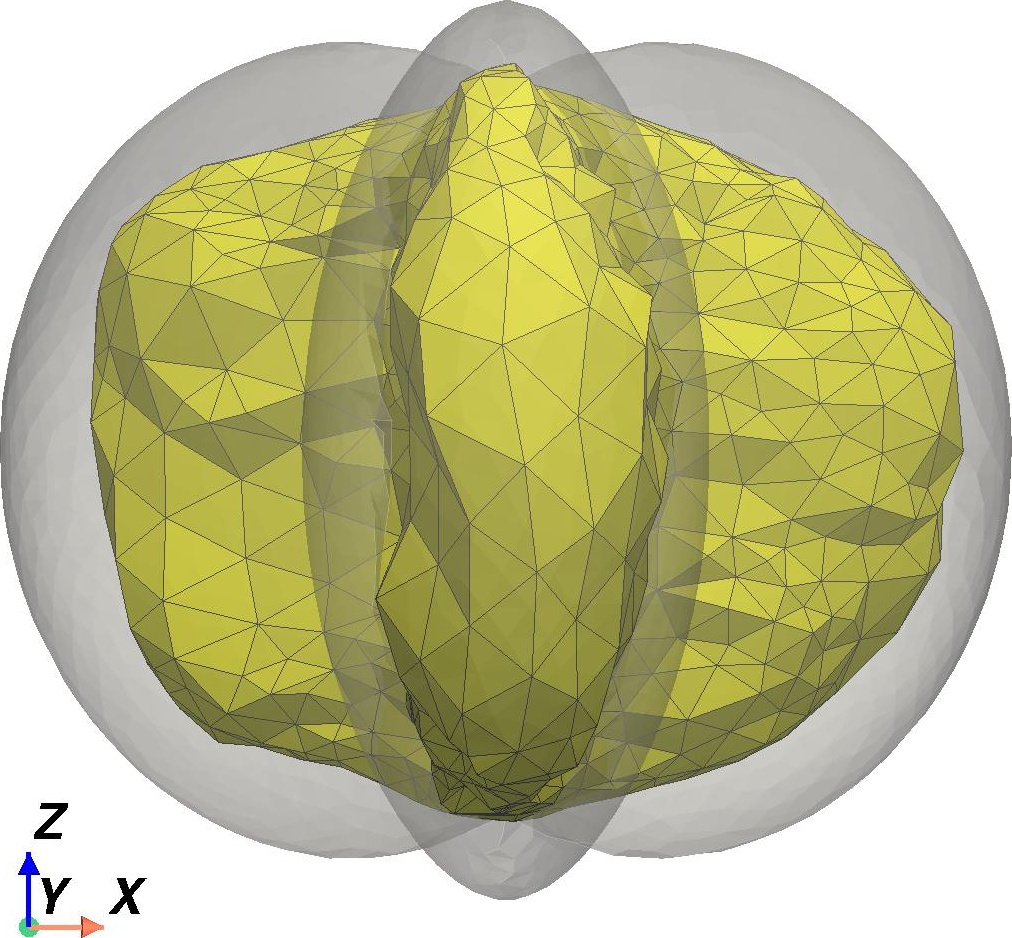}}
 \\[1em] 
\resizebox{0.16\linewidth}{!}{\includegraphics{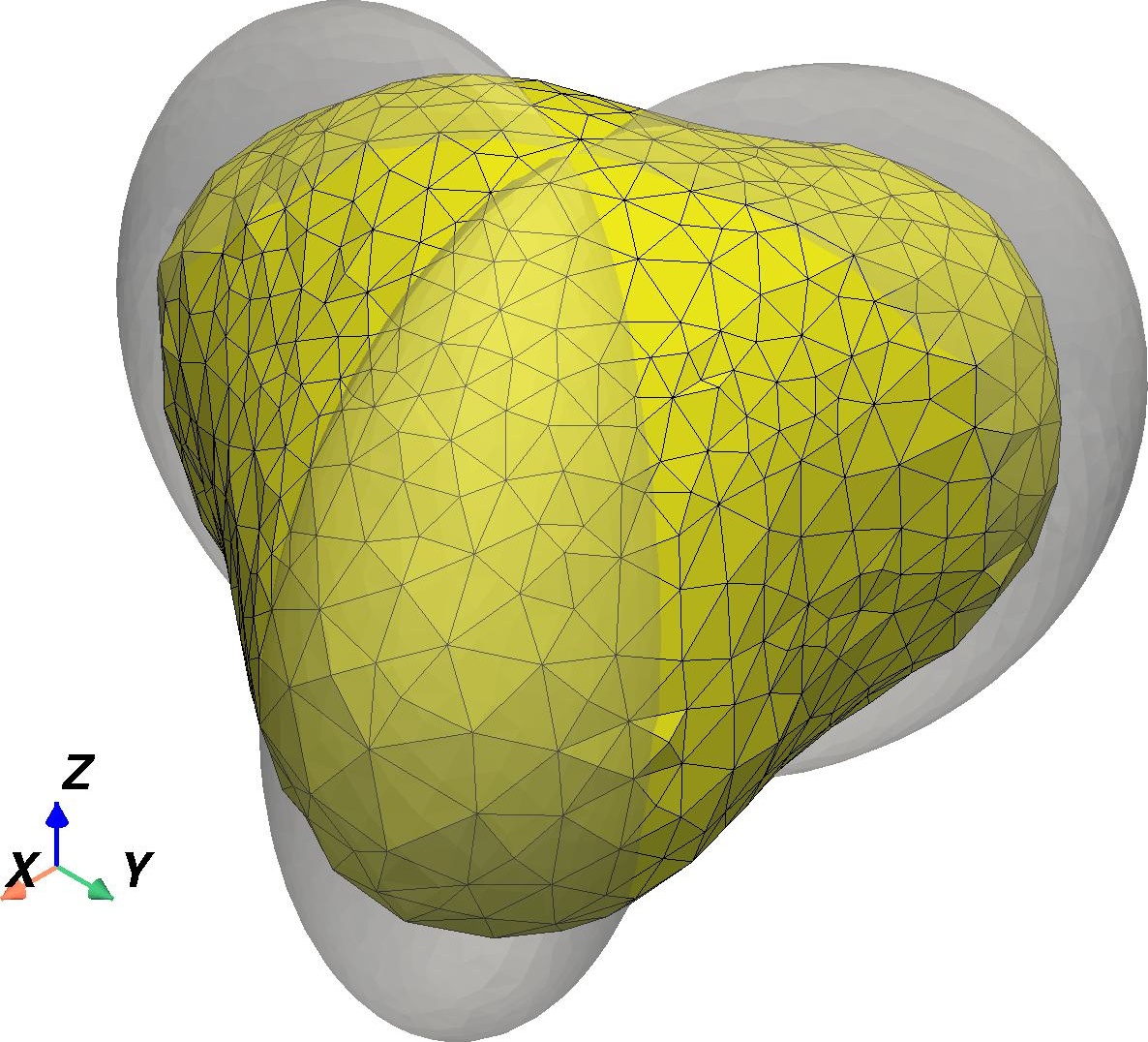}} 
\resizebox{0.16\linewidth}{!}{\includegraphics{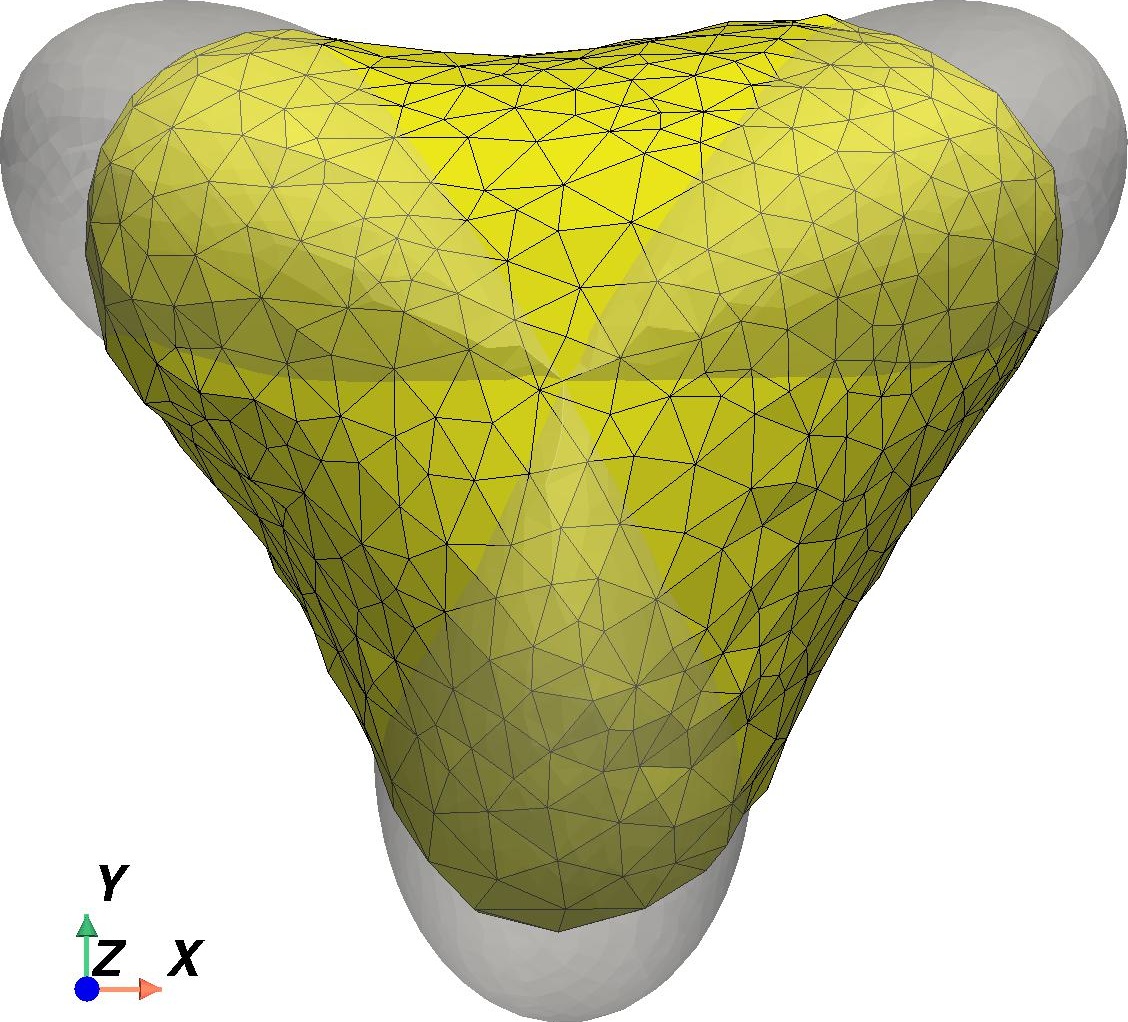}} 
\resizebox{0.16\linewidth}{!}{\includegraphics{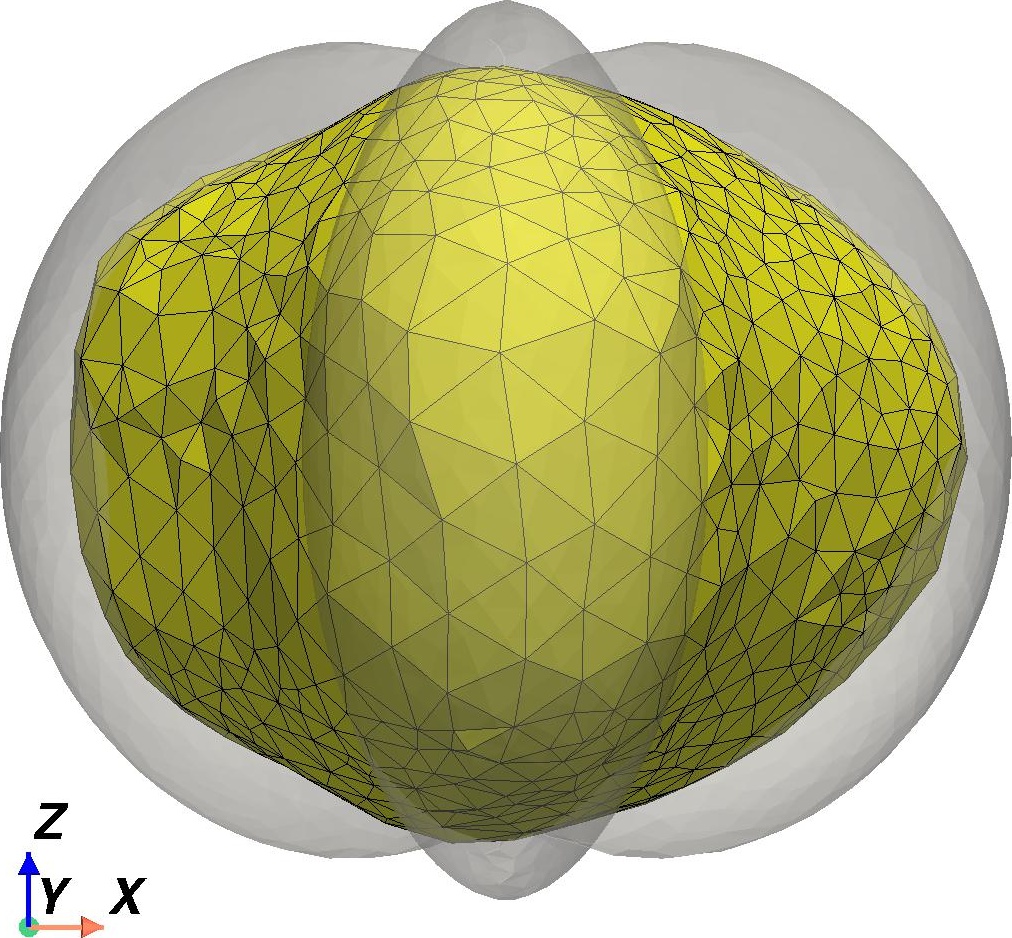}} \hfill
\resizebox{0.16\linewidth}{!}{\includegraphics{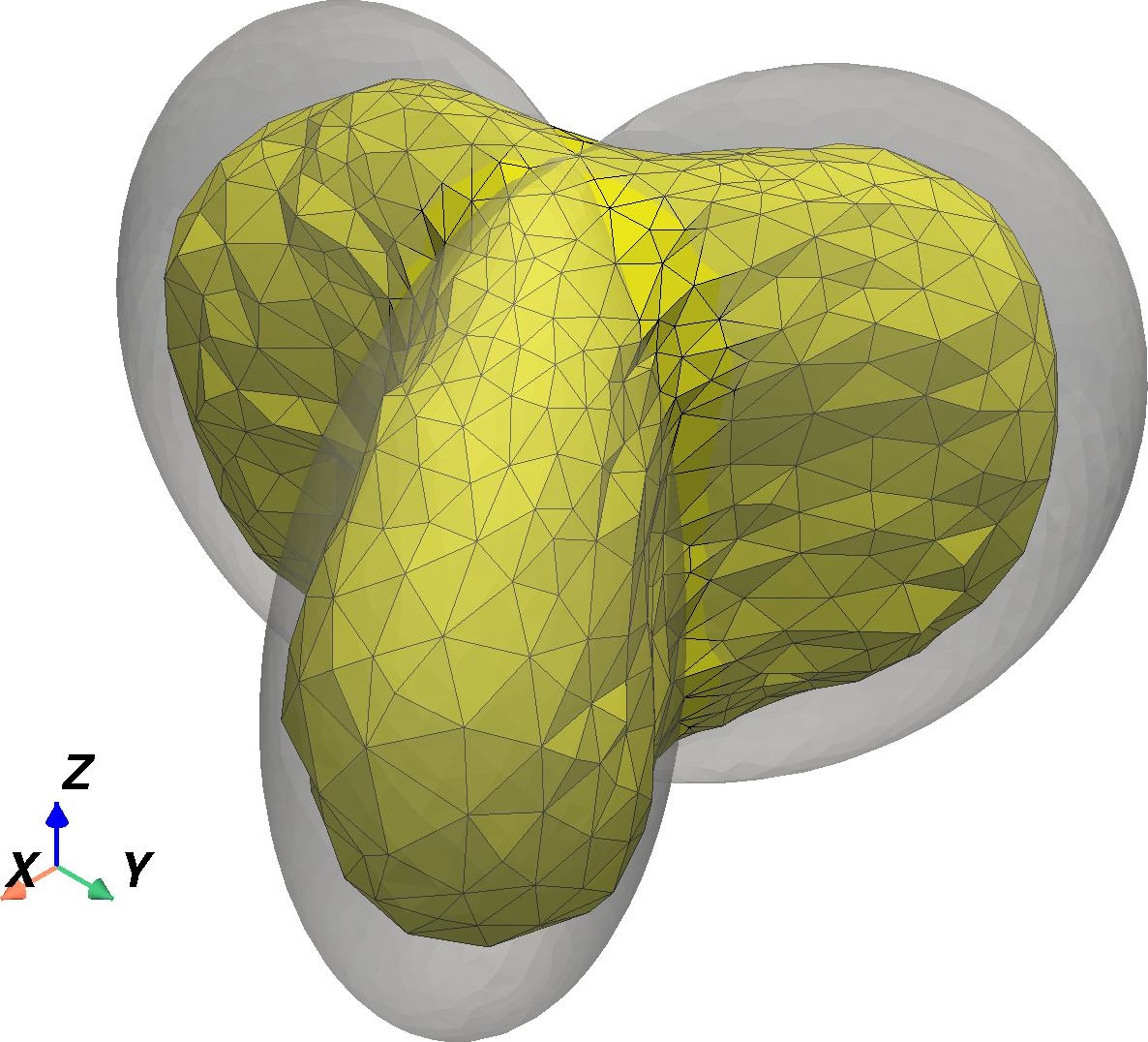}} 
\resizebox{0.16\linewidth}{!}{\includegraphics{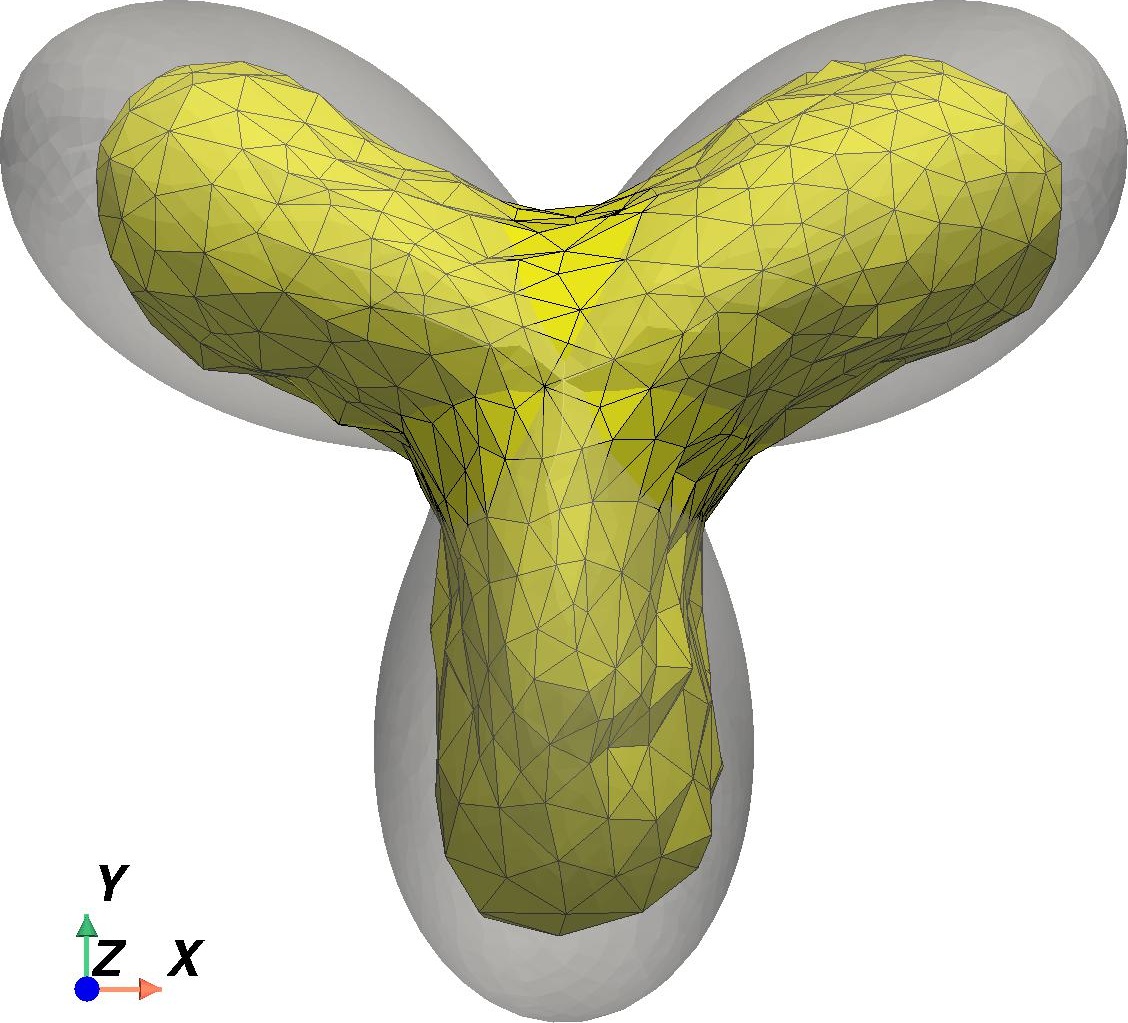}} 
\resizebox{0.16\linewidth}{!}{\includegraphics{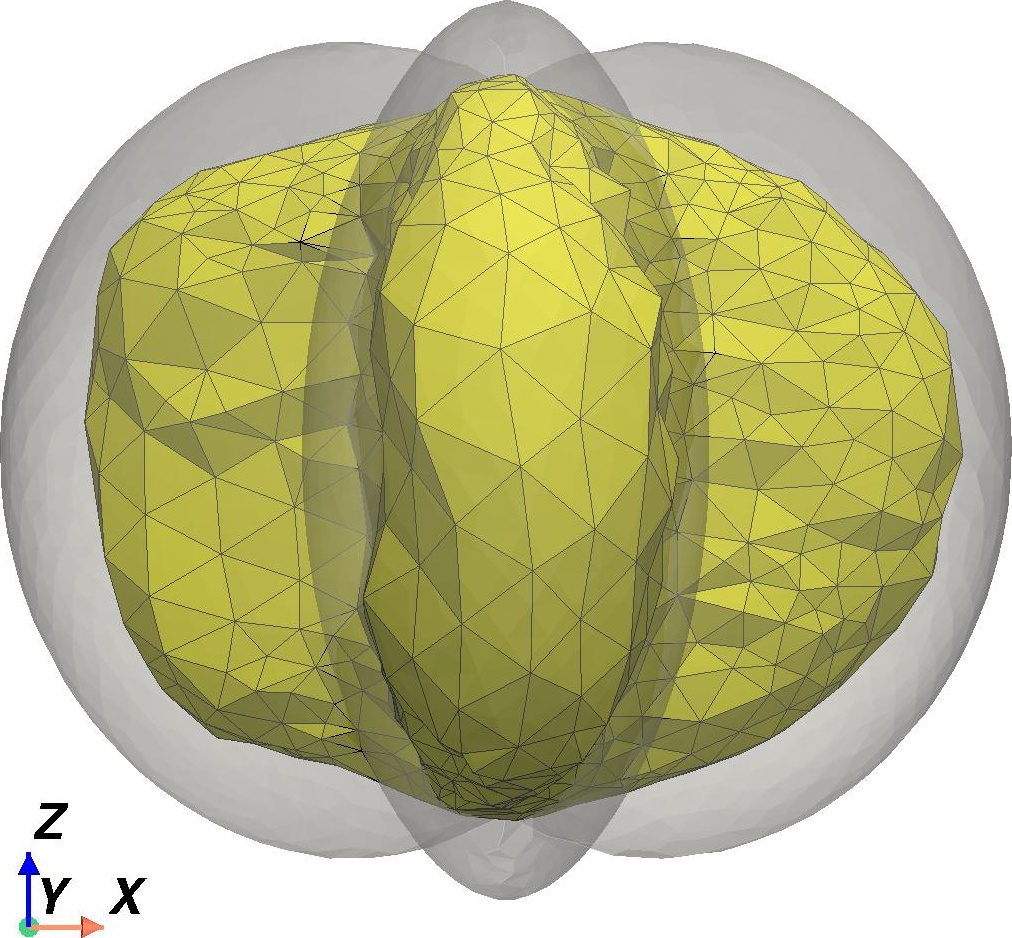}}
\caption{Reconstructions via SO (first three figures from the left) and ADMM (last three figures on the right) are presented with three different scenarios: exact data (first/top row), noisy data with a noise level of $\delta=15\%$ (second/middle row), and without regularization; as well as with regularization using $\gamma = 0.003$ (third/bottom row)}
\label{fig:figure3b}
\end{figure}

\begin{figure}[htp!]
\centering
\resizebox{0.24\linewidth}{!}{\includegraphics{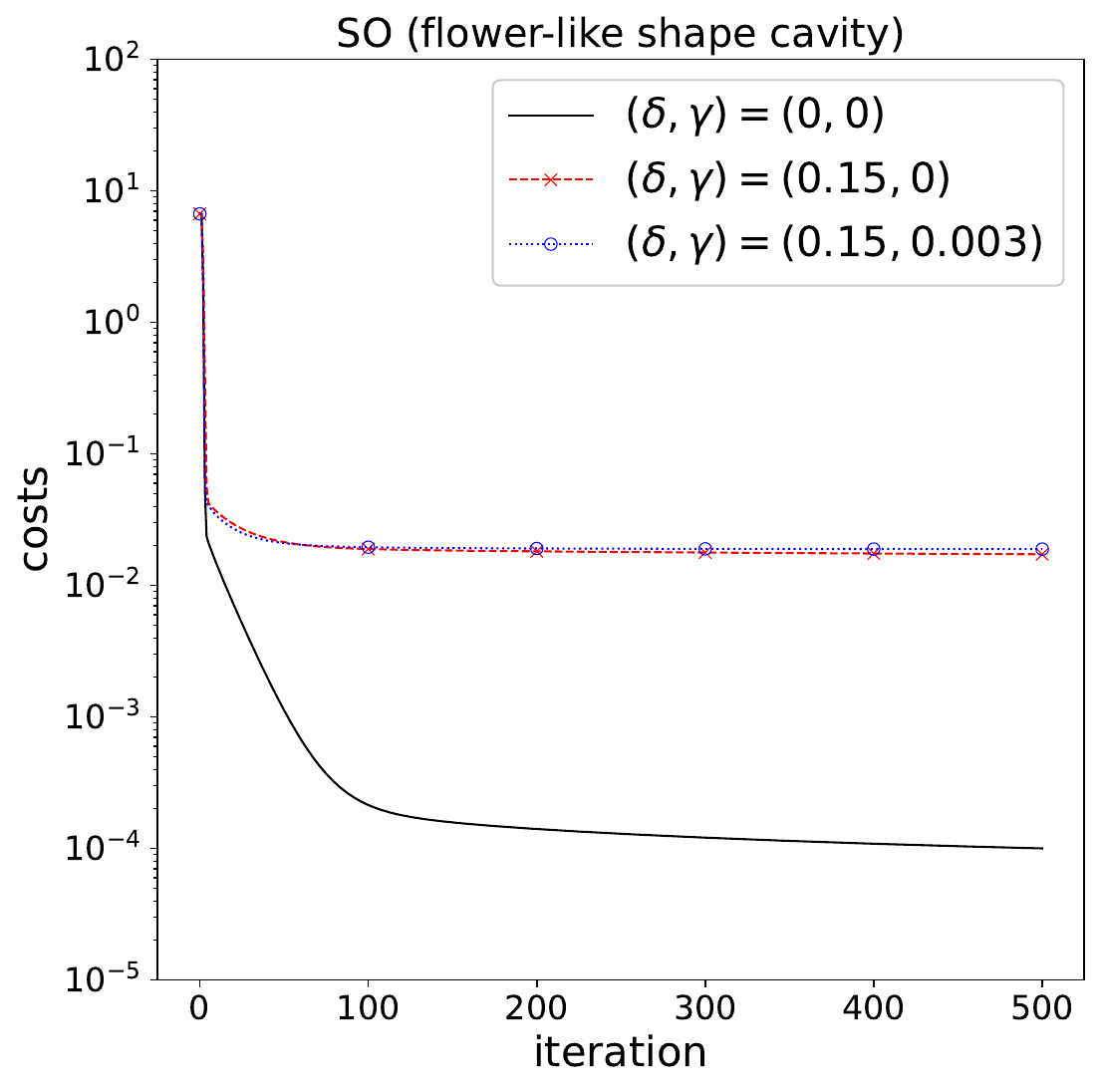}} \hfill
\resizebox{0.24\linewidth}{!}{\includegraphics{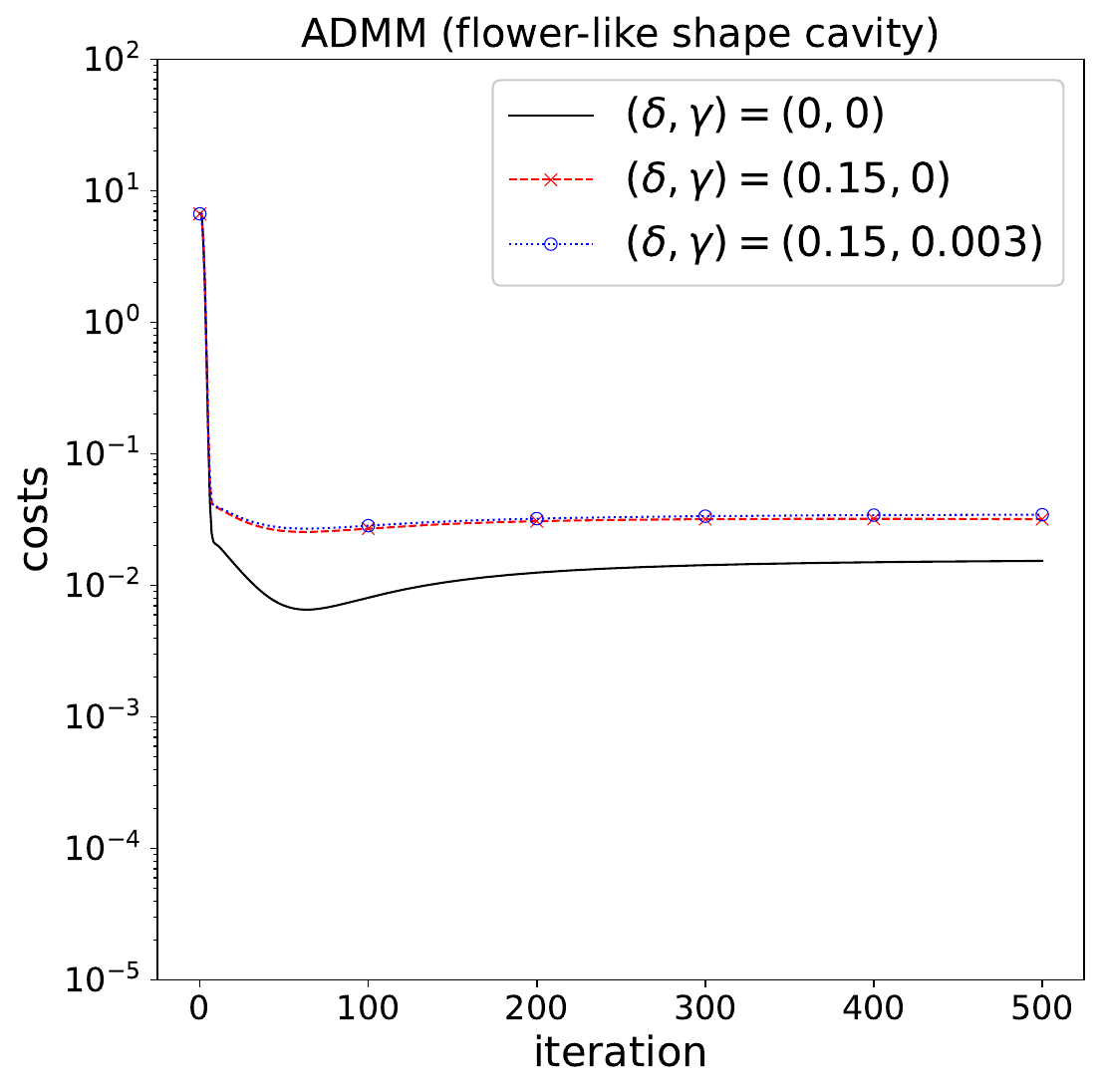}} \hfill
\resizebox{0.24\linewidth}{!}{\includegraphics{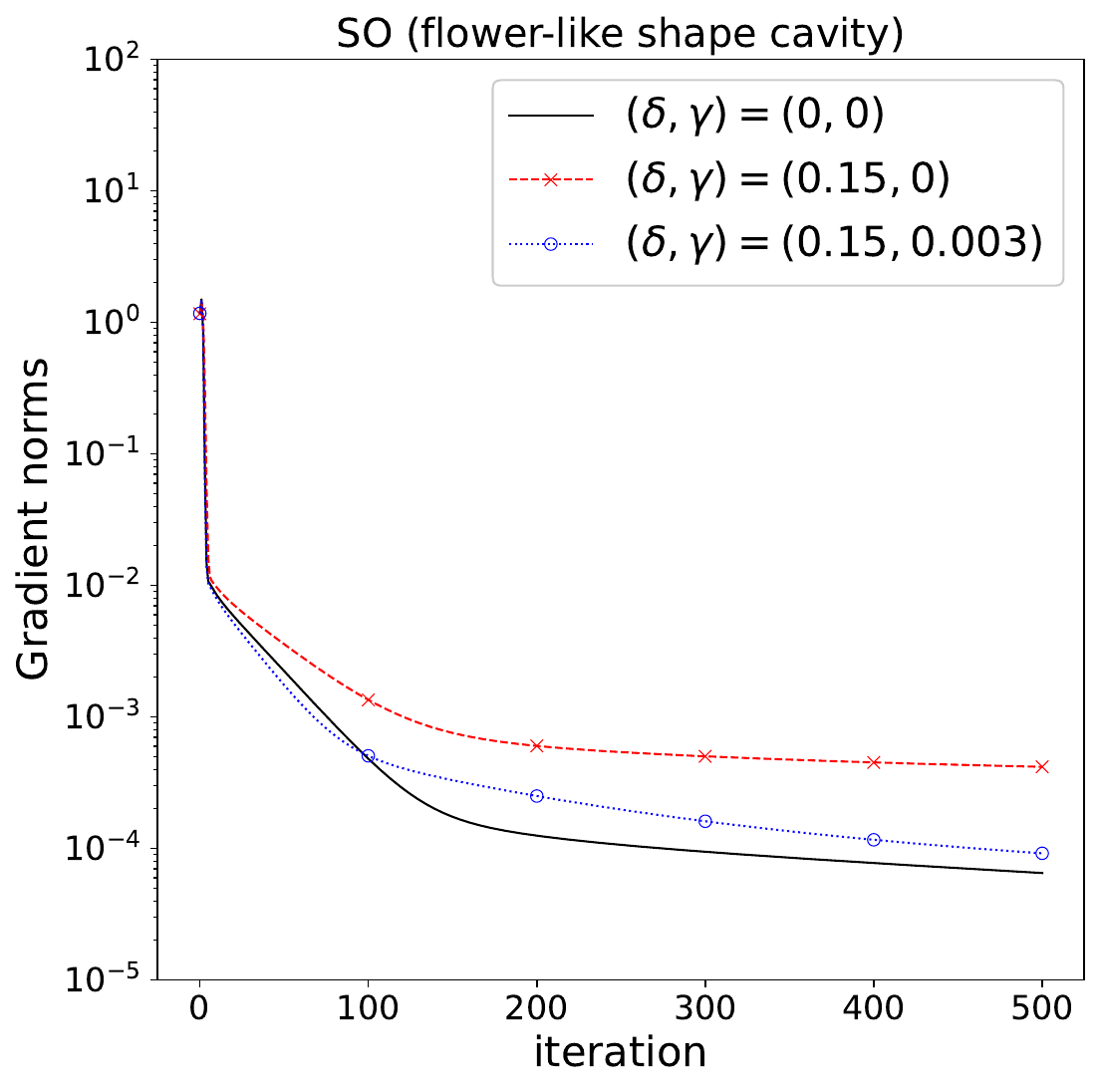}} \hfill
\resizebox{0.24\linewidth}{!}{\includegraphics{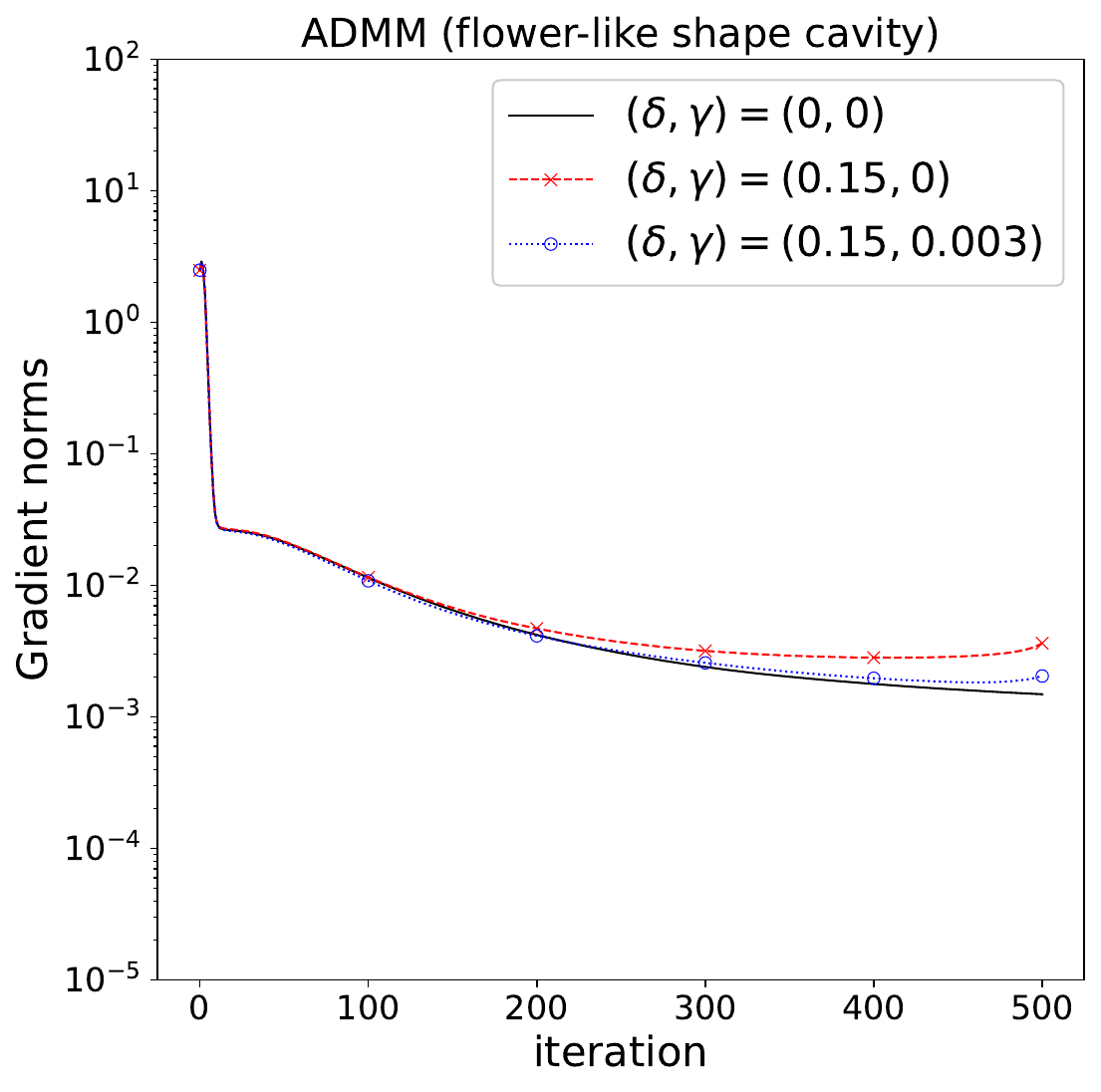}} \hfill
\caption{Histories of costs and gradient norms via SO and via ADMM}
\label{fig:figure3c}
\end{figure}

The numerical results are summarized in Figure \ref{fig:figure2a} through Figure \ref{fig:figure3c}, and the main observations align with those from the experiments in the 2D cases. Specifically, the results from ADMM are more accurate than those from SO, as ADMM can reconstruct the exact cavity with good accuracy even in the presence of noisy data, especially for large-sized cavities.

For the case of large cavities, refer to Figure \ref{fig:figure2a} for the reconstruction with exact data, Figure \ref{fig:figure2b} for the reconstruction with noisy data at a noise level of $\delta = 15\%$ without regularization, and Figure \ref{fig:figure2c} for the reconstruction with noisy data at a noise level of $\delta = 15\%$, now with regularization. Cross comparisons of the exact and computed shapes for both SO and ADMM for the last-mentioned test are shown in Figure \ref{fig:figure2d} and \ref{fig:figure2e}, and their corresponding histories of cost values and gradient norms are plotted in Figure \ref{fig:figure2f}.

From the experiments, the discernible effect of regularization becomes evident. Indeed, with the regularization parameter set to $\gamma=0.003$, the reconstruction of the exact cavity at $15\%$ noisy data becomes less rough compared to the case when no regularization is applied, as expected (compare the plots in Figure \ref{fig:figure2b} and Figure \ref{fig:figure2c}).

Meanwhile, as anticipated, reconstructing smaller-shaped cavities (e.g., see Figure \ref{fig:figure3a}) poses greater challenges, as illustrated by the reconstructed shapes in Figure \ref{fig:figure3b}. 
This difficulty is expected when the cavity is small in size or when its bounding surface is distant from the exterior boundary, where measurements are taken. 
It is worth noting that this challenge has been observed in previous studies.
Nevertheless, even in the case of smaller cavities, ADMM provides a fair and faster reconstruction, significantly outperforming SO. 
This superiority is evident in the illustrated computed shapes shown in Figure \ref{fig:figure3b}. For further clarification, we present the histories of cost values and gradient norms in Figure \ref{fig:figure3c}.
We observed, similar to the 2D case, a convergence in the cost values obtained via ADMM. 
It is important to point out here, however, that, for the case of noisy data, a better stopping condition (e.g., early termination of the algorithm) provides a more accurate reconstruction of the shape.
Additionally, we mention that, although not shown here, we tested other values of $\beta$, and from our experiments, we observed that, similar to previous experiments in 2D, setting $\beta$ to very small or large values is ineffective. 
It appears that the drawback of ADMM in the shape optimization framework, in general, lies in the necessity to calibrate the value of $\beta$ to achieve an effective implementation of the method.

In summary, with the appropriate choice of $\beta$, the proposed ADMM in the shape optimization setting offers a significant improvement in the results obtained from shape optimization (SO) for cavity detection. 
However, in general, we observe a computational time increase of about $15\%$ to $30\%$ when using ADMM compared to the conventional shape optimization method. (Refer to Table \ref{tab:tableCPU} for recorded computational times corresponding to the last test case with a smaller cavity.)
%
%
%
%
%
%
%
%
%
%
%
%
%
%
%
%
%
%
%
%
%
%
%
\begin{table}
\begin{center}
{\renewcommand{\arraystretch}{1.2} 
\begin{tabular}{||c c c c||} 
 \hline
Computational setup & SO & ADMM & Percentage Increase \\ [0.5ex] 
 \hline\hline
 exact data ($\delta = 0\%$) & 430 & 575 & $\approx 34\%$ \\ 
 \hline
 noisy data ($\delta = 15\%$) without regularization ($\gamma = 0$) & 442 & 575 & $\approx 30\%$ \\
 \hline
 \hline
 noisy data ($\delta = 15\%$) with regularization ($\gamma = 0.003$) & 492 & 573 & $\approx 17\%$ \\ [1ex] 
 \hline
\end{tabular}
}
\end{center}
\caption{Computational time in seconds}
\label{tab:tableCPU}
\end{table}
\section{Concluding Remarks and Future Work}\label{sec:conclusion}
We have proposed a novel application of the Alternating Direction Method of Multipliers (ADMM) to formulate a PDE-constrained shape optimal control problem with inequality constraints, addressing a shape inverse problem. 
The problem involves identifying a perfectly conducting inclusion within a bounded domain based on boundary measurements. As demonstrated here, the ADMM approach decouples the shape problem from the PDE constraint at each iteration. 
The resulting modified ADMM algorithm has proven to be easily implementable and numerically efficient in solving the shape inverse problem. 
Our proposed strategy improves the accuracy of the identification while maintaining a convergence rate comparable to classical shape optimization methods. 
However, it is important to note that ADMM introduces additional computational costs (approximately $15\%$-$35\%$ of the conventional shape optimization method) due to the steps involved in the approximation procedure. 
This drawback is expected, given that ADMM requires extra computational steps. Nevertheless, the improvements are found to be significant, as ADMM was able to reconstruct the concavities of the exact unknown interior boundary with high accuracy. 
In our future work, we will delve deeper--both theoretically and numerically--into investigating the application of ADMM for solving more complex shape inverse problems within the shape optimization framework.

\section*{Acknowledgement}
The work of JFTR is partially supported by the Japan Society for the Promotion of Science (JSPS) Grant-in-Aid for Early-Career Scientists under Japan Grant Number JP23K13012 and by the Japan Science and Technology Agency (JST) CREST Grant Number JPMJCR2014.




\vspace{-5pt}
\bibliographystyle{elsarticle-num.bst}
\bibliography{main.bib}







\end{document}